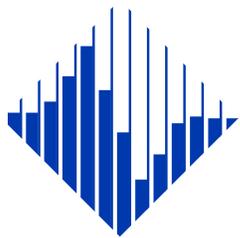

# PACIFIC EARTHQUAKE ENGINEERING RESEARCH CENTER

# Resilience of Critical Structures, Infrastructure, and Communities

**Gian Paolo Cimellaro**
**Ali Zamani Noori**
**Omar Kammouh**
Department of Structural, Building and Geotechnical Engineering
Politecnico di Torino, Torino

**Vesna Terzic**
Department of Civil Engineering and Construction
Engineering Management
California State University Long Beach

**Stephen A. Mahin**
Pacific Earthquake Engineering Research Center
University of California, Berkeley

# Resilience of Critical Structures, Infrastructure, and Communities

**Gian Paolo Cimellaro**

**Ali Zamani Noori**

**Omar Kammouh**

Department of Structural, Building and Geotechnical Engineering,
Politecnico di Torino, Torino

**Vesna Terzic**

Department of Civil Engineering and Construction
Engineering Management
California State University Long Beach

**Stephen A. Mahin**

Pacific Earthquake Engineering Research Center
University of California, Berkeley

# ABSTRACT

In recent years, the concept of resilience has been introduced to the field of engineering as it relates to disaster mitigation and management. However, the built environment is only one element that supports community functionality. Maintaining community functionality during and after a disaster, defined as resilience, is influenced by multiple components. This report summarizes the research activities of the first two years of an ongoing collaboration between the Politecnico di Torino and the University of California, Berkeley, in the field of disaster resilience.

    Chapter 1 focuses on the economic dimension of disaster resilience with an application to the San Francisco Bay Area; Chapter 2 analyzes the option of using base-isolation systems to improve the resilience of hospitals and school buildings; Chapter 3 investigates the possibility to adopt discrete event simulation models and a meta-model to measure the resilience of the emergency department of a hospital; Chapter 4 applies the meta-model developed in Chapter 3 to the hospital network in the San Francisco Bay Area, showing the potential of the model for design purposes Chapter 5 uses a questionnaire combined with factorial analysis to evaluate the resilience of a hospital; Chapter 6 applies the concept of agent-based models to analyze the performance of socio-technical networks during an emergency. Two applications are shown: a museum and a train station; Chapter 7 defines restoration fragility functions as tools to measure uncertainties in the restoration process; and Chapter 8 focuses on modeling infrastructure interdependencies using temporal networks at different spatial scales.

# ACKNOWLEDGMENTS

The work presented in this report summarizes the joint research activities performed at PEER from the Master students of the Politecnico di Torino. The work of each student is represented by each chapter in the report. Any opinions, findings, and conclusions or recommendations expressed in this material are those of the authors and do not necessarily reflect those of the Pacific Earthquake Engineering Research Center (PEER).

The research leading to these results has received funding from the European Research Council under the Grant Agreement no. ERC_IDEal reSCUE_637842 of the project IDEAL RESCUE - Integrated DEsign and control of Sustainable CommUnities during Emergencies and from the European Community's Seventh Framework Programme - Marie Curie International Outgoing Fellowship (IOF) Actions-FP7/2007-2013 under the Grant Agreement n°PIOF-GA-2012-329871 of the project IRUSAT— Improving Resilience of Urban Societies through Advanced Technologies.

The project has also received partial funding by the project "Internationalization" of the Politecnico di Torino sponsored by the Compagnia di San Paolo Foundation.

Any opinions, findings, and conclusions or recommendations expressed in this material are those of the author(s) and do not necessarily reflect those of the European Community's Seventh Framework Programme.

vi

# CONTENTS

# LIST OF TABLES

# LIST OF FIGURES

# 1  Modeling Disaster Resilience and Interdependencies of Physical Infrastructure and Economic Sectors

## 1.1  INTRODUCTION

Today modern communities and their economic sectors are highly interdependent, making them more vulnerable to natural and human-caused disruptive events. While interdependencies are considered to enhance a community in normal operating conditions as they promote a greater economic growth, they have serious drawbacks in the aftermath of natural or man-made disasters. After a catastrophic event such as an earthquake, the damaged region is affected by different types of losses. This depends on the level of interdependency among different economic sectors, their business downtime, and restoration. Disruptive events generate a chain reaction of production cutbacks among successive rounds of customers (demand reductions) and suppliers spreading through the entire regional economy.

In the last decade, several studies have addressed the behavior of communities in the aftermath of a disaster event from a global perspective. Cimellaro et al. [2016] defined seven dimensions of community resilience based on an index they derived that can characterize the behavior of a region. Among these dimensions, the economic effects post-disaster are certainly the most important. The study demonstrated that measuring the economic changes at the regional level triggered by a disaster is a crucial step towards disaster risk reduction. Recent discussions within the engineering and economics communities have focused on defining the economic structure of local and regional communities and the connections between productive sectors and consumers that will maximize the economic benefits.

The economic structure needs to be defined by identifying the weak points of the economic system that are to be "protected" to enhance economic resiliency in response to natural and man-made disasters. Webb et al. [2002] analyzed business characteristics that influence the long-term recovery after a catastrophic event. Rose and Liao [2005] focused on the estimation of indirect economic losses within a region stemming from a disruption of water service. Wasileski et al. [2011] examined how physical damage to the infrastructure, lifeline disruption, and business characteristics, among other factors, impact business closure and relocation following major disasters. Pant et al. [2013] developed a specific approach for the evaluation of interdependencies among multiple infrastructures able to support decision-making and resource allocation. Although the studies mentioned above deal with the economic effects of extreme events in the communities, they all focus on specific aspects of the problem.

In this chapter, a methodology is proposed for estimating the economic dimension of resilience that encompasses all types of losses that should be taken into account to predict the effects of natural disaster on a regional economy. The proposed methodology uses an economic framework that is an extension of HAZUS [2005] framework. The proposed framework divides losses in two main categories: direct and indirect. The direct losses include economic losses caused by physical damage to buildings and utilities. Economic losses generated by physical damage of buildings are based on the model by Terzic et al. [2014a]. Within the proposed methodology, a further step is made to develop a correlation between physical damage to utilities, their downtime, and losses. To account for the indirect losses that stem from the interdependence between different economic sectors, the structural growth model (SGM) introduced by Li [2010] is utilized as described by Martinelli et al. [2014].

The proposed methodology is demonstrated in a case study for the San Francisco Bay Area, a region with strong initiative in reducing earthquake risk by identifying performance goals that are to be achieved through design to improve resiliency of the region as described by Poland et al. [2009]. For this case study, the economic resiliency of the region is based on combination of the existing and simulated data. Simulated data were only used if the real data were not available. To take into account the uncertainty of the data used in the analysis, a sensitivity analysis was performed to identify preventive measures that could be facilitated to improve economic resiliency. Finally, an economic performance index, named economic resilience index $R$, was used as a measure of the economic ability of a region to withstand catastrophic events.

## 1.2   DESCRIPTION OF THE METHODOLOGY

Natural disasters may generate significant economic losses at both the local (regional) and global level. Given that regional losses are significantly higher than global losses, they will be the focus of this study. To estimate the total economic losses of a region struck by a natural disaster, the losses are disaggregated into direct and indirect losses. The direct economic losses are associated with the business-interruption cost due to physical damage of structures (buildings and lifelines), and indirect losses are associated with the disrupted inter-industry transactions.

### 1.2.1   Economic Loss Framework

The proposed framework for calculating economic losses due to hazard events, schematically presented in Figure 1.1, is an extension of HAZUS [2005]. HAZUS is software developed by the U.S. Federal Emergency Management Agency (FEMA) to estimate different types of losses generated by a natural hazard. The framework divides the losses in direct and indirect losses. The direct losses stem from building but also from utility damage (note, this kind of interdependency is disregarded by HAZUS) and are associated to the cost of reconstruction and business interruption. The indirect losses in the methodology are estimated as the general equilibrium effects of a disrupted inter-industry economy instead of being computed through the traditional Input-Output model of HAZUS.

Three main modifications of the framework shown in Figure 1.1 are applied to capture all the types of possible losses. The first is represented by the analysis of the industry loss of function due to the disruption of utilities. The second is given by a new method that is able to determine the probabilistic distribution of the time-dependent direct losses that affect a specific region of interest. Finally, the SGM is applied instead of the usual Input-Output model to quantify the indirect effects that arise when a cascade effect due to the business interdependencies occurs.

### 1.2.2 Direct Time-Dependent Losses

The proposed methodology refines the analysis of the time-dependent direct losses related to the building physical damage. Note that these losses do not include repair costs. The basic model assumes, as in HAZUS, that relocation occurs if the damage state of the building is greater or equal to moderate; in that case, the losses are given by relocation expenses, RE, rental income losses, RIL, and loss of income, LI.

Otherwise, the time-dependent direct losses are given only by the LI due to the loss of functionality that could arise even with slight damage to the building. Because the goal of the chapter is to quantify the global economic effects of a disaster on a specific region of interest, it should be taken into account that relocation may occur in different ways that influence the losses. To accomplish this goal, Equations (1.1), (1.2) and (1.3) have been implemented. The differences in the methodology adopted with respect to the HAZUS approach stem from two observations.

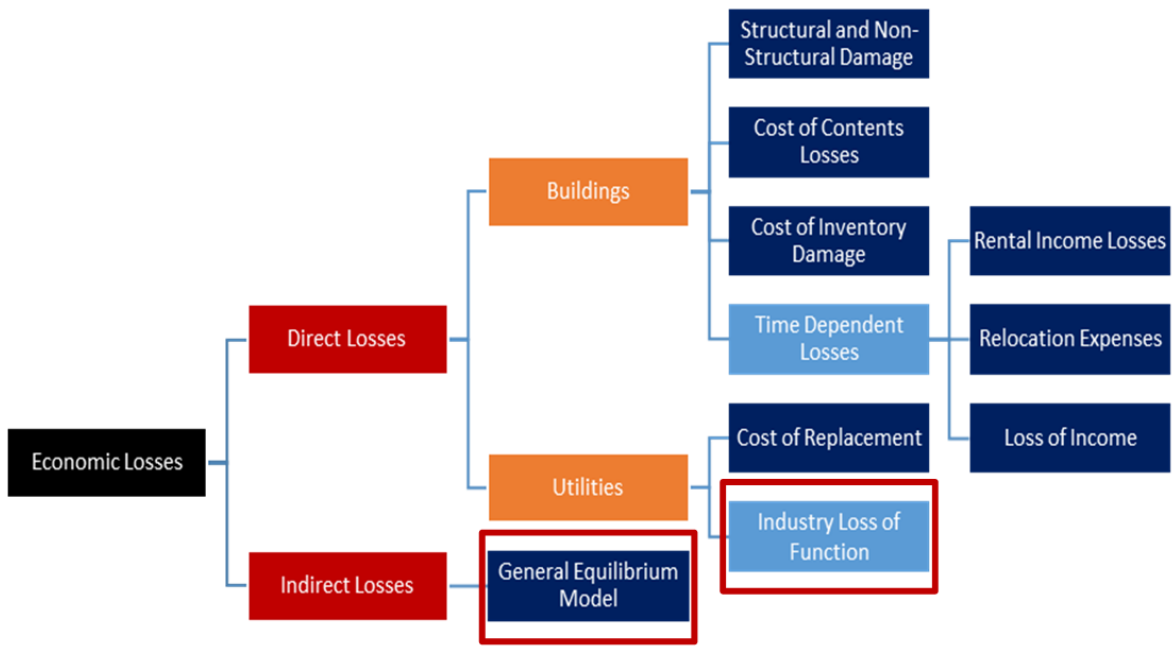

**Figure 1.1    Economic loss framework.**

$$RE_i = \left\{ (1-\%OO_i) \bullet \sum_{DS=3}^{5} \left( POSTR_{DS,i} \bullet DC_i \right) + \%OO_i \right.$$
$$\left. \bullet \sum_{DS=3}^{5} \left[ POSTR_{DS,i} \bullet DC_i (+RENT_i \bullet RT_{DS}) \right] \right\} \qquad (1.1)$$

$$LI_i = (1-RF_i) \bullet FA_i \bullet INC_i \bullet \left( \sum_{DS=1}^{5} POSTR_{DS,i} \bullet t_{DS} \right) \qquad (1.2)$$

$$RIL_i = (1-\%OO_i) \bullet FA_i \bullet RENT_i \bullet \left( \sum_{DS=3}^{5} POSTR_{DS,i} \bullet t_{DS} \right) \qquad (1.3)$$

where $\%OO_i$ is the percent owner occupied for occupancy $I$; $POSTR_{DS,I}$ is the probability of occupancy $i$ being in structural damage state $DS$; $DC$ is the disruption costs for occupancy $I$; $RENT_i$ is the rental cost for occupancy $I$: $RT_{DS}$ is the recovery time for the damage state $D$; $RF_i$ is the recapture factor for occupancy $I$; $INC_i$ is the income per day per square foot for occupancy $i$; and $t_{DS}$ is the period of time that depends on the $DS$, business property, and place of relocation (Figure 1.2).

First, since the goal of the research is to estimate the losses of a specific region, it is important to distinguish between inside and outside relocation. Second, HAZUS does not take into account the possibility that the industries forced to relocate own extra space in which move the activity, and that this space may be again inside or outside the region of interest.

The implemented algorithm takes into account these different possibilities by choosing different time windows to compute LI and RIL, depending on the location of the relocation space (HAZUS computes losses using the *loss of function* and the *recovery time*, respectively) and considering or not new rental costs and rental losses depending on if the property is owned by the business.

As described in HAZUS [2005], the lengths of the *recovery time* estimates of the median time for actual clean-up and repair, or construction. These estimates are extended to account for delays in decision-making, financing, inspection, and mobilization, and represent estimates of the median time for recovery of building functions.

The *recovery time* is then translated into a *loss of function* through multipliers that take into account the possibility that businesses rents alternative space or uses spare industrial capacity elsewhere. Once all the direct time-dependent losses for all the different occupancies in the different damage states are computed, the total direct time-dependent cost estimate is found by summing the total relocation expenses, rental income losses, and output losses. The flowchart of the method that refers to businesses that are owner occupied is represented in Figure 1.2.

The yellow blocks in the flowchart (Figure 1.2) are the decision blocks. Depending on these blocks, the algorithm takes different paths. Due to the scarcity of data, it is very difficult to obtain exact data for these blocks. For this reason a probabilistic approach has been adopted to take into account the uncertainty of the decision variables. However, if more data regarding the decision blocks become available, they could easily be substituted in the method to obtain outcomes that are more reliable. The probabilistic approach gives probability distributions for the different types of losses and offers the possibility to develop "sensitivity analyses" to identify the crucial factors that influence the losses.

The methodology implemented in this chapter is based on three assumptions: (1) that the greater the size of the business, the higher the possibility that the businesses own vacant space to relocate the activity; (2) that the probability that the vacant space is located within the region of interest is equal to the percentage of vacant buildings in the region; and (3) it is assumed that the longer the recovery time, the higher the probability that the external relocation will be permanent. These three assumptions have been translated into numeric values within the methodology. For example, the first assumption determines the decision block "*vacant space owned by the industry*" and is calibrated by assuming that a medium size business has 50% probability of owning vacant space or has excess capacity in other branches. Then, the probability density function of the business size for each building's occupancy in all fifty states is found considering the average number of employees per industry and assuming a normal distribution as shown in Figure 1.3. Finally, the cumulative distribution function is derived from the probability density function; see Figure 1.4. The cumulative distribution function is used to determine the values of the probabilities of the different industries to own extra space.

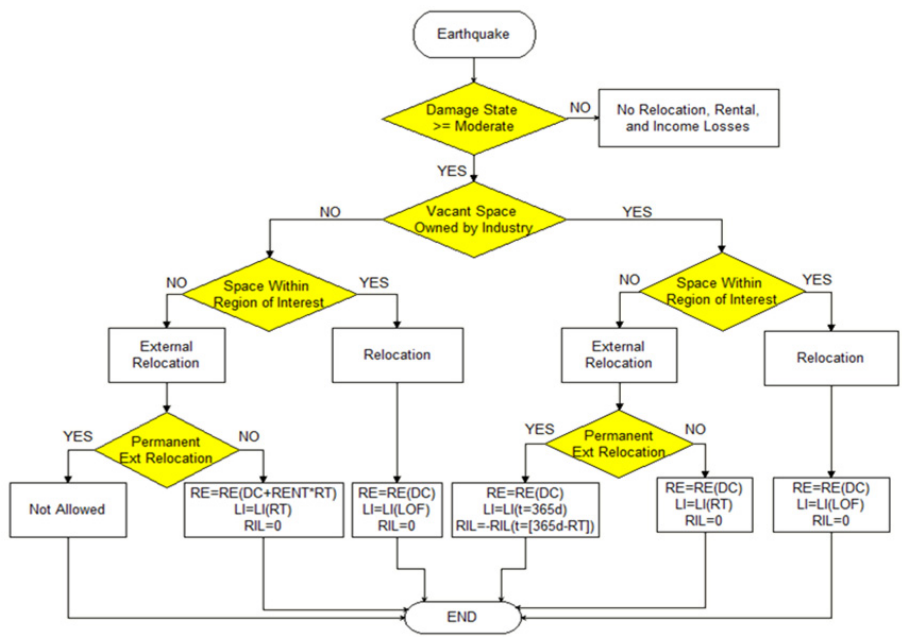

**Figure 1.2**     **Time-dependent losses algorithm for owner-occupied businesses.**

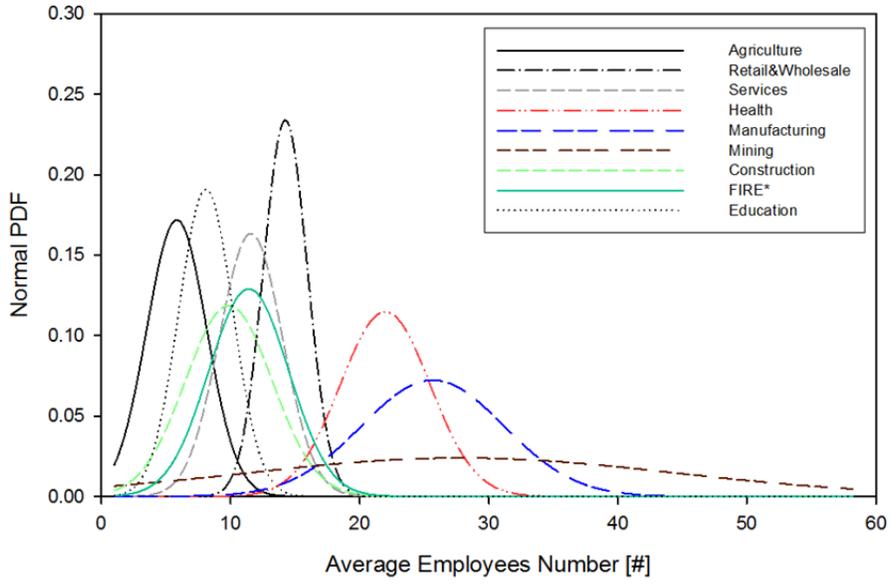

**Figure 1.3**     Probability density functions of the businesses sizes.

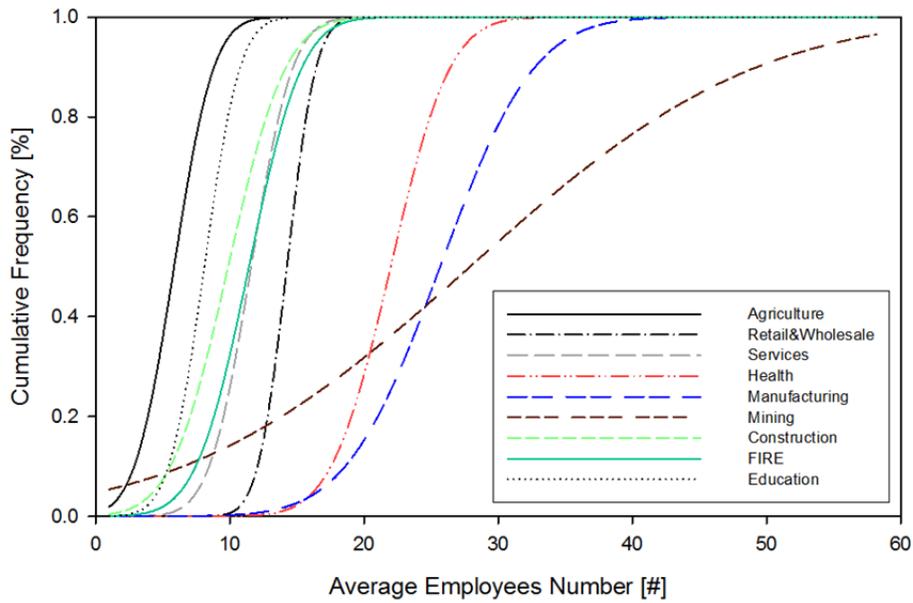

**Figure 1.4**     Cumulative distribution functions of the size of the businesses.

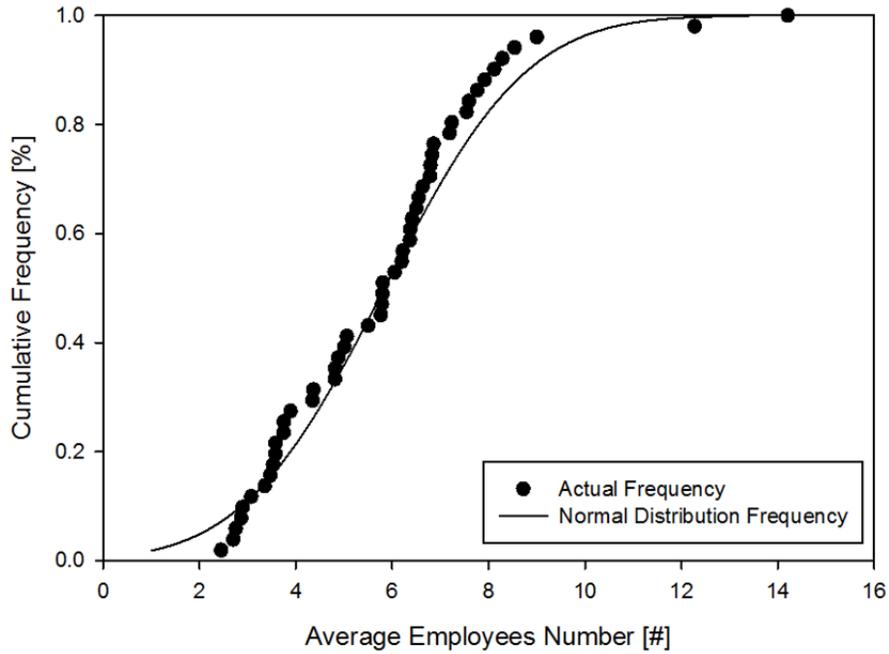

**Figure 1.5**     Validation of the normal distribution assumption.

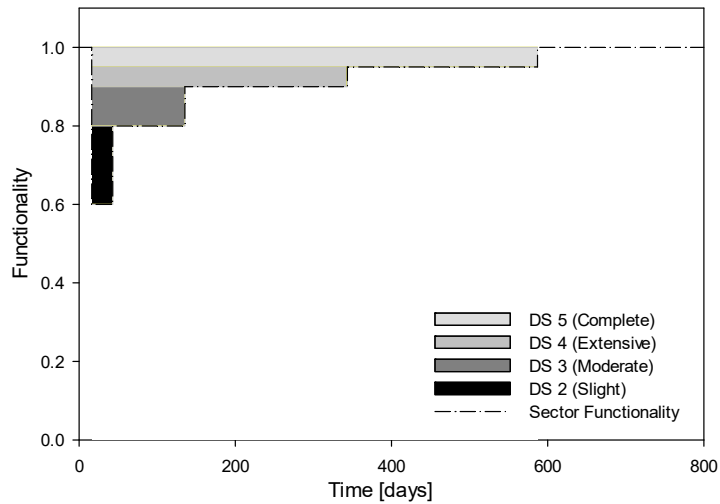

**Figure 1.6**     Loss of functionality for the educational sector for different damage states.

In other words, the mean value of the probability density function represents the typical size of a business which has 50% chance of owning vacant space while the value of the cumulative distribution function in correspondence of the average business size in the region of interest represents the probability of the sector to own vacant space in the probabilistic methodology adopted. The probability that the vacant space is located within the region of interest has been set equal to the percentage of vacant buildings in the region, which has been taken directly from HAZUS database.

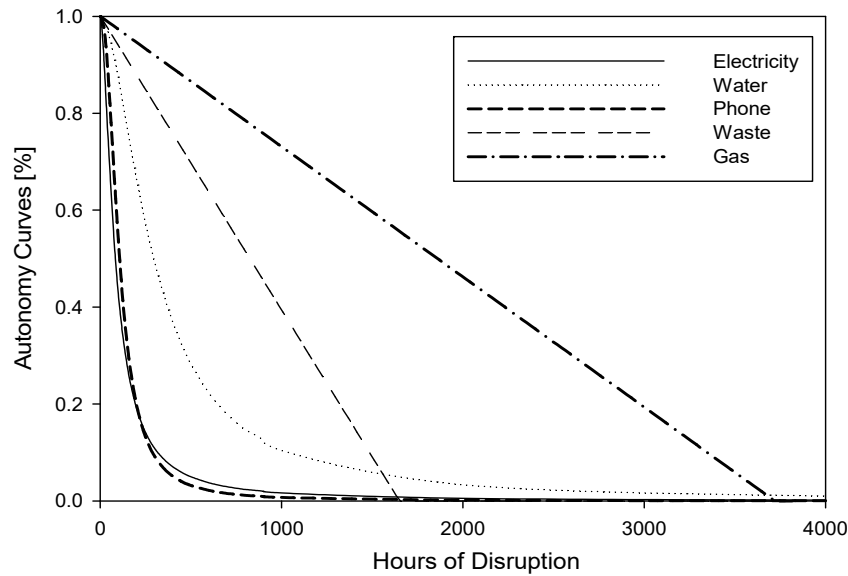

**Figure 1.7**  Autonomy factor curves of different lifelines for the retail and wholesale sector.

It has been assumed that after a specific number of days of relocation (which has been verified by sensitivities analyses), the relocation becomes permanent. This assumption has been implemented only for occupancies that are not of primary importance, like hospitals and schools, since a permanent relocation of these facilities is unlikely.

The data used to derive the business size distributions have been taken from the Economic Census (2007–2012), and the assumption of the normal distribution has been verified for each sector, as shown in Figure 1.5 for the Agricultural sector in the Bay Area case study. Note that the income losses considered in this section refer to the output losses suffered by the industries, which eventually represent the loss of functionality of each sector.

The cost of business interruption due to the physical damage of buildings is represented through a graph that shows the normalized output losses as a step function, where the different steps shown in Figure 1.6 for the Educational sector represents the number of damage states that contribute to the loss of business functionality.

After computing the building damage losses, the business losses due to lifelines disruptions are also taken into account in the proposed methodology. However, due to the scarcity of data, a hybrid approach has been adopted where both simulated and real data have been used. For example, the lifelines functionalities after the event are obtained by using the simulated data given by HAZUS. The real data is represented by the probability of business closure due to lifeline disruption. They have been derived using a procedure similar to the one explained by Chang et al. [2002] using data collected from surveys conducted on two natural disasters (the Northridge earthquake and Des Moines flood) described in the works of Tierney [1995] and the simulated results given by Rose et al. [2007]. In particular, a new function called an autonomy curve, which corresponds to the probability of business closure for a given lifeline, is derived using Equations (1.5) and (1.6). These autonomy curves represent the ability of each economic sector to withstand a utility outage of a different entity without losing functionality.

$$AF_i = 1 - P_{BI,i} \tag{1.4}$$

$$P_{BI,i} = P_{BC,j} \cdot P_{UO,j} \cdot \alpha_i \tag{1.5}$$

where $P_{BI,i}$ is the probability of business interruption due to utility $i$ outage; $P_{BC,j}$ is the probability of business closure for occupancy $j$; $P_{UO,ij}$ is percentage of business with utility $i$ outage for occupancy $j$; $\alpha_i$ is the average percentage of businesses that closed due to utility $i$ outage; and $AF_i$ is the autonomy factor of the sector on utility $i$.

The autonomy factor curves have been calibrated using the known temporal lifeline outage in the case study considered, while a different type of curve has been selected depending on the type of utility considered. For example, when analyzing the Retail and Wholesale sector for the electricity, water, and phone networks, a four-parameter logistic function was chosen; for the waste and gas system, a multi-linear curve was selected as shown Figure 1.7. All the figures shown in this chapter refer to the San Francisco Bay Area case study. The influence of each utility disruption on the economic sector functionalities is modeled applying the autonomy curves (AF) to determine the new sector functionalities using the following equation:

$$f_{sector}(t) = f_{utility}(t) + \left[1 - f_{utility}(t)\right] \cdot AF_{utility}(t) \tag{1.6}$$

where $f_{sector}$ is the functionality of the economic sector, and $f_{utility}$ is the functionality of the utility.

The limitation of Equation (1.6) is that the normal operating condition after lifeline disruption is reached is the same for both the Lifeline and the Economic sector, as shown in Figure 1.8, which considers the example of the water service. In reality, a lag exists between Economic sector and Lifeline recovery. Therefore, Equation (1.6) can be used until the Economic sector begins to recover. Then a lag factor, $\theta$, is introduced to take into account the delay of functionality of the other. The mathematical formulation for the lag factor is given by:

$$\begin{cases} \theta = 0 & t < t_r \\ \theta = \dfrac{t - t_r}{X_{gg}} & t_r < t < t_{fr} \end{cases} \tag{1.7}$$

$$f_{sector}(t + \theta \bullet t_r) = f_{utility}(t) + \left[1 - f_{utility}(t)\right] \bullet AF_{utility}(t) \tag{1.8}$$

where $t_r$ is the time instant when the recovery of the economic sector starts using Equation (1.6); $t_{fr}$ is the time instant when the recovery of the economic sector ends using Equation (1.6); $X_{gg}$ is the lag time of the economic sector with respect to the utility; and $AF_{utility}$ is the autonomy curves of the Economic sector with respect to the utility.

As a first approximation, the lag time needs to be calibrated. The lag time $\theta$ for the Economic sector is assumed as a fraction of the utility restoration time. Once all the autonomy curves—which describe the interdependencies between the Economic sector and the different lifelines—are determined, they are combined with the Economic sector functionality for determining the effect of all the different utilities.

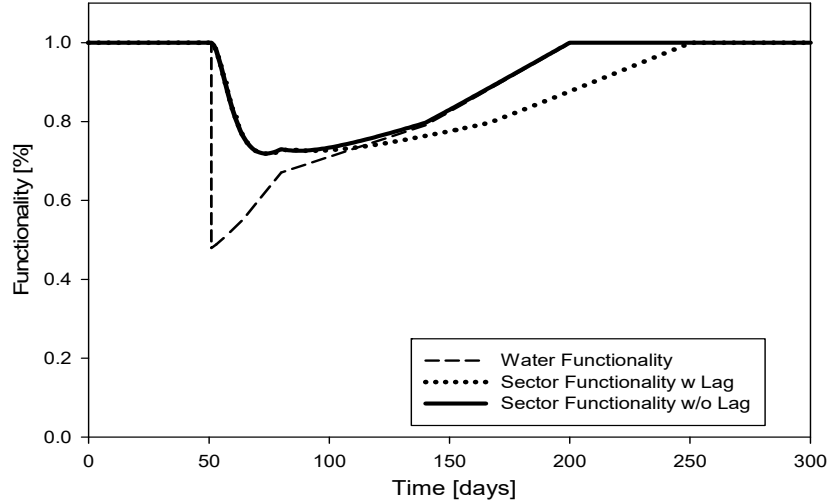

**Figure 1.8**     **Sector functionality influenced by the water service restoration.**

The new updated functionality curves are then combined to determine a single functionality curve for each Economic sector that captures the interdependencies between each lifeline. Note that the methodology overestimates losses due to utility disruption since it has been assumed that the businesses were affected separately by the utilities, which affects sector functionality most. Moreover, interdependencies are considered separately one by one, and it has not been taken into account the possibility that businesses that are forced to close due to utility disruption can reduce their losses by interacting with other utilities or by making up production at different times. To reduce this overestimation, the recapture factors provided by HAZUS have been used to decrease the losses.

Finally, the losses due to utilities disruption for each economic sector were summed with the output losses due to building damage, and a loss range was determined. The lower bound of this loss range is represented by the envelope of the two functionality curves affected separately by physical damage and the utilities disruption. The upper bound is represented by the sum of the two functionality losses. Then, depending on the conditional probability for a business to be affected simultaneously by building physical damage and utility disruption, it has been found a probable value within this range. Equation (1.9) is adopted to compute the global functionality.

$$F_{tot,\text{sec}} = \min(F_{\text{sec},utilities}; F_{\text{sec},building}) - P(BD \cap UO) \cdot \left[1 - \max\left(F_{\text{sec},utilities}; F_{\text{sec},building}\right)\right] \quad (1.9)$$

where $F_{sec,utilities}$ is the functionality of the sector influenced by the utilities; $F_{sec,building}$ is the functionality of the sector influenced by the building damage; and $P(BD \cap UO)$ is the probability that business is simultaneously affected by building damage; and ($BD$) and utility outage ($UO$).

### 1.2.3 Indirect Losses

After estimating the direct effect of the disaster event on each economic sector, the methodology applies the SGM to the scenario of interest, as described in Cimellaro et al. [2014], to estimate the indirect effects that stems from the interdependence between the sectors. In other words, the model applies an initial perturbation to the business functionalities that corresponds to the direct damages experienced by the sectors and then evaluates the recovery process, which is controlled by the price-adjustment velocity and by the depreciation factors of the goods. At the end of the analysis, it is possible to obtain the graph shown in Figure 1.9, which depicts the general equilibrium effects and from which the monetary losses due to the business interdependences can be derived.

### 1.2.3.1 The Structural Dynamic Growth Model

Li [2010] developed the structural growth model from the classical growth framework. This growth model can be used to compute general equilibrium effects. The model represents the production processes in the economy through two matrices: the input and the output coefficient matrices, respectively, **A** and **B**. For example, for the economy described in Equation (1.10), the two matrices are given by Equation (1.11):

$$\begin{cases} 280 \text{quarterswheat} + 12 \text{tonsiron} \Rightarrow 575 \text{quarterswheat} \\ 120 \text{quarterswheat} + 8 \text{tonsiron} \Rightarrow 20 \text{tonsiron} \end{cases} \quad (1.10)$$

$$\mathbf{A} = \begin{bmatrix} 56/115 & 6 \\ 12/575 & 2/5 \end{bmatrix}; \mathbf{B} = \mathbf{I} \quad (1.11)$$

where the $i$th column in matrix **A** represents the standard input bundle of agent $i$. In the classical economic growth framework, the equilibrium price vectors and equilibrium output vectors are the left and right P-F eigenvectors of **A**. The Structural Dynamic Growth model tries to integrate the market mechanism into the classical growth model by embedding in it an exchange process, which is represented by an exchange vector in order to reach equilibrium.

### 1.2.3.2 Exchange Process

The exchange process considers the economy as a discrete-time dynamic system and supposes economic activities such as price adjustment, exchange, production, etc., occur in turn in each period. With reference to the previous economic system, **S** in Equation (1.12) denotes the ($n \times m$) supply matrix, and **s** denotes the supply vector in the initial period.

$$\mathbf{S} = \begin{bmatrix} 575 & 0 \\ 0 & 20 \end{bmatrix}; \mathbf{s} = \begin{bmatrix} 575 \\ 20 \end{bmatrix} \quad (1.12)$$

Let **z** denote the vector consisting of purchase amounts of $m$ agents, where **z** is called the purchase vector or exchange vector (of standard input bundles), **Az** is called the sales vector of goods. It's possible to derive Equation (1.13) where $\hat{\mathbf{s}}$ represents the diagonal matrix with the vector **s** as the main diagonal and **u** the $n$-dimensional sales rate vector indicating the sales rates of $n$ goods.

$$\mathbf{u} \equiv \hat{\mathbf{s}}^{-1} \mathbf{A} \mathbf{z} \quad (1.13)$$

Under given price vector **p**, the purchase and sales values of $m$ agents are $\mathbf{p}'\mathbf{A}\hat{\mathbf{z}}$ and $\mathbf{p}'\hat{\mathbf{u}}\mathbf{S}$, respectively. It is assumed that the values of each agent purchases must be equal the value it sells, as in CGE income balance, thus obtaining Equation (1.14):

$$\mathbf{p}'\mathbf{A}\hat{\mathbf{z}} = \mathbf{p}'\hat{\mathbf{u}}\mathbf{S} \equiv \widehat{\mathbf{p}'\hat{\mathbf{s}}^{-1}\mathbf{A}\mathbf{z}}\mathbf{S} \quad (1.14)$$

If Equation (1.14) holds, there exists a unique normalized exchange vector that can be found by the following steps, which stands for the outcome of the exchange process:

- Step 1. Compute the matrix $\mathbf{Z} = \mathbf{p}\widehat{\mathbf{s}^{-1}\mathbf{A}\mathbf{z}}\mathbf{S}$
- Step 2. Find the normalized right P-F eigenvector of $\mathbf{Z}$, denoted by $\mathbf{x}$;
- Step 3. Find the minimal component of $\widehat{\mathbf{A}\mathbf{x}^{-1}\mathbf{s}}$, denoted by $\xi$;
- Step 4. Compute the exchange vector $\mathbf{z} = \xi\mathbf{s}$.

The economic system in period $t$ is represented by the following variables: $\mathbf{p}(t)$ = price vector; $\mathbf{S}(t)$ = supply matrix; $\mathbf{u}(t)$ = sales rate vector; and $\mathbf{z}(t)$ = exchange vector and production intensity vector. The market mechanism is embedded considering that in period $t+1$, the economy runs as in Equation (1.15) until the time where the system reaches the equilibrium.

$$\begin{cases} \mathbf{p}(t+1) = P[\mathbf{p}(t), \mathbf{u}(t)] \\ \mathbf{S}(t+1) = \mathbf{B}\mathbf{z}(t) + Q[e - u(t)S(t)] \\ [u(t+1), z(t+1)] = Z[A, p(t+1)S(t+1)] \end{cases} \quad (1.15)$$

where $P$ represent price adjustment process, $Q$ is the inventory depreciation function and stands for the depreciation process of inventories, and $Z$ is the exchange function depicted above.

### 1.2.3.3 General Equilibrium Effects After a Disruptive Event

Starting from the I-O matrices representative of the economy in normal operating condition, a shock that simulate earthquakes or other disasters is considered by modifying the exchange vector, which is the driver of the equilibrium process consistent with the direct damages experienced by the sectors. After the shock, the value of the loss due to the industries interdependencies can be estimated from the restoration curves of the system in Figure 1.9. To obtain data consistent with the economy of interest, the model requires calibration of two input parameters by which the model can be adapted to different regions. The first parameter is the price adjustment velocity; the second parameter is the depreciation factor of the goods.

### 1.2.4 Economic Resilience Index ($R_{EC}$)

Finally, the methodology evaluates the economic behavior of the analyzed region using a comprehensive resilience index, $R_{EC}$ determined according to the PEOPLES framework [Cimellaro et al. 2016]. $R_{EC}$ is the area under the function that is the sum of the direct losses and of the indirect losses normalized with respect to the value of the business functionality over the same control period; see Figure 1.10.

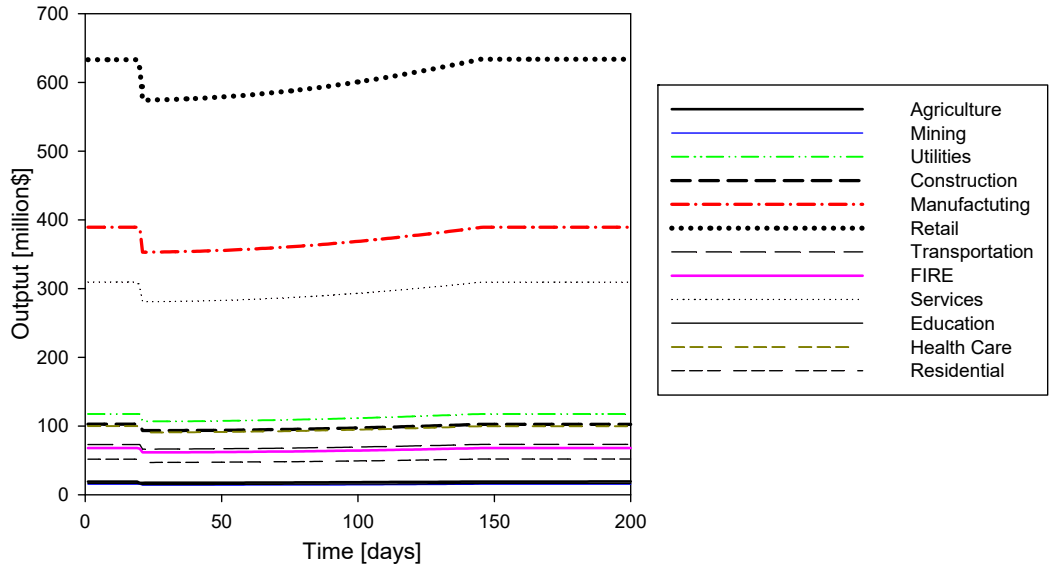

Figure 1.9  Indirect general equilibrium effects for the different sectors.

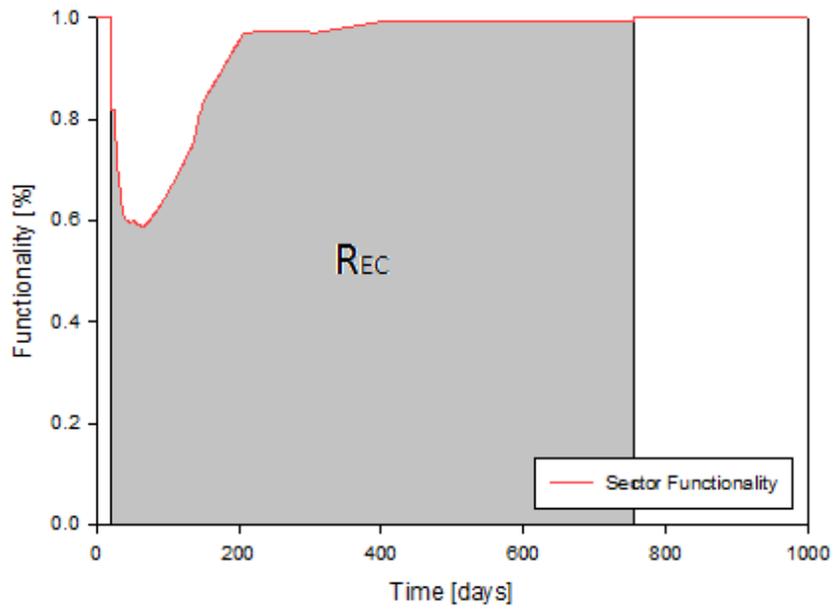

Figure 1.10  Area representing the economic resilience of the sector.

## 1.3  The San Francisco Bay Area Case Study

The San Francisco Bay Area shown in Figure 1.11 is considered as a case study to demonstrate implementation of the methodology described above. The U.S. Geological Survey (USGS) has estimated the maximum probability of 30% for a M > 6.7 in the Hayward Fault; therefore, the baseline scenario chosen is a M6.9 earthquake close to Oakland on the Hayward fault.

The structural/non-structural losses and the utilities functionalities have been derived from HAZUS after having loaded the soil map and the liquefaction susceptibility map of the region. Then the methodology described above is implemented to estimate the direct time-dependent losses. To do that, the values of *INCi* in Equation (1.2) have been updated coherently, with the output data of each sector published by the economic census.

The loss distributions for the economy in the region obtained by the methodology are shown in Figure 1.12. The estimated relocation expenses are $1.5 and 2.1 billion using the proposed methodology or HAZUS, respectively. The estimated rental income losses are $1.06 and $1.21 billion using the proposed methodology or HAZUS, respectively. Figure 1.13(a) represents the mean and dispersion of the direct output losses for each sector due to building damage; Figure 1.13b shows the losses taking into account disruption of utilities. The contribution of utility disruption to business loss of function has been computed assuming that the mean number of utilities that lost their businesses is 2.5, and that a business has about 50% of probability of being both affected by building damage and utility disruption.

The results show that for the Bay Area, the sectors experiencing the greatest losses are the Retail Wholesale, the Residential, and the Services & Government sectors; the Mining and Agriculture sectors experienced smaller losses. For the Retail & Wholesale and the Services & Government sectors, the greatest losses stem from the interdependencies, which affect the business interruption. For the Residential sector, a significant contribution is due to relocation expenses. The small losses of the Mining and Agriculture sectors are mainly due to the relatively small volume of business.

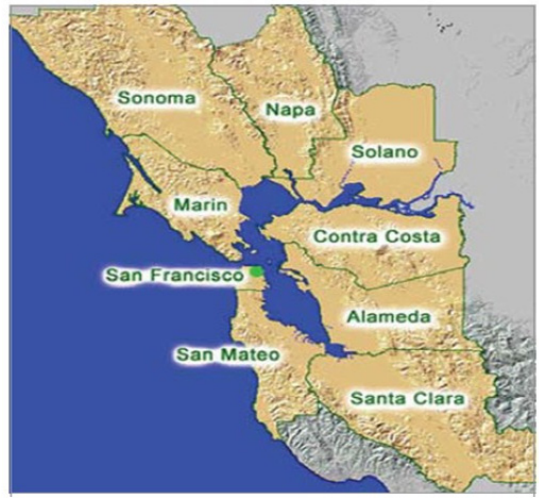

**Figure 1.11**    Region considered in the San Francisco Bay Area case study.

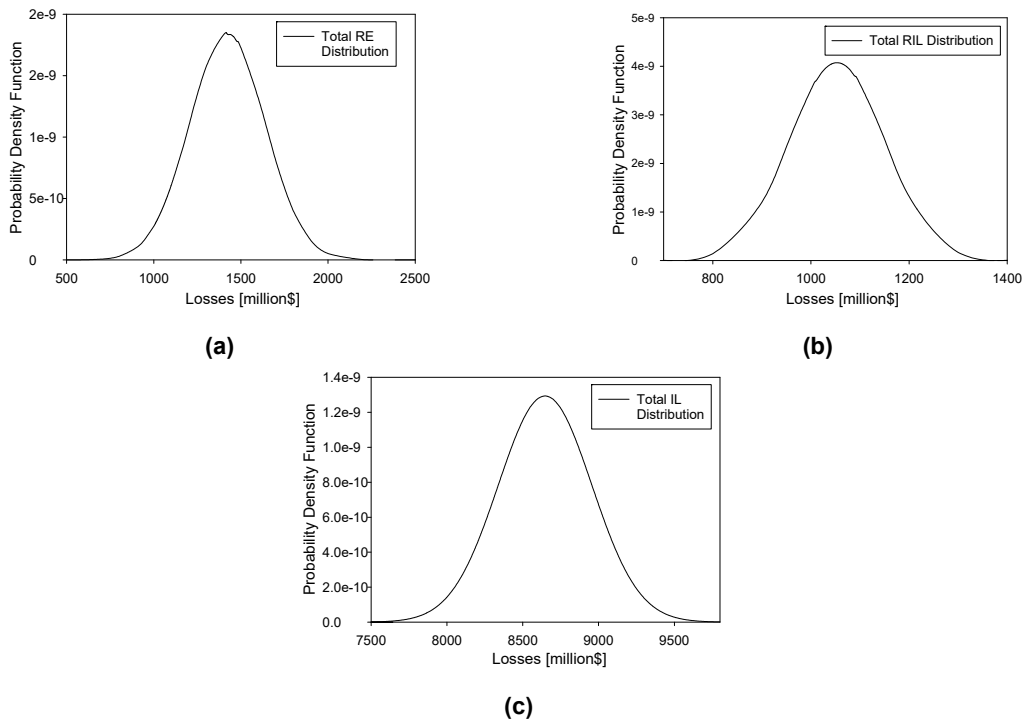

**Figure 1.12** Relocation expenses (a) RE, (b) loss of output IL, and (c) rental income losses RIL.

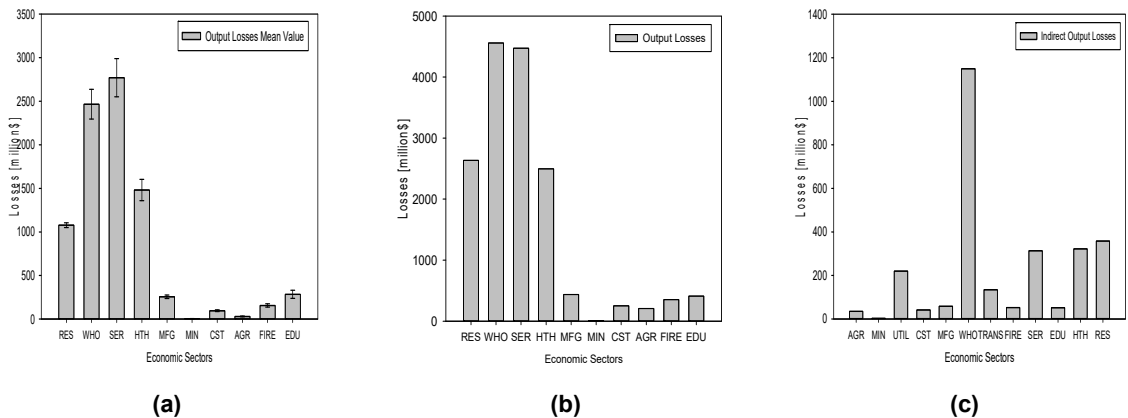

**Figure 1.13** Relocation expenses (a) RE (a), (b) loss of output IL, and (c) rental income losses RIL.

To estimate indirect losses, the SGM has been applied. The proposed method first computes the Input-Output matrix of the region of interest using the procedure explained by Chamberlain [2011] from the Make and Use matrices provided by the Bureau of Labor Statistics (1990–2010). Since public data are available only at the national level, the San Francisco Bay Area Input-Output matrix has been derived assuming a scaling factor based on the GDP value, which was applied to the U.S. Input-Output matrix.

The final indirect losses are represented in Figure 1.13(c). These losses have been derived by using depreciation factors for the different goods consistent with the values provided by the BEA. Similar to what has been done to estimate the output losses due to building damage and utility disruption, to take into account for the ability of business to make up production at different times indirect losses were reduced using the recapture factors provided by HAZUS. Note that the direct output losses were not been represented for the Utilities and the Transportation sector due to the unavailability of appropriate data for the methodology; however, direct output losses were taken into account in the total loss analysis considering the data provided in HAZUS. Finally, in the specific case study, the indirect losses represented approximately 15% of the direct losses.

### 1.3.1 Sensitivity Analysis

Performing a sensitivity analysis is a useful tool to demonstrate how the total losses are influenced by different parameters. Table 1.1 reports the results of different analyses performed for the total direct-time dependent losses; each analysis is distinguished by a specific assumption. Figure 1.14 summarizes the different outcomes. As is clear from Figure 1.14(b), a critical factor in minimizing loss is to avoid is permanent external relocation of businesses: this quadruples the loss of productivity of the sectors in the region; moreover, it does not allow the economy of the region to bounce back the pre-event levels of productivity since part of the functionality is lost forever. Indeed, the smaller losses are found assuming a high probability of vacant space within the region. Though it seems difficult to reach this condition in the reality, this observation can be taken as a guideline for the preventive measures implementation of the individual industries.

To show the uncertainty stemming from the unknown magnitude of the earthquake, Figure 1.15 represents the differences for the case of three different earthquakes. Finally, Table 1.2 summarizes the outcomes of the baseline case study from the $R_{EC}$ index point of view.

Table 1.1    Assumptions implemented in the sensitivity analysis.

| | |
|---|---|
| 1 | No excess capacity or vacant space owned by the business which relocate |
| 2 | Probability of vacant space within the region equal to 50% |
| 3 | Probability of vacant space within the region equal to 97% |
| 4 | High probability of permanent external relocation |
| 5 | Probability of excess capacity or vacant space owned by the industry close to 100% |
| 6 | Medium probability of permanent external relocation |

**Table 1.2**  Summary of total economic losses and resilience indices for the case study.

|  | M6.9 | M7.3 | M7.5 |
|---|---|---|---|
| Total Relocation Expenses (million$) | 1498 | 2129 | 2258 |
| Total Rental Income Losses (million$) | 1057 | 1511 | 1709 |
| Total Direct Output Losses (million$) | 18407 | 29094 | 25267 |
| Total Indirect Output Losses (million$) | 2555 | 4198 | 2359 |
| Total Structural/Non-Structural Losses (million$) | 29388 | 42804 | 49108 |
| $R_{EC}$ | 0.96 | 0.938 | 0.955 |

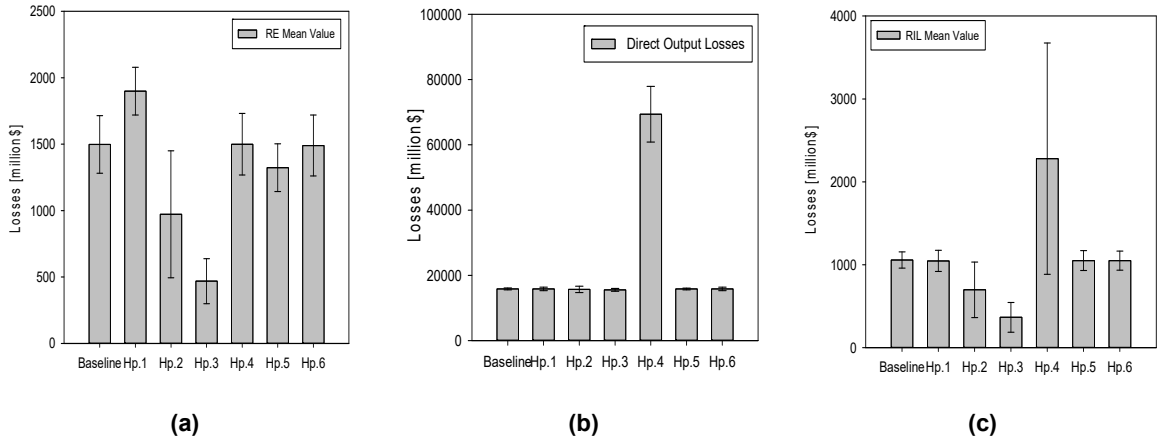

**Figure 1.14** (a) Time-dependent relocation expenses, (b) direct output losses, and (c) rental income losses for the different scenarios.

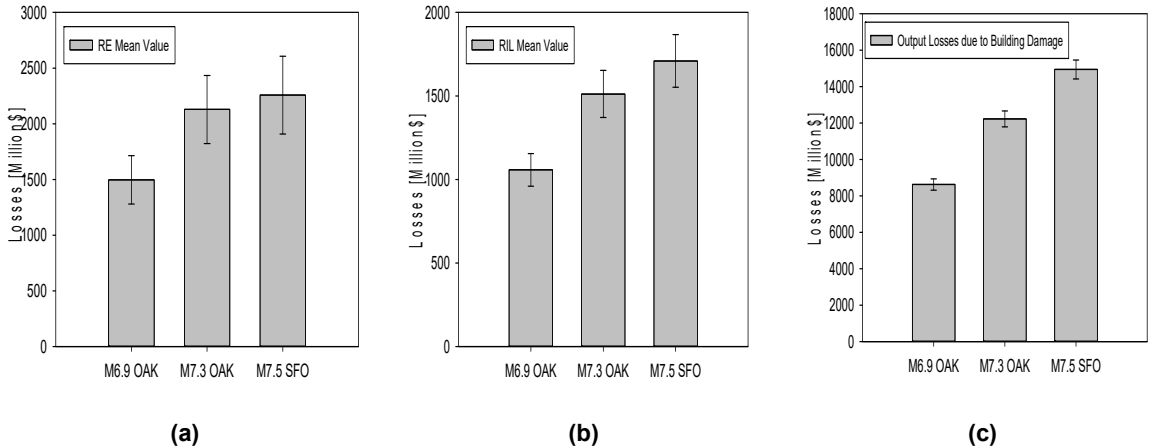

**Figure 1.15** Variation of the total direct time-dependent economic losses for the different assumed earthquake magnitudes.

## 1.4 REMARKS AND CONCLUSIONS

This chapter proposes a new methodology to evaluate the economic losses following a natural disaster. A new probabilistic framework to estimate time-dependent direct losses has been presented where the indirect losses have been estimated using the SGM, while interdependencies between the different economic sectors and lifelines during disruption are modeled using autonomy curves defined by the authors. The influence of different parameters uncertainties were analyzed and a final global economic resilience index $R_{EC}$ obtained, which can be used in more general community resilience frameworks (e.g., PEOPLES) to estimate the effects of the economic dimension.

Autonomy curves were derived using the probabilities of business closure collected from business surveys and simulation conducted mainly in California; thus, they are representative of the case study presented herein. These autonomy curves represent the main finding of the study.

As shown by the sensitivity analysis that simulates different earthquake magnitudes, the M7.5 earthquake in San Francisco causes less direct and indirect output losses compared to a M7.3 earthquake in Oakland, even though its magnitude is higher. This is because a M7.5 earthquake in San Francisco, considering the HAZUS approximation, will result in fewer utility losses; therefore, costs due to business interruption will be smaller. Further research will focus on addressing the limitations of the current methodology.

# 2 Utilizing Base-Isolation System to Increase Earthquake Resiliency

## 2.1 INTRODUCTION

In the event of structure failure, healthcare facilities represent a substantial hazard to human life (occupancy category III, per ICC IBC, [2012]). Therefore, they are designed following more stringent design requirements than buildings with residential and commercial occupancy. In the recent large earthquakes in Chile, New Zealand, and Japan, healthcare facilities generally performed well. However, healthcare closures due to extensive structural and nonstructural damage that resulted in the loss of their function have occurred [Miranda et al. 2012)]. As a result, increased attention is being placed on strategies to design facilities that are both safe and damage resistant. It is often presumed that such an approach increases costs to an unacceptable level. However, the cost-effectiveness of alternative design choices can be assessed using performance-based earthquake engineering (PBEE) methodology [Miranda and Aslani 2003] that quantifies expected future costs associated with damage repair, loss of functionality, casualties, etc.

This chapter presents the results of a study that compares the repair costs and repair times considering two designs for a three-story steel building: a high-performance special moment resisting frame (HP-SMRF) and a damage-resistant base-isolated intermediate moment resisting frame (BI-IMRF). The design of both systems comply with the occupancy category III [ICC 2012], allowing the building to serve as a healthcare facility. To aid in understanding the relative performance of these two systems, key engineering demand parameters (EDPs) (i.e., median values of maximum and residual story drifts and floor accelerations), repair costs, and repair times were compared at five hazard levels. These results were then used to estimate the resilience of the two systems. The value of PBEE analysis in identifying cost-effective seismic design strategies that produce more resilient, damage-resistant structures is discussed. Finally, two different ways to account for healthcare facilities equipment were considered to show the impact in terms of repair costs.

## 2.2 PBEE METHODOLOGY

The Pacific Earthquake Engineering Center (PEER) developed a methodology to probabilistically assess seismic performances of buildings, bridges, and other facilities. The purpose of the methodology is to provide stakeholders with the necessary information to make conscious decisions based on life-cycle considerations rather than on costs alone. The traditional way to measure a structure's performance is the quantitative evaluation of forces and deformations; however, these variables do not give a direct estimation of losses. Performance assessments provide information that can be immediately related to losses as needed by decision-makers and are assessed at a specified hazard level.

The performance-based process starts with the statement of performance objectives that take into account the site location and characterization, the configuration and occupancy of the building, the structural system, the non-structural components, and content location. Given the structure's characteristics, the process flowchart evaluates the building response after ground shaking and calculates the probable damage that can affect components as well as the associated losses.

The process is made of four logical progression steps that involve four generalized variables [Zareian et al. 2006]. The uncertainties in performance assessment can be calculated as outcome at each stage. Those variables are useful to quantitatively relate seismic hazard and other engineering parameters to structural performances. The four variables that are considered in the process are: Intensity Measure (IM), Engineering Demand Parameter (EDP), Damage Measure (DM), and Decision Variable (DV). Each quantity is defined as the conditional probability of exceedance that a certain demand does not exceed a fixed value.

The process starts with a hazard analysis whose aim is to describe the intensity of ground shaking. Earthquake hazard is the quantification of ground-motion intensity and, given the location $O$ and the design characteristics $D$, it can be expressed as the conditional probability that level of shaking is exceeded: $P(\text{im} \geq \text{IM} \mid o = O, d = D)$. In this study, an intensity-based approach was used to consider the response to a single acceleration response spectrum. Instead of a single value of spectral acceleration, a set of ground motions was selected and then scaled to match the response spectrum. Spectral acceleration corresponding to 5% damping of the first-mode period of the structure $S_a(T_1)$ can be used as ground-motion IM. The approach for the hazard analysis can be deterministic or probabilistic. Deterministic seismic hazard analysis (DSHA) considers a particular scenario; probabilistic seismic hazard analysis (PSHA) evaluates ground-motion uncertainties. In recent years, PSHA has become the preferred method for hazard analysis and was used herein. The PSHA outcome is the mean seismic hazard curve that relates the intensity measure $S_a(T_1)$ with its mean annual frequency of exceedance $\lambda(\text{IM})$, which varies in function of the site and of the vibration period of the building.

Given the ground motions exciting the system, a structural analysis is done to calculate the most representative EDPs describing the system. The most common parameters used to characterize the structural system are peak and residual interstory drifts, floor accelerations, component forces, and deformations. Structural response can be evaluated using nonlinear time history analysis or a simplified procedure based on equivalent static forces. Nonlinear analysis can be used to simulate the response of any building, while simplified analysis can be applied just for regular, low, and mid-rise structures with limited nonlinear response. In this study, a nonlinear response history analyses was performed to develop the IM–EDP relationships. The analysis was performed for multiple intensities and the outcome was the conditional probability that a certain demand parameter exceeded a defined EDP value for a given intensity measure: $P(\text{edp} \geq \text{EDP} \mid \text{im} = \text{IM})$. The estimation of EDPs accounts for uncertainties related to the structural mode, material properties, strength and deformation characteristics of structural components, variation in dead loads, and seismic mass. This study used the open-source software OpenSees (Open System for Earthquake Engineering Simulation; http://opensees.berkeley.edu) to perform the nonlinear analyses and calculate response parameters. The response quantities obtained from the analysis can be associated with the damage of structural/non-structural components and contents.

The next step in the process is to perform a damage analysis whose purpose is to correlate EDPs with DMs of all vulnerable components. Once the damage states for each element are defined, the conditional probability that a certain DS is exceeded for a given EDP [$P(dm \geq DM \mid edp = EDP)$] is evaluated. Component DSs are usually defined as damage levels corresponding to repair measures needed to restore components to their undamaged conditions. The relationships between DM and EDPs are usually expressed in the form of fragility curves.

The last stage of the process involves a loss analysis whose outcome is the definition of decision variables (DVs) that can be used by stakeholders to make informed decisions. In this step, damage information obtained by damage analysis is converted into final decision variables. The DVs involved in the process can be expressed in terms of economic losses, fatalities, duration of repairs or retrofit, or injuries. As seen in the previous steps, the results of this analysis can be represented as conditional probabilities that a DV overcomes a certain value given a specific DS [$P(dv = DV \mid dm = DM)$]. In contrast to DSs that are related to single elements, DVs are assigned at the system or building level. The process stages are graphically shown in Figure 2.1.

The four aforementioned PBEE steps can be analytically summarized in terms of a triple integral based on the total probability theorem. The integral expresses the probability that a decision variable exceeds a given value given as input the intensity measure of ground motion $P(dv \geq DV \mid im = IM)$. The combination of analysis stages is then given in accordance with the total probability by:

$$G(DV|IM) = \int_{allEDPs} \int_{allDMs} G(DV|DM) |dG(DM|EDP)| dG(EDP|IM) | d\lambda(IM) \quad (2.1)$$

As a result of this probabilistic analysis, it can be stated that DV values are not deterministic and vary with probability due to different sources of uncertainty.

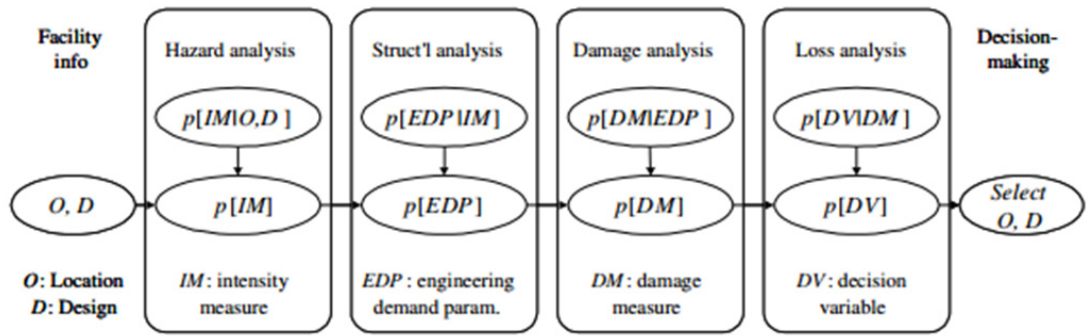

Figure 2.1    Underlying probabilistic framework [Moehle et al. 2003].

Healthcare facilities are critical buildings whose performance goals cannot be limited to merely reducing repair costs and guaranteed short downtimes. Energy generators and water suppliers have to remain functional in order to ensure the minimum required services. Finally, wall partitions and ceiling elements should maintain their isolating characteristics to help to prevent infection and fire propagation. The flow chart for the life-cycle cost analysis that was followed herein is shown in Figure 2.2. The scheme shows the various steps that were followed in the analysis to characterize the resilient features of both fixed-base and base-isolated buildings to facilitate the decision-making process.

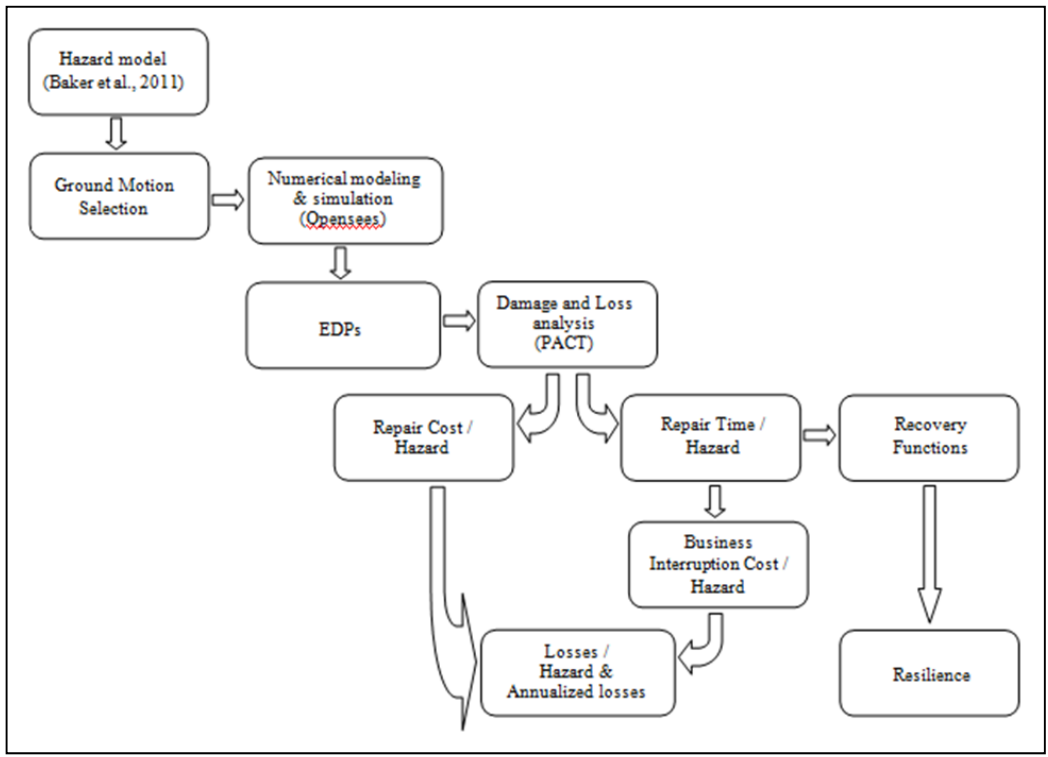

**Figure 2.2    Process flow chart.**

## 2.3    DESCRIPTION OF BUILDING

The study considered a three-story steel building located in Oakland, California, a site representative of the high-seismic hazard characteristics of western North America. The basic building plan dimensions are 120 ft (36.5 m) × 180 ft (54.9 m) with a bay spacing of 30 ft (9.1 m) in each direction. The building's plan is shown in Figure 2.3. The structure is located on relatively stiff soil having a shear-wave velocity on the top 30 m $V_{s30}$ equal to 360 m/sec. This velocity corresponds to the NEHRP site class C/D boundary with reference shear-wave velocity = 180 to 360 m/sec.

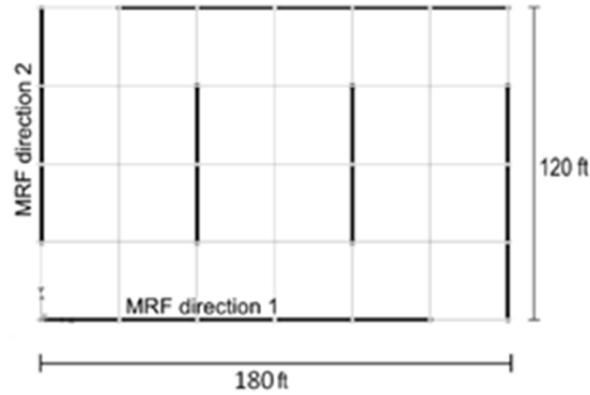

Figure 2.3    Building's plan [Mayencourt 2013].

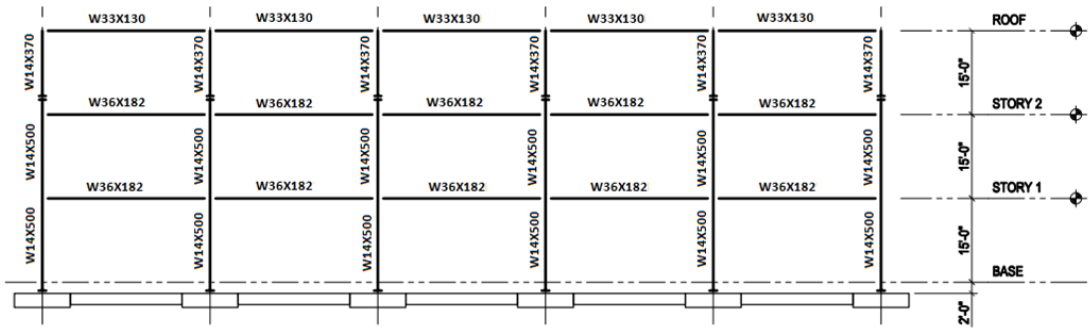

Figure 2.4    HP-SMRF lateral force resisting systems configuration [NEES 2008].

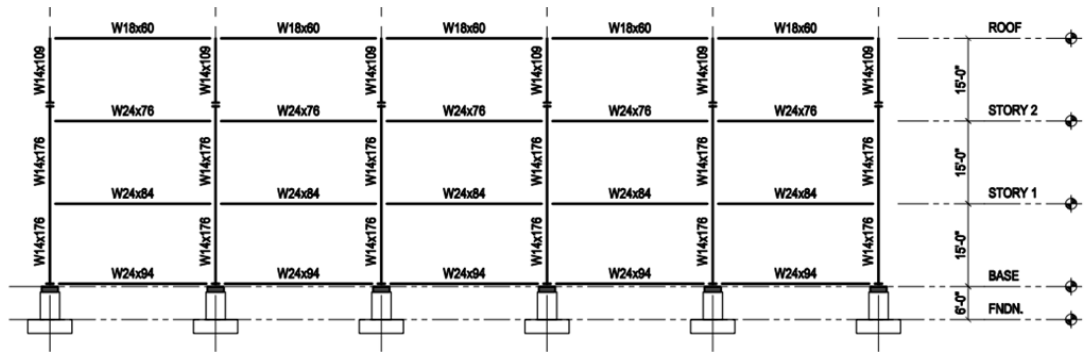

Figure 2.5    BI-IMRF lateral force resisting systems configurations [NEES 2008].

For the fixed base structure, special moment resisting frames (SMRF) were used and designed with higher performance objectives than the minimum code requirements. For the base-isolated configuration, intermediate moment resisting frames (IMRF) were used. The designs of the two considered systems, fixed-base and base-isolated moment resisting frames, are consistent with what might be used by many engineers and are compliant with the code standards for design according to the Equivalent Lateral Force Method [ASCE 2010].

The HP-SMRF was designed with a force reduction factor ($R/I_e$) of 6.4 (8/1.25), an interstory drift limit of 1.0% (more stringent than 2% required by code [ASCE 2010]; the reduced beam section (RBS) was used to ensure a ductile behavior. These connections are the only prequalified welded connections permitted for SMRF by AISC 341-05 [2005]. This resulted in fundamental period of the fixed-base system of 0.67 sec. Compared to the HP-SMRF, the BI-IMRF (Figure 2.5) was designed utilizing lower $R/I_e$ factor [1.69=(3/8) × (4.5/1)] and the same drift limit (1.0%). The IMRF uses welded unreinforced flange–welded web (WUF-W) connections. Because these connections do not contain flanges in the plastic hinge zone, they have higher rotational stiffness compared to the RBS. This higher stiffness is appropriate for isolated buildings where period separation is desirable. Additionally, the ductility demand for isolated buildings is expected to be low, eliminating the need for a highly-ductile connection such as the RBS. The fundamental period for the isolated system is 1.404 sec.

The isolation system was designed to have a maximum displacement of 30 in. under the maximum capable earthquake (MCE) event. It utilized triple friction pendulum bearings (TFPB) with the friction coefficients of the four sliding surfaces of 0.01, 0.01, 0.03, and 0.06, and effective pendulum lengths of 20, 122, and 122 in. Under an MCE event, this bearing had an effective period of 4.35 sec and effective damping of 15.1% (Table 2.1). More details on design of these two systems can be found in Mayencourt [2013] and Terzic et al. [2014a]. The hysteretic behavior of the bearings is represented using tri-linear models and for the considered isolator; see Figures 2.6 and 2.7.

Table 2.1    Isolation system parameters.

|  | DBE | MCE |
|---|---|---|
| Effective period | 3.95 sec | 4.35 sec |
| Effective damping | 22.9 % | 15.1 % |
| Isolator displacement | 16.1 in. | 30 in. |

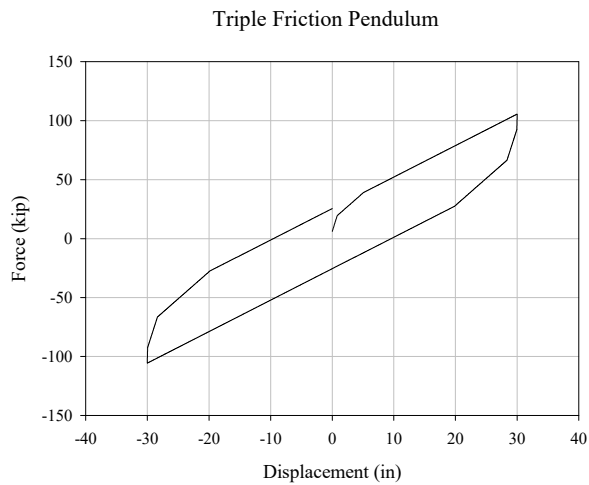

Figure 2.6    Hysteretic behavior (force versus horizontal displacement) for triple-friction-pendulum bearing.

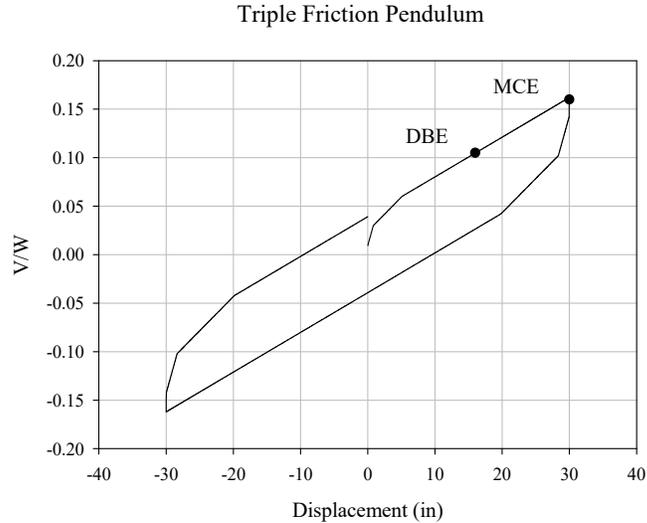

**Figure 2.7** Hysteretic behavior (normalized force versus horizontal displacement) for triple-friction-pendulum bearing.

Table 2.2 compares the member sizes between the code-minimum performance design and the considered HP-SMRF and BI-IMRF structures. As can be seen in Table 2.1, the high-performance system has much bigger frame sections than code-minimum members. This led to a 40% increase in the weight of the steel, which justifies the use of RBS to ensure a controlled ductile behavior.

The RBS connections (Figure 2.8) were used to protect the SMRF configurations by forcing the plastic hinge in the beam to form away from the column face (between $a$ and $a+b/2$ from the column face). To ensure a strong column–weak beam behavior, the RBS connections reduced the moment capacity of the beam locally, resulting in large strains occurring in the imposed location. The part of the beam that is between the RBS and the column face is called the panel zone, which must remain elastic. That zone should be free from damage so that damaged beams will be easier to replace after an earthquake. The floor weights of the building are shown in Table 2.3. The fixed-base structure doesn't have a base floor; instead, the base level of the isolated structure has been assumed to have the same weight of the first floor.

The construction material was ASTM A992 steel with 0.4901 kips/ft$^3$ weight per unit volume, 50 ksi (345 MPa) tensile yield strength, and 65 ksi (450 MPa) tensile ultimate strength. Material ductility of the considered material is defined by the maximum yield-to-tensile strength ratio equal to 0.85. This material was chosen because the AISC states that this is the preferred steel for wide-flange shapes.

Code spectral accelerations were determined from PSHA by the USGS for the given location. Healthcare facilities are assigned as Risk Category III, corresponding to structures whose failure could pose a substantial risk to human life [ASCE 2010]. The spectral accelerations for the MCE spectrum were selected to be $S_s = 2.2g$ for short periods and $S_1 = 0.74g$ at a period of 1 sec, which are representative of many locations in California. The long-period transition period $T_L$ was set equal to 8 sec for the Oakland site. Site coefficients $F_a$ and $F_v$ were chosen knowing that the site class of the building is $E$, and their values are 1.0 and 1.5, respectively. The response modification coefficient $R$ was set equal to 8 for steel SMRFs and equal to 4.5 for intermediate moment-resisting frames; the seismic importance factor $I_e$ for Risk Category III is 1.25. The design basis earthquake (DBE) spectrum was found by dividing per 1.5 the MCE spectrum. The code spectra are shown in Figure 2.9.

Table 2.2  Beams and columns sizes.

|  |  | Code minimum | HP-SMRF | BI-IMRF |
|---|---|---|---|---|
| Beams | Third story | W27x102 | W33x130 | W18x60 |
|  | Second story | W33x130 | W36x182 | W24x76 |
|  | First story | W33x141 | W36x182 | W24x84 |
|  | Base | - | - | W24x94 |
| Columns | Third story | W14x211 | W14x370 | W14x109 |
|  | Second story | W14x370 | W14x500 | W14x176 |
|  | First story | W14x370 | W14x500 | W14x176 |

Table 2.3  Floor weights.

| Story | Weight (kips) |
|---|---|
| Penthouse | 219 |
| Third story | 1773 |
| Second story | 1924 |
| First story | 1924 |
| Base | 1924 |

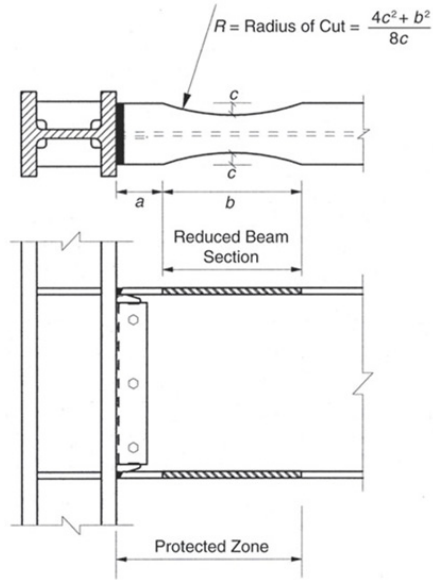

Figure 2.8  Reduced beam section [Chou and Uang 2003].

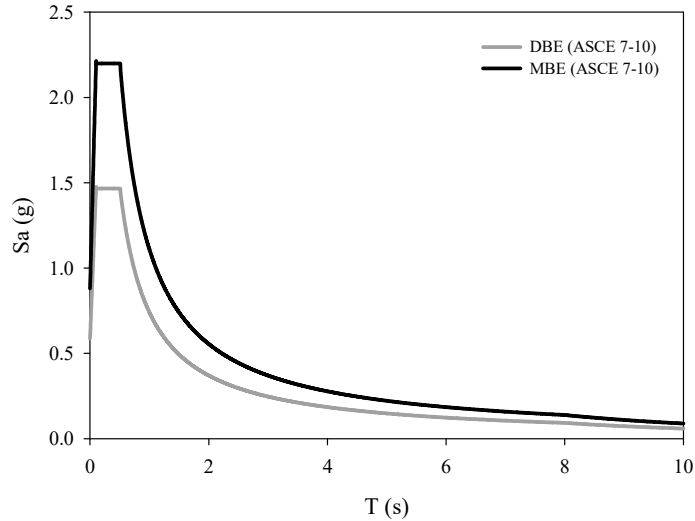

**Figure 2.9**   Code DBE and MCE spectra.

## 2.4   GROUND-MOTION SELECTION

The set of ground motions used in the analysis was selected to be representative of the hazard at the site of the I-880 viaduct in Oakland, California. For the hazard analysis, a location of 37.803N × 122.287W was used. The site is located between the Hayward and San Andreas faults as shown in Figure 2.10.

The site-specific ground motions were selected to match the uniform hazard spectrum [UHS] [USGS 2013]. The UHS for the different hazard levels were derived using the 2008 USGS hazard maps and interactive deaggregation tools [USGS 2013], and are shown in Figure 2.11. In addition, the mean magnitude, distance, and $\varepsilon$ values associated with occurrence of each spectral value were determined.

The deaggregation plots for $S_a$ values exceeded with 2% probability in 50 years at periods of 0.1 and 1 sec as shown in Figure 2.12 and Figure 2.13 to better understand the relative contributions of the various sources to the total seismic hazard. A comparison of the two plots shows that the contribution of the Hayward fault is higher for the lower periods, while the contribution of San Andreas fault increases as the period increases. Table 2.8 shows $\varepsilon$ values close to zero, implying that the spectral acceleration value is approximately equal to the median spectral acceleration value. As the probability of exceedance decreases, it shows positive $\varepsilon$ values, which are representative of spectral accelerations higher than the median ground motions. The median magnitudes increase with lower probabilities of exceedance, representing the extreme events.

Forty three-component ground-motion records were selected to represent the ground-motion hazard at each of three hazard levels: 2%, 10%, and 50% probabilities of exceedance in 50 years. More information on these motions can be found in Baker et al. [2011]. To better characterize the seismic hazard at the site, two additional sets of ground-motion records representative of hazard levels at 5% and 20% probabilities of exceedance were also used in the analysis. Those additional sets were selected from the PEER Ground-Motion Database by inputting the UHS at the different hazard levels.

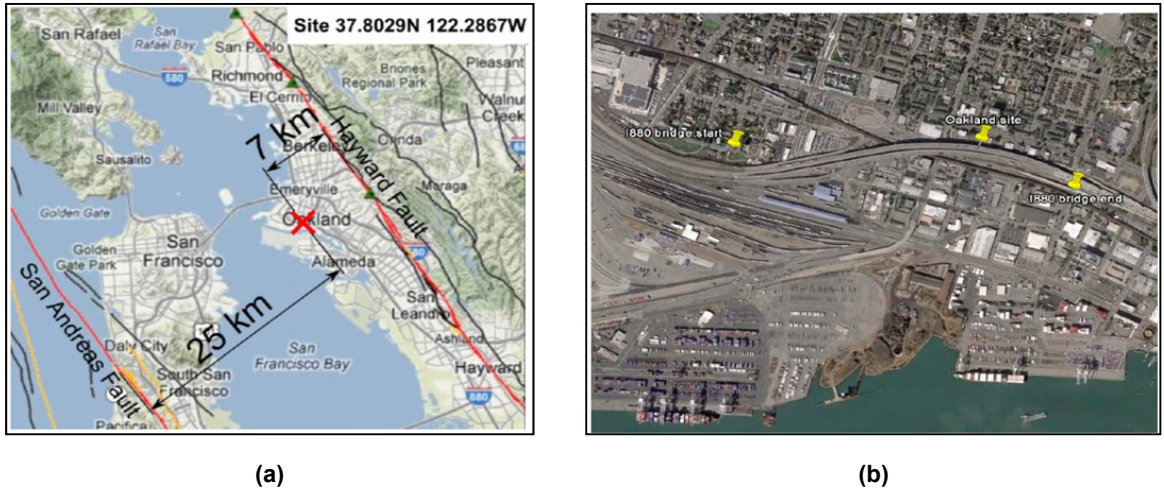

(a)                                      (b)

**Figure 2.10**     (a) Position of Hayward and San Andreas fault respect to the site of interest; and (b) location of I-880 bridge viaduct and the Oakland site. Aerial imagery from Google Earth (http://earth.google.com).

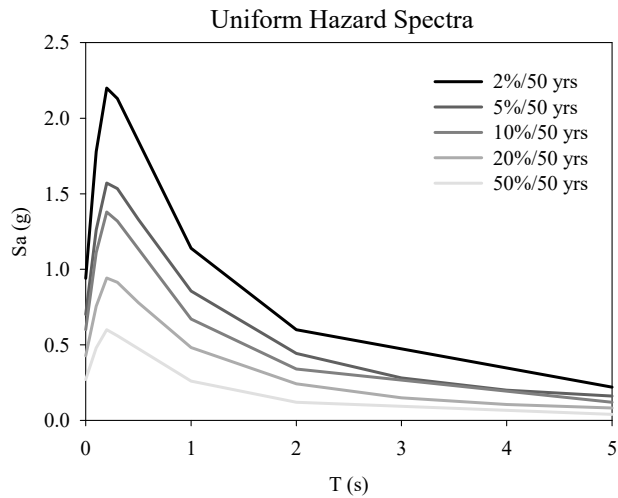

**Figure 2.11**     Uniform Hazard Spectra for the Oakland site.

Table 2.4    Uniform hazard spectrum and mean deaggregation values of
             distance, magnitude, and $\varepsilon$ for the Oakland site, with a 2%
             probability of exceedance in 50 years.

| | | 2% in fifty years | | |
|---|---|---|---|---|
| $T$ (sec) | Sa(g) | $R$ (km) | M | $\varepsilon$ |
| 0.0 | 0.94 | 8.8 | 6.78 | 1.70 |
| 0.1 | 1.78 | 8.4 | 6.73 | 1.76 |
| 0.2 | 2.20 | 8.4 | 6.77 | 1.74 |
| 0.3 | 2.13 | 8.5 | 6.81 | 1.73 |
| 1.0 | 1.14 | 9.9 | 7.00 | 1.74 |
| 2.0 | 0.60 | 13.4 | 7.20 | 1.74 |
| 5.0 | 0.22 | 16.0 | 7.43 | 1.64 |

Table 2.5    Uniform hazard spectrum and mean deaggregation values of
             distance, magnitude, and $\varepsilon$ for the Oakland site, with a 5%
             probability of exceedance in 50 years.

| | | 5% in fifty years | | |
|---|---|---|---|---|
| $T$ (sec) | Sa(g) | $R$ (km) | M | $\varepsilon$ |
| 0.0 | 0.70 | 9.9 | 6.80 | 1.27 |
| 0.1 | 1.26 | 9.8 | 6.73 | 1.31 |
| 0.2 | 1.57 | 9.8 | 6.77 | 1.30 |
| 0.3 | 1.53 | 9.8 | 6.82 | 1.31 |
| 1.0 | 0.86 | 11.1 | 7.00 | 1.40 |
| 2.0 | 0.44 | 15.3 | 7.18 | 1.40 |
| 5.0 | 0.16 | 16.5 | 7.36 | 1.31 |

Table 2.6    Uniform hazard spectrum and mean deaggregation values of
             distance, magnitude, and $\varepsilon$ for the Oakland site, with a 10%
             probability of exceedance in 50 years.

| | | 10% in Fifty Years | | |
|---|---|---|---|---|
| $T$ (sec) | Sa(g) | $R$ (km) | M | $\varepsilon$ |
| 0.0 | 0.60 | 10.1 | 6.80 | 1.05 |
| 0.1 | 1.11 | 10.0 | 6.75 | 1.10 |
| 0.2 | 0.38 | 10.0 | 6.78 | 1.10 |
| 0.3 | 1.32 | 10.2 | 6.82 | 1.09 |
| 1.0 | 0.67 | 11.8 | 7.00 | 1.09 |
| 2.0 | 0.34 | 15.6 | 7.15 | 1.09 |
| 5.0 | 0.12 | 16.9 | 7.313 | 1.01 |

**Table 2.7** Uniform hazard spectrum and mean deaggregation values of distance, magnitude, and $\varepsilon$ for the Oakland site, with a 20% probability of exceedance in 50 years.

| 20% in Fifty Years | | | | |
|---|---|---|---|---|
| $T$ (sec) | Sa(g) | $R$ (km) | M | $\varepsilon$ |
| 0.0 | 0.43 | 12.0 | 6.80 | 0.64 |
| 0.1 | 0.76 | 12.1 | 6.73 | 0.70 |
| 0.2 | 0.94 | 12.3 | 6.76 | 0.70 |
| 0.3 | 0.91 | 12.3 | 6.81 | 0.70 |
| 1.0 | 0.48 | 14.1 | 6.99 | 0.72 |
| 2.0 | 0.24 | 18.9 | 7.13 | 0.72 |
| 5.0 | 0.08 | 18.3 | 7.25 | 0.64 |

**Table 2.8** Uniform hazard spectrum and mean deaggregation values of distance, magnitude and $\varepsilon$ for the Oakland site, with a 50% probability of exceedance in 50 years.

| 50% in Fifty Years | | | | |
|---|---|---|---|---|
| $T$ (sec) | Sa(g) | $R$ (km) | M | $\varepsilon$ |
| 0.0 | 0.27 | 15.1 | 6.79 | 0.00 |
| 0.1 | 0.48 | 15.0 | 6.73 | 0.10 |
| 0.2 | 0.60 | 15..7 | 6.76 | 0.11 |
| 0.3 | 0.56 | 16.2 | 6.80 | 0.10 |
| 1.0 | 0.26 | 19.3 | 6.96 | 0.04 |
| 2.0 | 0.12 | 24.2 | 7.06 | 0.02 |
| 5.0 | 0.04 | 24.2 | 7.13 | -0.02 |

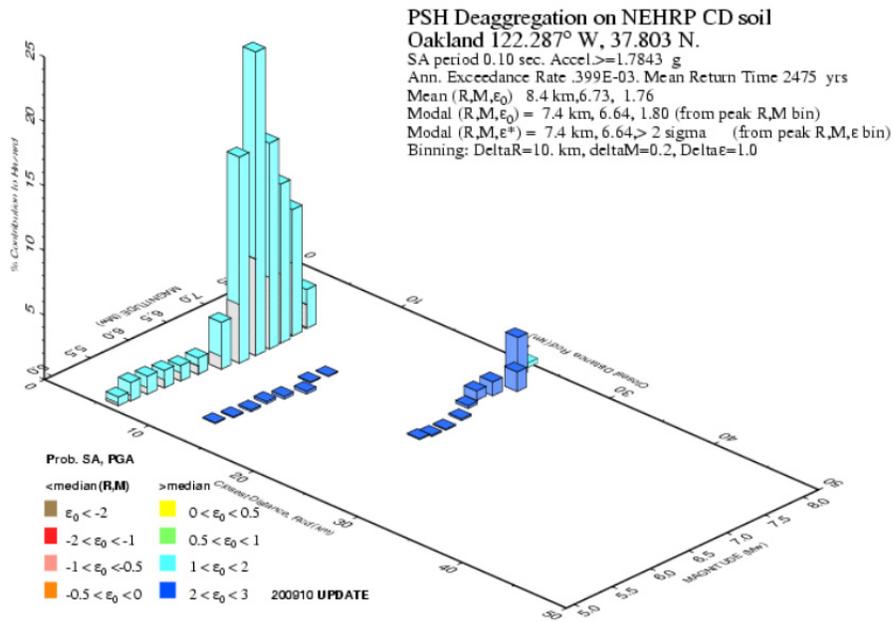

**Figure 2.12**   Deaggregation plot for *Sa* (0.1 sec) exceeded, with 2% probability in 50 years [USGS 2008].

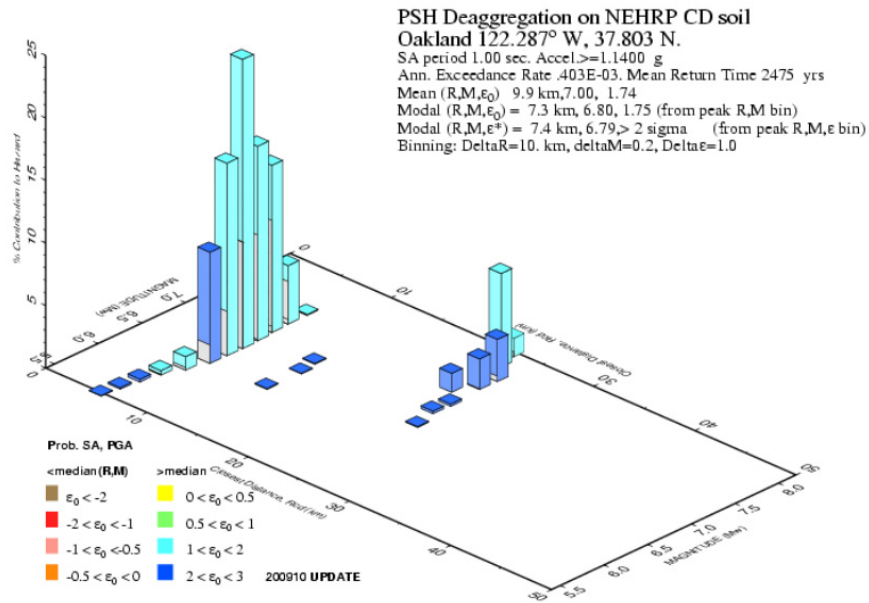

**Figure 2.13**   Deaggregation plot for $S_a$ (1 sec) exceeded, with 2% probability in 50 years [USGS 2008].

Table 2.9    Ground motions selected for the 20%/50 years hazard level.

| Record | NGA # | Event | Year | Station | SF | Mag | $V_{s30}$(m/sec) |
|---|---|---|---|---|---|---|---|
| 1 | 95 | Managua- Nicaragua-01 | 1972 | Managua- ESSO | 2.120 | 6.24 | 288.8 |
| 2 | 175 | Imperial Valley-06 | 1979 | El Centro Array #12 | 2.350 | 6.53 | 196.9 |
| 3 | 179 | Imperial Valley-06 | 1979 | El Centro Array #4 | 0.891 | 6.53 | 208.9 |
| 4 | 182 | Imperial Valley-06 | 1979 | El Centro Array #7 | 0.698 | 6.53 | 210.5 |
| 5 | 184 | Imperial Valley-06 | 1979 | El Centro Differential Array | 0.968 | 6.53 | 202.3 |
| 6 | 231 | Mammoth Lakes-01 | 1980 | Long Valley Dam (Upr L Abut) | 4.118 | 6.06 | 345.4 |
| 7 | 266 | Victoria- Mexico | 1980 | Chihuahua | 2.020 | 6.33 | 274.5 |
| 8 | 316 | Westmorland | 1981 | Parachute Test Site | 1.265 | 5.9 | 348.7 |
| 9 | 549 | Chalfant Valley-02 | 1986 | Bishop - LADWP South St | 2.657 | 6.19 | 271.4 |
| 10 | 587 | New Zealand-02 | 1987 | Matahina Dam | 2.914 | 6.6 | 424.8 |
| 11 | 718 | Superstition Hills-01 | 1987 | Wildlife Liquef. Array | 3.408 | 6.22 | 207.5 |
| 12 | 721 | Superstition Hills-02 | 1987 | El Centro Imp. Co. Cent | 1.293 | 6.54 | 192.1 |
| 13 | 728 | Superstition Hills-02 | 1987 | Westmorland Fire Sta | 1.449 | 6.54 | 193.7 |
| 14 | 767 | Loma Prieta | 1989 | Gilroy Array #3 | 1.434 | 6.93 | 349.9 |
| 15 | 779 | Loma Prieta | 1989 | LGPC | 0.514 | 6.93 | 477.7 |
| 16 | 802 | Loma Prieta | 1989 | Saratoga - Aloha Ave | 1.105 | 6.93 | 370.8 |
| 17 | 825 | Cape Mendocino | 1992 | Cape Mendocino | 0.984 | 7.01 | 513.7 |
| 18 | 827 | Cape Mendocino | 1992 | Fortuna - Fortuna Blvd | 1.846 | 7.01 | 457.1 |
| 19 | 983 | Northridge-01 | 1994 | Jensen Filter Plant Generator | 0.670 | 6.69 | 525.8 |
| 20 | 1042 | Northridge-01 | 1994 | N Hollywood - Coldwater Cyn | 1.660 | 6.69 | 446 |
| 21 | 1045 | Northridge-01 | 1994 | Newhall - W Pico Canyon Rd. | 0.587 | 6.69 | 285.9 |
| 22 | 1605 | Duzce- Turkey | 1999 | Duzce | 0.688 | 7.14 | 276 |
| 23 | 1611 | Duzce- Turkey | 1999 | Lamont 1058 | 2.968 | 7.14 | 424.8 |
| 24 | 2655 | Chi-Chi- Taiwan-03 | 1999 | TCU122 | 2.421 | 6.2 | 475.5 |
| 25 | 2699 | Chi-Chi- Taiwan-04 | 1999 | CHY024 | 4.342 | 6.2 | 427.7 |

The following criteria were used to better represent the hazard at the site for both 5% and 20% probability of exceedance in 50 years:

- magnitude values between 5.9 and 7.3 to be comprehensive of the range of probable magnitudes
- closest distance to the fault rupture between 0 and 20 km
- scale factor of all three components of ground motion between 0.25 and 4, which is the recommended limit on scaling
- shear-wave velocity at 30 m less than 550 m/sec
- periods between 0 and 5 sec

In addition to these criteria, ground motions were selected limiting to five the records coming from the same event in order to avoid the influence of a particular earthquake with many records (e.g., Northridge) on the set of ground motions. Each of the two additional sets of ground motions had 25 three-component ground motion records that were scaled by multiplying the spectral accelerations for the scale factor that provided the best match between the record and the UHS. Figure 2.15 compares UHS with the median pseudo-acceleration response spectra for the selected ground motions at all five considered hazard levels. There is a good agreement between UHS and median pseudo-acceleration response spectra for the range of periods of the considered structural systems.

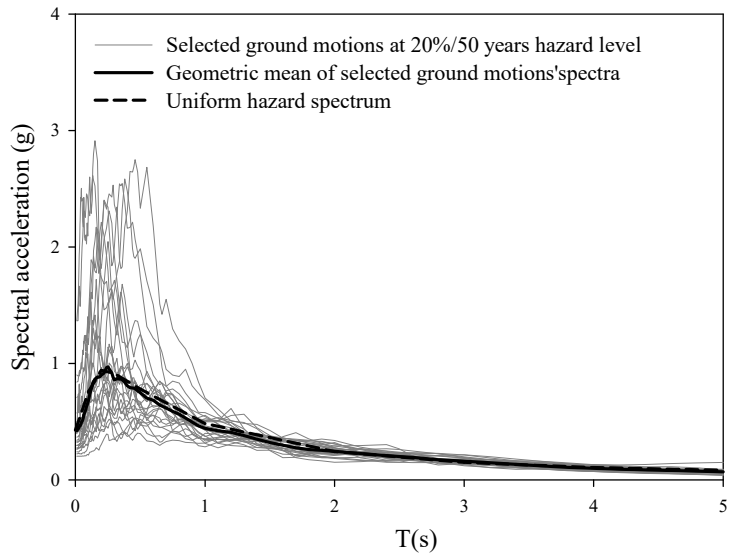

**Figure 2.14** Uniform hazard spectrum at the 20% in 50 years hazard level and response spectra of the selected ground motions.

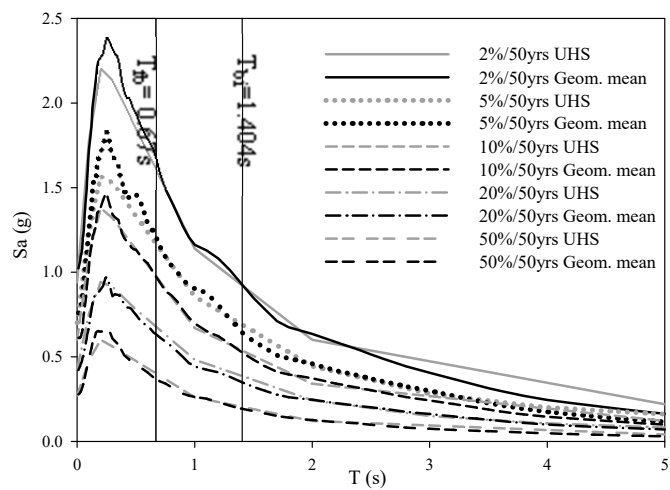

**Figure 2.15** Comparison of UHS with the median pseudo-acceleration response spectra for the selected ground motions at five considered hazard levels.

## 2.5  ANALYSIS MODEL AND METHOD USED

A nonlinear time-history analysis was performed for both fixed-base and base-isolated systems. To avoid complications related to a full three-dimensional (3D) model, the analysis was simplified by modeling an appropriate two-dimensional (2D) frame using OpenSees [McKenna et al. 2004]. This simplification is valid as the lateral load-resisting frames are located only on the perimeter of the building and do not have common elements. Gravity-load-only type connections were used elsewhere in the structure.

### 2.5.1  Modeling Assumptions

To achieve an acceptable level of accuracy using the simplified 2D analysis, numerical modeling assumptions were made for the two structural systems. Further details can be found in Terzic at el. [2014a]. The five main assumptions are listed below:

- Floor slabs were assumed to be axially inextensible. All nodes on the same floor were constrained to the same displacement.

- The 2D frame considered in the analysis carried half of the horizontal floor mass, equally distributed among the six nodes of the same floor. The vertical masses were assigned to the same nodes by using the principle of tributary area. The external nodes carried half of the interior nodes vertical mass.

- The frames were subjected to horizontal and vertical components of ground motion.

- The effects of large deformations of beam and column elements were accounted for using $P$-$\Delta$ nonlinear geometric transformation. The large deformations of the columns increased the distance between the initial centerline of elements and the point of application of the load. That implies an increase of the bending moment due to second-order effects. A negative stiffness was then added to the global system.

- Vertical gravity loads acting not only on the moment frame but on the entire structure, were included in the analysis since they act on the displaced location of the joints, giving birth to destabilizing $P$-$\Delta$ effects. $P$-$\Delta$ effects from the gravity columns were taken into account by using a single leaning column with applied gravity loads that relied on the seismic force-resisting system for lateral stability [Gupta and Kunnath 2000]. It was composed of elastic elements and connected to the 2D frame using a rigid truss element with an infinite axial stiffness, which was obtained by increasing the section up to an area of 10,000 in.². The column was composed of vertical elastic elements whose section was equal to half of the sum of the area of all gravity columns and of the moment-resisting frames in the opposite direction of the building. The connection at each floor intersection was modeled by using concentrated plastic hinges. The capacity of each hinge was equal to the flexural capacity of half of the sum of the capacity of all gravity columns and half of the columns of the frame in the opposite direction. The load applied on the column was equal to the half of the gravity load per floor reduced by the gravity load acting on the columns of the lateral load-resisting frame.

Table 2.10   Mass input at each node.

|  | Vertical mass (kip·sec²/in.) | Horizontal mass (kip·sec²/in.) |
|---|---|---|
| Second floor | 0.105 | 0.419 |
| Third floor | 0.105 | 0.419 |
| Roof | 0.096 | 0.430 |

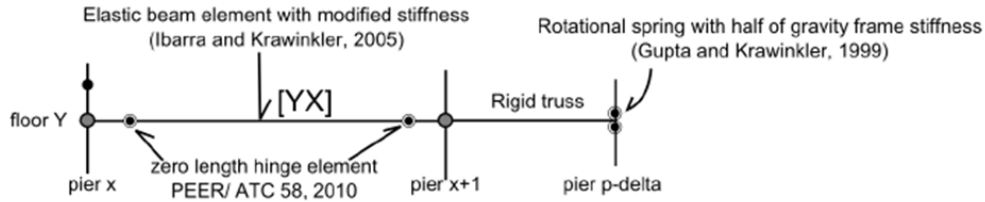

Figure 2.16   Beam assembly for a regular bay and *P*-delta bay [Mayencourt 2013].

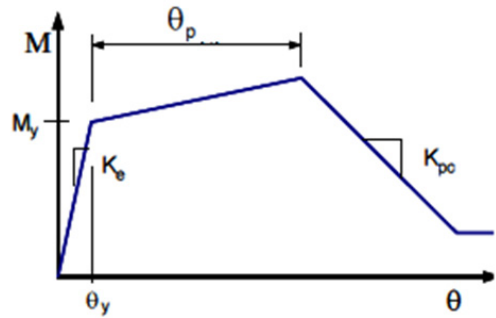

Figure 2.17   BS hinges backbone curve [PEER/ATC 2010]. $M_y$: effective yield strength, $K_e$: initial stiffness, $θ_y$: yield rotation, $θ_y + θ_p$: capping strength rotation, $K_{pc}$: post-capping stiffness.

To take into account the possibility of a nonlinear behavior of the system at high hazard levels, the beams of the HP-SMRF system were modeled as lumped plasticity elements. They were characterized by elastic elements with concentrated plastic hinges at the ends to recreate the lumped plasticity of the RBS. Plastic hinges were reproduced by rotational springs with deterioration rules based on experimental data.

Rotational springs were modeled in Open Sees by means of nonlinear zero-length elements placed at the center of the RBS. To maintain the same stiffness of the initial frame, the inertia of the elastic portion of the beam was increased by multiplying the initial inertia by 1.1 [Ibarra and Krawinkler 2005]. The hinge–elastic-beam–hinge assembly was used to model the lumped plasticity of the reduced beam section. The moment-rotation relationships for RBS connections were developed based on recommendations from PEER/ATC [2010] using *Hysteretic Material* of OpenSees. The initial backbone curve (Figure 2.17) defined the boundaries of the hysteretic response, which degraded (deterioration in strength and stiffness) as a function of damage and energy dissipated during the cyclic model.

Input values for the RBS are listed in Table 2.11. The effective yield strength $M_y$ defines the first point in the initial backbone curve. For beams with RBS connections, a mean value can be estimated as 1.06 $M_p$ [Lignos 2008]. For the same beams capping strength, $M_c$ can be estimated as the mean value of the ratio of capping strength to effective yield strength, $M_c/M_y$ = 1.09; Lignos suggests a value of 0.4 × $M_y$ to estimate the residual strength $M_r$. The panel zone at the beam and column intersection was modeled as a rectangle composed of four very stiff elastic beam elements to ensure that the plastic hinges form away from the columns.

Columns of the HP-SMRF and beams and columns of the BI-IMRF were both modeled using force-based beam–column elements of OpenSees. Thus, the spread of the plasticity along the element is taken into account. Sections were represented by fibers, and a steel material with a fatigue behavior was considered using the Menegotto-Pinto hysteretic model. Low-cycle fatigue failure of beams and columns was accounted for by using an OpenSees fatigue model. A global representation of a regular bay and leaning column is shown in Figure 2.18.

The isolated system was designed with an isolator beneath each column of the structural frame. The isolators were modeled with horizontal springs and tri-linear uniaxial material representative of the hysteretic behavior of triple-pendulum friction bearings. The capacity of each bearing of the considered frame is equal to half of the sum of the properties of all the isolators of the building divided by the six columns of the frame. Vertical displacements and rotations at the top and the bottom of isolators were assumed to be fixed. The leaning column was modeled with a roller at the base and was connected to the structural frame with a rigid truss element as at the other floors.

Following the recommendations of PEER/ATC [2010], damping was assigned to the frames in a way that represents the energy dissipation during the analysis. Since the analysis was performed by using nonlinear elements, most of the dissipation was already captured by the hysteretic response. In any case, in this study a damping ratio of 3% was taken for both structural systems. For the fixed-base building, mass and tangent stiffness proportional to Rayleigh coefficients were calculated based on two periods. For the 50% in 50-years hazard level, the first ($T_1$) and third ($T_3$) periods were selected. For the fixed-base building, mass and tangent stiffness proportional to Rayleigh coefficients were calculated based on two periods. The first ($T_1$) and third ($T_3$) periods were selected for the 50%/50-years and 20%/50-years hazard levels, while 1.5$T_1$ and $T_3$ were selected for the 10%/50-years, 5%/50years and 2%/50-years hazards. The first period was elongated 1.5 times to account for the change in period due to the nonlinear deformations of the system.

For the isolated building, the damping was assumed to be proportional only to the tangent stiffness of the structure. The stiffness-proportional damping was calculated from the fundamental period of the structure $T_1$ for the 50%/50-years hazard, and $T_{eff}$ for the 20%/50-years, 10%/50-years, 5%/50-years, and 2%/50-year hazards, where $T_{eff}$ is the effective fundamental period of the isolation system.

**Table 2.11    RBS hinges characteristics.**

| Section size | Effective yield strength $M_y$ (kip in.) | Yield rotation (rad) | Capping strength $M_c$ (kip in.) | Capping strength rotation (rad) | Residual strength $M_r$ (kip in.) | Residual strength (rad) |
|---|---|---|---|---|---|---|
| W33 × 130 | 19,695.49 | 0.000488 | 21,665.04 | 0.021709 | 7,878.196 | 0.120495 |
| W36 × 182 | 30,997.81 | 0.000429 | 34,097.59 | 0.021291 | 12,399.124 | 0.156850 |
| W36 × 182 | 30,997.81 | 0.000429 | 34,097.59 | 0.021291 | 12,399.124 | 0.156850 |

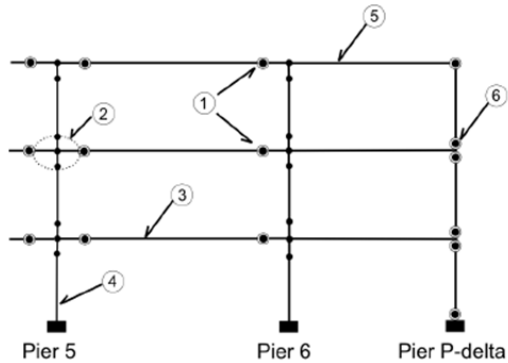

**Figure 2.18    Structural scheme of a generic bay and the *P*-Δ leaning column. (1) RBS connections, (2) panel zone, (3) elastic beam-column elements, (4) columns, (5) rigid truss, (6) *P*-Δ floor connections [Mayencourt 2013].**

### 2.5.2  Comparison of Structural Response

While numerous parameters need to be considered to fully evaluate structural response, it is common practice to correlate performances to EDPs as story drifts, floor accelerations, and residual drifts. By comparing the EDP average peak values for the five considered hazard levels, the performance characteristics of the system can be assessed. The severity of damage to various structural and nonstructural components associated with EDPs can be quantified using fragility relations from FEMA P-58 [FEMA 2012]. Losses associated with damage will be evaluated in the next section.

The base-isolated moment frame substantially reduced accelerations and drifts compared to the fixed-base frame; see Figure 2.19 and Figure 2.20. While the effectiveness of the isolation system in reducing story drift increased by increasing the ground shaking intensity (ranging from 20% to 62% with an average of 49%), the reduction of acceleration was consistently high at all hazard levels (ranging from 84% to 90% with an average of 88%). The BI-IMRF, with the uniform acceleration profile over the height of the building and the peak median value, which reached 0.22$g$ at the 2% in 50-years hazard level, most likely will not trigger any damage of the acceleration sensitive components (e.g., ceiling, MEP, contents).

Figure 2.21 shows that isolation system was also effective in eliminating residual drifts of the moment frame at the highest hazard levels. At the 50%/50-years hazard level, the HP-SMRF developed maximum median drift of 0.46%, 20% larger than maximum median drift of the BI-IMRF of 0.37%; see Figure 2.19. Since both moment frames are expected to yield at drift ratios slightly larger than 1%, elastic structural behavior was anticipated at this hazard level. The damage to interior partitions was expected for both the HP-SMRF and BI-IMRF system since the median drift associated with initiation of damage to partition walls commonly used in healthcare facilities is 0.21% [FEMA 2012]. Median horizontal accelerations in the HP-SMRF ranged from 0.26$g$ to 0.67$g$ over the height of the building—see Figure 2.20a—likely triggering damage to piping, electronic, and medical equipment in the upper levels [FEMA 2012].

At the 20%/50-years hazard level, greater differences in story-drift demands were observed between the two systems; see Figure 2.19. Compared to the BI-IMRF, the fixed-base HP-SMRF had about two times larger drift ratio at each level, with the peak median value reaching 0.84%. This would likely result in a greater damage to partition walls and initiation of damage to staircases, which initiates at a drift of 0.5%, per FEMA [2012]. At this hazard level, damage to structural elements was not anticipated. Median horizontal accelerations in the HP-SMRF ranged from 0.37$g$ to 1.13$g$ over the height of the building; see Figure 2.20(b). These accelerations extended the regions of the building that typically undergo acceleration-related damage, triggering additional damage to ceilings, chillers, fire sprinkler drops, bookcases, and filing cabinets [FEMA 2012]. At the 10%/50-years hazard level. Figure 2.19(c) shows even greater differences in story-drift demands between the two systems.

The fixed-base HP-SMRF had a peak median drift ratio of 1.24%, which suggests initiation of yielding of the system and probable extensive damage to wall partitions and moderate damage to staircases. The BI-IMRF, with a peak median drift ratio of 0.57%, is anticipated to remain elastic with slight damage to wall partitions and staircases. Median horizontal accelerations in the HP-SMRF ranged from 0.59$g$ to 1.54$g$ over the height of the building; see Figure 2.20(c). These accelerations extended the regions of the building that typically undergo acceleration-related damage observed at lower hazard levels, triggering additional damage to lightening, cooling tower, HVAC ducts, and air handling units.

At the fixed-base HP-SMRF had a median peak story drift of 1.57% (5%/50-years hazard level) and 2.24% (2%/50-years hazard level)—see Figure 2.19(d) and (e) —suggesting damage to both structural and nonstructural components thus requiring substantial repair. The BI-IMRF, with the peak median drift ratios of 0.62% (5%/50-years hazard level) and 0.83% (2%/50-years hazard level) is anticipated to remain elastic with slight non-structural damage. Median horizontal accelerations in the HP-SMRF ranged from 0.78$g$ to 1.61$g$ for the 5%/50-year hazard level and from 0.97$g$ to 1.81$g$ for the 2%/50-year hazard level—see Figure 2.20(d) and (e)—causing damage to all acceleration sensitive non-structural components and content except for electrical systems and components.

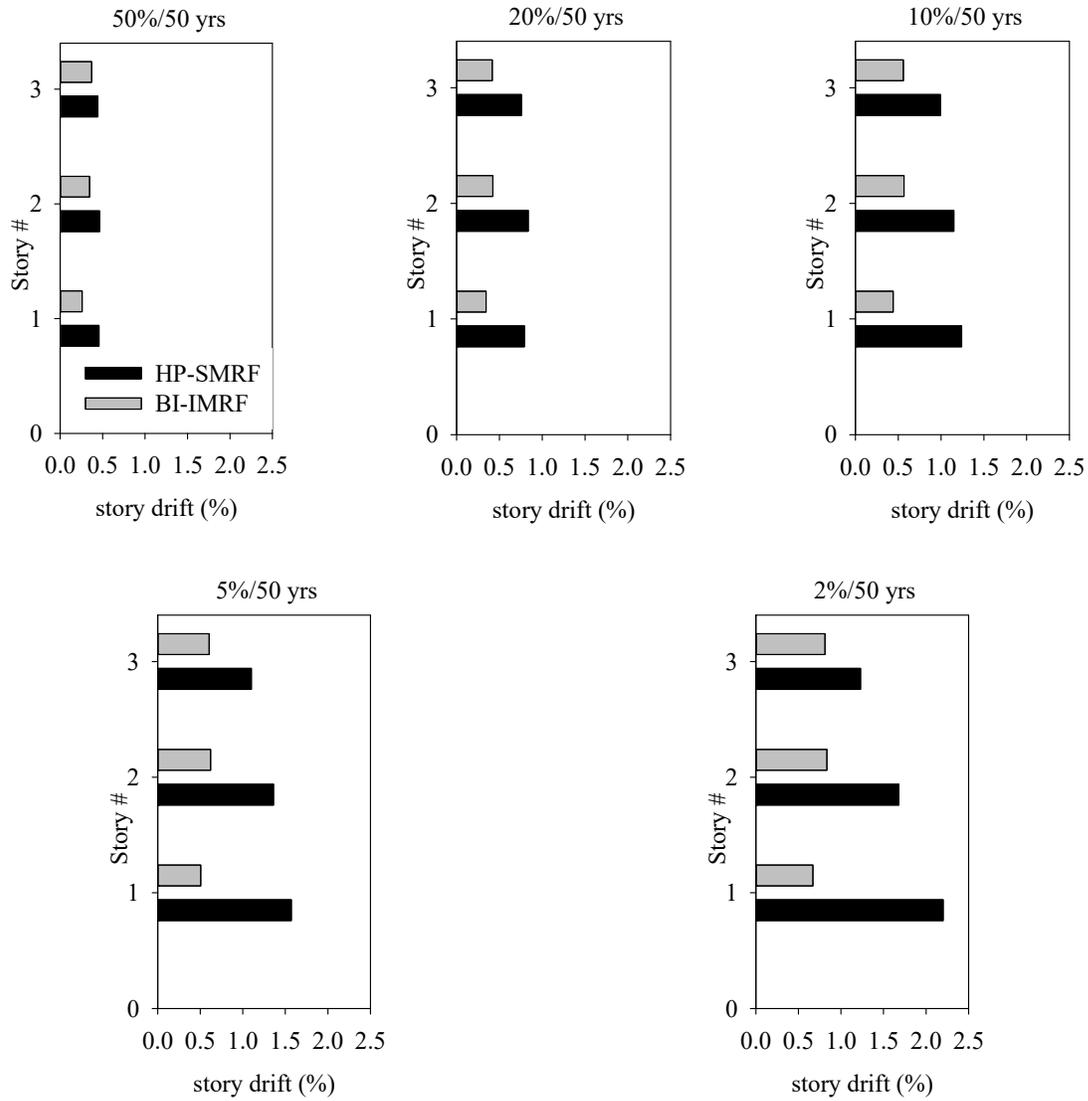

**Figure 2.19** Median story drifts of the HP-SMRF and the BI-IMRF on TFPBs for five hazard levels: 2%, 5%, 10%, 20%, and 50% probabilities of exceedance in 50 years.

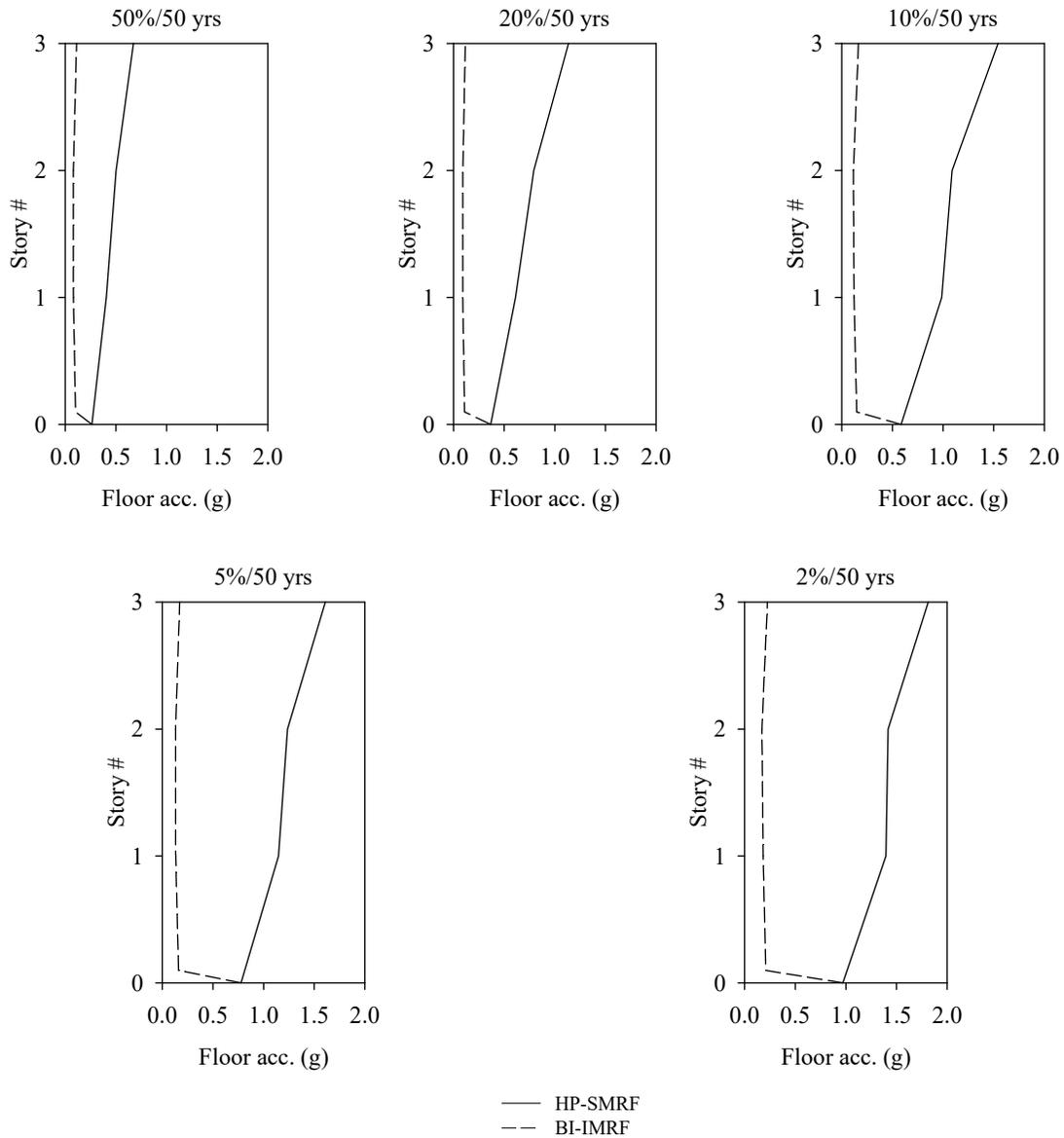

**Figure 2.20    Median absolute floor accelerations of the HP-SMRF and the BI-IMRF on TFPBs for five hazard levels: 2%, 5%, 10%, 20%, and 50% probabilities of exceedance in 50 years.**

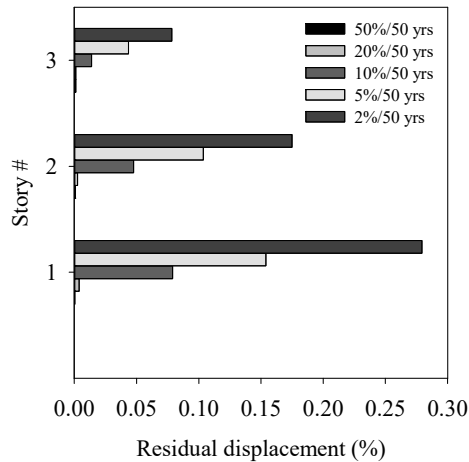

**Figure 2.21** Residual story drifts of the HP-SMRF at five hazard levels. Superstructure of the BI-IMRF had no residual drifts.

### 2.5.3 Hysteretic Cycle Response

The behaviour of the HP-SMRF from serviceability limit states to near-collapse limit states is shown in terms of the hysteretic cycle response of the system at all considered hazard levels. Each graph shows the base shear and the first floor displacement of the ground-motion record that better represents the median structural response of the building at each hazard level. As obvious, this specific ground motion cannot summarize the effect of a set of ground motions, but it can give an idea of the median behaviour of the system at an increasing the level of hazard. An analysis of the hysteretic cycle response allows for insight into the strength and stiffness degradation of the structure under random loading histories as well as the modes (slow or rapid) in which this deterioration can occur.

At 50%/50 years and 20%/50 years hazard levels, the behaviour of the system is kept in the elastic range, reaching maximum first-floor displacements of 0.89 and 1.63 in., respectively (Figure 2.22 and Figure 2.23), which correspond to story drift ratios of 0.43% and 0.8%. As previously explained, both moment frames were expected to yield at drift ratios slightly larger than 1%. This happens at 10%/50-years hazard level when the displacement reached 2.55 in., corresponding to a first-floor story drift ratio of 1.25%; see Figure 2.24. The hysteretic cycle was symmetric, with a maximum positive and negative base shear slightly bigger than 2500 kips. At 5%/50-years hazard level, the shape of the hysteretic cycle is similar, but there is more evidence of the nonlinear behaviour of the system; see Figure 2.25. The force did not increase significantly (maximum base shear of 2700 kips), but a big increase in displacement can be observed (3.27 in.), with a maximum first-floor story drift ratio of 1.6%. This effect is even more noticeable at 2%/50-years hazard level, when a further slight increase of force led to story drift ratios that reached 2.1% for maximum displacements of 4.25 in.; see Figure 2.26. The hysteretic cycle is no longer symmetric, and the degradation in stiffness appears more evident than in the previous case, but it is not conspicuous even at this hazard level.

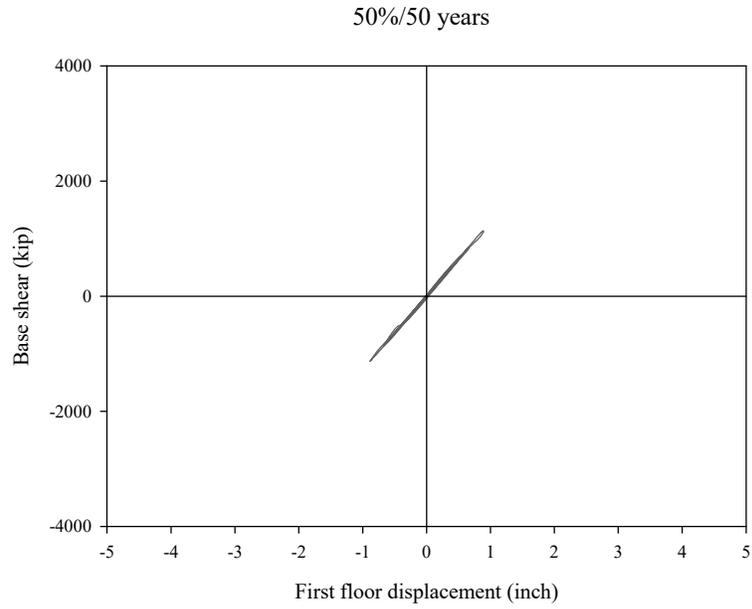

**Figure 2.22**     Hysteretic cycle response of the system at 50%/50-years hazard level.

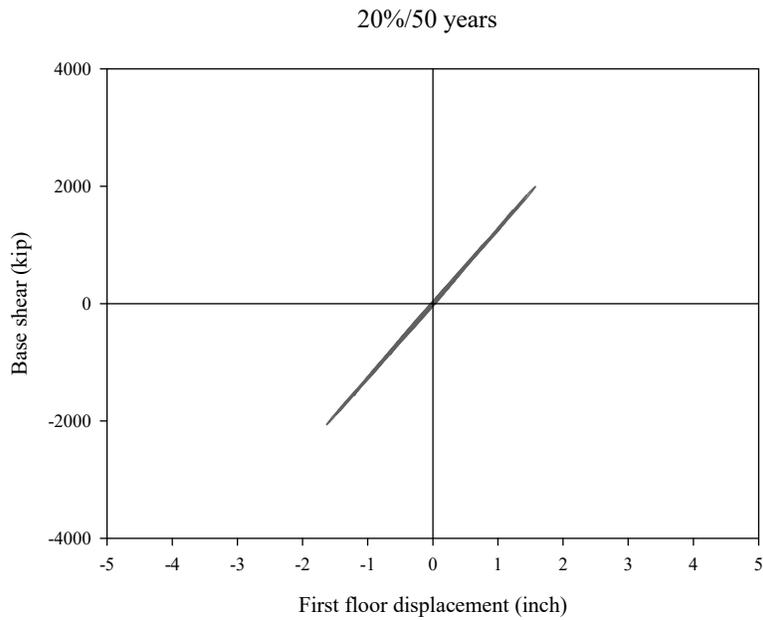

**Figure 2.23**     Hysteretic cycle response of the system at 20%/50-years hazard level.

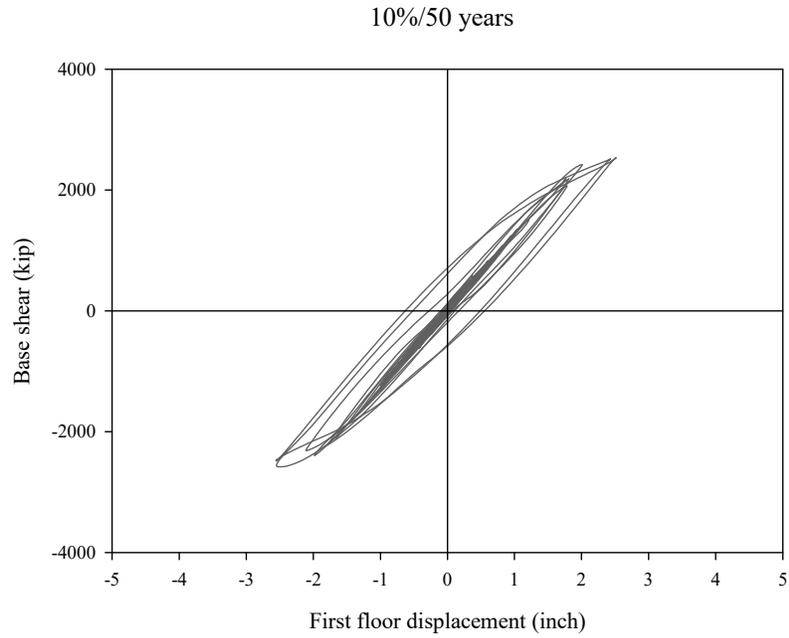

**Figure 2.24**    Hysteretic cycle response of the system at 10%/50-years hazard level.

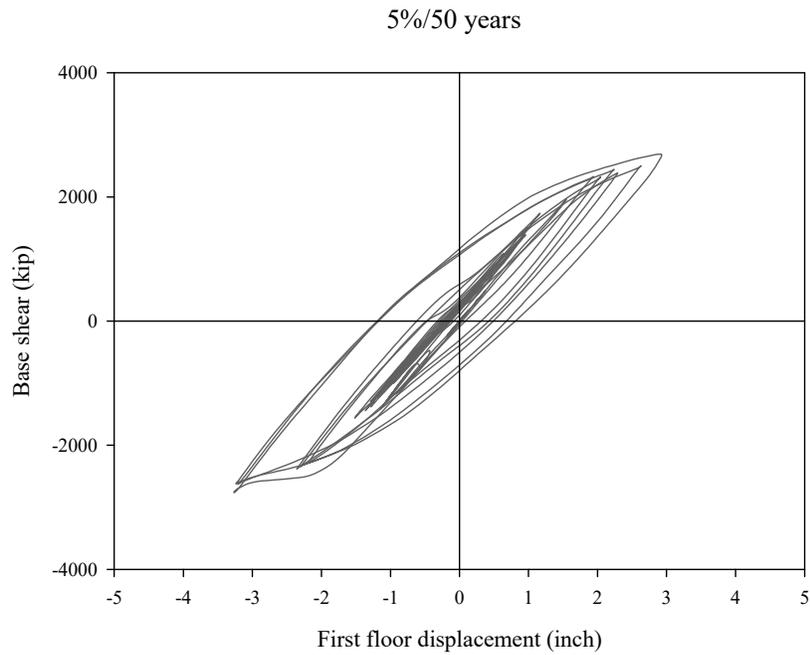

**Figure 2.25**    Hysteretic cycle response of the system at 5%/50-years hazard level.

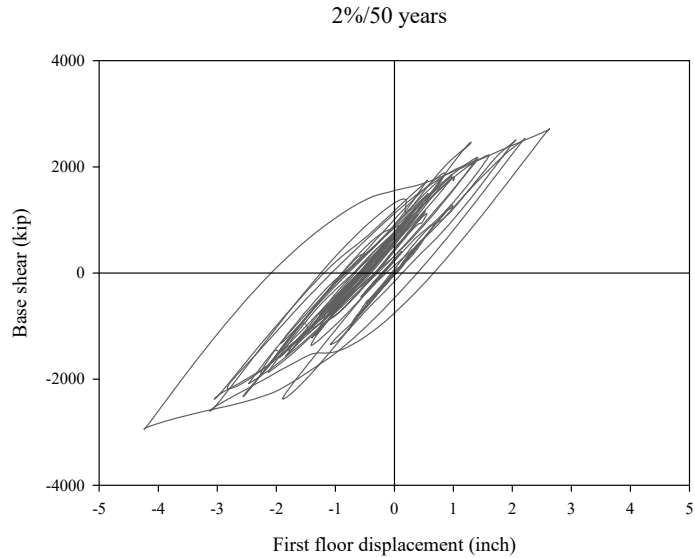

**Figure 2.26**     **Hysteretic cycle response of the system at 2%/50-years hazard level.**

## 2.6   LOSS ANALYSIS: HEALTHCARE FACILITY

Two loss metrics were used to estimate effectiveness of isolation system in reducing the total financial losses: (1) financial losses associated with the cost required to implement repairs and (2) repair time. A seismic life-cycle assessment to calculate repair costs and repair times for the two systems (fixed-base and base-isolated moment frames) and healthcare occupancy, at each of five considered hazard levels was performed using the computer software Performance Assessment Calculation Tool (PACT) [ATC 2012].

### 2.6.1  Introduction to PACT

The seismic performance assessment methodology depicts building performance after earthquake shaking. The consequences resulting from earthquake damage are expressed in terms of direct economic losses and indirect economic losses. The direct economic losses are identified by the building repair costs that represent the cost to restore the building to its pre-earthquake condition. The indirect economic losses are related to the business interruption activities and with building downtime. These performance measures make it possible to quantify earthquake consequences in terms that are meaningful to decision-makers and more directly useful in the decision-making process compared to the traditional discrete performance levels that have been used in the past.

Many input data are required in PACT to perform the analysis and predict the performance of the building under earthquake shaking. The user must input the basic building information and a description of the population model that includes the distribution of people inside the building and the variability of this distribution over the time. Then a detailed description of type and quantities of structural/non-structural elements and contents has to be added. Several types of occupancies are available in PACT to model the population and to define components types and quantities without performing a specific inventory of the building.

In the final stage, the EDPs computed with the nonlinear response history analysis (maximum story drifts, maximum absolute floor horizontal accelerations, and peak residual story drifts) are inputted for each ground motion at each hazard level. Building components are generally sensitive to story drifts or floor accelerations. In PACT, each building component is associated with a fragility curve that correlates EDPs to the probability of that item reaching a particular damage state. The component's damage is then related to a loss (e.g., repair cost or repair time) utilizing consequence functions. The total loss at each hazard level is then estimated by integrating losses over all components of a system. To account for the uncertainties affecting calculation of seismic performance, FEMA P-58 methodology uses a Monte Carlo procedure to perform loss calculations [FEMA 2012].

### 2.6.2 Basic Building Data

Basic building information includes the size of the building (number of stories, story height, and floor area at each level) and the replacement cost. Building replacement costs are equal to the initial construction cost increased by 20% to include cost allowances for demolition and site clearance [FEMA 2012]. It is used as input in PACT to compute the cost associated with damage that makes the building irreparable. It can occur if the repair cost of the structure exceeds 0.4 times the replacement cost. This threshold is suggested by the *FEMA P-58* based on past studies that show that a building is likely to be replaced rather than repaired if the ratio between repair cost and replacement cost is bigger than 0.4.

Construction cost of healthcare facilities cannot be estimated without including the huge cost of medical equipment. For this reason, it was calculated using the metric of $597.7/ft$^2$ estimated by M. Phipps per Mayencourt [2013]. Considering the footprint of the three-story building, the initial construction cost of the healthcare is estimated to be $38,730,960; this is the same for the two considered structural systems. Construction and replacement costs for the considered systems are shown in Table 2.12.

Table 2.12    Construction and replacement costs for healthcare facilities.

|  | **HP-SMRF** | **BI-IMRF** |
|---|---|---|
| Construction cost | $38,730,960 | $38,730,960 |
| Replacement cost | $46,477,152 | $46,477,152 |

### 2.6.3 Population Models

The building population model defines the number of occupants for 1000 ft$^2$ of building floor space. PACT provides information for the peak population for each occupancy and the time of the day when this peak is expected to occur. The peak population for healthcare facilities is 5 occupants for 1000 ft$^2$, which is reached at 3:00 pm. The PACT manual also contains data about the variation in the population as a percentage of the peak population over the course of a 24-hr. period; see Figure 2.27.

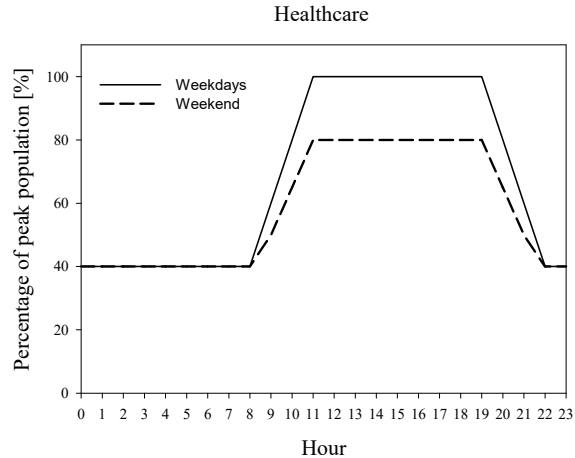

**Figure 2.27** Population distribution over the course of 24 hours for healthcare facilities [FEMA 2012].

### 2.6.4 Fragility Groups

To determinate the damage, all structural/non-structural components and contents considered vulnerable to earthquake shaking are described by using fragility functions and consequence functions. Each fragility function specifies DS probability for a certain demand parameter. Because of the uncertainty in evaluating the type and the extent of damage, fragility functions are lognormal distributions. Each component is generally defined by three DSs (slight, moderate, or extensive), and a unique fragility function is required for each of them. Fragility curves are defined by a median demand value ($\lambda$) at which there is the 50% of probability that the DS will initiate and by the dispersion ($\beta$), which takes into account the uncertainty that the DS will initiate at that particular demand value. A repair cost is associated to the DS by means of consequence functions that indicate the probable loss associated to a certain level of damage. Fragility information is contained in PACT database (*FEMA P-58 Fragility Specifications*). This database may not include all components. In case of unavailability of data, user-defined fragility curves can be developed and added to PACT to suit the needs of a specific building. Figure 2.28 shows a sample fragility function for steel column base plates given in the database. In this sample curve for a story drift ratio of 0.06 rad, there is 10% probability the component is damaged more than the third state (complete fraction of the column), 25% probability the component reaches an intermediate DS between the DS2 (propagation of brittle crack into column) and DS3, 50% probability the damage is between DS1 (initiation of crack) and DS2, and 15% probability the component is undamaged.

All components are categorized in PACT by using fragility groups. A fragility group is a collection of components with similar construction characteristics, similar potential model of damage and probability of incurring damage, and similar potential consequences resulting from the damage. Each fragility group is identified by a specific classification number. The first letter of this number indicates the macro-category of each group:

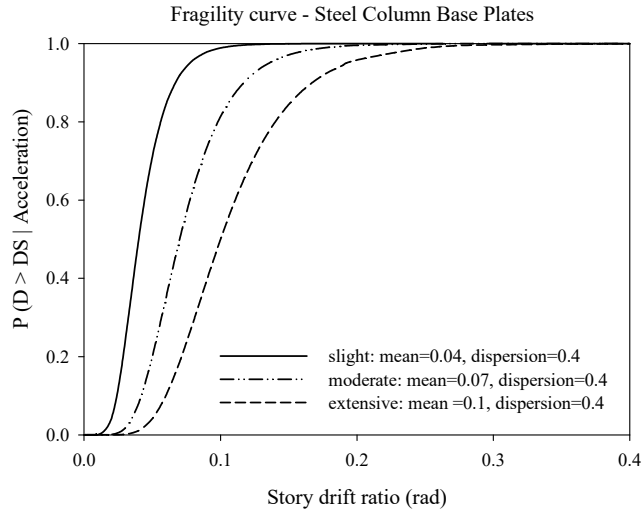

**Figure 2.28** Steel column base plates fragility curves.

- A – Substructure
- B – Shell
- C – Interiors
- D – Services
- E – Equipment and Furnishings

The type and quantities of each component must be clearly defined. For most non-structural components and contents used in the loss analysis, normative quantities recommended by *FEMA P-58* [FEMA 2012] were used. Structural/non-structural components and building contents considered in the loss analysis are given in Table 2.13. They include: (i) structural: moment connections, shear tab gravity connections, base plates, and column splices; (ii) non-structural: partition walls, curtain walls, cladding, ceiling, lighting, stairs, elevators, and MEP components, and (iii) content: bookcases, filing cabinets, computers, servers, and medical equipment. Isolator devices and utilities at the isolation level are not included in the loss model due to unavailability of their fragility functions in PACT. Healthcare occupancy implies higher quantities of most of the components due to higher standards required as well as specific types of non-structural components.

Shell components include the superstructure (B10), the exterior enclosure (B20), and the roof (B30). Type and quantities of structural components (such as column base plates or RBS connections) were estimated directly from design drawings of the building. The only differences between fixed-base and isolated-systems components are in the superstructure components. Interiors (C) include wall partitions, stairs, and ceilings. Services (D) are the components that make the building operational and include plumbing, HVAC, and electrical systems. As an example, Table 2.14 lists components that form the HP-SMRF superstructure.

Equipment and furnishings (D) are strongly related to the type of occupancy. The PACT database contains fragility specification for file cabinets, bookcases, and electronic equipment.

Healthcare facilities contain heavy and costly-to-repair equipment that could highly affect the overall repair cost of the building. Fragility functions for medical equipment are not available in PACT nor are consequence functions. In this study, fragility functions for medical equipment were developed using two methods: a global fragility curve representing all medical equipment [Yao and Tu 2012] and multiple curves representing the single medical items of a generic hospital inventory [ATC 63].

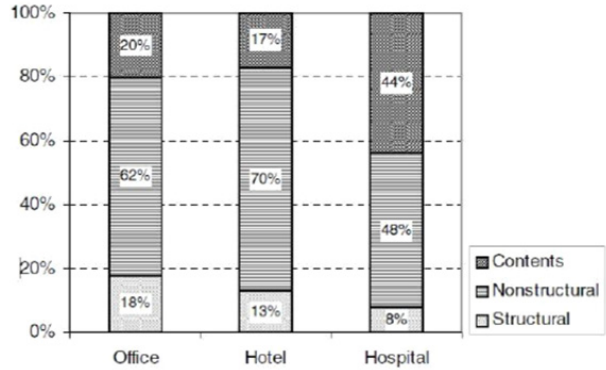

**Figure 2.29    Cost breakdown of office buildings, hotels, and hospitals.**

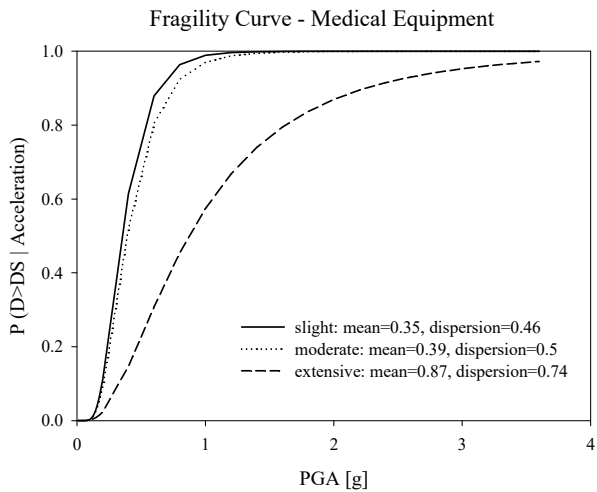

**Figure 2.30    Fragility curve for medical equipment for slight, moderate, and extensive damage levels [Yao and Tu 2012].**

First, the simplified method was used to make a comparison between fixed-base and base-isolated designs. The Yao and Tu medical equipment curve was derived by investigating 41 healthcare buildings after the 1999 Chi-Chi, Taiwan, earthquake. Figure 2.30 shows the median and dispersion values for slight, moderate, and extensive levels of damage. Because hat item is sensitive to peak ground acceleration, it was placed at the building ground floor. Per Miranda and Taghavi [2003], the consequence functions relating the component damage to the repair costs for an extensive level of damage were based on assuming that medical equipment cost would be 44% of the total building cost. The consequence functions relating damage of medical equipment to repair time were not developed due to unavailability of data.

Table 2.13    Type of components considered in the analysis.

| Type | Structural | Building Content | Non-structural | | | |
|---|---|---|---|---|---|---|
| Component | Bolted shear tab gravity connections | Tall File Cabinet | Curtain Walls | Precast Concrete Panels | Chiller | Control Panel |
| | Steel Column Base Plates | Unanchored Bookcase | Wall Partition, Type: Gypsum with metal studs | Cold Water Piping | Cooling Tower | Fire Sprinkler Drop Standard Threaded Steel |
| | Column Splices | Desktop Electronics | Wall Partition, Type: Gypsum + Ceramic Tile | Hot Water Piping - Small Diameter Threaded Steel | Air Handling Unit | Transformer/primary service |
| | Post-Northridge RBS connection with welded web | Medical Equipment | Suspended Ceiling | Hot Water Piping - Large Diameter Welded Steel | HVAC Galvanized Sheet Metal Ducting | Motor Control Center |
| | Post-Northridge welded steel IMRF connection other than RBS | | Recessed lighting in suspended ceiling | Sanitary Waste Piping | HVAC Drops | Low Voltage Switchgear |
| | | | Prefabricated steel stairs | Steam Piping - Small Diameter Threaded Steel | Variable Air Volume (VAV) box | Distribution Panel |
| | | | Hydraulic Elevator | Steam Piping - Large Diameter Welded Steel | Concrete tile roof | |

Table 2.14    HP-SMRF superstructure fragility groups.

| PACT ID | Component type | EDP | Unit | Floor | Dir. 1 | Dir. 2 | Non dir. |
|---|---|---|---|---|---|---|---|
| B1031.001 | Bolted shear tab gravity connections | Story drift ratio | Ea. | 1 | 136 | 84 | - |
|  |  |  |  | 2 | 136 | 84 | - |
|  |  |  |  | 3 | 136 | 84 | - |
| B1031.011a | Steel column base plates, column W < 150 plf | Story drift ratio | Ea. | 1 | 5 | 4 | - |
| B1031.011c | Steel column base plates, column W > 300 plf | Story drift ratio | Ea. | 1 | 12 | 14 | - |
| B1031.021a | Welded column splices, column W < 150 plf | Story drift ratio | Ea. | 3 | 5 | 4 | - |
| B1031.021c | Welded column splices, column W > 300 plf | Story drift ratio | Ea. | 3 | 12 | 14 | - |
| B1035.002 | Post-Northridge RBS connection with welded web, beam one side of column only, beam depth >= W30 | Story drift ratio | Ea. | 1 | 4 | 8 | - |
|  |  |  |  | 2 | 4 | 8 | - |
|  |  |  |  | 3 | 4 | 8 | - |
| B1035.012 | Post-Northridge RBS connection with welded web, beams both sides of column, beam depth >= W30 | Story drift ratio | Ea. | 1 | 8 | 6 | - |
|  |  |  |  | 2 | 8 | 6 | - |
|  |  |  |  | 3 | 8 | 6 | - |

### 2.6.5 Repair Costs and Loss Ratio

Repair cost estimates can provide the design engineer with valuable insights regarding the desirability and cost-effectiveness of enhancements to the structural system. Table 2.15 and Figure 2.31 show the median repair costs for the fixed-base and base-isolated moment frames for the five considered hazard levels. The diagram clearly shows effectiveness of base-isolated systems in mitigating damage. Reductions in the cost of repairs are consistently high at all hazard levels for the base-isolated systems. The reduction in repair costs ranged from 76% to 88%, with an average of 85%. While the fixed-base system generates disproportionally greater losses for the most expensive facility, the base-isolated system generates proportionally greater losses.

To identify the major contributors to the losses, the partial contributions of structural/non-structural components and contents to repair costs is shown in Figure 2.31. Non-structural components and contents dominate the losses for healthcare facility. In the case of the fixed-base healthcare facility, non-structural components dominate the losses (72% contribution) at the lower hazard levels (50% and 20% in 50 years). At the 10%/ and 5%/50-years hazard levels, non-structural components and content have an almost equal contribution to the total repair costs. At the 2%/50-years hazard level, damage to the medical equipment, which is the primary source of the content damage, dominates the losses. The cost of repairs of structural components, although minor for the fixed-base system at higher hazard levels, is decreased through utilization of base isolation. Although total repair costs points to the improvement provided by the isolation system, the partial contribution of structural/non-structural elements and contents needs to be considered; see Table 2.16. The results obtained are graphically shown in Table 2.17 in which the fixed-base and base-isolated systems are compared.

To facilitate decisions on whether to repair or replace a building damaged after an earthquake, repair costs can be expressed in terms of a loss ratio, which FEMA P-58 [FEMA 2012] defines as the necessary repair costs divided by the building's replacement costs. The building's replacement cost in this case is based on the initial construction cost associated with structural and nonstructural components; the contents are excluded. According to FEMA P-58, building owners typically prefer to replace a building rather than repair it when the loss ratio exceeds 40%. Figure 2.32 plots loss ratios for each system at the five considered hazard levels. Although the fixed-base healthcare buildings have significantly higher loss ratios than the base-isolated buildings at all hazard levels, the highest loss ratio of the fixed-base system of 0.26 is significantly smaller than the FEMA P-58 replacement threshold of 0.4. Therefore, none of the buildings will require replacement even for a very rare earthquake with the 2% probability of exceedance in 50 years.

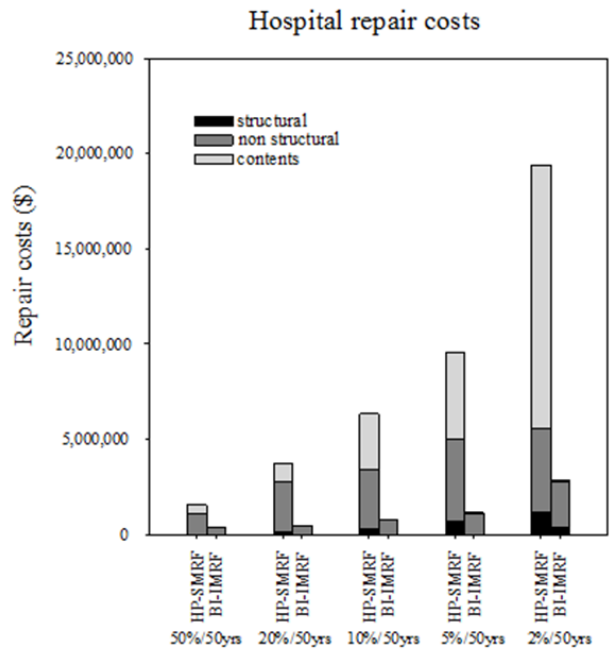

**Figure 2.31** Median repair costs for the HP-SMRF and the BI-IMRF for five hazard levels: 2%, 5%, 10%, 20%, and 50% probabilities of exceedance in 50 years.

**Table 2.15  Median repair costs at all hazard levels for all systems.**

| | Repair costs ($) | | | | |
|---|---|---|---|---|---|
| | **50%/50yrs** | **20%/50yrs** | **10%/50yrs** | **5%/50yrs** | **2%/50yrs** |
| HP-SMRF | 1,551,724 | 3,723,077 | 6,333,333 | 9,500,000 | 19,333,333 |
| BI-SMRF | 370,656 | 475,229 | 796,178 | 1,135,762 | 2,852,941 |

**Table 2.16  Median repair cost of structural components, non-structural components, and contents.**

| | HP-SMRF repair cost ($) | | | | BI-IMRF repair cost ($) | | |
|---|---|---|---|---|---|---|---|
| | **Structural** | **Non structural** | **Contents** | | **Structural** | **Non structural** | **Contents** |
| 50%/50yrs | 0 | 1,131,536 | 417,189 | 50%/50yrs | 0 | 369,057 | 1,912 |
| 20%/50yrs | 136,775 | 2,662,720 | 922,603 | 20%/50yrs | 0 | 472,315 | 3,099 |
| 10%/50yrs | 303,006 | 3,105,441 | 2,922,484 | 10%/50yrs | 0 | 794,419 | 1,915 |
| 5%/50yrs | 804,520 | 4,229,748 | 4,553,780 | 5%/50yrs | 11,702 | 1,089,125 | 38,883 |
| 2%/50yrs | 1,140,311 | 4,413,423 | 13,830,020 | 2%/50yrs | 386,178 | 2,380,857 | 80,162 |

**Table 2.17  Repair cost contribution at 2%/50-years hazard level.**

| Repair cost contribution | | |
|---|---|---|
| | **HP-SMRF** | **BI-IMRF** |
| Medical equipment | 73.29% | 0.49% |
| Wall Partition | 7.42% | 52.72% |
| Suspended Ceiling | 4.39% | 0.90% |
| Precast Concrete Panels | 3.91% | 26.99% |
| Concrete tile roof | 2.40% | 0.26% |
| Post-Northridge connections | 1.39% | 4.46% |
| Bolted shear tab gravity connections | 1.03% | 8.81% |
| Desktop electronics | 0.60% | 2.23% |
| | 94.43% | 96.88% |

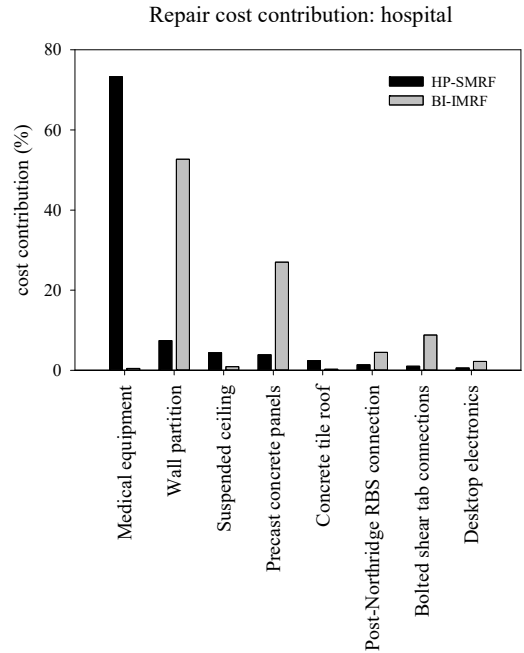

**Figure 2.32    Repair cost contribution at 2%/50-years hazard level.**

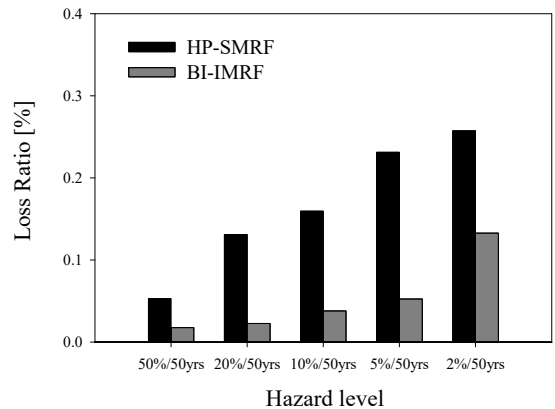

**Figure 2.33    Median loss ratio for the HP-SMRF and the BI-IMRF for five hazard levels and healthcare occupancy.**

### 2.6.6  Repair Time

To estimate the resilience of the system and the revenue losses resulting from the business interruption following an earthquake, business downtime as a function of time needs to be characterized. Business downtime should include the time required to: (1) identify damage, design repairs or upgrades, obtain permits and financing, and mobilize supplies and manpower; and (2) make the repairs necessary to restart operations. If the structural/non-structural damage is such that the building is not deemed economical to repair, which corresponds to a loss ratio bigger than 40%, downtime is assumed equal to the replacement time for the building. If there is a minor damage that doesn't affect the operability of the system, the downtime is set equal to zero. If there is major non-structural damage that affects a building's functionality, the owner can simultaneously begin site preparation, building clean up, mobilizing a general contractor, and securing necessary funding before initiating repairs. In the case of extensive non-structural damage or major structural damage, a detailed inspection is needed. Site preparation and clean up can occur simultaneously with the inspection. Although business models exist for commercial occupancy types (e.g., Terzic et al. [2014a]), a similar type of model could not be found for a healthcare facility. Therefore, the study presented herein will use repair time as a metric for comparing the two systems and two occupancy types.

Estimating the time required to repair a structure is difficult without specific information about the availability of workers and material. To calculate repair time, a number of assumptions were made. PACT estimates repair time considering the number of labor hours associated with the required repair at each DS. One of the input parameters in PACT is the number of workers that can occupy the building at the same time, which is dependent on whether the building is occupied during construction. In this study, it was assumed that supplies and workers were available to permit the necessary work. A high density of workers (one worker per 500 $ft^2$) was used, assuming that the building would be unoccupied during the repair of damaged building components. The repair time was calculated considering two repair schemes: (1) the *parallel* scheme assumes simultaneous repair at all three floors, and (2) the *serial* scheme assumes sequential repair at three floor levels [FEMA 2012]. Both repair schemes assume sequential repair of all damaged components within one floor level.

These repair schemes are not optimal but provide a good estimate of repair time of the lower and upper bound for the chosen density of workers. While the assumptions made may be feasible for the systems with the smaller level damage (i.e., the isolated system), they may be hard to achieve for the systems with more extensive damage (i.e., the fixed-base system). Therefore, these assumptions are advantageous for the HP-SMRF relative to the base-isolated system as they reduce the relative benefits of the isolated system. Figure 2.35 shows the median repair times for the HP-SMRF and the BI-IMRF for all five hazard levels considering the two repair strategies: parallel and serial. The base-isolated systems are again very effective in reducing repair time, thus significantly reducing the total downtime for the isolated buildings. Upper (serial) and lower (parallel) bounds of the repair times are both several magnitudes smaller for the isolated buildings relative to the fixed-base buildings. For the 50%/50-years hazard level, the repair times of the base-isolated buildings are 2–3 times smaller than for the fixed-base buildings. For the higher hazard levels, 20%/, 10%/, and 5%/50-years hazard level, base isolation is even more effective, resulting in 4–6 times smaller repair time. For the 2%/50 years hazard level, the reduction in repair time is 3–4 times, which is still significant. Repair times would have been even higher if repair time of medical equipment was included in the consequence function.

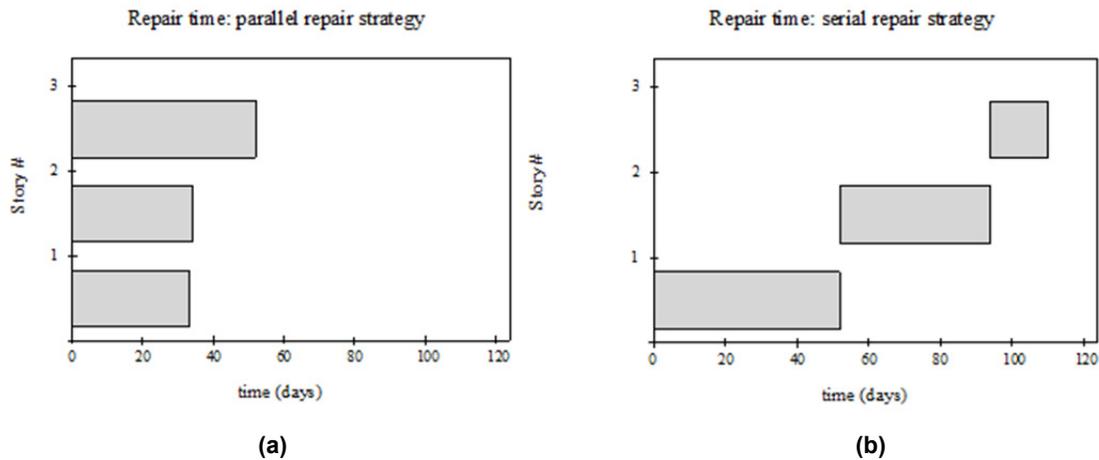

(a)                                                   (b)

**Figure 2.34**      Comparison between parallel (a) and serial (b) repair strategy.

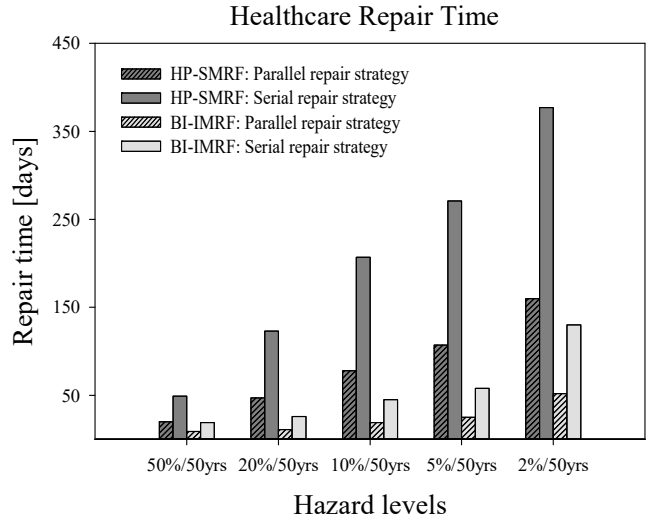

**Figure 2.35**      Median repair times for the HP-SMRF and the BI-IMRF for five hazard levels considering two repair strategies: parallel and serial.

### 2.6.7 Indirect Losses

Indirect losses were calculated as the hospital's daily loss in income due to business interruption. By averaging the square footage of four California hospitals located in the Bay Area, the number of beds available for a 64,800 ft$^2$ hospital was estimated to be 92. Since the average revenue per inpatient bed is approximately $1 million per year for California hospitals [Meade et al. 2002], the hospitals' daily income is $252,055.

In the case of the serial repair strategy, a floor is considered functional after being repaired even if the other floors are not yet functional. Assuming that the population is equally distributed over three floors, the daily income after the first floor is repaired is a third of the total income under normal conditions. The same approximation can be applied to the other floors. In the case of the parallel repair strategy, although the functionality is recovered at once, there will be no daily income until the total recovery of the structure is completed. Analyzing the damage states in PACT, it was found that several structural components were seriously damaged after an earthquake with 2%/50-years probability of exceedance occurs. Assuming that an inoperative building will correspond to a drop to zero in terms of its functionality is realistic. Figure 2.36 compares the two repair strategies after the drop of functionality for the buildings that have experienced an earthquake with 2%/50-years probability of exceedance.

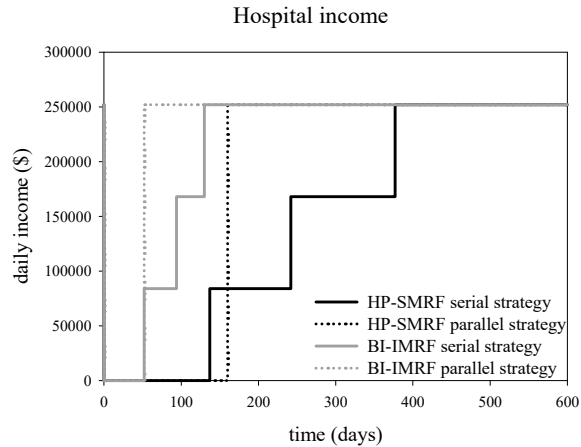

**Figure 2.36**  Comparison between serial and parallel repair strategy for the structures considered.

### 2.6.8 Total Annualized Losses

Annualized losses (due to both repair costs and business interruption costs) were calculated by integrating direct and indirect losses at each hazard level weighted by the probability that the hazard was exceeded. Then the total losses were obtained by adding the direct and indirect losses. Table 2.18 shows the results for serial and parallel repair strategy. Figure 2.37 shows the comparison between the considered systems.

As can be seen in Table 2.18, isolating the buildings led to a significant reduction of the total annualized losses, resulting in an average decrease of more than 70% for the healthcare facility. For both structural systems, the contribution of indirect losses is significantly bigger than the contribution of the direct losses. In both cases, the serial repair strategy resulted in bigger indirect losses when considering repair times. The chart shows that even if the assumptions are conservative, the performance of the base-isolated structures is superior compared to fixed-base structures.

**Table 2.18**  Annualized losses for all systems considered.

|  | Serial repair strategy | | | Parallel repair strategy | | |
|---|---|---|---|---|---|---|
|  | **Direct loss** | **Indirect loss** | **Total** | **Direct loss** | **Indirect loss** | **Total** |
| HP-SMRF | 69,061 | 328,377 | 397,438 | 69,061 | 208,248 | 277,309 |
| BI_SMRF | 11,263 | 95,885 | 107,148 | 11,263 | 65,725 | 76,989 |

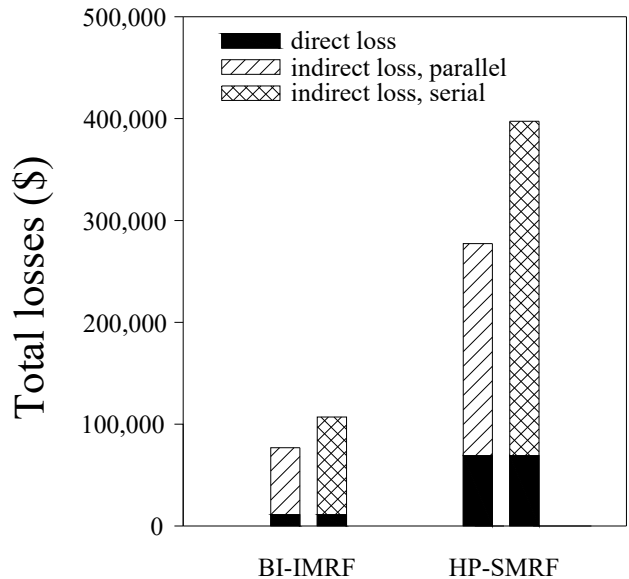

Figure 2.37    Annualized losses for all systems considered.

## 2.7  LOSS ANALYSIS: SCHOOL

Next, we performed a seismic life-cycle assessment using the same methodology—PACT [ATC 2012] —to calculate repair costs and repair times for the two systems (fixed-base and base-isolated moment frames) for a school at each of five considered hazard levels. Two loss metrics were used to estimate the effectiveness of isolation system in reducing the total financial losses: (1) financial losses associated with the costs required to implement repairs and (2) repair time.

### 2.7.1  Basic Building Data

Basic building information includes the size of the building (number of stories, story height, and floor area at each level) and the replacement cost. Replacement costs for the buildings are equal to the initial construction cost increased by 20% to include cost allowances for demolition and site clearance [FEMA 2012]. It is used as input in PACT to compute the cost associated with damage if the building is deemed irreparable, which occurs if the repair cost of the structure exceeds 0.4 times the replacement cost. This threshold is suggested by the FEMA P-58 based on past studies that show that a building is likely to be replaced rather than repaired if the ratio between repair cost and replacement cost is bigger than 0.4.

The initial construction cost of the school is estimated to be $17,823,000 for the HP-SMRF and $17,408,000 for the BI-IMRF, which is the same as if it were a commercial building [Terzic et al., 2014a; Ryan et al. 2010]. The HP-SMRF required 40% more steel, which increased the overall cost of the building compared to the code-minimum design. Because of the cost of the isolators, the base-isolated structure shows an increase in cost when compared to the code-minimum design; see Table 2.19.

Table 2.19   Construction and replacement costs for schools.

|  | Code minimum | HP-SMRF | BI-IMRF |
|---|---|---|---|
| Construction cost | $16,800,000 | $17,823,000 | $17,408,000 |
| Replacement cost | $20,160,000 | $21,387,600 | $20,889,600 |

### 2.7.2 Population Model

The building population model defines the number of occupants for 1000 ft$^2$ of building floor space. PACT provides information for the peak population for several types of occupancy (including schools) and the time of the day in which this peak is expected to occur; see Table 2.20. It also contains data about the variation of population as a percentage of peak population over a 24-hour period; see Figure 2.38.

Table 2.20   Peak population model.

| Occupancy description | Peak population model (occupant per 1000 sq. ft.) | Peak population model (time of day) |
|---|---|---|
| Education(K-12): Elementary Schools | 14 | Daytime |

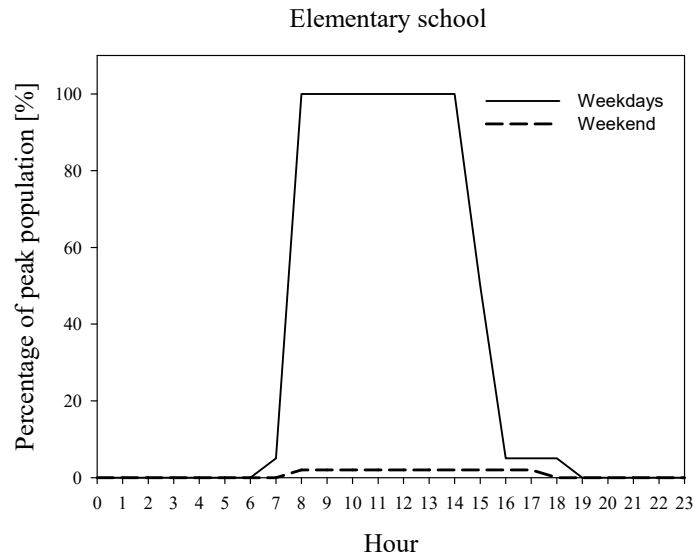

Figure 2.38   Population distribution of occupancy for a school over a 24-hour period [FEMA 2012].

### 2.7.3 Fragility Groups

As described in Section 2.6.4, all structural/non-structural elements and contents considered vulnerable to earthquake shaking are described using fragility functions and consequence functions; see Table 2.1.3 and Figure 2.2.8. The shell components include the superstructure (B10), the exterior enclosure (B20), and the roof (B30). Type and quantities of structural components (such as column base plates or RBS connections) were estimated directly from design drawings of the building. Fixed-base and isolated-system components differed only for the superstructure. Interiors included wall partitions, stairs, and ceilings. Services considered are only those components that make the building operational. They included plumbing, HVAC, and electrical systems. Equipment and furnishing were strongly related to the type of occupancy studied. The PACT database contains fragility specifications for file cabinets, bookcases, and electronic equipment, which are considered adequate representative of standard contents for a school.

### 2.7.4 Repair Costs and Loss Ratio

Repair cost estimates provide the design engineer with valuable insights regarding the desirability and cost-effectiveness of enhancements to the structural system. Table 2.21 and Figure 2.39 show the median repair costs for the fixed-base and base-isolated moment frames for the five considered hazard levels and occupancy type. It clearly shows effectiveness of base-isolated system in mitigating damage. Reduction in cost of damage repair is consistently high at all hazard levels. It ranged from 66% to 82%, with an average of 76%.

To identify the major contributors to the losses, Figure 2.39 shows the partial contributions of structural/non-structural elements and contents to the cost of repairs. For both structural systems, non-structural components dominated the losses (for a contribution greater than 73%); however, the contribution for the base-isolated building was smaller. The cost of repairing structural components, although minor for the fixed-base system at higher hazard levels (up to 23%), was completely diminished if base-isolation was used. Although the total repair costs highlights the improved performance provided by the isolation system, the partial contribution of structural/non-structural elements and contents needs to be considered; see Figure 2.39 and Table 2.22.

Table 2.21    Median repair costs at all hazard levels for all systems.

|  | Repair costs ($) | | | | |
| --- | --- | --- | --- | --- | --- |
|  | 50%/50yrs | 20%/50yrs | 10%/50yrs | 5%/50yrs | 2%/50yrs |
| HP-SMRF | 592,967 | 1,415,217 | 2,415,217 | 3,242,424 | 4,605,555 |
| BI-SMRF | 199,355 | 259,655 | 438,710 | 567,568 | 1,312,500 |

Table 2.22 Median repair cost of structural components, non-structural components and contents for school structures.

| | HP-SMRF school repair costs ($) | | | | BI-IMRF school repair costs ($) | | |
|---|---|---|---|---|---|---|---|
| | structural | non structural | contents | | structural | non structural | contents |
| 50%/50yrs | 0 | 503,005 | 88,959 | 50%/50yrs | 0 | 196,327 | 3137 |
| 20%/50yrs | 7304 | 1,245,003 | 161,546 | 20%/50yrs | 0 | 255,669 | 4170 |
| 10%/50yrs | 68,701 | 2,178,899 | 171,154 | 10%/50yrs | 0 | 426,509 | 11,494 |
| 5%/50yrs | 364554 | 2,700,564 | 181,969 | 5%/50yrs | 2,937 | 549,598 | 14,474 |
| 2%/50yrs | 1,073,583 | 3,343,834 | 187,953 | 2%/50yrs | 100,838 | 1,175,365 | 33,166 |

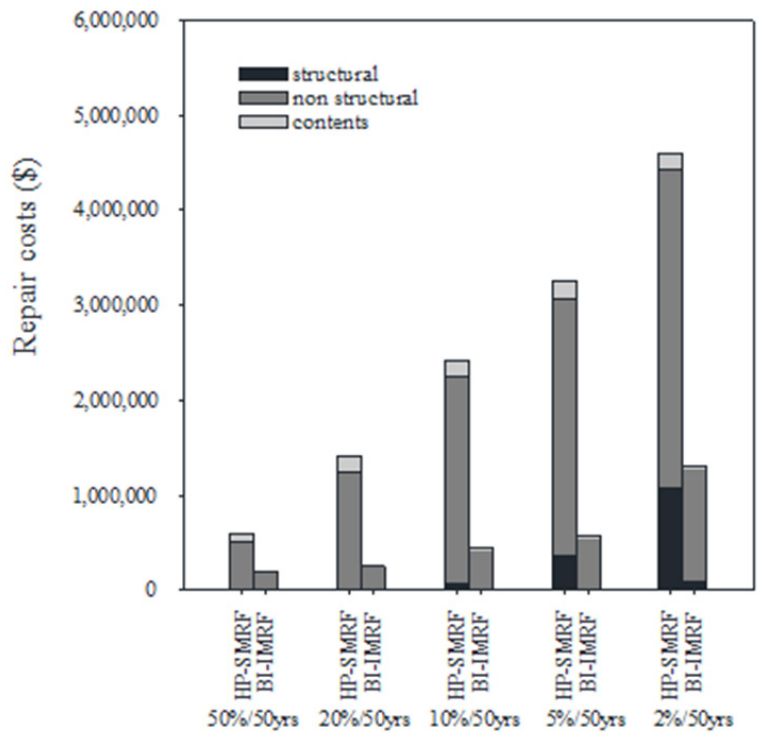

Figure 2.39 Median repair costs for the HP-SMRF and the BI-IMRF for five hazard levels: 2%, 5%, 10%, 20%, and 50% probabilities of exceedance in 50 years.

A more detailed analysis was performed to understand the contribution of the various elements in terms of repair costs. The percentages shown in Table 2.23 are related to the 2%/50-years hazard level, the level with the highest accelerations, which produces the highest repair costs. In the isolated building, a bigger decrease in accelerations with respect to drifts led to acceleration values that do not reach the first DS of the acceleration-sensitive components; i.e., the ceilings and roof. Although these components do not affect total repair cost, the drift-sensitive components contribute the most. Conversely, for the fixed-base building, acceleration-sensitive elements account for more than 30% of the total repair costs. Drift-sensitive non-structural components such as wall partitions and panels significantly affect the repair costs for both systems. Desktop electronics (including computers, monitors, stereos, etc.) account for a small contribution. The results are compared graphically in Figure 2.40.

To facilitate decisions on whether to repair or replace a building damaged after an earthquake, repair costs can be expressed in terms of loss ratio per FEMA P-58 [FEMA, 2012]. Figure 2.41 plots loss ratios for each system at the five considered hazard levels. Although the fixed-base building has significantly higher loss ratios than the base-isolated building at all hazard levels, the highest loss ratio of the fixed-base system of 0.21 is significantly smaller than the FEMA P-58 replacement threshold of 0.4. Therefore, none of the buildings would require replacement even for a very rare earthquake with the 2% probability of exceedance in 50 years.

Table 2.23    Repair cost contribution at 2%/50-years hazard level.

|  | HP-SMRF | BI-IMRF |
|---|---|---|
| Wall Partition | 21.09% | 54.75% |
| Suspended Ceiling | 21.46% | 0.02% |
| Precast Concrete Panels | 19.26% | 33.63% |
| Concrete tile roof | 10.39% | 0% |
| Post-Northridge connections | 5.80% | 2.67% |
| Bolted shear tab gravity connections | 8.74% | 6.41% |
| Desktop electronics | 2.53% | 1.45% |
|  | 89.27% | 97.47% |

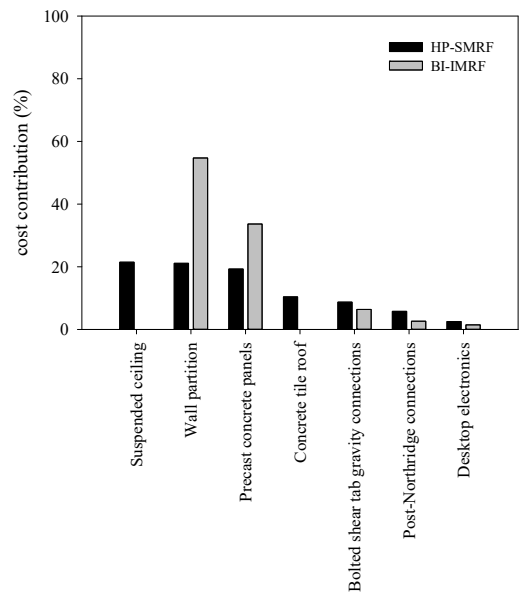

Figure 2.40    Repair cost contribution at 2%/50-years hazard level.

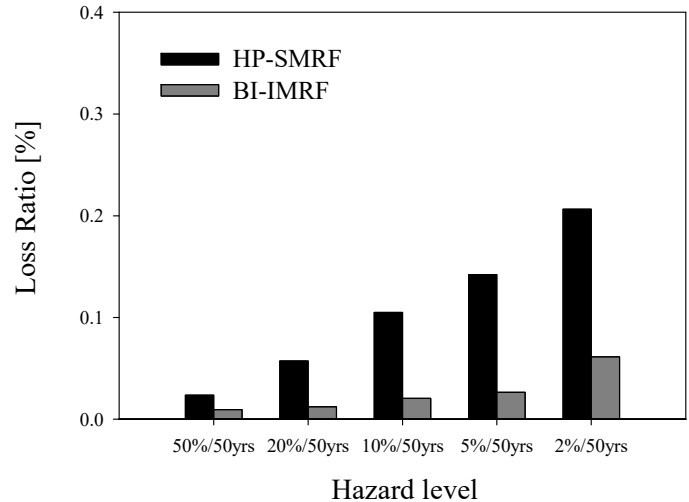

**Figure 2.41    Median loss ratio for the HP-SMRF and the BI-IMRF for five hazard levels.**

## 2.7.5  Repair Time

Although business models exist for the commercial occupancy-type structures (e.g., Terzic et al. [2014a]), such a model could not be found for school occupancy. Therefore, the study presented herein will use the repair time as a metric for comparing the two systems. Similar assumptions made for hospital occupancy considering availability of workers and materials were made for the school. A high density of workers (one worker per 500 ft$^2$) was used, and it was assumed that the building was not occupied during the repair of damaged building components. As with the hospital, the repair time was calculated considering the *parallel* scheme that assumes simultaneous repair at all three floors, and the *serial* scheme that assumes sequential repair at three floor levels [FEMA 2012]. Both repair schemes assumed sequential repair of all damaged components within one floor level; see Figure 2.42.

Figure 2.43 shows the median repair times for the HP-SMRF and the BI-IMRF for five hazard levels considering two repair strategies: parallel and serial. Base isolation was again very effective in reducing the repair time, resulting in a significantly smaller overall downtime for the isolated buildings. The upper (serial) and lower (parallel) bounds of the repair times are both several magnitudes smaller for the isolated buildings relative to the fixed-base buildings. For the 50%/50-years hazard level, the repair times of the base-isolated buildings were 2–3 times smaller compared to the fixed-base buildings. For the higher hazard levels, 20%/, 10%/, and 5%/50-years hazard levels, the base-isolation was even more effective, resulting in 4–6 times smaller repair times. For the 2%/50-years hazard level, the reduction in repair time was 3–4 times which is still significant.

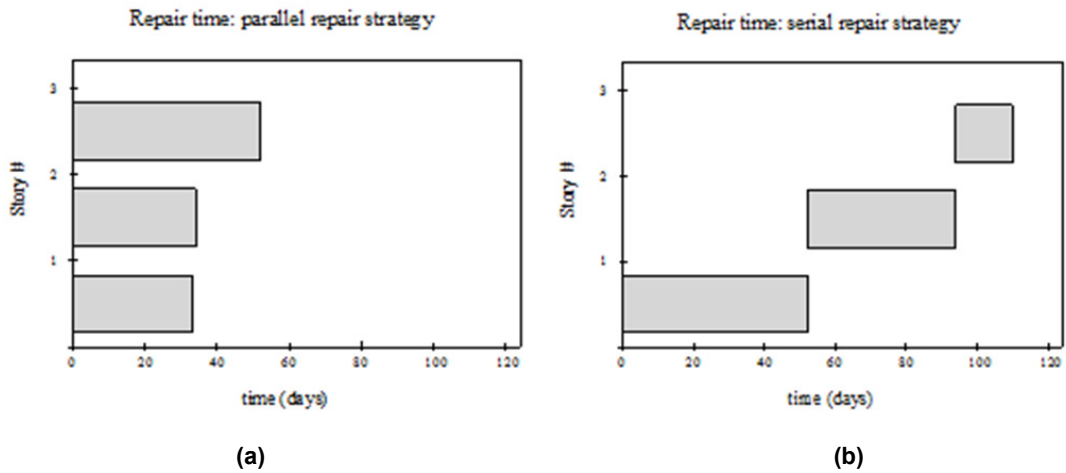

(a)                         (b)

**Figure 2.42**     Comparison between parallel (a) and serial (b) repair strategy.

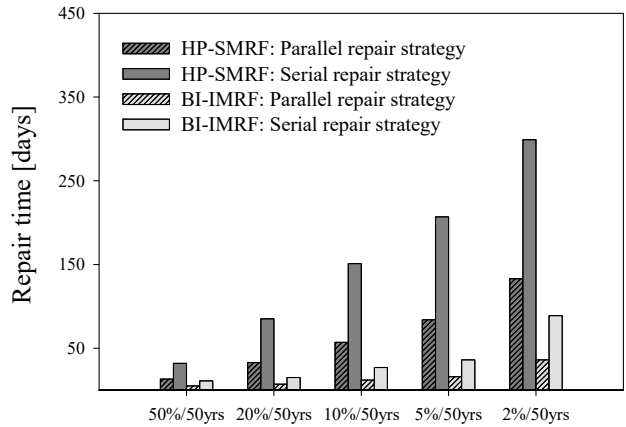

**Figure 2.43**     Median repair times for the HP-SMRF and the BI-IMRF for five hazard levels considering parallel and serial repair strategies.

### 2.7.6 Indirect losses

Indirect losses were calculated as the daily loss in income of schools due to business interruption. By averaging the square footage of four California schools, the number of students for a 64,800ft$^2$ school was estimated to be 290. Assuming the school is a private elementary school, the annual income per student is $6733 [National Centre for Education Statistics 2011]; therefore, the annual earning of the school is $1,952,570 with a daily income of $5350. The rational for both the serial and parallel repair strategy was discussed in Section 2.6.7. Figure 2.44 shows the comparison between the two repair strategies after a drop in functionality for buildings experiencing an earthquake with 2% in 50 years probability of exceedance.

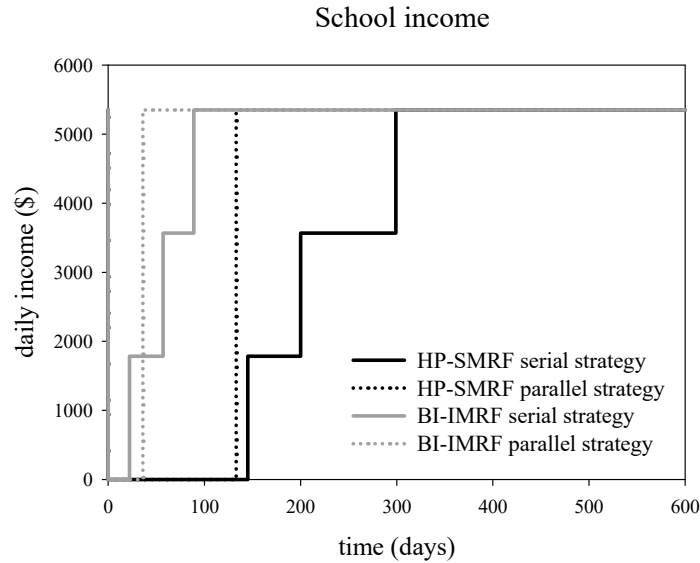

**Figure 2.44**    Comparison between serial and parallel repair strategy for the structures considered.

### 2.7.7 Total losses

Annualized losses (due to both repair costs and business interruption costs) were calculated by integrating direct and indirect losses at each hazard level weighted by the probability that the hazard was exceeded. Then the total losses were also found by adding the direct and indirect losses contribute. Table 2.24 shows the results for serial and parallel repair strategy. Figure 2.45 shows the comparison between all the considered systems.

As can be seen in Table 2.24, isolating the buildings led to a significant reduction of the total annualized losses, resulting in an average decrease of 76%. When considering school buildings and given that the daily loss in income is small, the contribution of direct losses is bigger than the contribution of indirect losses. In regards to repair time, the serial repair strategy resulted in bigger indirect losses. The chart shows that even if the assumptions are conservative for the base-isolated structures, the system performs better compared to the fixed base structure.

**Table 2.24**    Annualized losses for all considered systems.

|  | Serial repair strategy | | | Parallel repair strategy | | |
| --- | --- | --- | --- | --- | --- | --- |
|  | Direct loss | Indirect loss | Total | Direct loss | Indirect loss | Total |
| HP-SMRF | 24,825 | 5,205 | 30,029 | 24,825 | 3,144 | 27,969 |
| BI-SMRF | 5,964 | 1,120 | 7,084 | 5,964 | 839 | 6,803 |

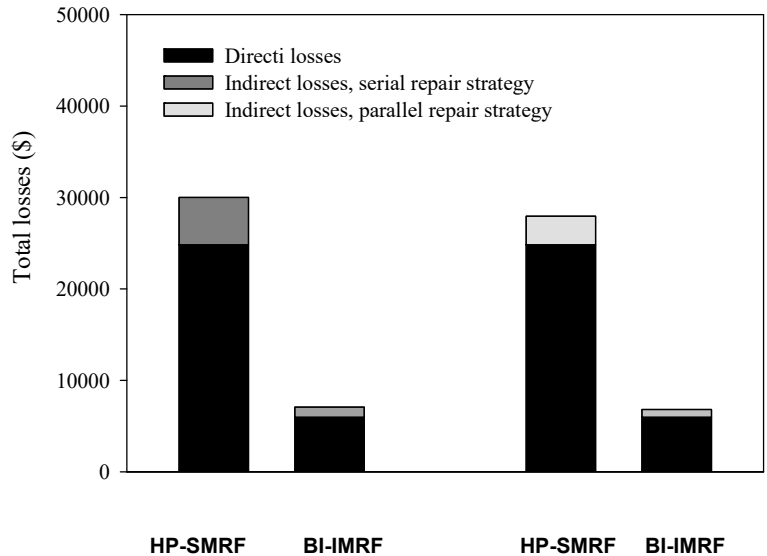

**Figure 2.45    Annualized losses for all considered systems.**

## 2.8  RESILIENCY

Resiliency is the ability of a system to sustain a level of functionality or performance or to re-establish its functionality following a hazard event. The level of resiliency is measured by integrating the recovery function of the system within a certain period of time [Cimellaro et al, 2010a, b]. To quantify the resiliency of the considered building, it is necessary to determine the recovery functions. Recovery functions describe how a given system returns to the same level of functionality before the extreme event. Functionality ranges from 0 to 100% where 100% means no reduction in performance, while 0% means total loss. In this study, two different ways to recover functionality of the system were examined. In the first recovery model, functionality is recovered all at once over a period of time that corresponds to the repair time, i.e., the parallel repair strategy used in PACT. The second recovery model is a step function that is associated with the serial repair strategy where total functionality is recovered step-wise as the repairs of each floor of the building are completed sequentially.

### 2.8.1  Healthcare Facility

When considering the healthcare facility,, it is obvious that the base-isolated buildings are more resilient than the fixed-base buildings as they have significantly reduced repair times and thus their return to functionality will be faster. To better quantify resilience, resilience functions at a very rare earthquake with a 2% probability of exceedance in 50 years were determined. For this hazard level, it was assumed that both the fixed-base and the base-isolated system incurred enough damage to trigger closure of the buildings. The probable lower and upper bounds for the recovery and therefore resiliency were established based on the lower (parallel scheme) and upper (serial scheme) bounds of repair times. Resilience functions as well as recovery functions for all the considered systems are shown in Figure 2.46.

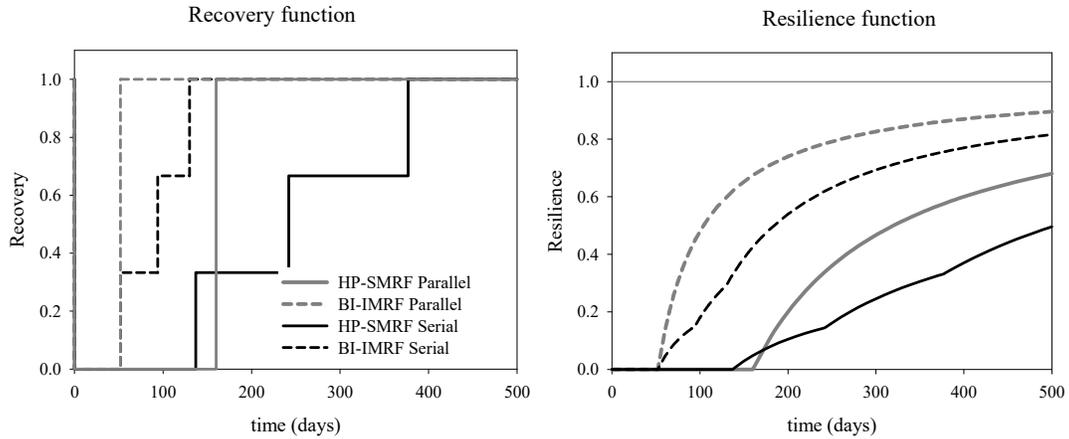

**Figure 2.46**  Healthcare facility: recovery (a) and resilience functions (b) of all considered systems at 2% in 50-year hazard level.

### 2.8.2 School Facility

Figure 2.46 clearly shows significantly greater resilience of the base-isolated system compared to the fixed-base system. Note that following a rare event, the base-isolated system starts to recover its functionality after 52 days compared to the fixed-base system where functionality only begins between 137 and 160 days. Resiliency functions are much steeper for the base-isolated building, indicating faster recovery times. Table 2.25 shows possible expressions for resilience functions over the period of time for the given systems; the standard deviation is reported to show the good match between resiliency functions and experimentally calculated data.

It is also possible to calculate the resilience factor corresponding to a specific period of time. In this study, 400 days was chosen as reference period of time because by that point, all the systems have recovered their full functionality. As can be seen in

Table 2.26, in all the considered cases, the isolated systems have higher resilience than the fixed-base systems. In particular, resilience factors for isolated systems considering the serial repair strategy are more than double with respect to the corresponding fixed-base systems.

When considering the school building, Figure 2.47 clearly shows significantly greater resilience of the base-isolated system relative to the fixed-base system. The base-isolated system starts to recover its functionality between 22 and 36 days. The fixed-base system starts to re-establish its function between 133 and 145 days following a very rare earthquake Resiliency functions are much steeper for the base-isolated system indicating faster recovery times. Table 2.27 shows possible expressions for resilience functions over the period of time for the given systems. where the standard deviation is reported to show a good match between resiliency functions and experimental calculated data.

It is also possible to calculate the resilience factor corresponding to a specific period of time. In this study, 400 days was chosen as reference period of time because by that point, all the systems had recovered their full functionality. As can be seen in Table 2.28, in all cases consider, the base-isolated systems had higher resilience than the fixed-base systems. In particular, the resilience factor for isolated system considering the serial repair strategy is more than double with respect to the corresponding fixed-base system.

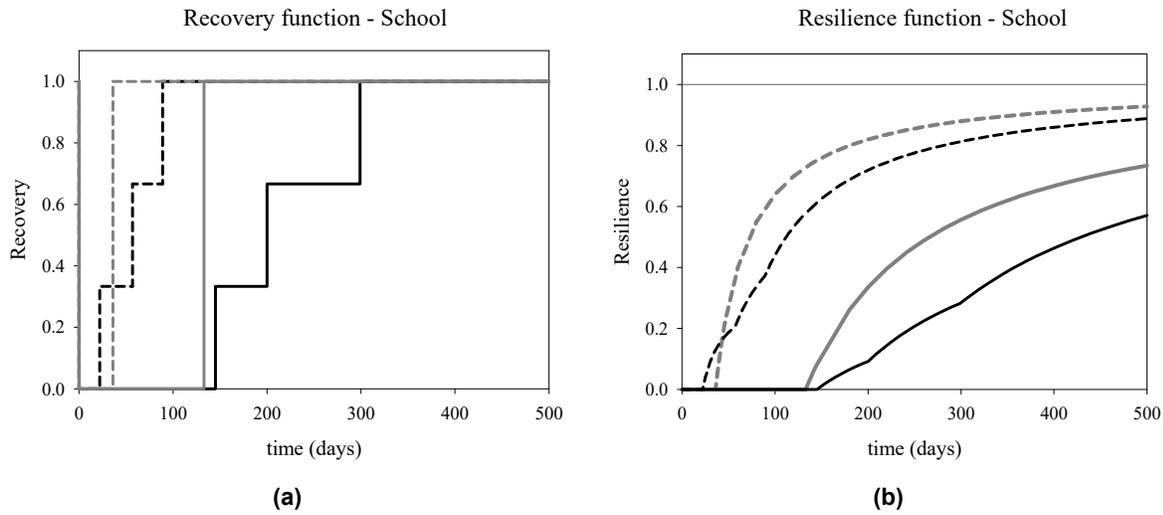

**Figure 2.47** School facility: recovery (a) and resilience functions (b) of all considered systems at 2% in 50-year hazard level.

**Table 2.25** Resilience functions for healthcare buildings.

| Structural system | Resilience function | $R^2$ |
|---|---|---|
| HP-SMRF PARALLEL | $R=1\cdot10^{-8}t^3-2\cdot10^{-5}\,t^2+0.09\,t-1.076$ | 0.998 |
| HP-SMRF SERIAL | $R=-2\cdot10^{-9}\,t^3+2\cdot10^{-6}\,t^2+0.000\,t-0.131$ | 0.998 |
| BI-IMRF PARALLEL | $R=2\cdot10^{-8}\,t^3-2\cdot10^{-5}\,t^2+0.008\,t-0.238$ | 0.971 |
| BI-IMRF SERIAL | $R=8\cdot10^{-9}\,t^3-1\cdot10^{-5}\,t^2+0.006\,t-0.311$ | 0.997 |

**Table 2.26** Resilience factors at $t$ = 400 days for healthcare buildings.

| Resilience factors | | | |
|---|---|---|---|
| Parallel repair strategy | | Serial repair strategy | |
| HP-SMRF | BI-IMRF | HP-SMRF | BI-IMRF |
| 0.6000 | 0.8700 | 0.3700 | 0.7700 |

**Table 2.27** Resilience functions for the school facility.

| Structural system | Resilience function | $R^2$ |
|---|---|---|
| HP-SMRF PARALLEL | $R=1\cdot10^{-8}t^3-2\cdot10^{-5}\,t^2+0.010\,t-0.996$ | 0.998 |
| HP-SMRF SERIAL | $R=-3\cdot10^{-9}\,t^3+2\cdot10^{-6}\,t^2+0.001\,t-0.269$ | 0.999 |
| BI-IMRF PARALLEL | $R=2\cdot10^{-8}\,t^3-3\cdot10^{-5}\,t^2+0.009\,t-0.2$ | 0.972 |
| BI-IMRF SERIAL | $R=1\cdot10^{-8}\,t^3-2\cdot10^{-5}\,t^2+0.007\,t-0.116$ | 0.996 |

Table 2.28    Resilience factors at *t* = 400 days for school buildings.

| School | | | |
|---|---|---|---|
| Parallel repair strategy | | Serial repair strategy | |
| HP-SMRF | BI-IMRF | HP-SMRF | BI-IMRF |
| 0.6675 | 0.9100 | 0.4633 | 0.8600 |

## 2.9    NONLINEAR ANALYSIS USING ATC 63 FRAGILITY CURVES FOR MEDICAL EQUIPMENT

In order to better evaluate the healthcare facility building performance after an earthquake, further analyses were made. The aim of these new analyses is to move on from the simplified analysis with a single fragility curve representing all of the medical equipment as proposed by Yao and Tu [2012]. For this reason, new fragility curves modeling the equipment inventory from a typical healthcare facility were used to replace the single curve previously used. These new fragilities were derived as part of the Applied Technology Council project (ATC 63) in which a methodology was developed to assess seismic design provisions for building systems based on the FEMA P-58 methodology. The objective of the study is to evaluate how ASCE 7-10 design provisions can affect seismic performance by impacting the building losses.

The FEMA P-58 methodology is a powerful tool used to evaluate building-specific seismic performance in terms that can directly be used by decision-makers. The ATC 63-2/3 projects [FEMA 2013] provides benchmark data for various types of buildings analyzed using the FEMA P-58 methodology. For this reason, 18 representative building archetypes designed according to the requirements of ASCE 7-10 were investigated. Different levels of ground shaking as well as different geometries (two or five stories buildings), construction materials (wood, steel, and concrete), and occupancies (residence, office, emergency operator center, medical office building, and an acute care hospital) were considered to cover a wide range of cases. Building models for each archetype were built and implemented in the PACT software system. As part of the project, a list of possible components that can cause losses for all archetypes was performed and the DS at which the building would be "barely functional" determined.

The buildings components were chosen using the Normative Quantity Estimation Tool. Although the FEMA P-58 fragility library is extensive, it is not comprehensive of all the fragility specifications needed to model the archetypes. As part of the project results, the fragility database was expanded to include new medical equipment fragility curves. An inventory of major equipment items was generated based on a representative 70,000 ft$^2$ facility and data obtained from a major healthcare provider. Thirty-nine additional hospital items fragility curves were developed to represent most major medical components. For each item, different curves were determined in terms of the risk category (II or IV) and the component height in the building divided by the total height $x/h$.

The medical items included in the list are both fixed and mobile. For mobile items, the possibility to slide versus tip was analyzed based on the dimensions of the unit. A slick floor surface was used, which implies that all items were prone to slipping as opposed to tipping over.

Mobile items are velocity sensitive, and two DSs were derived to represent the possible consequences. For DS1, the demand parameter is the floor velocity capable of making the item slide 6 in. for a total of, 22 in./sec. DS1 has minor impact on the unit, resulting in a still functional element. For DS2, the demand parameter is the floor velocity capable of making the item slide 12 in. for a total of, 31 in./sec. The DS has significant impact, leading to a non-functional element. Fixed components are acceleration sensitive and their anchorage is divided into categories: Overhead Mount and Floor Mount. These elements have just one DS, and they are divided in special elements; that equipment that is required to maintain basic post-earthquake functionality, and that equipment that has experienced anchorage failure and is damaged beyond repair.

As part of that study, the number of units and floor location for a 2- and 5-story building was given. This study considered a 3-story building; a graphic showing where the departments are assumed to be located is shown in Figure 2.48. A list of the most typical items housed in each department of the healthcare facility is given in Table 2.29.

Since the mobile items are velocity sensitive, it was necessary to input in PACT the absolute velocities at each floor and for each hazard level. Absolute velocities were derived by adding the relative velocities to the velocities generated by the ground motion. The first contribution was found by running OpenSees and recording the relative velocities at each node of the considered structure. To find the ground-motion velocities, the software Seismosignal was used. By inputting the accelerograms from PEER Strong Motion Database, velocity time series that were found for all considered hazard levels. The relative and ground-motion contributions were added to determine the absolute velocities for both fixed base and base isolated structures; see Figure 2.49.

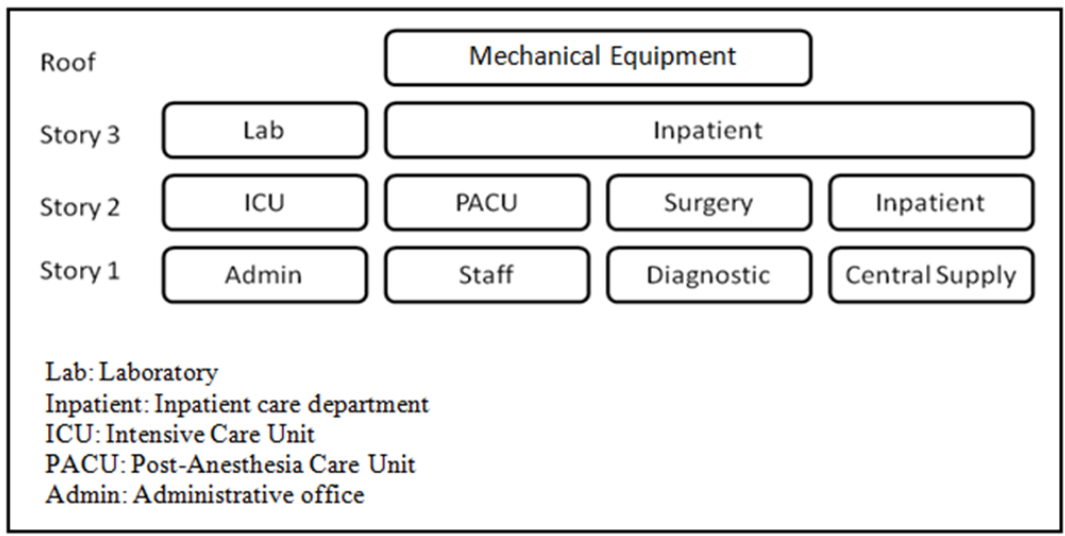

**Figure 2.48**   Location of various departments in the healthcare facility.

**Table 2.29     Medical equipment components grouped by department.**

| Surgery | | Diagnostic | ICU/PACU | Inpatient | Lab | Central supply |
|---|---|---|---|---|---|---|
| Anesthesia machine | Surgical table | Catheter cabinet | Patient bed | Patient bed | Freezer | Shelving machine |
| Anesthesia boom | Sterilizer | Breathing circuit dryer | Refrigerator | Refrigerator | Laminar flow hood | Washer / disinfector |
| Hypothermia system | Warming cabinet | Endoscope system | Headwall | Refrigerator under counter | Chemistry analyzer | |
| Ice slushier | Dual surgical light | Medstation | Ice dispenser | Headwall | Cryostat | |
| Equipment boom | Blood recovery system | Ultrasound unit | | Ice dispenser | Shelving system | |
| Balloon pump | Light system | Examination light | | | Refrigerator-lab | |
| Cath lab system | | Decontamination washer | | | Bio safety hood | |

As can be seen in the figures above, the velocity trend is similar to the acceleration trend, and a significant decrease is noticeable in the base-isolated system compared to the fixed-base structure. The previous loss analysis only considered a single fragility curve as representing the whole equipment. In the PACT model, all new components were added at the floor corresponding to their location within the building. The quantities were chosen taking into account the information provided in the ATC-63 study and were consistent with the department's position inside the building. The risk category was assumed to be IV, which is the category for healthcare facilities. As an example, see Table 2.30 for a listing of the medical equipment components at the first floor. After implementing the PACT model including fragility curves from ATC 63 project, a more detailed loss assessment was done that reflected the impact on repair costs. Total repair costs for the fixed base HP-SMRF structure are shown in Table 2.31.

As can be seen in Table 2.31, total repair costs found that although the two different procedures are comparable at the lower hazard levels, there is a huge difference at the 2%/50-years hazard level. To understand the reason for this difference, the partial contribution to total repair costs of structural/non-structural elements and contents must be determined.

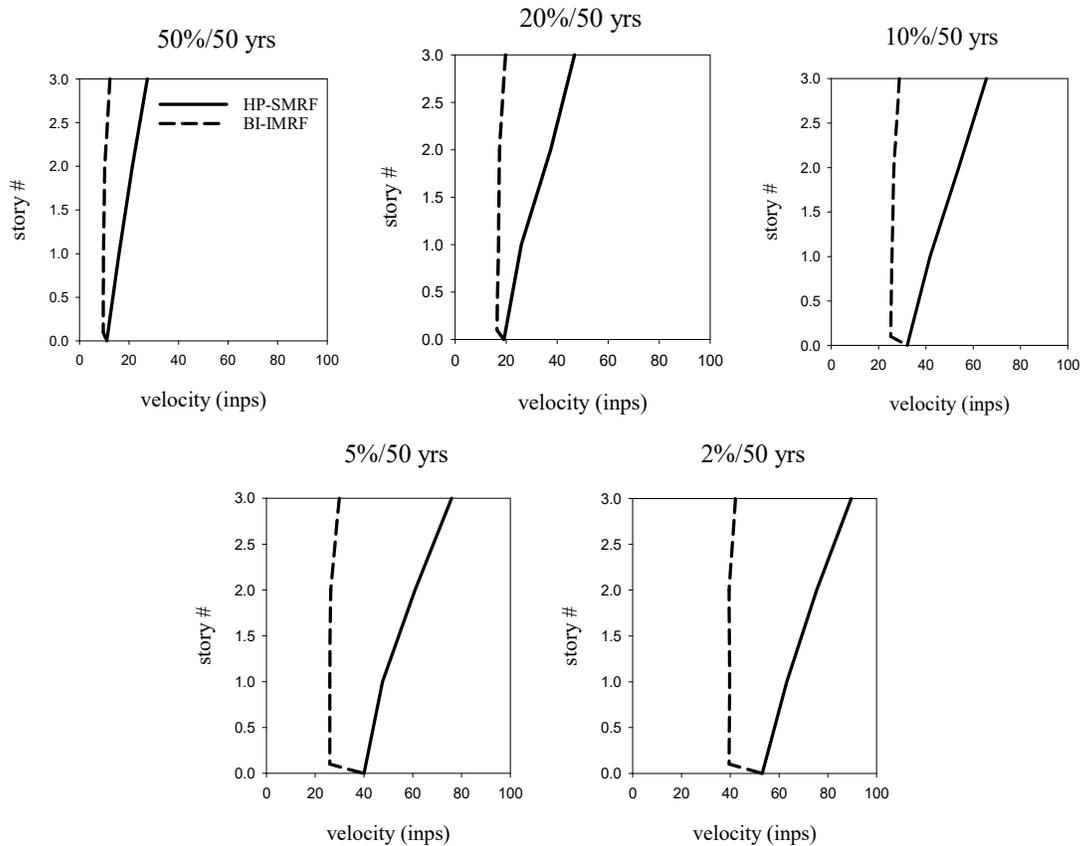

**Figure 2.49** Median absolute floor velocities of the HP-SMRF and the BI-IMRF on TFPBs for five hazard levels: 2%, 5%, 10%, 20%, and 50% probabilities of exceedance in 50 years.

A comparison of the two approaches shows a good match between structural/non-structural repair costs because those elements are the same in both two analyses. The slight differences are due to the probabilistic analysis in PACT. For the lower hazard levels of exceedance, the content repair cost associated with the new curves is larger than costs associated with medical equipment. At 5% and 2%, the costs related to the new curves increases. This large increase in content repair cost using Yao and Tu [2012] approach is that input accelerations are greater than accelerations related to the third DS, whose repair cost is $16,663,662. This means there is a high probability of reaching content costs close to the DS3-level costs. In the case of the ATC 63 approach, at 2% hazard level, input velocities and accelerations are larger than the worst possible DS. For this reason, a limit case in which all the curves reach their maximum damage was considered, and the associated repair cost is $6,072,882, which is close to the value obtained from analysis at 2%/50-years hazard level ($4,897,670). With the first simplified method, content repair costs are underestimated if compared with the more exact analysis at the lower hazard levels while the content cost at 2%/50 years is significantly overestimated. The same comparison was done for the base-isolated building. The total repair costs considering the two approaches are shown in Table 2.33.

For the BI-IMRF system, repair costs calculated using the two methods at each hazard level are comparable. A slight increase in costs at all hazards can be observed using the Yao/Tu approach but the difference is never more than 5%. The contribution of the model elements is shown in Figure 2.51 and Table 2.35.

For the isolated structure, accelerations are such that the third DS of Yao/Tu medical equipment fragility curve is never reached. Using the ATC 63 approach, there is a low probability to incur damage since input accelerations are always smaller than the mean acceleration related to the DS, because the highest acceleration at 2% hazard level is equal to 0.2$g$ while the mean DS1 acceleration is 0.92$g$. For this reason, medical equipment that is acceleration sensitive has a low probability of being damaged. The same can be said for velocity sensitive components since maximum velocities exciting the system are around 50 in./sec, and the DS1 mean velocity is 266 in./sec. Therefore, content evaluated using the ATC 63 method results in a minimal impact on repair costs in contrast to the simplified approach.

Table 2.30    Medical equipment components at floor 1.

| PACT ID | Component type | EDP | Unit | Number of non dir. FCs |
|---|---|---|---|---|
| E1028.003a | Catheter Cabinet, unanchored laterally | Peak Floor Velocity | Ea. | 1 |
| E1028.021a | Endoscopy System, unanchored laterally | Peak Floor Velocity | Ea. | 2 |
| E1028.023a | Ultrasound Unit, unanchored laterally | Peak Floor Velocity | Ea. | 2 |
| E1028.108c | Washer/Disinfector - Risk Cat IV - x/h=0 to 0.2 - anchorage fragility | Acceleration | Ea. | 1 |
| E1028.205b | Decontamination Washer - Risk Cat IV - x/h=0 to 0.4 anchorage fragility | Acceleration | Ea. | 1 |
| E1028.212b | Breathing Circuit Dryer - Risk Cat IV - x/h=0 to 0.4 anchorage fragility | Acceleration | Ea. | 1 |
| E1028.342a | Medstation- Risk Cat IV - x/h=0 to 0.4 anchorage fragility | Acceleration | Ea. | 1 |
| E1028.502a | Storage Shelving - Risk Cat IV - x/h=0 to 0.4 anchorage fragility | Acceleration | Ea. | 7 |

Table 2.31    Comparison between the two approaches of median repair costs at all hazard levels for fixed-base system.

| | Total repair cost ($) | | | | |
|---|---|---|---|---|---|
| | 50%/50yrs | 20%/50yrs | 10%/50yrs | 5%/50yrs | 2%/50yrs |
| HP-SMRF (Yao and Tu, 2012) | 1,551,724 | 3,723,077 | 6,333,333 | 9,500,000 | 19,333,333 |
| HP-SMRF (ATC 63, 2012) | 1,481,123 | 4,990,779 | 7,489,839 | 8,746,201 | 10,466,415 |

**Table 2.32** Structural, non-structural and content contribution to repair cost for fixed-base system.

| | Repair costs ($) | | | | | |
|---|---|---|---|---|---|---|
| | HP-SMRF [Yao and Tu 2012] | | | HP-SMRF [ATC 62 2012] | | |
| | Structural | Non-structural | Content | Structural | Non-structural | Content |
| 50%/50 yrs | 0 | 1,131,536 | 417,189 | 0 | 1,130,712 | 350,411 |
| 20%/50 yrs | 136,775 | 2,662,720 | 922,603 | 135,393 | 2,565,367 | 2,290,019 |
| 10%/50 yrs | 303,006 | 3,105,441 | 2,922,484 | 381,825 | 3,121,835 | 3,986,179 |
| 5%/50 yrs | 804,520 | 4,229,748 | 4,553,780 | 699,581 | 4,024,186 | 4,022,434 |
| 2% /50 yrs | 1,140,311 | 4,413,423 | 13,830,020 | 1,132,779 | 4,435,966 | 4,897,670 |

**Table 2.33** Comparison between the two approaches of median repair costs at all hazard levels for base-isolation system.

| | Total repair cost ($) | | | | |
|---|---|---|---|---|---|
| | 50%/50yrs | 20%/50yrs | 10%/50yrs | 5%/50yrs | 2%/50yrs |
| BI-IMRF (Yao and Tu, 2012) | 370,656 | 475,229 | 796,178 | 1,135,762 | 2,852,941 |
| BI-IMRF (ATC 63, 2012) | 369,137 | 471,812 | 794,391 | 1,102,934 | 2,701,409 |

**Table 2.34** Structural, non-structural and content contribution to repair cost for base-isolation system.

| | Repair cost ($) | | | | | |
|---|---|---|---|---|---|---|
| | BI-IMRF [Yao and Tu 2012] | | | BI-IMRF [ATC 63 2012] | | |
| | Structural | Non-structural | Content | Structural | Non-structural | Content |
| 50%/50yrs | 0 | 369,057 | 1,912 | 0 | 369,137 | 0 |
| 20%/50yrs | 0 | 472,315 | 3,099 | 0 | 471,812 | 0 |
| 10%/50yrs | 0 | 794,419 | 1,915 | 0 | 794,111 | 280 |
| 5%/50yrs | 11,702 | 1,089,125 | 38,883 | 10,055 | 1,091,276 | 1,604 |
| 2%/50yrs | 386,178 | 2,380,857 | 80,162 | 334,590 | 2,362,515 | 4,304 |

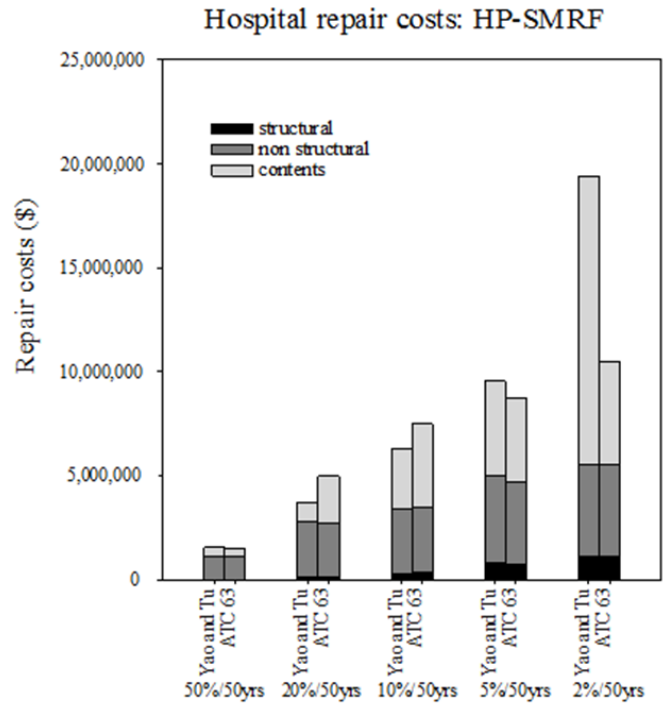

Figure 2.50   Median repair costs for HP-SMRF with the Yao and Tu and ATC 63 methods for five hazard levels: 2%, 5%, 10%, 20%, and 50% probabilities of exceedance in 50 years.

A more detailed analysis was conducted to understand the relative contribution of single components to repair costs. Table 2.35 shows that for the HP-SMRF system, repair costs associated with medical equipment will be high: their contribution is almost 44% of the total. The Cath lab system will sustain the greatest damage, with the median acceleration at the second floor being 1.42$g$ and DS acceleration is 0.92$g$. For the base-isolated system, medical components do not contribute significantly to repair costs. The single-components contribution is shown in Figure 2.52; the bar chart clearly shows that medical items do not have any impact on the repair costs for the base-isolated system.

Because the loss ratio does not take into account the repair costs of the contents, the results obtained using Yao/Tu approach versus the ATC 63 method are the same. Medical equipment fragility curves developed with both methods do not have repair times associated with damage. For this reason, indirect losses associated with business interruption remain the same with respect to the Yao/Tu approach. Since direct losses are different, total annualized losses can be recalculated to see the impact of the new set of fragilities integrated at each hazard level and weighted by the probability that the hazard was exceeded. Table 2.36 shows there is a slight difference between the total annualized losses due to the lack of information about repair times made it impossible to properly evaluate indirect losses. As a result of this lack of data, the recovery and resilience functions are equal in both cases. It is suggested that additional studies be conducted to associate repair times and replacement costs to medical fragility curves to better estimate indirect losses contribution related to healthcare facilities equipment.

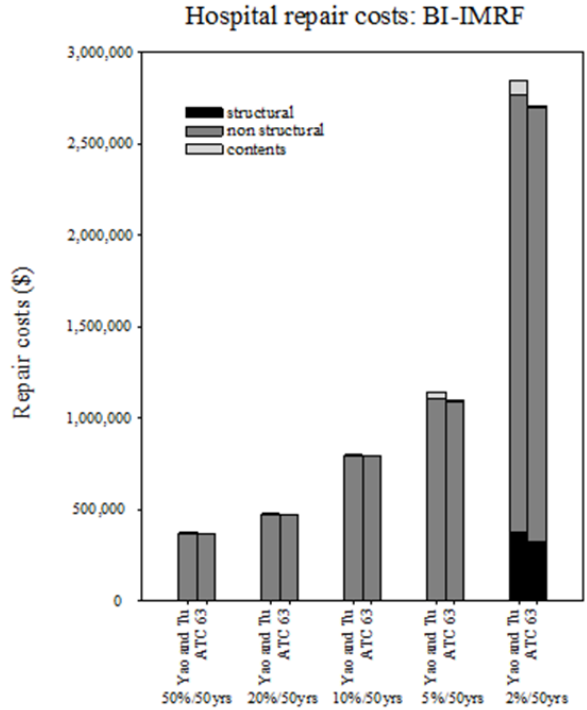

**Figure 2.51** Median repair costs for BI-IMRF comparing the Yao/Tu and ATC 63 methods for five hazard levels: 2%, 5%, 10%, 20%, and 50% probabilities of exceedance in 50 years.

**Table 2.35** Repair cost contribution at 2% in 50 years hazard level using the ATC 63 fragility curves.

|  | HP-SMRF | BI-IMRF |
|---|---|---|
| Cath Lab System | 31.93% | 0% |
| Wall Partition | 16.15% | 60.82% |
| Precast Concrete Panels | 11.87% | 30.43% |
| Suspended Ceiling | 8.04% | 0.02% |
| Headwall | 4.17% | 0% |
| Chemistry Analyzer | 4.16% | 0% |
| Post-Northridge RBS connections | 3.42% | 2.99% |
| Bolted shear tab gravity connections | 2.49% | 4.81% |
| HVAC Drops | 2.46% | 0% |
| Variable Air Volume (VAV) | 2.31% | 0% |
| Sterilizer | 1.71% | 0% |
| Washer/Disinfector | 1.28% | 0% |
| Concrete tile roof | 1.17% | 0% |
| TOTAL | 91.17% | 99.08% |

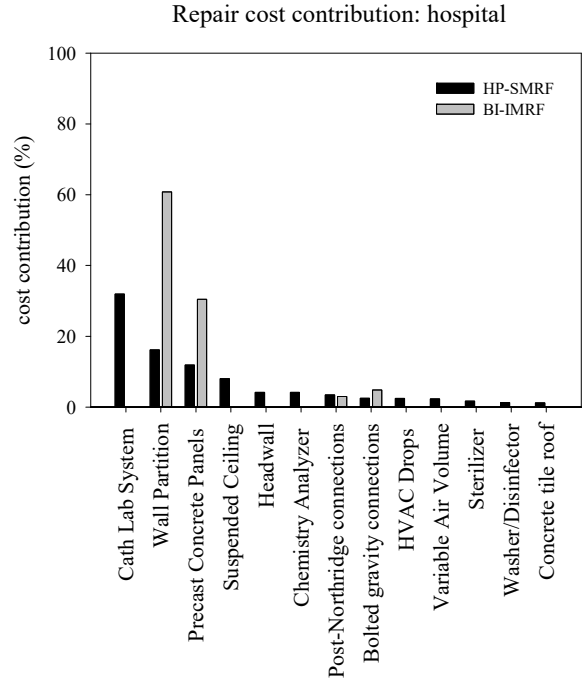

**Figure 2.52**  Repair cost contribution at 2% in 50 years hazard level using ATC 63 fragility curves.

**Table 2.36**  Comparison between total annualized losses calculated using the Yao/Tu and ATC 63 approaches.

| | **Total annualized losses ($)** | |
| --- | --- | --- |
| | **Yao and Tu** | **ATC 63** |
| HP-SMRF serial | 397,438 | 400,225 |
| HP-SMRF parallel | 277,309 | 280,096 |
| BI-IMRF serial | 107,148 | 106,435 |
| BI-IMRF parallel | 76,989 | 76,275 |

## 2.9.1  Sensitivity Analysis

The purpose of a sensitivity analysis is to study how the different sources of uncertainty in a model can affect the outputs. This type of analysis helps to investigate the robustness of the results and the correlation existing between inputs and outputs. Once the most significant sources of uncertainty are determined, it is easier to find possible errors in the model. After the sensitivity analysis is done, Tornado Diagrams can be used to plot the results, providing a way to identify those factors whose uncertainty has the largest impact. To build a Tornado Diagram, the possible source of uncertainty is modeled as a variable, and the impact of variation of that quantity is investigated with all the other quantities maintained at a baseline value. In this study, median EDP values of the fragility curves that provided the biggest contribution to repair costs were chosen as variable. The EDPs were varied by adding and decreasing by 30% the median values of a single curve at a time. Repair costs were collected for each case and at each hazard level. Tornado Diagrams were built considering the repair cost without any change to fragility curves as reference value (Figure 2.53 to Figure 2.57).

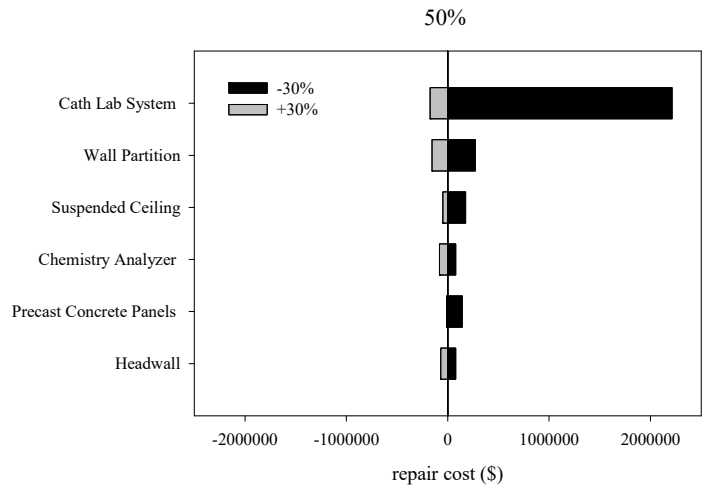

**Figure 2.53** Tornado diagram at 50% hazard level.

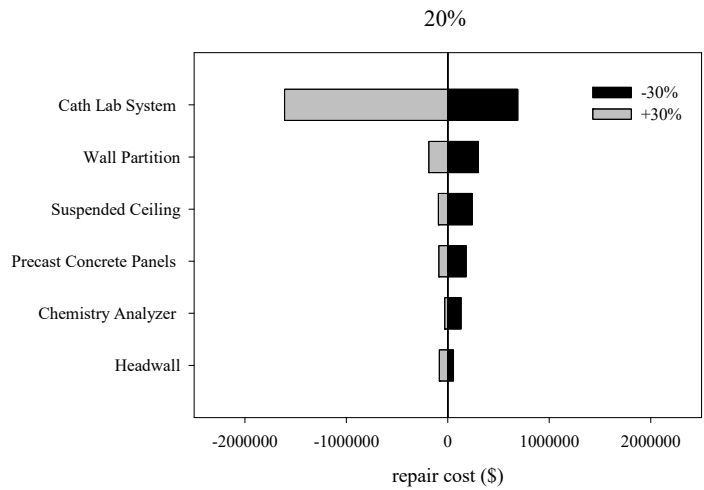

**Figure 2.54** Tornado diagram at 20% hazard level.

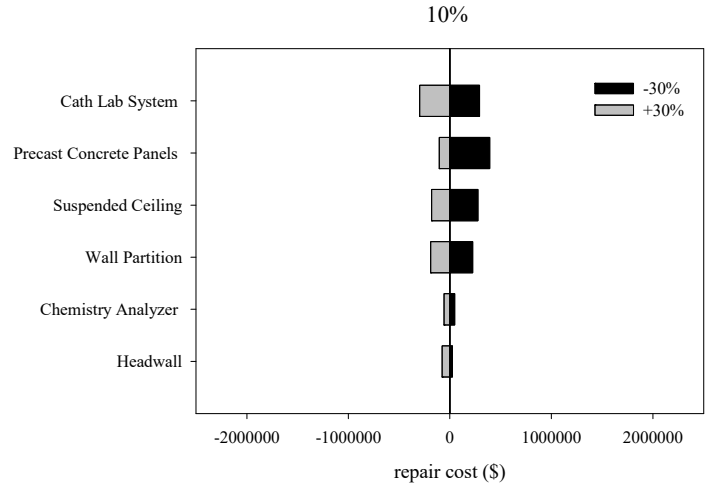

**Figure 2.55**   Tornado diagram at 10% hazard level.

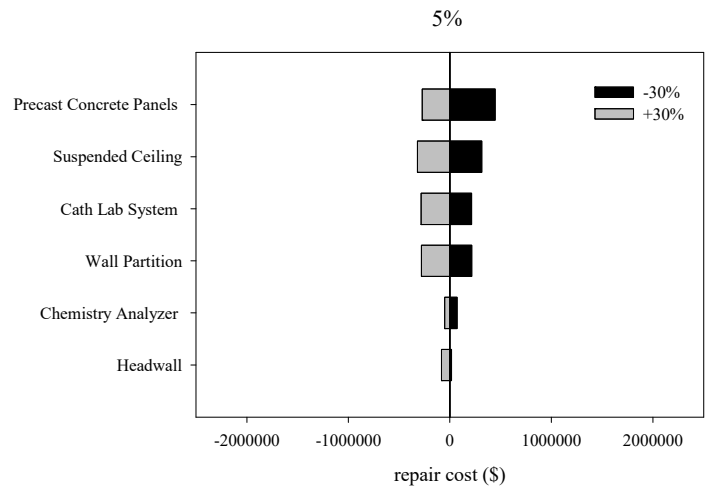

**Figure 2.56**   Tornado diagram at 5% hazard level.

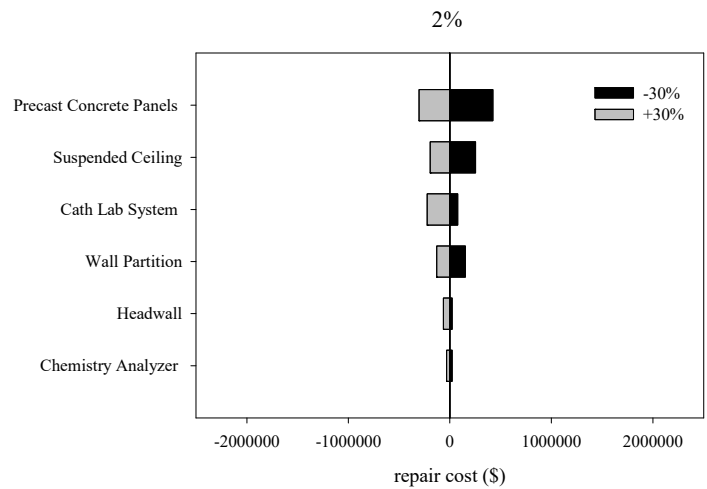

**Figure 2.57**   Tornado diagram at 2% hazard level.

As is known, repair costs are mostly influenced by non-structural components and contents. The charts show that Catheterization Laboratory (Cath Lab) systems are the main source of uncertainty at low-hazard levels; at higher levels of hazard precast concrete panels are the highest variability. The considered medical components have acceleration sensitive fragility curves whose median value is 0.92I; input acceleration at the 50% hazard level ranges from 0.26$g$ to 0.67$g$. The consequence of decreasing the median EDP by 30%–from 0.92$g$ to 0.64$g$–is that the increased probability of those curves being damaged at low-hazard levels. For all hazard levels; the Cath Lab systems had the highest impact on the uncertainty, implying that particular attention should be paid to identifying the curve parameters.

## 2.10 SIMPLIFIED ANALYSIS OF SCHOOL FACILITY

For buildings with less than 15 stories, it is common to follow the simplified loss analysis procedure of FEMA P-58. The simplified procedure uses linear mathematical structural models to estimate median values of response parameters. This study applied the simplified analysis to the three-story building to compare the results of nonlinear and linear analysis to determine structural response. The analysis was performed only for the HP-SMRF as nonlinear analysis has demonstrated wide excursions in plastic field for high-hazard levels.

The simplified approach can be used only if the building complies with certain pre-requisites, such as regular shape in plan and elevation, low or medium height, and limited nonlinear response. It is based on the equivalent lateral force method. Several assumptions were made: (1) the building was assumed to have independent translational response along each axes so two separate analysis are required for each orthogonal direction; (2) $P$-$\Delta$ effects were not included and, as consequence, story drift ratios are limited to 4%; and (3). story drift ratios did not have to exceed four times the yield drift ratio in order to avoid an excessive degradation of the structure.

### 2.10.1 Pseudo-Lateral Forces

The simplified analysis uses pseudo lateral forces to simulate the ground motion, which are. applied separately along each direction of the structure. Since the overall stiffness of the frames is the same in each direction, it adequate to apply these forces in one direction. The total load $V$ can be calculated using Equation (2.1):

$$V = C_1 C_2 S_a(T_1) W_1 \tag{2.1}$$

where $C_1$ is an adjustment factor for inelastic displacements; $C_2$ is an adjustment factor for cyclic degradation; $S_a(T_1)$ is the 5% damped spectral acceleration at the fundamental period of the building in the direction under consideration (from the uniform hazard spectrum for the selected level of ground shaking); and $W_1$ is the building's first modal weight in the direction under consideration but cannot taken as less than 80% of the total weight, $W$. $W_1$ can be calculated as $C_{sm}W$ where $C_{sm}$ is defined in ASCE/SEI 41-06 and is equal to 0.9. For the fundamental period of the building, the adjustment factors $C_1$ and $C_2$ are given by Equations (2.2) and (2.3):

$$C_1 = 1 + \frac{(S-1)}{(aT_1^2)} \tag{2.2}$$

$$C_2 = 1 + \frac{(S-1)^2}{(800T_1^2)} \tag{2.2}$$

where $T_1$ is the fundamental period of the building, and $a$ is a function of the soil site class (for class C, $a = 90$, ASCE/SEI 7-10). $S$ is a strength ratio and can be computed using Equation (2.3):

$$S = \frac{S_a(T_1)W}{V_y} \tag{2.3}$$

where $V_y$ is the estimated structure lateral yield strength at the first mode response. A simplified way to calculate $V_y$ is provided by ASCE 2010'; it is assumed that during the design phase it can be taken in the range given by Equation (2.4):

$$\frac{1.5 S_a(T)W}{R/I} \leq V_y \leq \frac{\Omega_0 S_a(T)W}{R/I} \tag{2.4}$$

where $S_a(T)$ is the structure's design spectral acceleration (DBE) at the fundamental period of the structure, as defined in ASCE/SEI 7-10; $R$ is the response modification coefficient per ASCE/SEI 7-10; $\Omega_0$ is the over-strength factor per ASCE/SEI 7-10; and $I$ is the importance factor per ASCE/SEI 7-10. $V_y$ was taken equal to 2278 kip, the mean value of the range above.

The values of the total load $V$ as well as the coefficients used to calculate it are shown in Table 2.37 for each hazard level. The distribution of the pseudo lateral force $V$ over the height of the building can be found by multiplying at each floor level the force for a vertical distribution factor given by Equation (2.5):

$$C_{vx} = \frac{w_x h_{(x-1)}^k}{\sum_{(j=2)}^{(N+1)} w_j h_{(j-1)}^k} \tag{2.5}$$

where $w_j$ is the lumped weight at floor level $j$; $h_{j-1}$ ($\underline{h}_{x-1}$) is the height above the effective base of the building to floor level $j$; and $k$ is equal to 2 for $T_1$ greater than 2.5 sec and equal to 1 for $T_1$ less or equal to 0.5 sec (linear interpolation used for intermediate periods). The distribution obtained is shown in Table 2.38 for each hazard level.

Table 2.37    Pseudo-lateral force $V$ at each hazard level

|  | $Sa(T_1)$ (g) | S | $C_1$ | $C_2$ | $V$ (kips) |
|---|---|---|---|---|---|
| 50%/50 yrs | 0.37 | 0.95 | 1 | 1 | 1941 |
| 20%/50 yrs | 0.63 | 1.63 | 1.02 | 1 | 3391 |
| 10%/50 yrs | 0.97 | 2.51 | 1.04 | 1.01 | 5362 |
| 5%/50 yrs | 1.21 | 3.11 | 1.05 | 1.01 | 6806 |
| 2%/50 yrs | 1.67 | 4.31 | 1.08 | 1.03 | 9845 |

**Table 2.38    Distribution of the pseudo-lateral force $V$ over the height.**

|  | 50%/50 yrs | 20%/50 yrs | 10%/50 yrs | 5%/50 yrs | 2%/50 yrs |
|---|---|---|---|---|---|
| $V_3$ (kips) | 987 | 1725 | 2727 | 3461 | 5007 |
| $V_2$ (kips) | 634 | 1108 | 1752 | 2224 | 3217 |
| $V_1$ (kips) | 319 | 558 | 883 | 1120 | 1620 |

## 2.10.2 Analysis: Model and Methods

### 2.10.2.1 Modeling and Assumptions

The building was modeled using SAP 2000. To better compare the results in term of linear and nonlinear analysis, the elastic model was built following some general assumptions. To be consistent with the previous model, a two-dimensional model was built using a linear analysis there wouldn't be complications related to a three-dimensional model. The moment frame chosen was the frame in direction 1 with the same distribution of masses previously explained. To account the stiffness of the gravity frame, a further column was added having the same dimensions of the leaning column used in OpenSees and connected to the frame by means of a rigid truss element. All nodes at the same floor have been constrained to move together using the internal constrain *diaphragm* of SAP 2000.

Since the linear static analysis maintains the use of a linear stress–strain relationship, the possibility of a nonlinear behavior of the system was not taken into account. As a consequence, the lumped plasticity of the beams or the spread of the plasticity along the columns were not included in the model.

### 2.10.2.2 Comparison of Structural Response

The pseudo-lateral forces were then applied to the frame, and the story drift ratios determined. The procedure proposed by FEMA P-58 applies a correction factor to the story drift ratios to account for inelastic action and higher mode effects. The correction is given by Equation (2.6):

$$\Delta_i^* = H_{\Delta i}(S, T_1; h_i, H)\Delta_i \tag{2.6}$$

where $\Delta_i$ is the story drift ratio obtained from the analysis, $\Delta_i^*$ is the corrected story drift ratio, and $H_{\Delta i}$ is the correction factor is given by Equation (2.7)

$$ln(H_{\Delta i}) = a_0 + a_1 T_1 + a_2 S + a_3\left(\frac{h_i}{H}\right) + a_4\left(\frac{h_i}{H}\right)^2 + a_5\left(\frac{h_i}{H}\right)^3 \tag{2.7}$$

where $H$ is the total height of the building, and $S$ is the strength ratio previously obtained. The values of coefficients $a_0$, $a_1$, $a_2$, $a_3$, $a_4$, and $a_5$ were derived from FEMA P-58 study for different structural systems that included the moment frame used in this study (Table 2.39).

Table 2.39   Correction factor for story drift ratio, floor velocity and floor accelerations for 2 to 9 story buildings [FEMA 2012].

|  | Frame type | $a_0$ | $a_1$ | $a_2$ | $a_3$ | $a_4$ | $a_5$ |
|---|---|---|---|---|---|---|---|
| Story drift ratio | Braced | 0.90 | -0.12 | 0.012 | -2.65 | 2.09 | 0 |
|  | Moment | 0.75 | -0.044 | -0.010 | -2.58 | 2.30 | 0 |
|  | Wall | 0.92 | -0.036 | -0.058 | -2.56 | 1.39 | 0 |
| Floor velocity | Braced | 0.15 | -0.10 | 0 | -0.408 | 0.47 | 0 |
|  | Moment | 0.025 | -0.068 | 0.032 | -0.53 | 0.54 | 0 |
|  | Wall | -0.033 | -0.085 | 0.055 | -0.52 | 0.47 | 0 |
| Floor acceleration | Braced | 0.66 | -0.27 | -0.089 | 0.075 | 0 | 0 |
|  | Moment | 0.66 | -0.25 | -0.080 | -0.039 | 0 | 0 |
|  | Wall | 0.66 | -0.15 | -0.084 | -0.26 | 0.57 | 0 |

Figure 2.58 shows the drift story ratios for the HP-SMRF using either linear or nonlinear analysis. The drift values of the nonlinear analysis are the median values of all the ground-motion records selected for a certain hazard level.

As can be seen from Figure 2.58, story drift ratios obtained from linear analysis are always larger than those obtained from nonlinear analysis, and the difference increases with the increase of the hazard level. This effect is on the conservative side since the linear static analysis is a simplified procedure to estimate response parameter, and the correction factor has the effect to further increase the gap between nonlinear and linear analysis, especially at the third floor. Figure 2.59 shows story drift ratios before and after the correction factor application. This correction is strongly dominated by the coefficients calibrated by FEMA P-58 and theoretically should reduce the difference between the linear and nonlinear analysis, rendering making the two procedures more or less comparable. With the moment frame consider in this case, however, the story drift ratios at the first and at the second floor are the same, significantly increasing the value at the third floor and reaching a peak of 5.4% at 2%/50-years hazard level; see Figure 2.59(e). This disqualifies use of the simplified procedure because the story drift exceeds 4% for the 2%/50-years hazard level.

Regarding peak floor accelerations, the simplified procedure of FEMA P-58 does not account for the results obtained from the linear model using SAP 2000. Indeed, peak floor accelerations are obtained from the PGA values by using the following Equation (2.8):

$$a_i^* = H_{ai}(S, T_1; h_i, H) \, pga_i \tag{2.8}$$

where $H_{ai}$ (S, $T_1$; $h_i$, H) is the acceleration correction factor given by Equation (2.9):

$$ln(H_{ai}) = a_0 + a_1 T_1 + a_2 S + a_3 \left(\frac{h_i}{H}\right) + a_4 \left(\frac{h_i}{H}\right)^2 + a_5 \left(\frac{h_i}{H}\right)^3 \tag{2.9}$$

The values of coefficients $a_0$, $a_1$, $a_2$, $a_3$, $a_4$, and $a_5$ are shown in Table 2.39. The results obtained following this procedure are shown in Figure 2.60 and compared with the values obtained from nonlinear analysis.

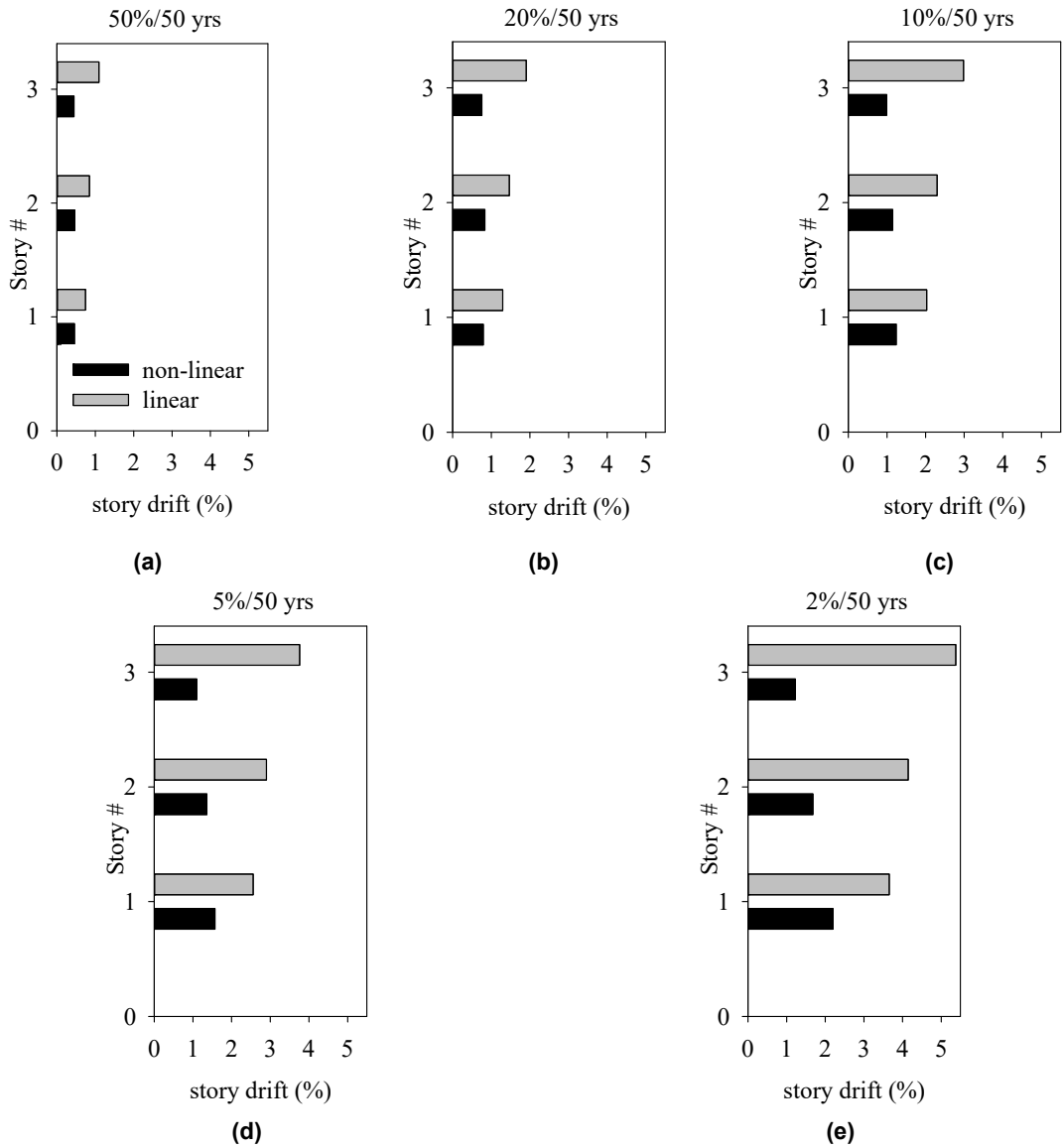

**Figure 2.58    Comparison of story drift ratios at each hazard level subjected to linear and nonlinear analysis.**

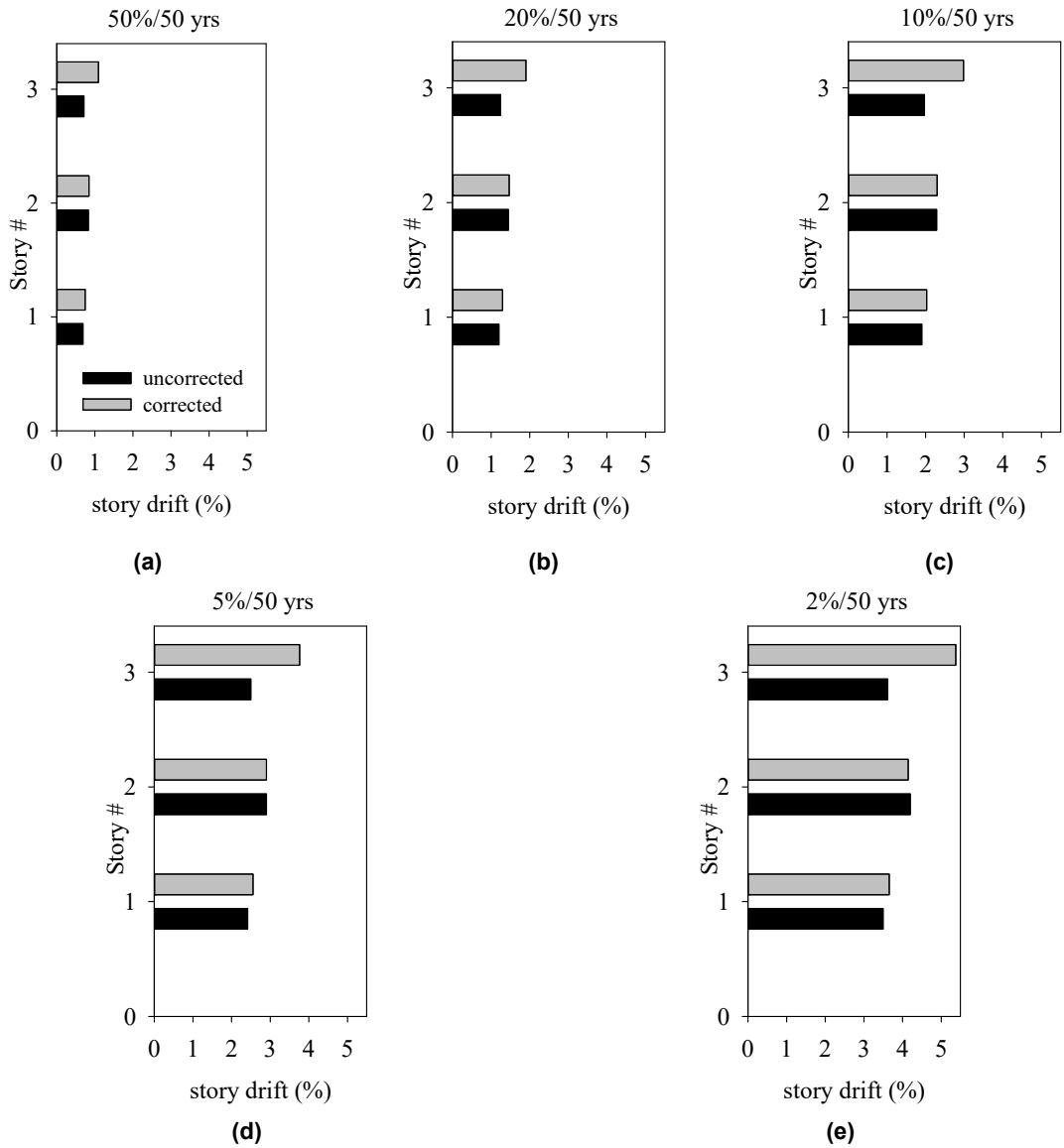

**Figure 2.59** Uncorrected and corrected story drift ratios at each hazard level using linear analysis.

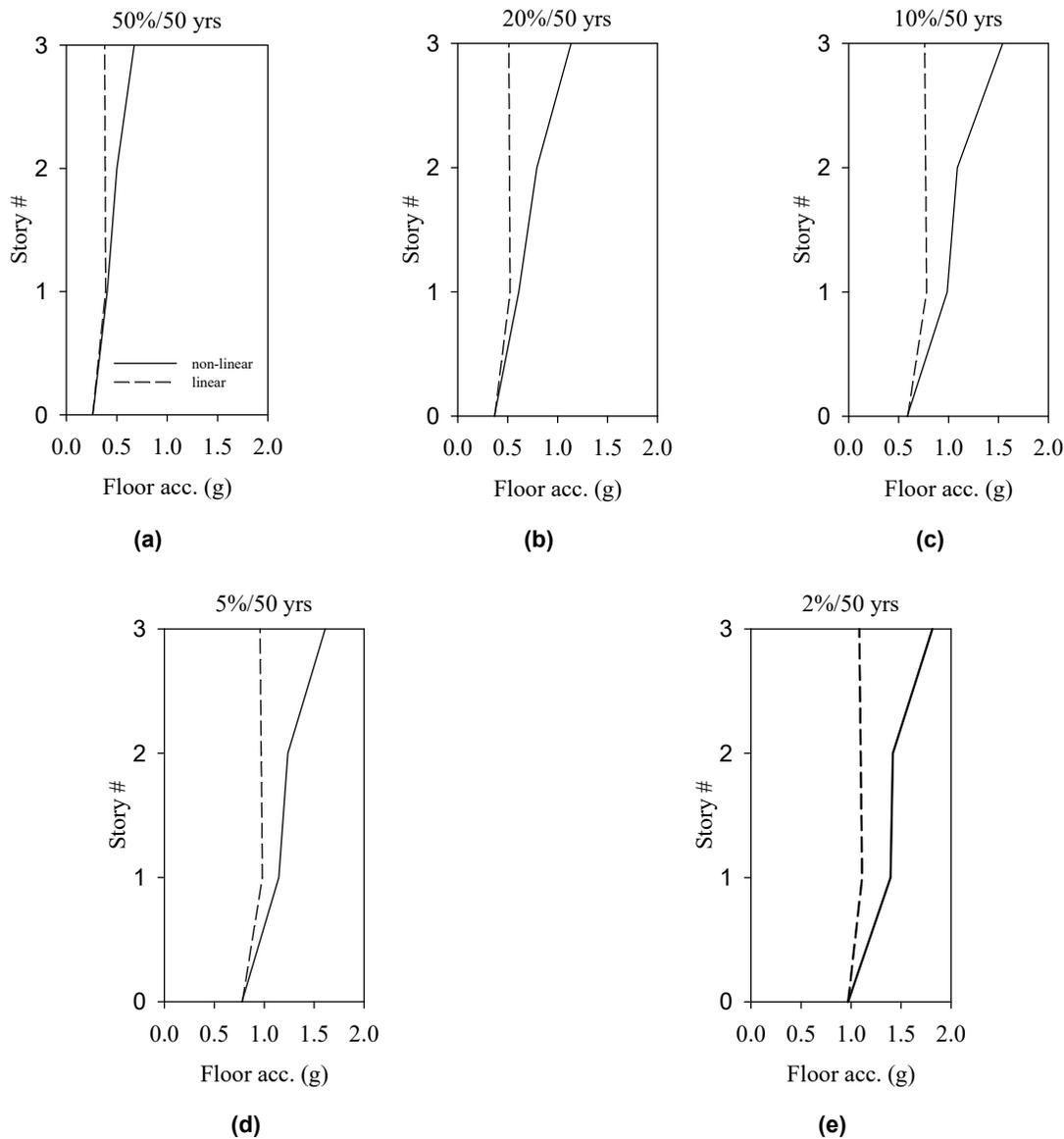

**Figure 2.60   Comparison of peak story accelerations at each hazard level using linear and nonlinear analysis.**

The results of peak story accelerations from linear analysis are always smaller than those obtained from the nonlinear analysis. Starting from the same value of PGA at the ground floor, the gap between nonlinear accelerations determined in OpenSees and linear accelerations calculated with this simplified procedure increases with every floor. This trend is even greater at high hazard levels and also for story drift ratios.

To sum up, the EDPs obtained by the simplified procedure of FEMA P-58 are comparable to those low floor levels found using the nonlinear analysis. At high floor levels, the discrepancy between the two procedures is marked. For the loss analysis, the greater story drift ratios will cause damage to structural components as all are sensitive to this EDP; however, lesser damage to the contents is expected because of the reduction in peak story accelerations.

Non-structural components are either drift or acceleration sensitive, making it difficult to predict their damage in advance.

### 2.10.2.3 Linear Time History Analysis

Applying pseudo-lateral forces to a structure is a convenient way to reproduce the effect of ground motion. To be on the safe side, the application of these forces over building height represents an overestimation of the seismic action itself. To better understand possible connections between these forces and the selected ground motions for a specific hazard level, multiple time history analysis were performed in SAP 2000 to determine the action of a single ground-motion record on the structure. Since the simplified procedure uses pseudo-lateral forces to predict story drift ratios only, the following comparisons will be made for this specific EDP.

Figure 2.61 compares the linear time history analysis obtained using SAP2000 and the nonlinear time history analysis using OpenSees at 5%/50-years hazard level for two ground-motion records. Neither analysis showed significant residual displacement. The comparison shows a good match between linear and nonlinear time history analysis. The story drift ratio at the first floor is slightly higher using the nonlinear analysis, proving that nonlinear behavior of the system at that floor causes an increase in story drifts in response to a small increase in intensity. Figure 2.64 presents an additional comparison for the ground motion used in Figure 2.62, resulting in a larger residual displacement. A high value of this EDP can be coupled with significant nonlinear behavior of the structure as shown in Figure 2.63, which plots the hysteretic cyclic response of the system subjected to this ground-motion record. The graph in Figure 2.63 computes the base shear versus the first-floor displacement. As can be seen, the structure starts to accumulate a huge plastic deformation immediately during the second reload cycle. As a consequence, story drift ratios computed by nonlinear analysis will be significantly larger than those obtained with the linear analysis (Figure 2.63) where it obvious that a big increase of displacement occurs in response to a small increase of force.

The time history analysis performed for the ground-motion 18 fault-normal also highlights the similarity between these story drift ratios obtained using the linear analysis and those obtained using the linear procedure of FEMA P-58. A comparison between the linear time history analysis of the strongest ground motion record and the simplified analysis of FEMA P-58 is shown in Figure 2.65 at all hazard levels reporting both uncorrected and corrected story drift ratios. Note that in Figure 2.65(a), the story drift ratios of the time history analysis are similar to those calculated applying pseudo-lateral forces *without* applying the correction factor. As shown in Figure 2.65(b), after the application of the correction factor, the match is still good for first and second floors, but the gap increases in the third floor. The simplified procedure was not used at 2%/50-years hazard level because at that level the story drift ratio exceeded 4%. In summary, there is little difference between the uncorrected/corrected story drift ratios obtained using the simplified procedure and the story drift ratios obtained from the linear time history analysis of the strongest ground-motion record at all hazard levels.

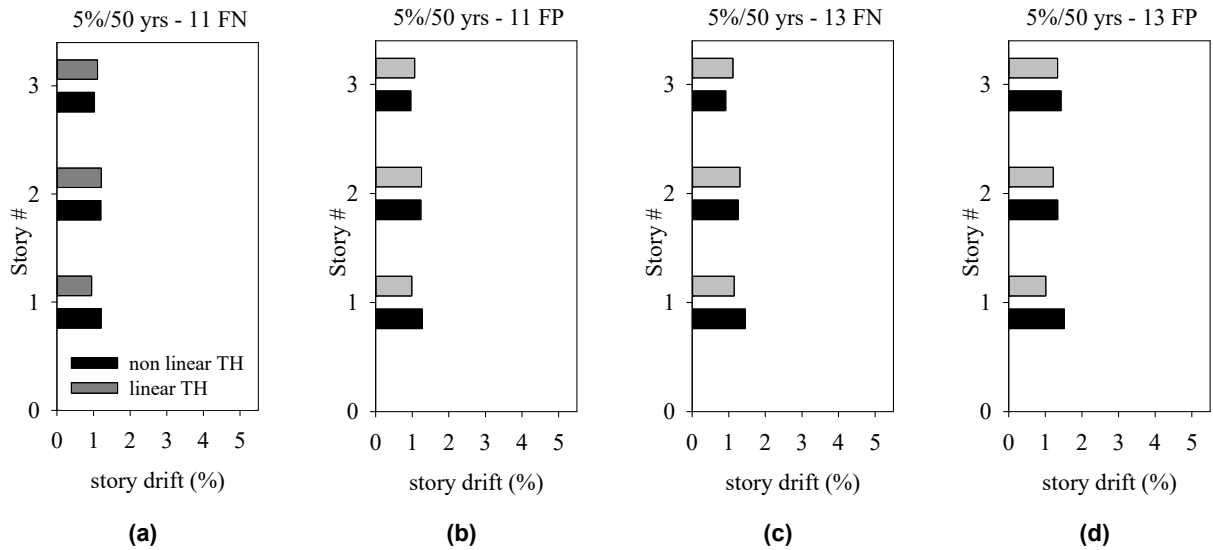

**Figure 2.61** Comparison between linear and nonlinear time history analysis for two ground-motion records: (a) 11 fault normal; (b) 11 fault parallel; (c) 13 fault normal; and (d) 13 fault parallel.

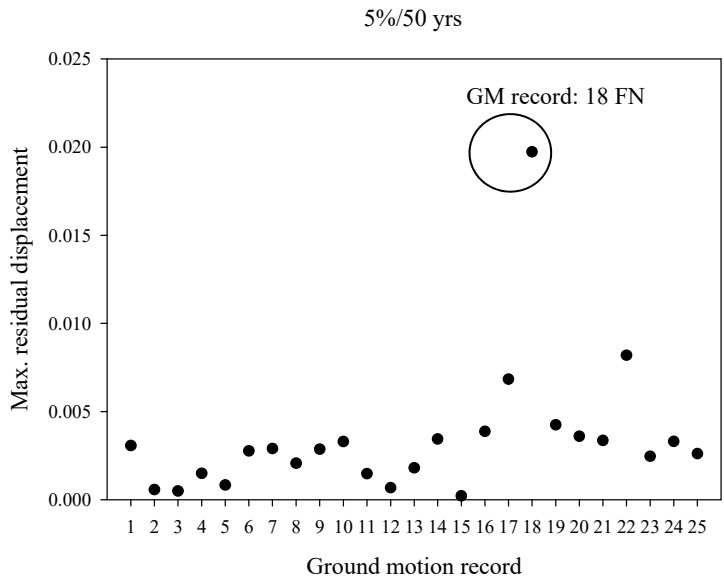

**Figure 2.62** Maximum residual displacement for all 25 ground motion records at 5%/50-years hazard level.

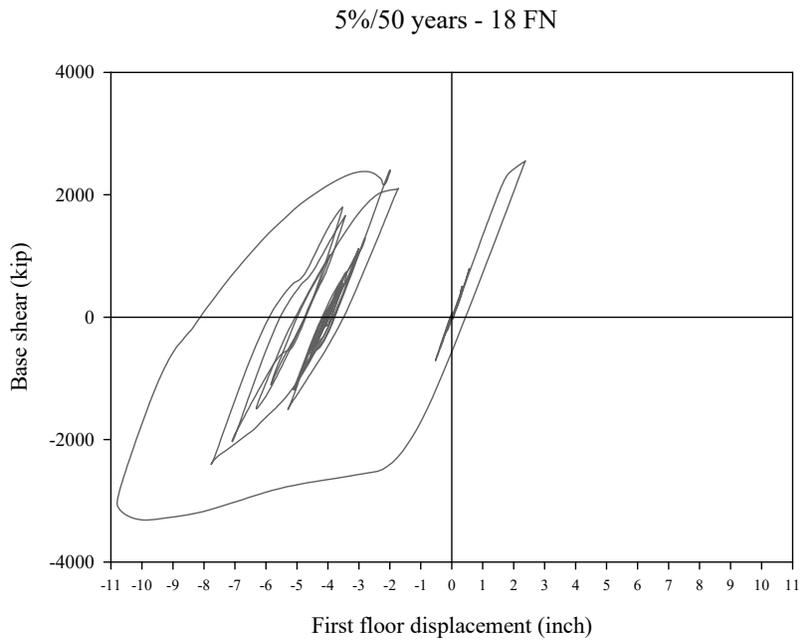

**Figure 2.63** Hysteretic cycle response of HP-SMRF subjected to the ground-motion record 18 fault normal.

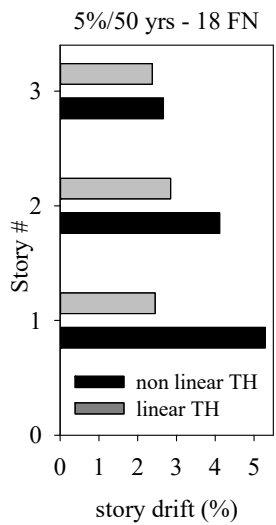

**Figure 2.64** Comparison between linear and nonlinear time history analysis for the ground-motion record 18 fault normal; note the increase in residual drift.

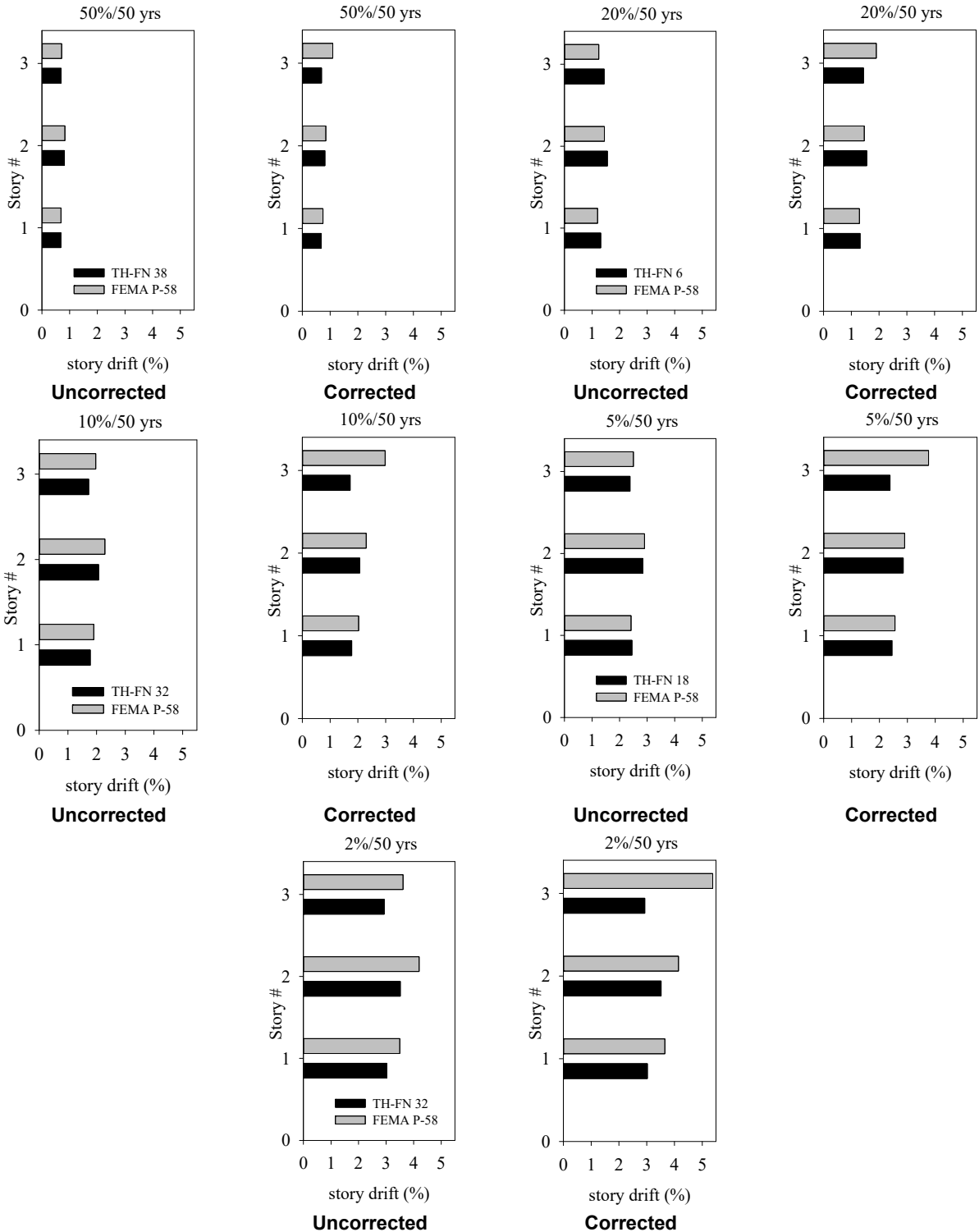

**Figure 2.65** Comparison between (a) the uncorrected and (b) the corrected story drift ratios obtained using the simplified procedure and the story drift ratios obtained from the linear time history analysis of the strongest ground motion record at all hazard levels.

### 2.10.3 Loss Analysis

#### *2.10.3.1 Repair Cost and Loss Ratio*

With the EDPs calculated above, a simplified loss analysis was performed using PACT. A single set of story drift ratios and floor accelerations were input for each direction at each hazard level. Dispersion values calculated in function of the fundamental period $T_1$ and the strength ratio $S$ (FEMA P-58) were added to account for the uncertainty of the distributions. The results in term of repair costs are shown in Table 2.40 by comparing linear and nonlinear analysis at each hazard level. The costs associated with repairing damage are consistently higher using the simplified procedure, ranging from 58% at 50%/50-years hazard level to a peak of 135% at 5%/50-years hazard level. Since the drift at 2%/50-years hazard level exceeds the limit, the estimated cost of repairs is not considered reliable.

    To identify the major contributors to the losses, the total repair cost was disaggregated into structural components, non-structural components, and contents; see Figure 2.66 and Table 2.41. As expected, the biggest difference between the two analyses is related to the repair cost of the structural components. Since these components are all drift-sensitive, the net increase of story drift ratios leads to an increase of this portion of the repair cost. Figure 2.66 shows that the cost of structural components rises up to a value of more than $6,000,000, which is six times bigger than the cost obtained using the nonlinear analysis. For the two higher hazard levels, it becomes the major contribution to the loss.

    The portion of structural repair cost given by each component cost is shown in Table 2.42 for the 2%/50-years hazard level. Bolted shear tab gravity connections provided most of contribution, comprising more than the half of the repair cost of all structural components. Indeed, by using the simplified analysis, drifts capable of damaging the connections are reached for high hazard levels. Since there are high quantities of connections inside the structure, they contributed significantly to the total repair cost. The repair costs of non-structural components increased as well, but were maintained in the same range as obtained in the previous analysis. A slight reduction in content repair costs is due to the reduction in accelerations, but it is not so pronounced since the repair cost of contents is only a small portion of the overall repair cost.

    Loss ratios are shown in Figure 2.67, which take into account only the repair costs of structural/non-structural components. The loss ratio at the 2%/50-years hazard level exceeds 40%, which means that the building is likely to be replaced rather than repaired. Given that this occurred at 2%/50-years hazard level, the results of the simplified procedure cannot considered accurate for this level of damage.

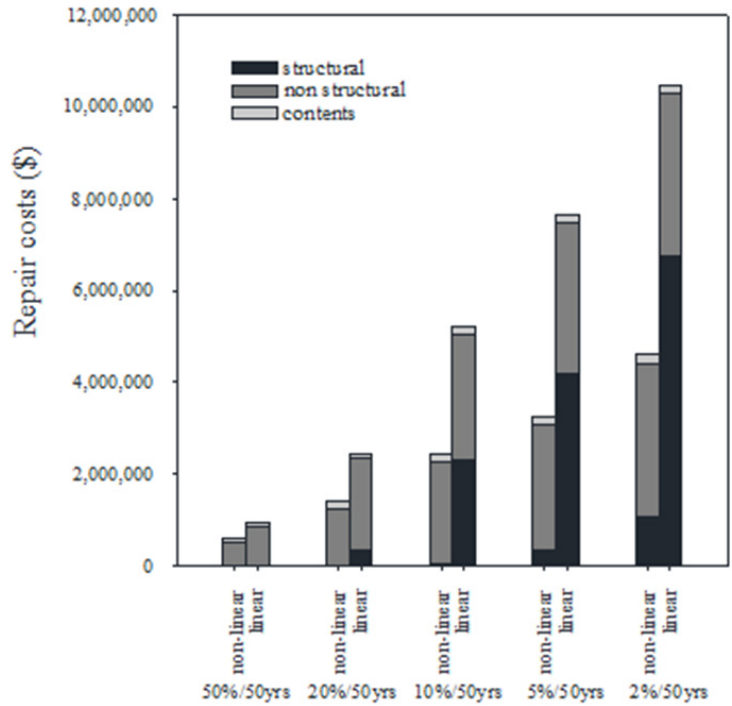

**Figure 2.66** Median repair costs with linear and nonlinear analysis for five hazard levels.

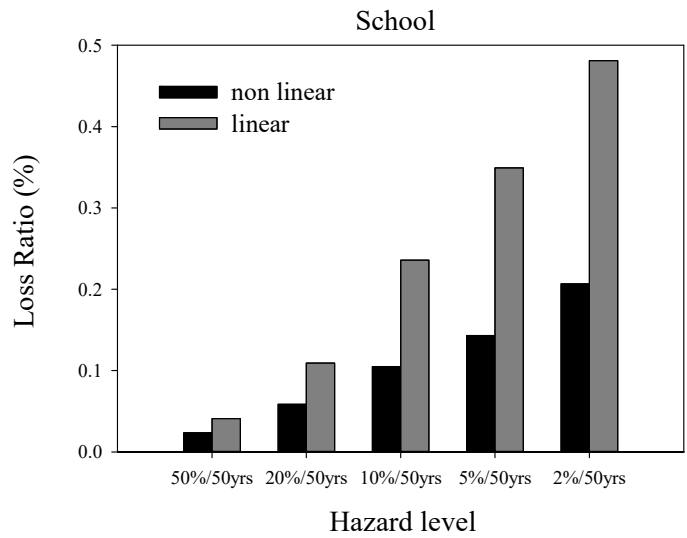

**Figure 2.67** Median loss ratio with linear and nonlinear analysis for five hazard levels.

Table 2.40   Median repair costs at all hazard levels with linear and non-linear analysis.

| | Repair costs [$] | | | | |
|---|---|---|---|---|---|
| | **50%/50yrs** | **20%/50yrs** | **10%/50yrs** | **5%/50yrs** | **2%/50yrs** |
| HP-SMRF non-linear | 592,967 | 1,415,217 | 2,415,217 | 3,242,424 | 4,605,555 |
| HP-SMRF linear | 938,333 | 2,450,000 | 5,200,000 | 7,633,333 | 10,471,428 |

Table 2.41   Structural, non-structural and content repair cost using either linear or nonlinear analysis.

| HP-SMRF nonlinear: Repair cost [$] | | | | HP-SMRF linear: Repair cost [$] | | | |
|---|---|---|---|---|---|---|---|
| | **Structural** | **Non-structural** | **Contents** | | **Structural** | **Nonstructural** | **Contents** |
| 50%/50yrs | 0 | 503,005 | 88,959 | 50%/50yrs | 17,026 | 855,051 | 63,898 |
| 20%/50yrs | 7,304 | 1,245,003 | 161,546 | 20%/50yrs | 344,206 | 1,991,583 | 116,465 |
| 10%/50yrs | 68,701 | 2,178,899 | 171,154 | 10%/50yrs | 2,320,736 | 2,718,631 | 157,328 |
| 5%/50yrs | 364,554 | 2,700,564 | 181,969 | 5%/50yrs | 4,196,128 | 3,296,791 | 160,783 |
| 2%/50yrs | 1,073,583 | 3,343,834 | 187,953 | 2%/50yrs | 6,751,189 | 3,549,207 | 179,499 |

Table 2.42   Percentage of structural repair cost given by each component at 2%/50-years hazard level.

| Structural components | Portion of structural repair cost |
|---|---|
| Bolted shear tab gravity connections | 59.3% |
| Post-Northridge connections | 29.3% |
| Steel Column Base Plates | 6% |
| Welded column splices | 5.4% |
| | 100% |

### 2.10.3.2 Repair time

Because the structural components are so damaged, repair times (Figure 2.68) are significantly bigger than those previously obtained using nonlinear analysis. Upper (serial) and lower (parallel) bounds of the repair times are both more or less doubled at all hazard levels. At the 2%/50-years hazard level, repairing the building takes more than two years if using the serial repair strategy.

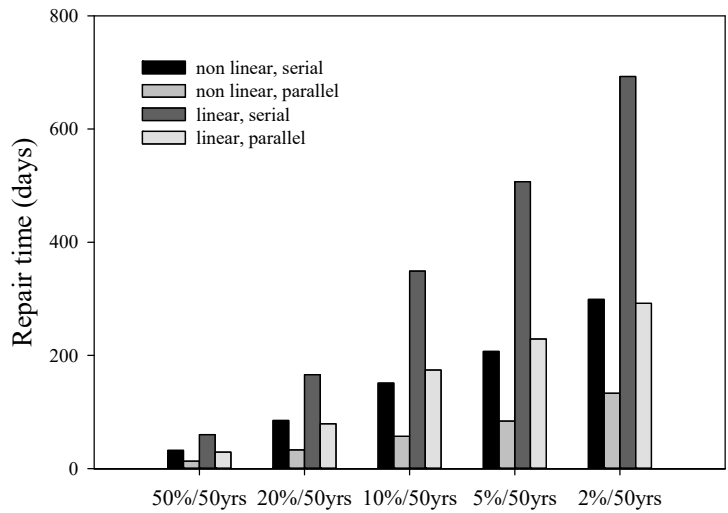

**Figure 2.68** Median repair times with linear and nonlinear analysis for five hazard levels, considering two repair strategies, parallel and serial.

### *2.10.3.3 Indirect losses*

Table 2.43 shows the annualized losses for serial and parallel repair strategy. Note that utilizing the simplified procedure of FEMA P-58 leads to a doubling of the total annualized losses because of damage incurred by structural components. Indirect losses are doubled as well, since repair of the building requires more time, thus increasing the time of the structure will be inoperable.

**Table 2.43** Annualized losses with linear and non-linear analysis.

|  | Serial repair strategy | | | Parallel repair strategy | | |
| --- | --- | --- | --- | --- | --- | --- |
|  | **Direct loss** | **Indirect loss** | **Total** | **Direct loss** | **Indirect loss** | **Total** |
| Nonlinear | 24,825 | 5,205 | 30,029 | 24,825 | 3,144 | 27,969 |
| Linear | 46,977 | 10,042 | 57,019 | 46,977 | 7,889 | 54,866 |

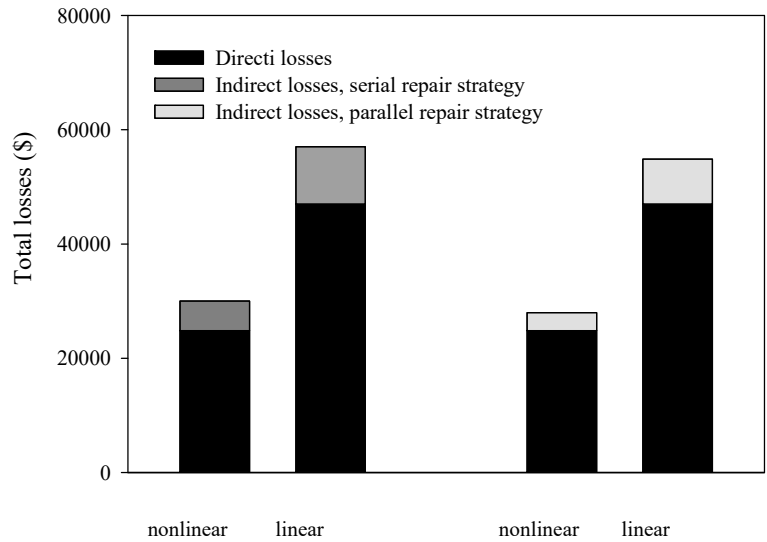

Figure 2.69    Annualized losses for all considered systems.

## 2.11  REMARKS AND CONCLUSIONS

Over the past two decades, performance-based earthquake evaluation (PBEE) has developed to a point where it can be effectively used to design structures. In particular, it can identify the contributions of different structural/nonstructural elements and contents to the total cost of ownership over the life span of a building, thereby enabling the designer to enhance a given design or choose alternative structural systems to improve performance and mitigate damage. This report applied the PBEE methodology to a healthcare facility and to a school building.

The healthcare facility was located in Northern California. The PBEE methodology was applied to evaluate its resiliency in the event of an earthquake. The loss assessment of the healthcare facility was performed by using two different approaches: (i) medical equipment was modeled using a simplified global fragility curve; and (ii) medical equipment was modeled by inserting fragility curves for each component. The first approach was used to compare the performance of a fixed-base and a base-isolated structure. The comparison shows a significant median damage-saving and repair-time reduction using the isolation system as opposed to the fixed-base system. The efficiency of the base-isolation system was widely proven, showing an 85% average reduction in repair costs and three to six times smaller repair times compared to the fixed-base system. This improved efficiency was verified using the resilience index, which had an increment of 63% when considering the base-isolation system with respect to the fixed-base system.

Next, the ATC 63 approach was used to model the medical equipment, and the results of this analysis was compared to the results obtained using global fragility curve. The comparison showed significant differences in the fixed-base system. For the isolated structure, the differences generated by using specific fragility curves for *content* is irrelevant. Although the ATC 63 method is more time-consuming than the simplified method, the results for the fixed-base system at higher hazard levels need to be taken into account because the first method overestimates the repair costs by up to almost three times at the 2%-50-years hazard level. For the base-isolated system, the differences between the two approaches are reduced due to the small contribution of contents to repair costs. Based on the results presented herein, the global fragility curve method could be used to evaluate the isolated structure, although it is not recommended.

For a school building located in Oakland, California, the base-isolated system provided significant median damage savings and reduction in repair time compared to the fixed-base system. This is due to the substantial reduction in accelerations, drifts, and residual drifts when considering the isolated building. The reduction in repair cost ranges from 66% to 82%, with an average of 76%. Such large reduction in cost of damage repairs of base-isolated systems comes primarily from preventing damage of structural components and minimizing the damage of non-structural components. Repair times are three to six times smaller for the isolated buildings relative to the fixed-base buildings. For the design-basis earthquake (a 10% probability of exceedance in 50 years) the repair time of the fixed-base building is expected to be in the range of 57 and 151 days compared to 12 and 27 days for the base-isolated building. Such a dramatic reduction in repair time implies significantly smaller downtime and higher resilience of the base-isolated buildings.

In addition, in application of the FEMA P-58 procedures, the pseudo-lateral forces applied on the structure can be assimilated to the envelope of the set of ground motions selected for a certain hazard level. This was found by comparing the story drift ratios obtained using the simplified procedure with the story drift ratios obtained by performing a linear time-history analysis of the strongest ground-motion record at each hazard level. For high levels of hazard, story drift ratios were much higher than those obtained by nonlinear analysis. At the 2%/50-years hazard level, they exceeded the 4% limit; therefore, the simplified analysis cannot be considered valid. The increase in repair cost ranged from 58% at the lower hazard level to 135% at 5%/50-years hazard level.

Comparing the results with the previous analysis, the damage to structural components was significant. Structural damage represented a higher portion of repair costs for probability of exceedance higher than 5% in 50 years. Such a dramatic increase in structural repair costs implies significantly larger downtimes and smaller resilience for the system studied using the simplified procedure. The procedure may be considered validated for four hazard levels out of the five even though it overestimates the loss. However, this study clearly demonstrates that nonlinear analysis is the optimum procedure to obtain a realistic estimation of damage, thus facilitating the decisions facing stakeholders in evaluating a building's structural integrity in the event of an earthquake.

# 3 Using Discrete Event Simulation to Evaluate Resilience of Emergency Departments

## 3.1 INTRODUCTION

Regardless of whether it is a large-scale area or a smaller-scale area, the ability of a geographical area to respond to an emergency situation is related to the proper functioning of its own infrastructure system. This becomes painfully evident when critical infrastructure systems fail during disasters, resulting in economic loss and human fatalities. Critical infrastructure comprising the assets, systems, and networks are necessary for the normal operation of cities, regions, and states; therefore, their incapacitation or destruction have a debilitating effect on security, public health, or any combination thereof. An effective way of measuring and analyzing how an infrastructure system will react in the event of a disaster is the resilience analysis.

The resilience of a complex system is defined as the capacity to prepare for and adapt to changing conditions, and restart operations after an extreme event. Improving the resilience of countries and communities means that structures will have the ability to resist, absorb, and accommodate to and recover from the effects of a hazard in a timely and efficient manner. Conceptually, resilience entails three interrelated concepts: (1) reduce the probability of failure; (2) limit and control consequences if a failure occurs; and (3) improve recovery time in the event of a failure [Chang 2009]. The emphasis on consequences and recovery regarding urban infrastructure systems in the event of a disaster comprises multiple disciplines. Improving resilience is not only a technical and structural problem: it includes societal and organizational dimensions. Therefore, the identification of all system vulnerabilities and the development of appropriate and innovative techniques to absorb shocks while maintaining functionality is the first step to ensure the protection of critical infrastructures and limit the effects of a disaster.

### 3.1.1 Resilience of Healthcare Facilities

The resilience of a healthcare facility can be defined as a hospital's ability to withstand the event, absorb the shock of disasters while addressing the surging medical demand in order to recover quickly to its original state, or develop strategies to adapt to a new one. Recent events have shown how systems (regions, communities, structures, etc.) are vulnerable to natural disasters of every type, such as earthquakes, hurricanes, explosions, and any other type of catastrophe regardless of whether the event is caused by a natural catastrophe, terrorism, or emerging infectious diseases. Hospitals have been recognized as critical components in disaster response. Their ability to supply essential health services and continued functionality when an emergency occurs is mandatory. Within a short time, hospitals must provide care to a large number of injured whose lives are at risk, and they must have the ability to expand their services quickly beyond normal operating conditions to meet an increased demand for medical care. Emergencies and disasters often occur without warning and, as in the case of an earthquake, may damage or destroy parts of the hospital. Under such circumstances, it is necessary that hospitals and other healthcare facilities must remain safe, accessible, and functioning at maximum capacity in order to provide critical services. Contingency plans for disasters should be in place and health personnel trained to maintain functionality.

Emergency Departments (ED) are the most affected areas in hospitals in the event of a disaster. They play a pivotal role in the delivery of acute ambulatory and inpatient care in providing immediate medical assistance [Morganti et al. 2013]. In the event of a disaster, the number of incoming patients and the severity of their injuries rise significantly. A change in patients' arrival rates entails an increase of crowding and prolonged injured waiting times for the injured, thus increasing the risk of aggravating patients' conditions. Considering all these aspects, EDs should have an emergency plan ready for implementation during catastrophic events. The Emergency Response Plan (ERP) consists of procedures designed to respond efficaciously to those situations in which normal operating procedures cannot provide essential health services. Such a plan should ensure continuation of patient case, availability of equipment and treatment material, and protocols for the appropriate interaction with other critical facilities in the event of an emergency. Generally, the ERP is activated when the number of ill or injured exceeds the normal capacity of the ED or the normal operations of multiple departments to provide the quality of care required. That said, testing the effectiveness of the ERP before a disaster occurs is problematic.

Discrete event simulation (DES) models are useful tools to test ERPs under a rapid increase in the volume of incoming patients. Using discrete-event Monte Carlo computer simulations, hospital administrators can model different scenarios of the hospital to see how they compare to desired performance [Morales 2011] and assists in the planning of the most effective use of hospital resources [Šteins 2010]. Emergency Response Plans require identification of what factors represent the quality of healthcare services and what can best describe the performance of an ED during a dramatic event. Different parameters can be used to evaluate the efficacy of ERPs. Among these parameters, the most representative is the patient wait time (PWT) until medical intervention occurs.

Patient wait times play an increasingly important role in measuring a hospital's ability to provide emergency care to all the injured in an extreme situation [Cimellaro et al. 2011]. The length of time patients wait to see a provider is considered a visible and significant indicator of an ED's resilience. Overcrowding in the ED is undesirable. It creates access issues and leads to delays in care, thus putting lives at risk as an injured person's condition may escalate as they wait for care. The ability to predict PWTs when a disaster occurs could be an effective way to solve the overcrowding problem and improve the ERP. Hospitals can achieve this goal by adopting operations-management techniques and related strategies to enhance efficiency, taking into account not only the internal organization of the hospital, but interaction and coordination with other healthcare facilities. Relocating the injured to a facility that is not as badly impacted is key in saving lives.

In this research, a simplified model has been developed in order to quantify behavior in an ED during emergencies. Patient wait times have been selected as the most representative parameter to evaluate hospital resilience under seismic events. A DES model has been built for the hospital's ED considering different scenarios. Then, a meta-model was developed from the results of the DES model. It provides PWTs as a function of the seismic input and the number of the available treatment rooms in the ED. Finally, a general meta-model is proposed to evaluate the resilience index for any considered hospital. The advantage is that the meta-model can be applied as a formula with a reduced number of parameters.

## 3.2 STATE-OF-THE-ART

Although the word "simulation" has different meanings, the most comprehensive one is "the representation of something real" [Webster 1993]. In particular, in the branch of the scientific analysis, it may be defined as "the action of performing experiments on a model of a given system" in which the word "system" refers to a collection of entities that act and interact together toward the accomplishment of some logical end [Schimdt and Taylor 1970]. A model may be taken as representative of that system even if it is an abstract of the system. The greater the level of accuracy of the model, the closer it will be to real conditions. When modeling a dynamic system that describes real-life scenarios—such as a hospital in this specific case—DES is normally used.

In 1960, Keith Douglas Tocher developed the first such simulation program named the General Simulation Program (GSP). Over the ensuing decades, other developers have improved steadily simulation programs and created systems to fit different real-world problems. By the mid-1990s, a significant step forward was made by developing simulation models that could be built very quickly, giving the possibility to study, experiment, and analyze the interactions of any system and its subsystems. For this reason, DES is currently considered as a powerful and versatile tool for the analysis of complex systems. Being able to create a simulation model of a system provides the user with many benefits including: (1) the ability to analyze a number of observations of a system; (2) improved system understanding for more rapid analysis; (3) provides an opportunity to test modifications that could improve system operation; and (4) is generally less costly than direct system study [Fishman 2001]. Thus, DES has become a useful tool to study systems in many fields, but it is especially applicable in engineering, health, management, social, and transportation sciences.

Discrete event simulation is increasingly used to analyze healthcare facilities because of its multifaceted structure. There are multiple interactions between patients, doctors, nurses, technicians, different departments, and circulation patterns. How do we evaluate each of these components and take into account all these multiple interactions that affect the whole system? For this reason, DES is widely used to model the medical field because it can be represented in a chronological sequence of events that occur at a definite instant in time and marks changes in the system. In this way, the end of each event marks the start of the next event as is typical in hospitals and medical clinics worldwide. In addition, a DES model could be considered as a useful communication tool between hospital administration and modelers that help enhance administrator's understanding of the main performance drivers of the healthcare processes [Curran at al. 2005]. There are a number of comprehensive literature reviews available [Günal and Pidd 2010]. Early reviews include England and Roberts [1978], which analyzed reports of 92 simulation models.

The majority of healthcare simulation studies focus primarily on PWTs. Dansky and Miles's research [1997] evaluated the efficiency of the healthcare facilities considering PWTs as the main response parameter and investigated the relationship between PWTs times and the satisfaction with the service received. They found that customer satisfaction is strictly related to PWTs. Using DES, the current number of patients waiting as well as the time that patients spend moving can help determine utilization of resources. Martin et al. [2003] studied how all these values can be used to evaluate hospital performance and improve healthcare response. Many research projects have studied strategies on how to decrease PWTs times. One of the earliest DES studies was conducted by Fetter and Thompson [1965]. They analyzed doctors' utilization rates with respect to PWTs times by using different input variables, such as patient load, patient arrival patterns, walk-in rates, and physician service. More recently, Yerravelli [2010] studied PWTs at the KCH ED.

The objective of the research was to evaluate hospital performance as well as identify strategies for reducing waiting times by developing a KCH ED model. Resources utilization was taken into account in order to determine the required staffing levels and minimize operating costs. Santibáñez et al. [2009] provided a framework of how to reduce wait time and improve resource utilization by developing a computer simulation model of an Ambulatory Care Unit (ACU). Duda [2011] conducted a similar study and analyzed the flow of patients, the time spent in the hospital through arrival, and service characteristics in order to identify which processes need to be changed to achieve alignment and which alternatives have to be taken into account to increase the effectiveness of the patient flow process and reduce waiting times. Takakuwa and Shiozaki [2004] proposed a procedure for planning emergency-room operations that reduced PWTs times by adding a more appropriate number of doctors and medical equipment. A similar study to assess the effect of some possible changes in the ED processes was performed by Mahapatra et al. [2003]. This study showed that the addition of an ACU improved average waiting times by at least 10%.

Other strategies to address PWTs times include the study conducted by Lau [2008], who analyzed three Orthopedic Clinics across Ontario to find solutions to long PWTs times and proposed a new scheduling algorithm to decrease the average waiting time. A DES model was developed for an existing clinic by Hu [2013] to study an optimal human resource allocation in order to reduce PWTs. Similarly, Aeenparast et al. [2013] used simulation model to predict changes in patients waiting time and physicians' idle time due to changes in system. Kirtland et al. [1995] identified three alternatives that saved, on average, thirty eight minutes of wait time per patient.

Generally, EDs are one of the most popular areas for the application of DES. Medeiros et al. [2008] is one example; they constructed a simulation model for one ED and in a subsequent step implemented a new approach known as PDQ (provider-directed-queuing) that reduced non-critical PWTs and increased room availability for critical patients. A DES model has been used also by Morgareidge et al. [2014] to optimize the care process and design the ED for a particular case study. A new approach was also proposed by Davies [2007] that developed a computer simulation model for an ED called "See" and "Treat "method, where the triage process is eliminated, and the patients are directed by a trained receptionist to the doctor or ENP (emergency nurse practitioner) based on the patient condition. The "See" and "Treat" process simplifies the service by eliminating queues between patients and the hospital's human resources, which reduced the unnecessary waiting time between these queues. Similarly, Samaha et al. [2003] developed an ED simulation model to reduce the length of stay of patients by evaluating different scenarios. Their results demonstrated that longer PWTs are process-related not resource-related, and that a triage process that included "fast track" area reduced a patient's length of stay. In 2005, Komashie and Mousavi conducted a what-if analysis for one considered ED that varied the number of beds, doctors, and nurses in the simulation model, thus reducing PWTs and improve patient follow through to discharge.

This research studied the ED of a hospital has been using a DES model. Different scenarios were analyzed assuming a catastrophic event resulting in structural damage in specific parts of the building and a variable patient arrival rate dependent on the seismic intensity. An analytical model is proposed in order to obtain PWTs without running a complex DES model several times.

## 3.3 METHODOLOGY

This study assessed the resilience of a hospital ED in the event of a catastrophic event. Specifically, it examined the ED's ability to mount a robust response to unforeseen, unpredicted, and unexpected demands by adjusting its functioning prior to, during, or following a catastrophic situation. Two different aspects were considered.

First, using a DES model numerical data of WTs were determined for one hospital under normal operating conditions and those results were compared when an ERP was applied. Determining trends in PWTs in both situations provides a tool to evaluate the effectiveness of the ERP. Using the results from this comparison, a mathematical model was developed for different disaster scenarios, including the seismic input and any structural damage to the hospital.

Second, hospital emergency networks were studied in terms of resilience concepts. During an emergency, interaction among healthcare facilities is essential to sustaining required operations, to cope with a dynamic and extended influx of patients, and to provide adequate patient care under unusual conditions. In a complex system such as an urban area, the management of an emergency is closely related to the efficiency of its infrastructure network, and the continual influx of patients is another dynamic to consider. Therefore, a network model was developed for one city assuming a catastrophic event occurs. A number of strategic nodes have been identified within the considered city, and a network was created to study how the system reacts to emergencies. All nodes were defined as functioning EDs in healthcare facilities. This network has been specifically designed to adapt the resources to potentially changing demands, which includes the possibility that patients are transported to another healthcare facility in the network to that they received care in the shortest possible time.

Creating a valid network framework requires consideration of a number of factors. Among all the parameters related to such complex scenario, the most representative parameters are: PWTs, the ED's resilience, and transportation considerations. In order to evaluate PWTs, the meta-model described above was used. An evaluation analytical model for assessing comprehensively hospital disaster resilience is proposed below.

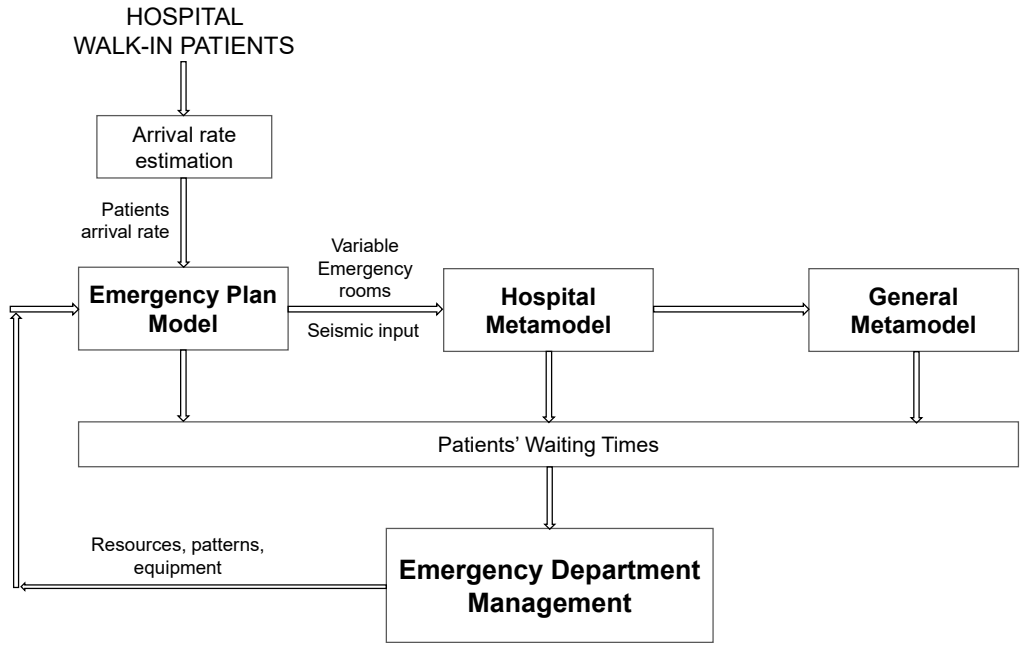

**Figure 3.1**     **Research logic framework.**

Figure 3.1 shows the research framework. First, a simulation model for the ED when the ERP plan is applied, using the estimated arrival rate for patients arrival rate as input data; see Section 3.4. Once the output parameters have been collected, the most significant PWT is chosen. The hospital meta-model is developed by varying the seismic input and the number of the available treatment rooms; see Section 3.5. Using this meta-model, a general meta-model that is applicable to all hospitals was developed; see Section 3.6. Based on PWTs, the general meta-model provides hospitals with metrics regarding its resilience, with the goal of improving its response to catastrophic events.

## 3.4    DISCRETE EVENT SIMULATION FOR EVALUATING AN EMERGENCY DEPARTMENT'S PERFORMANCE

A simulation model for the Mauriziano Hospital's ED is presented below; see Figure 3.2. First, the problem must be defined, with specific objectives and questions to be addressed. This model ties directly hospital performance to PWTs**.** In emergency situations the length of time patients spend in the ED from admittance to discharge is considerable; therefore, several steps have been considered. The entire process by which patients enter the ED, interact with medical staff, receive appropriate treatment, and finally are discharged has been considered.

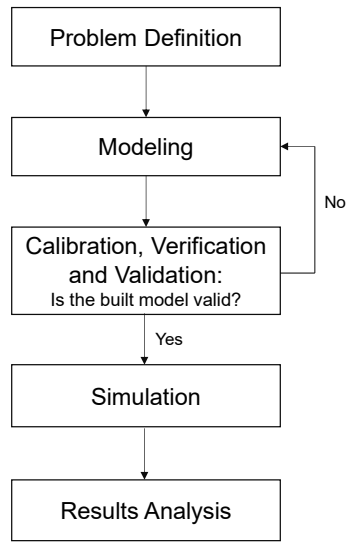

Figure 3.2        Steps followed to develop Mauriziano simulation model.

### 3.4.1  Description of the Case-Study

This research developed a computer simulation model for the Umberto I Mauriziano Hospital located in Turin, Italy. The hospital stands out in the landscape of healthcare facilities in the Piedmonte region because of its broad array of both basic case and specialization services. Built in 1881, the hospital is located in the southeast part of the city, roughly 3 km from the city center. It was bombed several times during World War II, necessitating rebuilding portions of the hospital. The addition of several other buildings over time has resulted in a large hospital complex. Presently it includes 17 units, covering an overall surface of 52,827 $m^2$. While developing the simulation model, only the ED located in building 17 were considered.

The ED consists of an entrance area where "triage" is performed; four macro areas correspond to four different color codes—red, yellow, green and white—that represent the severity of injury. Red (emergency) identifies those patients with compromised vital functions whose lives are at risk. Yellow (urgency) identifies those patients who are not in immediate danger of life but present a partial impairment of vital functions. Green (minor urgency) identifies those patients that need medical care but do not have a critical injury affecting vital functions. White (no urgency) identifies those patients whose medical issues could be addressed by a general doctor and who do not require emergency care.

When the ERP is in force, the number of color-coded areas is reduced to three, with those patients color-coded "white" being sent to another facility outside the ED. Under normal operating conditions, patients color-coded yellow and green codes share the same area, i.e., the treatment rooms. Under emergency conditions, the red-coded area is located immediately in front of the ambulance entrance and contains two rooms in which patients receive initial care. The yellow-coded area is parallel to this area, comprising three of the treatment rooms. Separate from this zone, the green-coded area is situated perpendicular to yellow- and red-coded areas and includes two treatment rooms. Each area is provided with waiting rooms for patients. Inside the ED, there are also a number of recovery rooms where patients can stay before being discharged or awaiting transport to another part of the hospital.

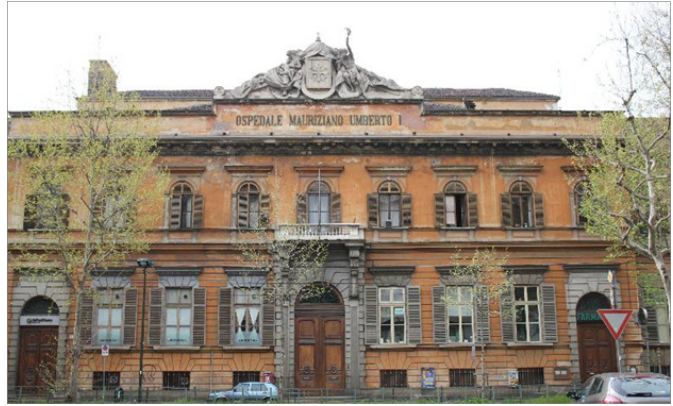

**Figure 3.3**    Umberto I Mauriziano hospital, Turin.

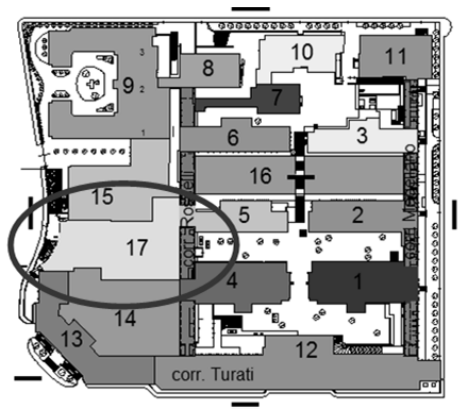

**Figure 3.4**    Location of Emergency Department building within the Mauriziano Hospital complex.

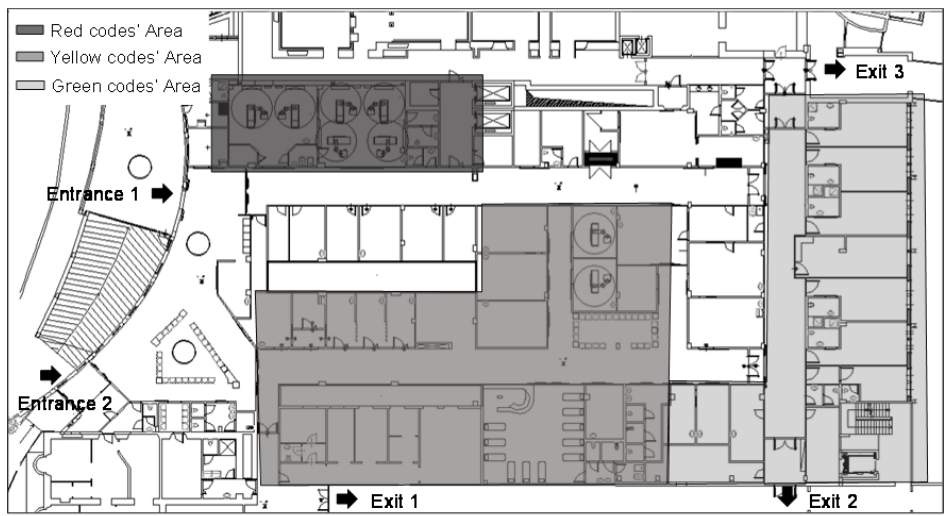

**Figure 3.5**    Emergency Department color-codes areas.

### 3.4.2 Emergency Department Simulation Model

The modeling phase involves all the details regarding data collection, data inputted to the simulation model, elements used in simulation model, and details about computer simulation model. The starting point is an empty system without any characteristics. In this research, an ED simulation model was developed using ProModel version 7.0, downloaded on February 15, 2014. A number of simplifying assumptions were made to describe the ED operation more concisely.

First, the model describes the hospital's ERP so that an emergency scenario, e.g., earthquake, can be considered. In order to represent an extreme situation, a seismic arrival rate was used. Unless specifically noted, in this case it was assumed that the hospital's structural/non-structural elements remained undamaged due to the earthquake. Therefore, the description of the ED operation is based on the ERP with no other parameters considered.

Second, a 13-day simulation was run considering that the seismic event occurs after two days of simulation. In order to simplify the model, it was assumed that the ERP was applied in the first day of the simulation even though the minimum required conditions to apply it were not satisfied.

Third, it was assumed that even after a disaster occurs, the system behavior does not change. This is an unrealistic assumption because in emergencies, several modifications may affect system characteristics but these changes were not considered herein. This hypothesis has been done to analyze how the system could respond to an emergency event such as an earthquake while maintaining the initial characteristics.

Fourth, in order to define the arrival rates, a division according to the injury color codes was considered from the moment patients first arrive at the hospital. In reality, the injury code is determined once patients are evaluated, i.e., care administrated during "triage." A patient's color code may change during their stay at the ED. For this simulation, it was assumed that those patients whose color code changed, all did so at the same point of their treatment.

In this case study, only four codes were considered: red, yellow, green, white. In reality, the ERP for this hospital considers blue and black codes, which represent "compromised vital functions" and "died," respectively. These two designations were not considered in this model because they had no influence on PWTs.

Finally, because the ERP has never been implemented and there is no data regarding its efficacy, the probability values entered for the construction of the model were obtained from interviews with the hospital's medical staff. In the case where the information gleaned from the interviews were deemed unsatisfactory, a probability of 50% has been considered.

### 3.4.3 Data Collection

When building the simulation model for the DES, the most significant data that describes all of the hospital's processes were identified. The ED can be characterized by the number of treatment rooms, the number of resources (doctors, nurses, and healthcare staff), and the procedures conducted inside the different rooms, as well as the circulation patterns and patient arrival rates. Three main methods were used to collect data.

First, the hospital's ERP was considered as a source of information regarding organizational aspects during emergencies, such as resources schedule, locations, and patient path networks.

Second, patient arrival rates were calculated using the hospital's register statistics. Information regarding the patient inflow, check-in and checkout times, and the time spent in each room as well as patients' movements from one location to another were obtained. Patient arrivals in the ED vary from hour to hour and, in order to determine the patient arrival distributions, an arrival cycle was defined using data from the hospital registers. The data provided by hospital's records were used to validate the simulation model. Figure 3.6 shows the considered arrival cycle, which represents the percentage of daily patients that arrive at a given instant within a twenty-four-hour period.

Third, researchers interviewed medical professionals who work in the ED to understand and describe patient flow. The product of this phase is a flow map, which was submitted to the hospital's personnel for review and approval. By the direct conversation with ED staff, all processes that take place in the ED during an emergency situation were examined.

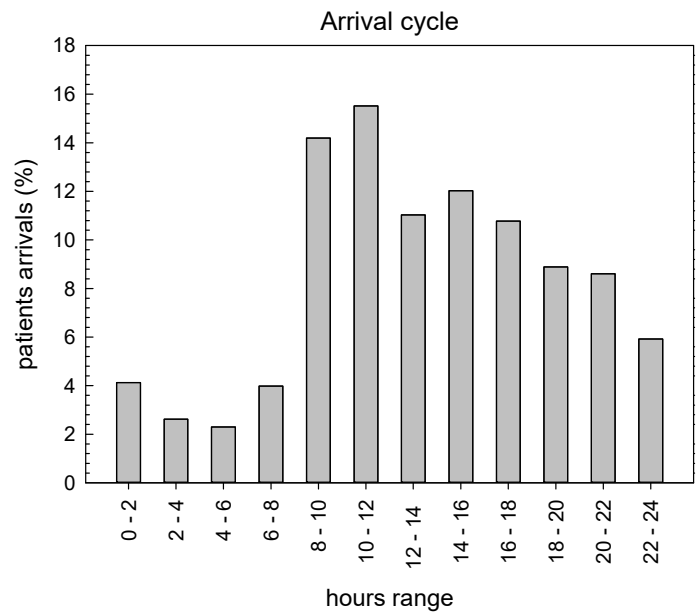

Figure 3.6     Percentage of patients arriving for each hour in a twenty-four-hour period.

### 3.4.4 Seismic Input

To take into account the increase of the patients flow due to a catastrophic event and the consequent overcrowding of the ED, a seismic event was considered. The data collected from a California hospital during 1994 Northridge, California, earthquake was used in the model to simulate the seismic event. The earthquake damaged a number of buildings and freeways, and many people were injured. Even though it is considered a moderate event, its impact was substantial because it occurred in a metropolitan area. The patient arrival rate in the aftermath of the Northridge earthquake was selected because it is the only event with documentation of patient arrival rates [Stratton et al. 1996; Peek-Asa et al. 1998; and McArthur et al. 2000].

The pattern of the arrival rates of patients in the aftermath of the Northridge earthquake is provided by Cimellaro et al. [2011]. For this study, patient arrival rates were scaled to correspond to the analyzed geographic area (Turin, Italy). An earthquake with a return period of 2500 years was considered, assuming a nominal life for a building of strategic importance of 100 years according to the Italian seismic standards [NTC-08 2008]. Initially, a scaling procedure based on the peak ground acceleration (PGA) was used, but because of its limitations, a second procedure based on the Modified Mercalli Intensity (MMI) scale was selected. Figure 3.7 plots the seismic input during the three-day period after the earthquake occurred. The patient arrival rates related to Northridge were scaled with respect to the ratio between the PGA and the MMI values. The seismic arrival rate shown in Figure 3.7 was divided in different color codes following a similar distribution proposed by Yi [2005].

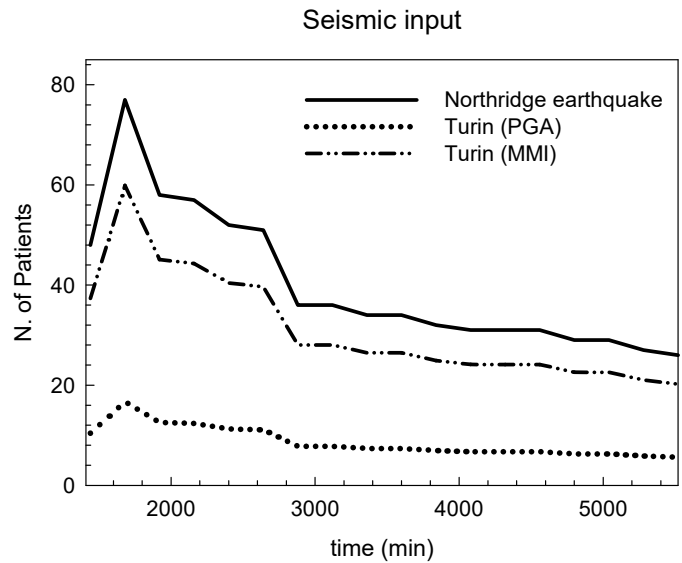

Figure 3.7    Patient arrival rates for Northridge earthquake and arrival rate scaled with respect to PGA and MMI.

### 3.4.5 Minimum Requirements for Application of the Emergency Response Plan

The ERP is applied when the flow of incoming patients exceeds the normal flow. According to the Mauriziano Hospital's protocol, this situation occurs when there is the simultaneous access (or within a short period) of 10 or more patients in critical condition, at which point the ERP is activated. Per the ERP, patients in critical condition are coded either red or yellow. Figure 3.8 shows the number of patients arriving at the ED during the three-day period post-earthquake.

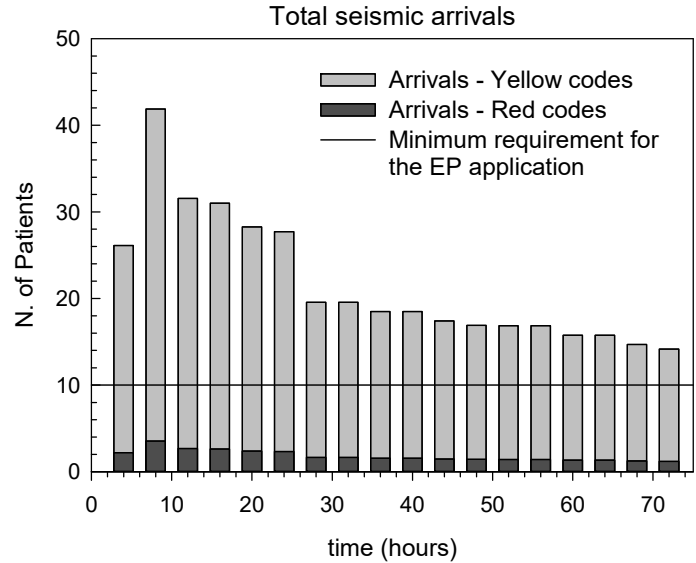

**Figure 3.8**   Total arrival rate during an emergency (red and yellow codes).

### 3.4.6  Model Architecture

For this study, the DES model was built using ProModel [2007] simulation software, which is considered the best choice to develop a complex hospital model; see Figure 3.9. The first step is define all elements of the model. These elements along with the simulation constructs are listed below:

*Location elements:* This element represents all the places in the ED where entities carry out their duties. The physical facilities consist of treatment rooms (color-coded given the condition of the patient), two waiting rooms, a triage room, an examination area, a critical area (one shock room and one intensive reanimation room), observation rooms, and some separate stations. There are two possible entrances to the ED; one is used by the ambulances only and the other is for patients and/or visitors. The first entrance is located in the northwest part of the ED near the red-coded area, and the second entrance is on the southwest side. Those patients whose mode of arrival is either by ambulance or car (assuming the patient is in critical condition) enter though the north entrance, which is closest to the shock and intensive care rooms. All other walk-in patients use the south entrance, which is nearest to the yellow- and green-coded areas. There are three exits located in the south, northeast, and southeast sides of the ED. Which exit patients use depends on their destination (others healthcare facilities, hospital wards, having been discharged, etc.). Each location has been assigned a capacity. Some locations, such as entrances, exits, and waiting rooms, have infinite capacity while others, like treatment rooms, the shock room, and intensive care room, have a defined number of patients who can be treated at the same time.

*Entities*: This element represents the active elements created in the system. They move within the system and are affected and processed by the status of the system. In this model, entities are patients visiting the ED that are categorized according to the severity of their ailments. Entities that flow through the simulation model have been divided into four categories corresponding to the four color codes: red, yellow, green, and white codes. For each type of patient, an entry and a path have been assigned. Each entity has a travel speed: 50 mpm for yellow-, green-, and white-coded patients, and 60 mpm for red-coded patients.

*Path network*: This element consists of patients, nurses, and doctors who proceed through the model on a path network. It consists of nodes connected by segments, which can be unidirectional or bidirectional. In this model, movement along these path segments connected to the path nodes are defined in terms of distance and speed. The path network created for the ED is shown Figure 3.10. Dotted lines indicate the path segments on which entities can move.

An explicit mapping for some destination nodes has been created for the specific branches that entities and resources must take when traveling through the model. For speed and distance networks, if there are multiple paths emanating from one node to another node, the default path selection is based on the shortest distance. In the ED in question, however, some paths can only be used by medical staff (as in the case of the passage from the red to yellow area). In these cases, the mapping definition has been used to solve the default problem.

*Resources:* This element is defined in the model as persons who assist entities in performing operations. In this particular case, resources represent the medical staff, including doctors, nurses, and healthcare operators. They are divided into two categories: those persons who provide service from a fixed station and those persons who travel through the ED. A work schedule has been considered when modeling resources. As summarized in Table 3.1, each color-coded area has its own staff team. In addition to these teams specific for each area, there are a number of supporting resources, including doctors, nurses, and healthcare staff, who perform triage, transport patients from one location to another, and assist with admittance to the hospital.

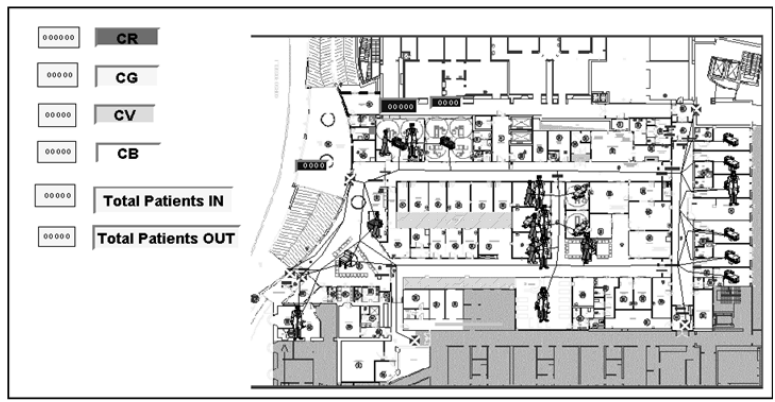

**Figure 3.9** DES model extract from Promodel software [2007].

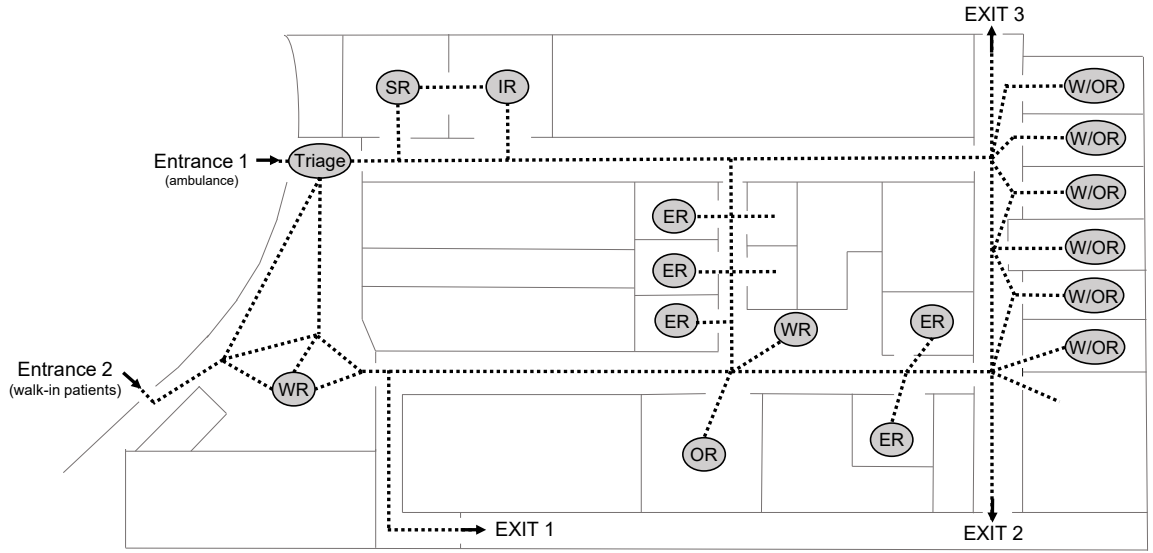

**SR**=Shock room   **IR**=Intensive care room   **ER**=Emergency room   **WR**=Waiting room   **W/OR**=Waiting and observation room   **OR**=Observation room

**Figure 3.10**   Network for patients to move through the Emergency Department.

Table 3.1        Resources definition.

| Color codes area | Work schedule | Resources |
|---|---|---|
| Red area | hours 8/20 | 2 doctors, 4 nurses |
|  | hours 20/8 | 2 doctors, 3 nurses |
| Yellow area | hours 8/20 | 5 doctors, 3 nurses |
|  | hours 8/20 | 5 doctors, 3 nurses |
| Green area | hours 8/20 | 3 doctors, 5 nurses |
|  | hours 8/20 | 2 doctors, 3 nurses |

*Processing tool:* This element defines the way in which entities move between locations or remain at a given location. To build a model that represents adequately the ED operation, all activity regarding patients from when they arrive until they leave the ED must be considered. It is necessary to take into account not only patients movements from one location to another, but also how much time they spend in each location and how and for how long they use a particular resource. According to the ERP, patients follow different paths and make diverse actions depending on their color code. According to data provided by the hospital staff and the ERP, the processing phase has been developed considering the following actions for each color code.

- Red Codes: red-coded patients generally arrive by ambulance at Entrance 1. As soon as they arrive, they are moved directly to the red area after triage is performed. Red-coded patients mainly use two rooms: the shock room and the intensive care room where the most critically ill or injured patients are treated immediately. After receive the first treatment in these two rooms, some patients are relocated to the yellow area, others are transferred to the appropriate hospital ward, and others leave the hospital, either being discharged or moved to another to another healthcare facility.

- Yellow Codes: yellow-coded patients may arrive at both Entrance One and Two. After receiving triage, they wait in the yellow-coded waiting room until one of the treatment rooms is available. Some patients are kept in the observation room where they receive the initial treatment. After being treated, some patients leave the hospital while others are sent to the exam room. Once examined, patients are sent back to either treatment rooms or to the green-coded area. From treatment rooms, some patients leave the ED (and are transferred to hospital wards or to relocated to other healthcare facilities), and some patients are sent back to the exam room until it is deemed appropriate to discharge them from the ED.

- Green Codes: green-coded patients normally arrive at Entrance 2. Considering their less serious condition, these patients are requested to wait until triage can be performed as priority is given to yellow- or red-coded patients. Once they are coded green, they are sent to the green area where they wait in several observation rooms. While waiting, patients who present less severe injuries are treated by an available nurse and then leave the hospital. The others wait until one of the green-coded treatment rooms is available. After receiving treatment, they leave the hospital or move to the exam room. Once examined, they leave the ED, either relocating to a hospital ward or are discharged.

- White codes: white-coded patients arrive at Entrance 2. These patients are requested to wait until triage can be performed as priority is given to yellow- or red-coded patients. According to the ERP, patients designated "white" are allowed to leave the hospital once triage is performed. They are not treated in the ED because their conditions do not warrant it. Once the triage is performed, they leave the ED. This is because white codes are for those patients that do not have any serious injuries or disease; for this reason, in emergency situations they are not treated in the ED.

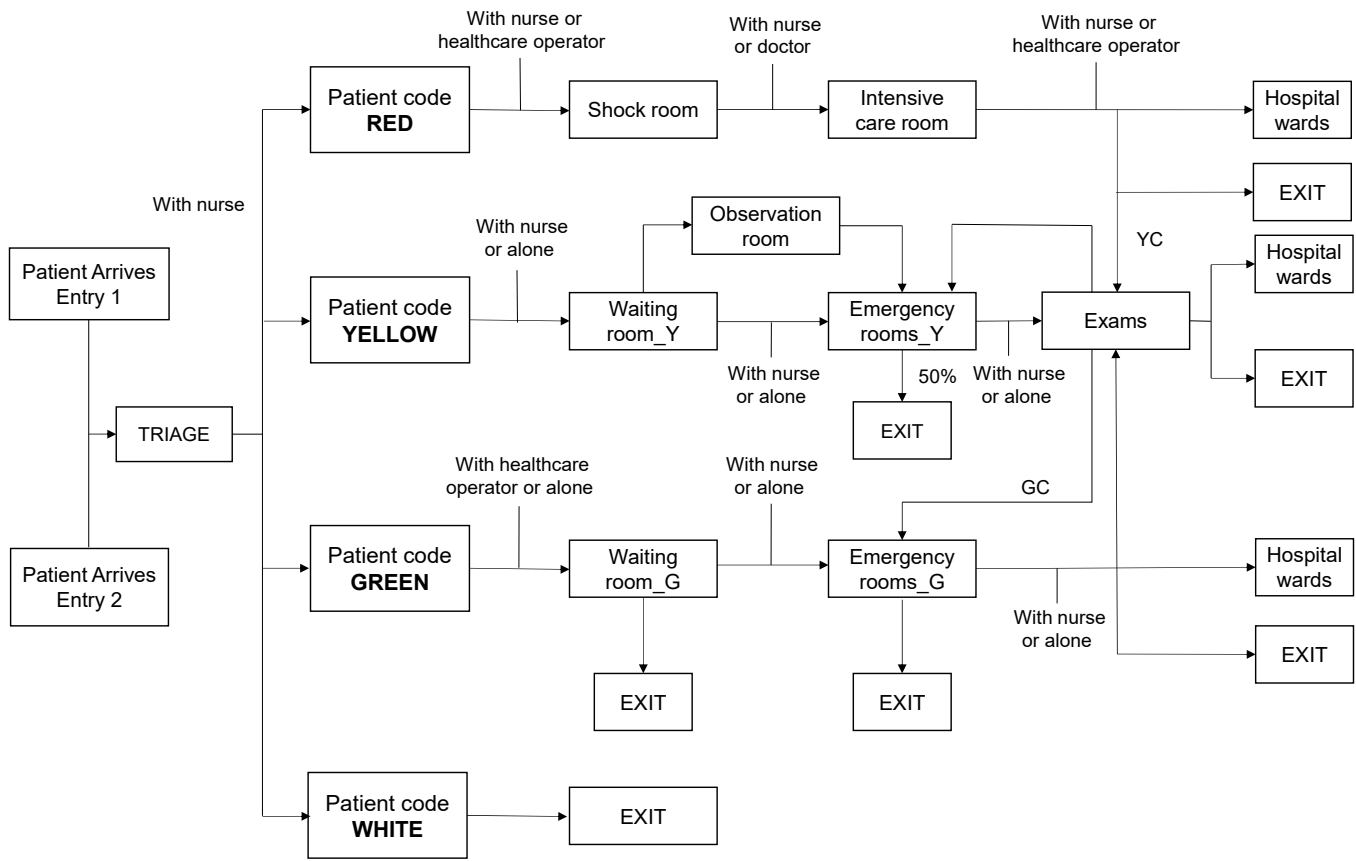

**Figure 3.11** Process map for the Emergency Department.

### 3.4.7 Verification, Validation, and Simulation

After building the model has been completed, it has to be verified and validated. Once the model was completed, it was verified to ensure that the computer programming and implementation of the conceptual model was correct. Patient routing, service times, and all assumptions were validated with ED staff. The process map representing the model logic was checked with the ERP director, and all required corrections and changes were made.

Once built the model and its accuracy verified, numerical simulations of the DES hospital model were performed. A simulation period of 13 days was run in each simulation, dividing the seismic input into three main parts. The simulation began with a three-day simulation period that followed normal operating conditions under extreme conditions, then run in emergency operating conditions determined using the scaled arrival rate assuming a seismic event, and another eight-day period was run in normal operating conditions because the system needs time to return to the pre-earthquake steady state; see Figure 3.12.

Different output variables were collected to obtain a general idea of how an increase in patient arrival rates can affect the operation of the ED. The thirteen-day simulation was run 100 times for each different scenario, taking into account the differences between different runs and correcting any errors related to the individual simulations.

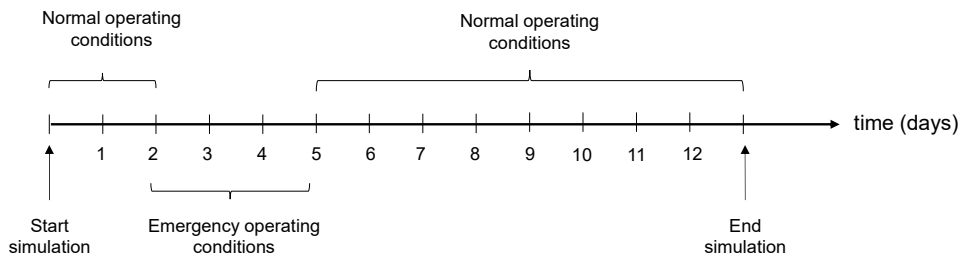

**Figure 3.12**   Simulation framework for the DES hospital model.

### 3.4.8 Analysis Results

Results demonstrate that both yellow- and green-coded patients experienced longer PWTs in normal operating conditions during an extreme situation. In particular, yellow-coded PWTs reached an average peak value of about 720 min while green-coded PWTs was roughly 750 min. When the ERP is implemented, PWTs reached an average peak value of about 30 min for yellow-coded patients and roughly 190 min for green-coded patients Obviously, PWTs increased in the aftermath of the earthquake and returned to normal after the emergency period was over.

The comparison between waiting times for both normal and emergency operating conditions demonstrates that simulated model results were consistent with expectations. The ERP application has a greater impact on reducing PWTs. Figure 3.13 illustrates that PWTs varied substantially when moving from normal operating conditions to emergency operating conditions. A 96% decrease was observed in PWTs for yellow-coded patients moving from normal to emergency operating conditions, while there was a 75% decrease for green-coded patients.

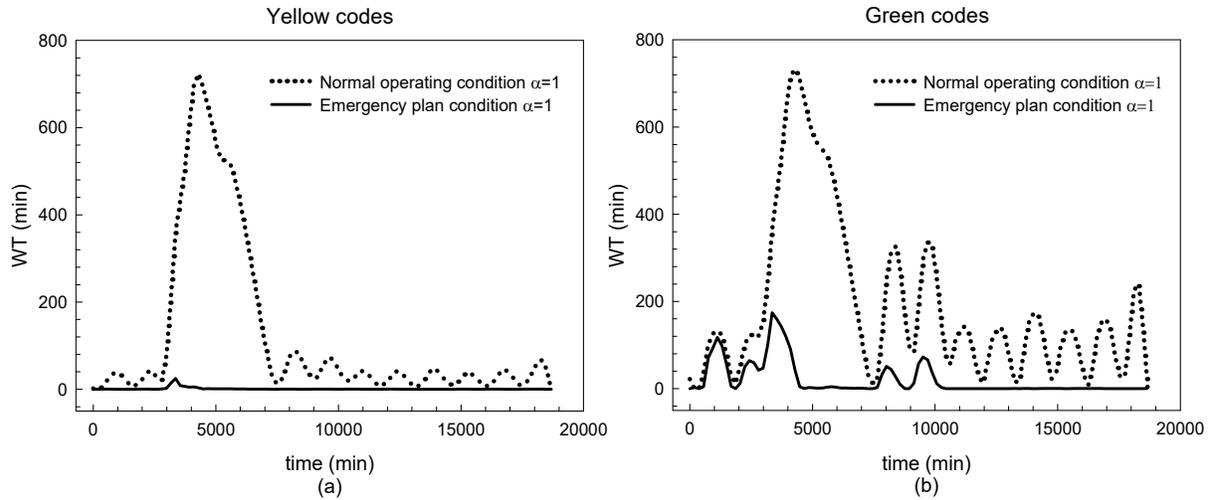

**Figure 3.13** Comparison between normal operating conditions and emergency conditions with $\alpha = 1$ for (a) yellow-coded patients and (b) green –coded patients.

### 3.4.9 Amplified Seismic Input

To study the impact of earthquake magnitude and how this affects PWTs, an amplified seismic input was considered. The magnitude was amplified in order to analyze the sensitivity of the operations of the ED in regards to the amplitude of the earthquake. Multiplicative scale factors ranging from 1.1 to 1.6 were used to amplify the input data. The factors ($\alpha$) used for the analysis are shown in Figure 3.14.

Once the amplified magnitude has been calculated using the different scale factor, new input data were introduced into the simulation model, and a Monte Carlo simulation was run to obtain the average curve of PWTs for the ED. Using the same amplified seismic input, the average curve of PWTs was obtained also for normal operating conditions. A comparison between the normal and the emergency operating conditions with a factor value of 1.6 is shown in Figure 3.15 for both yellow- and green-coded patients.

Figure 3.15 demonstrates that PWTs decreased substantially from normal operating conditions to emergency operating conditions depending on the seismic input ($\alpha = 1.6$). An average peak value of about 3200 min for yellow-coded patients and roughly 3250 min for green-coded patients was reached under normal operating conditions. When the ERP was activated, the PWTs reached an average peak value of about 300 min for yellow-coded patients and about 785 min for green-coded patients.

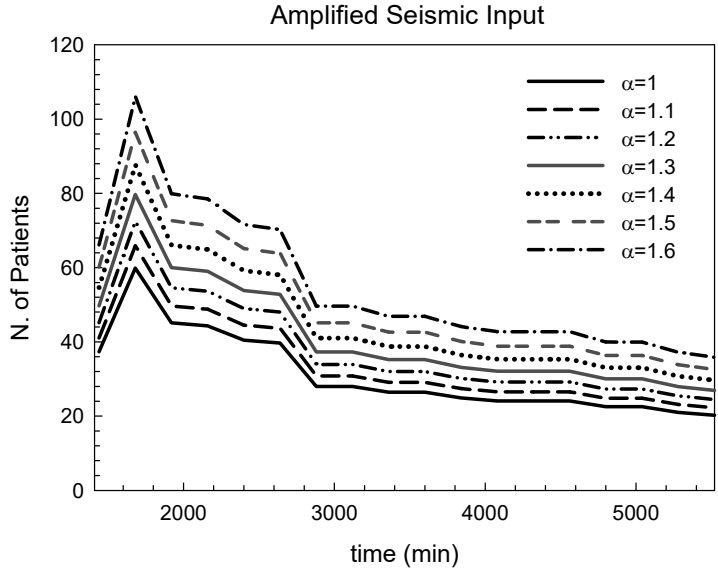

Figure 3.14    Amplified seismic input for different scale factors *α*.

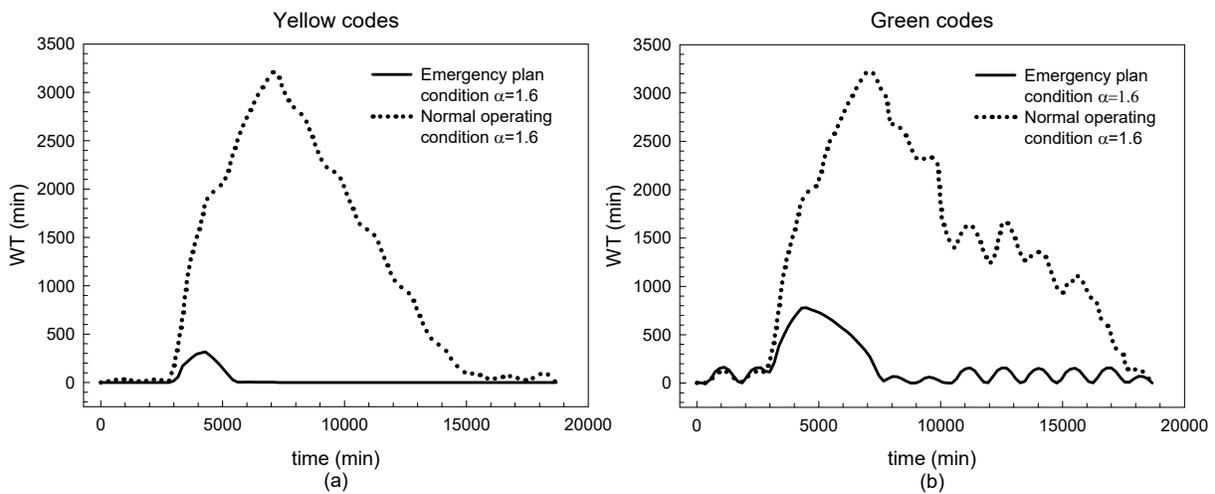

Figure 3.15    Comparison between normal and emergency operating condition with an amplified seismic input ($\sigma$=1.6) for (a) yellow-coded patients and (b) green-coded patients.

### 3.4.10 Evaluating Changes in the Emergency Response Plan

Changes in the ERP were analyzed to reduce PWTs. Green-coded patients were considered for the following reasons: (1) these patients face a lengthy wait in an emergency situation, reaching an average peak wait of about 800 min (13 hours) when the seismic input is amplified with a multiplicative scale factor of 1.6. Although green-coded patients must give priority to red and yellow codes, delayed diagnosis and treatment could affect a patient's condition, potentially affecting treatment and outcomes. This delay can also lead to complications, putting patients' lives and well-being in jeopardy; (2) these patients can receive treatments also outside of the treatment rooms. For this reason, an additional doctor could decrease PWTs without the need to add an additional treatment room.

This research studied several possible strategies to reduce PWTs. An analysis of the hospital's ERP took into account various options. First, an increase in doctors or number of treatment rooms was considered. Three different solutions were tested to study the system's sensitivity towards a change in the number of resources and locations.

First, all hospital characteristics were maintained, but the model was run considering one additional doctor for the green-coded area. Then, an extra location was considered, assuming that green-coded patients used the treatment room designated for white-coded patients. Finally, one supplementary resource and location were added and modeled. Monte Carlo simulations were performed for each case, and the results are shown in Figure 3.16. The graph illustrates that when one additional doctor is considered, the average peak of PWTs decreased roughly 39%. If a treatment room is added results in a reduction of about 74% compared with the existing configuration. Finally, considering both resource and location, a peak of about 90 min is reached, which corresponds to a 88% decrease.

Based on these results, the best solution to decrease PWTs for green-coded patients is to add one additional doctor and one additional treatment room. Closer examination of the results demonstrates that increasing the number of treatment rooms for green-coded patients is optimal. It would not increase costs, whereas the addition of one doctor and one treatment room would. Furthermore, the difference in the reduction in PWTs between the first solution and the second one is about 14%, corresponding to about 90 min for green-coded patients under emergency conditions. This decrease is not significant when compared with the additional incurred costs if an additional doctor and treatment room were added.

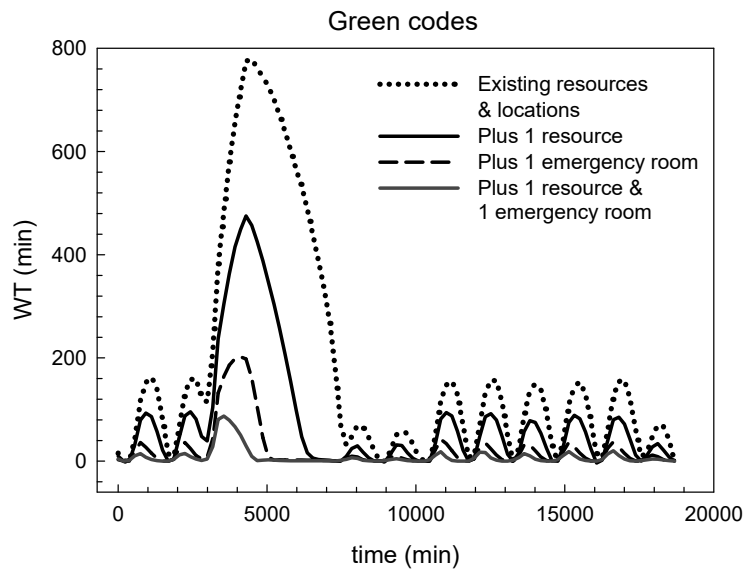

**Figure 3.16**     Extra resources and additional locations for green-coded patients under emergency operating conditions.

## 3.5 MAURIZIANO EMERGENCY DEPARTMENT META-MODEL

### 3.5.1 Motivation for a Meta-Model

This study used a simulation model of an ED to evaluate the ERP of a hospital and used the results to analyze ways to improve its efficiency. Although a simulation model is useful for predicting ED response, two different aspects make this tool unusable by the hospital staff during emergencies. First, the long run times necessary to explore the simulated scenario to retrieve statistically meaningful results render it impractical to use it for evaluating real-time hospital factors. Secondly, DES models produce a significant amount of complex data; thus, the results should only be interpreted by experts.

Therefore, an approximation of the simulation model was developed. This is a preferable strategy to study an environment such as an ED, where the stochastic nature of its operation demands that decisions be made quickly. This model, called a meta-model, simplifies the system, exposing more clearly the input–output relationship. The objective is to create a relatively simple functional relationship between the system behavior and selected variables that represent the surrounding conditions. The meta-model is easier to manage compared with DES models, and it is also useful in order to analyze model parameters without having to perform additional simulation runs. In this specific case, a set of equations were developed that approximate the simulation.

Several steps have been taken in order to develop the meta-model framework. First, the parameter that describes the system behavior was identified. In this case, PWTs were chosen as the most representative factor in determining the efficacy of the ERP. The meta-model uses PWTs as a function of specified variables correlated to the hospital's characteristics. This research considered two variables: the seismic input and the number of available treatment rooms. These variables describe a number of different scenarios that could represent all considered real situations. In this way, an equation can be used as a "real-time" decision aid to determine the best alternative when facing an emergency situation.

### 3.5.2 Methodology

This section details the structured approach used to building the meta-model, which consisted of defining the problem and then generation of the meta-model. Defining the goal of the model required identifying the meta-model input and response as well as determination of the most relevant characteristics of these data. The considered meta-model inputs are the earthquake magnitude ($a$) and the number of non-functioning treatment rooms ($n$) due to the earthquake. The main output parameter is PWTs.

Once the input and output data have been defined, then the meta-model itself can be generated. Herein, a sensitivity analysis was used. It measured how the system output varied with respect to a change in system parameters or inputs. In this research, the sensitivity analysis has been performed in respect to the two defined meta-model inputs. First, it simulated the closure of some treatment rooms one by one, considering that possible structural damage due to the earthquake has rendered these rooms unusable. Second, an amplified seismic input was considered using a number of scaling factors. Then, Monte Carlo simulations were run for all the considered scenarios with all possible *α* and *n* combinations. Using the data from the simulations, a nonlinear curve regression was used in order to find an equation that could describe PWTs. Equation (3.1) provides the mathematical correlation between independent and dependent variables, whose general form is shown below:

$$Y = f(x_1, x_2, ..., x_n) \tag{3.1}$$

where *Y* represents the independent variable that is the output of the meta-model, and $x_1, x_2, ..., x_n$ are the dependent variables or inputs. By putting only the input values in Equation (3.1), it is possible to calculate the output value without the need to perform further simulations. Sigmaplot 12.0 software was used to curve-fit the data.

### 3.5.3 Assumptions

The assumption during construction of the meta-model are as follows:

First, the proposed meta-model was built based on the DES model; all the parameters of the meta-model were calculated considering the simulations results. Nevertheless, these results are strongly related to the entire hypothesis made at the simulation level. Therefore, all the assumptions made during the DES model construction must be taken into account when building the meta-model.

Second, some approximations related to the mathematical structure of the meta-model were made. Thus, the construction of the analytical solution starts by postulating a specific form for the model and then testing its validity. It means that first, a trend for the output variable has been chosen and considers that all the examined scenarios can be described by the same equation; only the dependent parameters re varied. The average PWT resulting from the simulations has been considered in order to select this unique equation that represents the trend of the output variable.

Third, the meta-model was built only for yellow-coded scenarios. The procedure for a green-coded meta-model is the same developed for yellow-codec scenarios.

Finally, it was assumed that the surrounding conditions of the ED remain intact. During emergencies several modifications may affect system characteristics, but these changes have not been considered herein, i.e., resources (doctors, nurses, etc.), path networks, and locations remain the same, and the only two considered variables are the seismic input and the number of the available treatment rooms. This hypothesis has been done in order to analyze the sensitivity of the system with respect to a change in the number of treatment rooms and the intensity of the earthquake.

### 3.5.4 Meta-Model Architecture

The meta-model for the hospital's ED was constructed considering both normal and emergency operating conditions. It describes hospital response in real time during an extreme situation. It takes into account the influence of possible structural damage due to seismic input, rendering some treatment rooms non-functional. Patient wait times have been identified as the main parameter in order to describe hospital response under an emergency situation and the most representative indicator to evaluate the resilience of an ED.

The proposed meta-model has been developed by representing the simulation model results by a mathematical function whose general form is:

$$Y = f(x_1, x_2, ..., x_n, t) \tag{3.2}$$

where $Y$ is the system response, $x_1, x_2 ... x_n$ are the considered input variables, and $t$ represents the time. In this specific case, the system response is the PWTs, and the input variables are the seismic input and the number of non-functional waiting rooms.

$$WT = f(t, n, \alpha) \tag{3.3}$$

where $WT$ represents patient wait time, $n$ is the number of non-functional treatment rooms, $\alpha$ is the seismic input, and $t$ is the time in minutes.

To develop this equation, a specific form for the meta-model has been postulated. A simple scientific approach was considered. First a function for the model that may closely follow the output variable $WT$ was formulated. Then, the parameters of the selected model were estimated. Finally it was demonstrated that proposed meta-model adequately represented ED behavior. A log normal function was chosen, and its parameters were calibrated based on simulation model data. After fitting different equations, the following approximating function was selected as the best one in order to describe the trend of WTs for all the different scenarios considered.

$$WT = \frac{a}{t} * \exp\left[-0.5 * \left(\frac{\ln t/b}{c}\right)^2\right] \tag{3.4}$$

where $WT(t, n, \alpha)$ is the waiting time in minutes, $t$ the time in minutes, $n$ the number of non-functioning treatment rooms, $\alpha$ the seismic input, and $a, b, c$ are the calculated parameters. Both normal and emergency operating conditions were studied, and two different meta-models have been proposed.

### 3.5.5 Meta-Model Parameters under Normal Operating Conditions

Using the approach described above, a meta-model for the ED under normal operating conditions was created. As mentioned before, PWTs was used as a key parameter to describe hospital behavior in an emergency situation. It can be expressed mathematically by a lognormal equation in which $a, b, c$ parameters are calibrated for the yellow-coded patients. In order to find the value of these parameters, Monte Carlo simulations were performed for all the possible $\alpha$ and $n$ combinations; see Table 3.2.

First, the earthquake magnitude given in Figure 3.7 has been amplified proportionally using scaling factors shown in Figure 3.14. A Monte Carlo simulation was performed assuming a constant value of $n$ and a variable value of $\alpha$. The average $WT$ for each scaling factor was considered in order to calibrate the analytical model. The simulations results for each $n$ value are illustrated in Figure 3.17. The figure shows that by increasing the seismic magnitude, PWTs rise proportionally with the scaling factors.

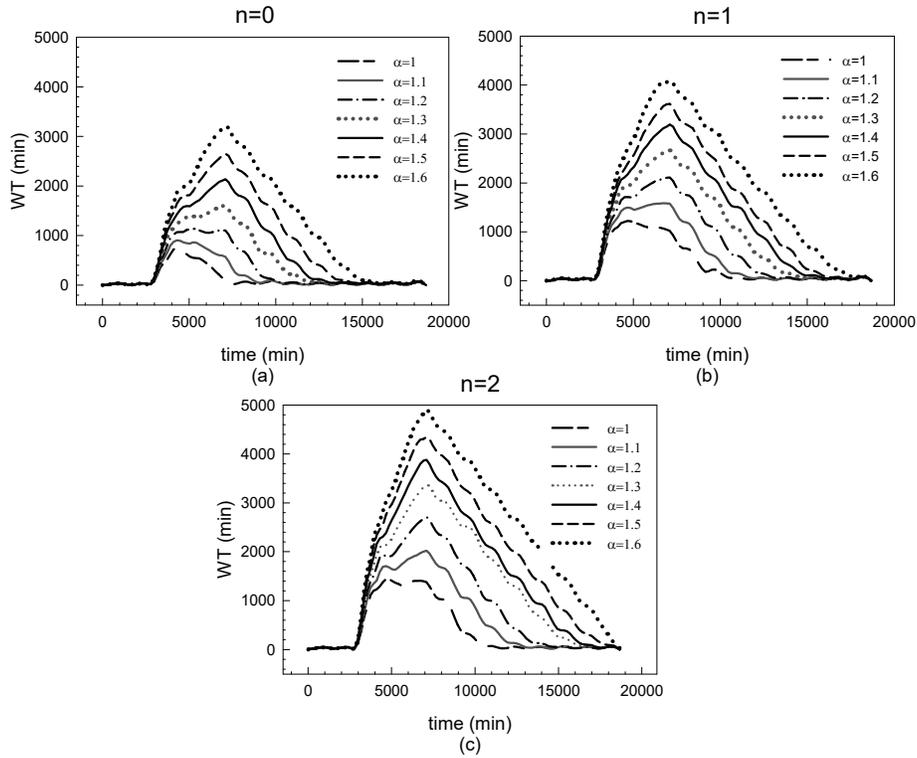

**Figure 3.17** Simulations results under normal operating conditions with a constant value of *n* and a variable value of *α* for (a) *n*=0; (b) *n*=1; and (c) *n*=2.

Table 3.2    All performed DES simulations.

|       | α=1 | α=1.1 | α=1.2 | α=1.3 | α=1.4 | α=1.5 | α=1.6 |
|-------|-----|-------|-------|-------|-------|-------|-------|
| N=0   | 0_1 | 0_1.1 | 0_1.2 | 0_1.3 | 0_1.4 | 0_1.5 | 0_1.6 |
| n=1   | 1_1 | 1_1.1 | 1_1.2 | 1_1.3 | 1_1.4 | 1_1.5 | 1_1.6 |
| n=2   | 2_1 | 2_1.1 | 2_1.2 | 2_1.3 | 2_1.4 | 2_1.5 | 0_1.6 |

Then, the effect of the hospital experiencing structural damage and its effect on PWTs was investigated. The simulation investigated closing treatment rooms (ER) one by one (*n*), assuming that the damage was extensive enough that they were non-functional. A Monte Carlo simulation was run considering a constant value of *α* and a variable value of *n*. The simulation results are shown in Figure 3.18 for three different *α* values.

As shown in Figure 3.18, it is possible to see that by closing some treatment rooms, PWTs increased drastically. In particular, when the *α* factor is 1.6 and two treatment rooms are closed, PWTs reach a peak of about 5000 min, corresponding to approximately 84 hours (three and a half days): the system is congested due to a patient volume that exceeds hospital capacity.

In order to build the meta-model, it has been observed that the trend of the graphs shown previously could be approximated with a bell-shaped curve. Thus, for any fixed value of $\alpha$ and $n$, the *WT* curve always presents a steady state before the seismic event, which is followed by a peak that represents a transient period in which the system experiences an increase in patient flow. A regression equation has been obtained to represent *WT*, where a lognormal function was selected to represent the *WT* trend, whose general form under normal operating conditions is:

$$WT(t,n,\alpha) = \frac{a_n}{t} * \exp\left[-0.5 * \left(\frac{\ln(t/b_n)}{c_n}\right)^2\right] \qquad (3.5)$$

where $WT(t, n, \alpha)$ is the waiting time in minutes, $t$ is the time in minutes, $n$ is the number of non-functioning emergency rooms, $\alpha$ is the seismic input, and $a_n$, $b_n$ and $c_n$ are three parameters dependent on the $\alpha$ and $n$ values calculated under normal operating conditions.

In order to determine $a_n$, $b$ and $c_n$, values, their dependence on parameter $\alpha$ must be determined. Therefore, it has been observed that these coefficients are quadratic functions of $\alpha$ as shown in Figure 3.19, which considers just one fixed $n$ value; the others are similar.

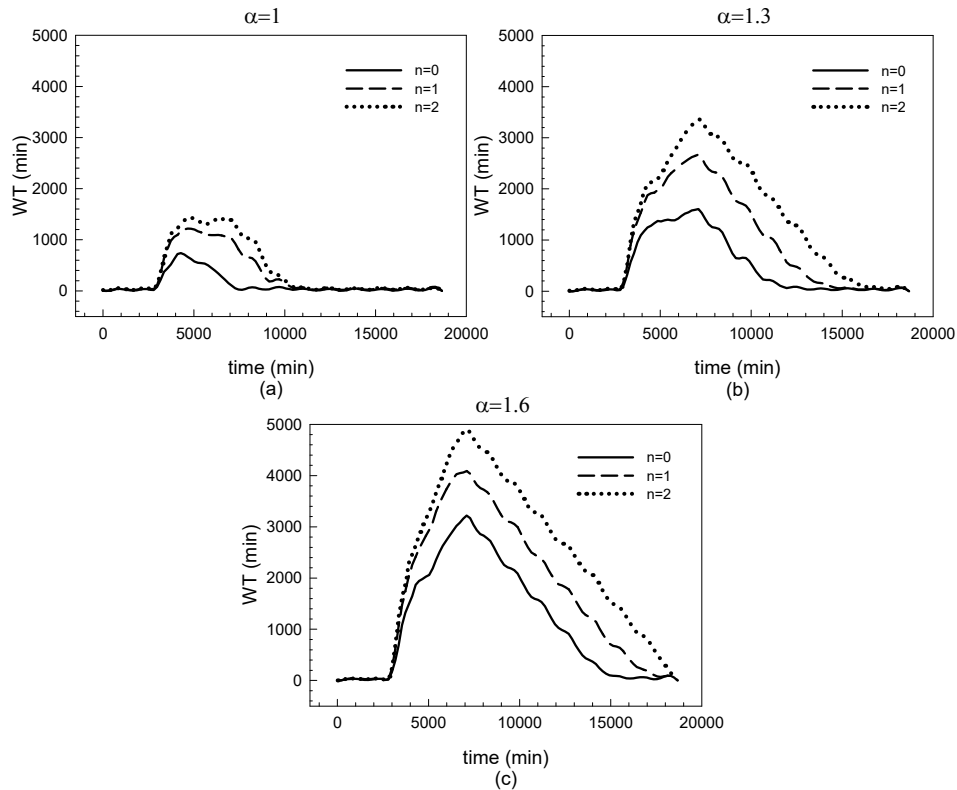

**Figure 3.18** Simulations results for the normal operating condition model with a constant value of $\alpha$ and a variable value of *n* for (a) $\alpha$ = 1; (b) $\alpha$ = 1.3; and (c) $\alpha$ =1.6.

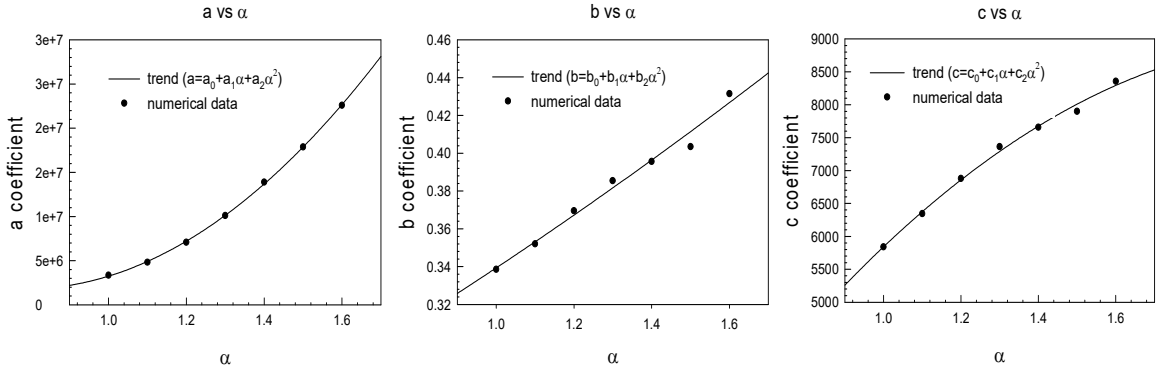

Figure 3.19    Quadratic interpolation, coefficients *a*, *b*, *c* versus $\alpha$ for *n* = 1.

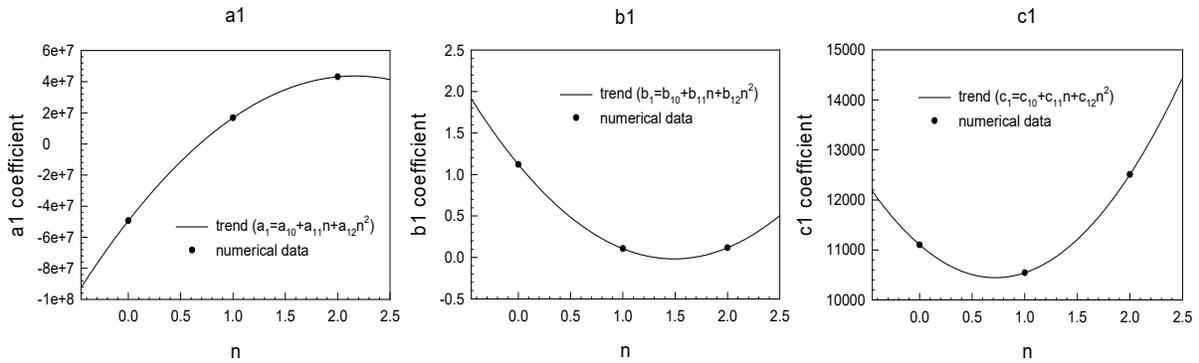

Figure 3.20    Quadratic interpolation, coefficients $a_1$, $b_1$, $c_1$ versus *n*.

Thus, it is possible to express the values of these coefficients as:

$$a_n(\alpha) = a_0 + a_1\alpha + a_2\alpha^2 \tag{3.6}$$

$$b_n(\alpha) = b_0 + b_1\alpha + b_2\alpha^2 \tag{3.7}$$

$$c_n(\alpha) = c_0 + c_1\alpha + c_2\alpha^2 \tag{3.8}$$

The dependence from the *n* parameter was studied next and a quadratic model considered in order to represent coefficients $a_0, a_1, a_2, b_0, b_1, b_3, c_0, c_1,$ and $c_2$ as a function of the number of non-functioning treatment rooms; see Figure 3.20.

Therefore, it is possible to express the values of these coefficients as:

$$\begin{cases} a_0(n) = 21178533,7 - 50687867,5n - 10938560,2n^2 \\ a_1(n) = -49405307,7 + 86079082,9n - 19905188,7n^2 \\ a_2(n) = 31467171,4 - 30777131,8n + 8057254,1n^2 \end{cases} \tag{3.9}$$

$$\begin{cases} b_0(n) = -0,5166 + 1,1094n - 0,3743n^2 \\ b_1(n) = 1,121 - 1,529n + 0,5132n^2 \\ b_2(n) = -0,3514 + 0,5445n - 0,1776n^2 \end{cases} \quad (3.10)$$

$$\begin{cases} c_0(n) = -3955,3 + 3131,5n - 1393,7n^2 \\ c_1(n) = 11100,9 - 1821,2n + 1262,6n^2 \\ c_2(n) = -2328,4 + 45,4n - 200,1n^2 \end{cases} \quad (3.11)$$

Thus, a relationship among PWTs, the time in minutes, and the considered variables—the seismic input ($\alpha$) and the number ($n$) of non-functioning treatment rooms—have been established, representing the proposed meta-model. Therefore:

$$WT(t,\alpha,n) = \frac{a_n(\alpha,n)}{t} * \exp\left\{-0.5 * \left[\frac{\ln\left(\frac{t}{b_n(\alpha,n)}\right)}{c_n(\alpha,n)}\right]\right\} \quad (3.12)$$

By applying Equation (3.9), it is possible to calculate the length of time patients must wait to see a doctor at in the ED at a given instant in time under normal operating conditions.

### 3.5.6 Meta-Model Parameters under Emergency Operating Conditions

Following the same approach used above, the meta-model for the ED under emergency operating conditions has been created as follows:

$$WT(t,n,\alpha) = \frac{a}{t} * \exp\left[-0.5 * \left(\frac{\ln(t/b)}{c}\right)^2\right] \quad (3.13)$$

Coefficients $a$, $b$, and $c$ have been calibrated for the ED working under emergency operating conditions, i.e., the ERP has been implemented. As done earlier, two steps were followed to build the meta-model. First, an amplified seismic input using the scaling factor $\alpha$ was considered. After running Monte Carlo simulations, the data was collected for each analyzed scenario. The results are shown in Figure 3.21 where an increase in patient flow results in a rise in PWTs proportional to the scaling factors. Second, it was assumed that some treatment rooms were rendered non-functional due to structural damage from the earthquake. The values resulting from Monte Carlo simulations are presented in Figure 3.22 for three different $\alpha$ values.

As shown in Figure 3.22, PWTs increase with the number of damaged treatment rooms. Note that when two of the available treatment rooms are not functional, the PWTs reached a peak of about 6000 min (4 days) for a scaling factor $\alpha=1.6$. This value of waiting times is higher with respect to the same conditions when the emergency plan is not applied.

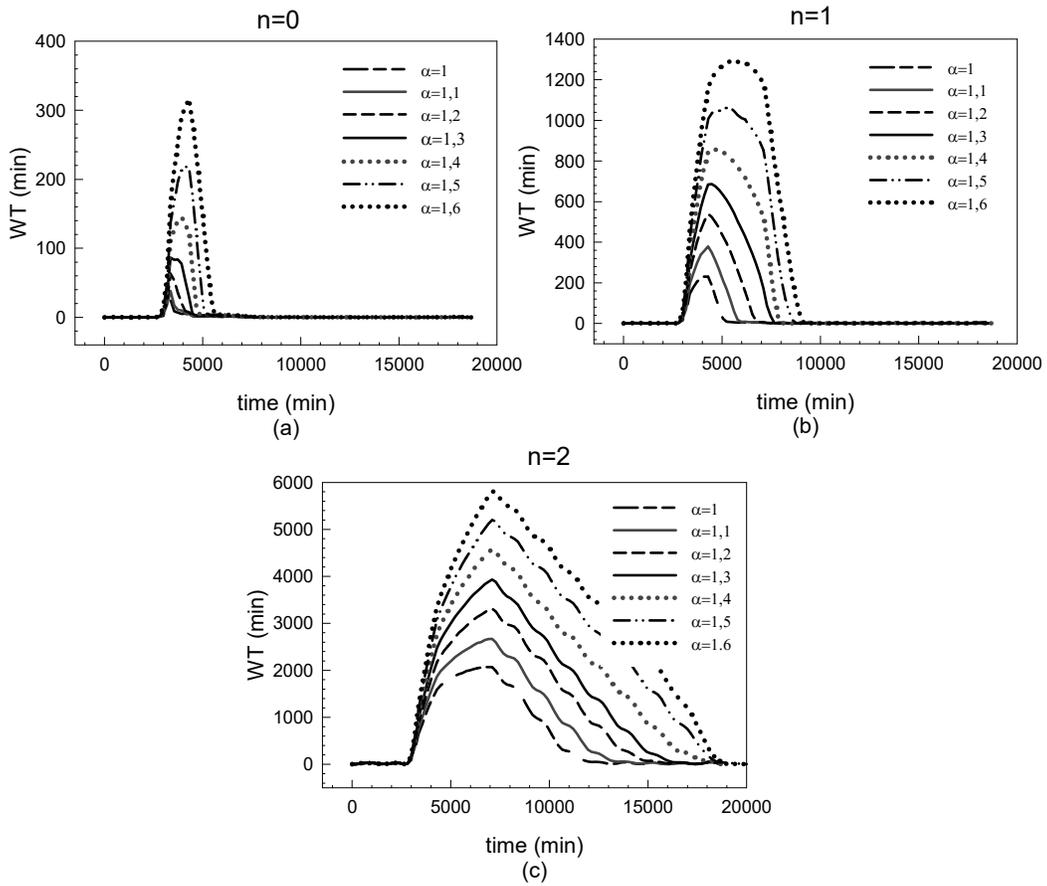

**Figure 3.21** Simulations results for the emergency condition model with a constant value of *n* and a variable value of $\alpha$ for (a) *n* = 0, (b) *n* = 1, and (c) *n* = 2.

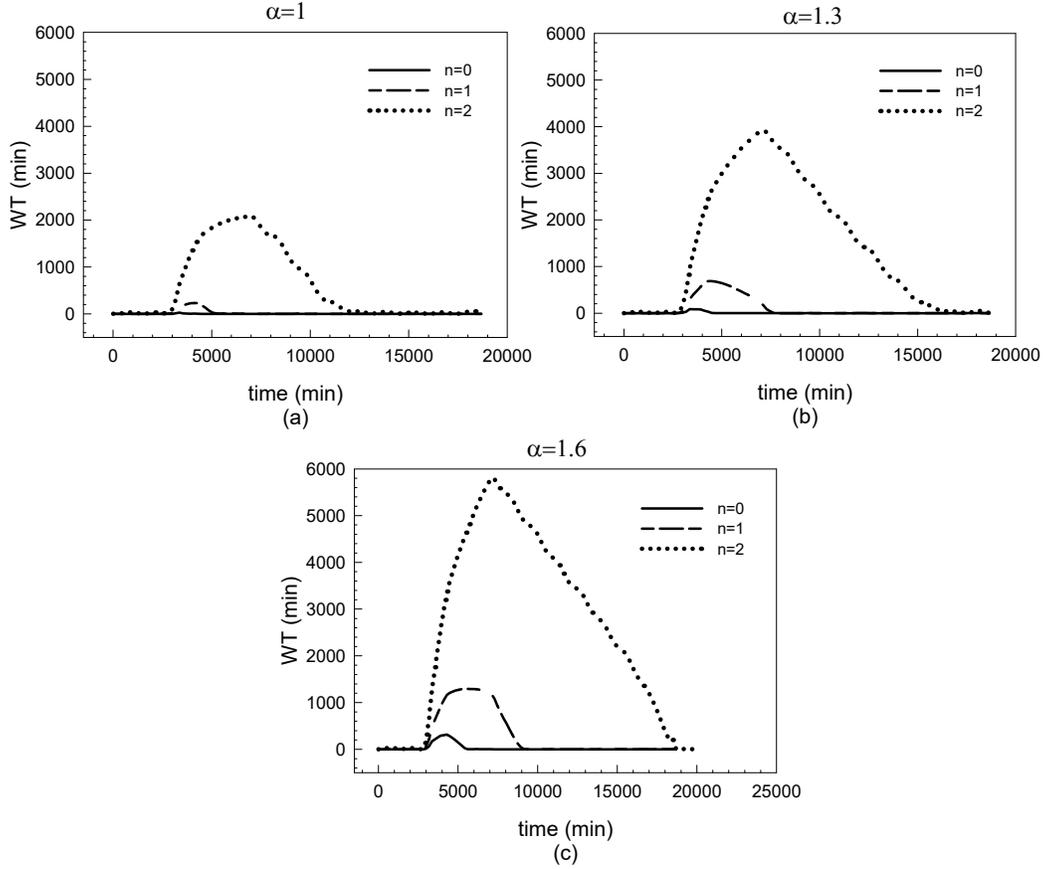

**Figure 3.22** Simulations results for the emergency condition model with a constant value of $\alpha$ and a variable value of $n$ for (a) $\alpha = 1$, (b) $\alpha = 1.3$, and (c) $\alpha = 1.6$.

Under normal operating conditions, there are five treatment rooms for patients coded both green and yellow. Under emergency conditions, there are three treatment rooms coded yellow, which are now located in a different area of the hospital. In addition, the model was built considering that yellow codes have priority over green codes. Thus, in normal operating conditions, when two treatment rooms are closed, yellow-coded patients share a total number of three treatment rooms with green-coded patients. Considering their priority, it could be assumed that yellow-coded patients can use two of these surgeries. In contrast, if two of the available treatment rooms are not functional under emergency operating conditions, yellow-coded patients can be treated only in one surgery. For this reason, when $n = 2$, PWTs under normal operating conditions are lower than PWTs under emergency operating conditions. A lognormal function has been selected to represent the *WT* trend whose general form under emergency operating conditions is:

$$WT(t,n,\alpha) = \frac{a_e}{t} * \exp\left[-0.5 * \left(\frac{\ln(t/b_e)}{c_e}\right)^2\right] \qquad (3.14)$$

where $WT(t, n, \alpha)$ is the waiting time in minutes, $t$ the time in minutes, $n$ the number of non-functioning treatment rooms, $\alpha$ the seismic input and $a_e$, $b_e$, and $c_e$ are three parameters dependent on the $\alpha$ and $n$ values calculated assuming emergency operating conditions.

As aforementioned, to determine $a_e$, $b_e$, and $c_e$, values, first we consider the dependence from the seismic input. These coefficients are quadratic functions of α, as shown in Figure 3.23.

Thus, it is possible to express the values of these coefficients as:

$$a_e(\alpha) = a_0 + a_1\alpha + a_2\alpha^2 \tag{3.15}$$

$$b_e(\alpha) = b_0 + b_1\alpha + b_2\alpha^2 \tag{3.16}$$

$$c_e(\alpha) = c_0 + c_1\alpha + c_2\alpha^2 \tag{3.17}$$

Thus, a quadratic model has been considered in order to represent coefficients $a_0, a_1, a_2, b_0, b_1, b_3, c_0, c_1,$ and $c_2$ as a function of $n$ as shown in Figure 3.24 to represent the number of non-functional treatment rooms.

It is possible to express the values of these coefficients as:

$$\begin{cases} a_0(n) = 4313145 + 13231212,6n - 9439291,9n^2 \\ a_1(n) = -8170064,6 - 25095914,1n - 14299370,7n^2 \\ a_2(n) = 3947395,5 + 6797542,2n + 1122876,7n^2 \end{cases} \tag{3.15}$$

$$\begin{cases} b_0(n) = -0,1195 - 1,099n + 0,6206n^2 \\ b_1(n) = 0,1625 + 1,728n - 0,8719n^2 \\ b_2(n) = 0,0033 - 0,61n + 0,3148n^2 \end{cases} \tag{3.18}$$

$$\begin{cases} c_0(n) = 3304,5 - 6345,4n + 3260,9n^2 \\ c_1(n) = -939,3 + 8878,9n - 3687n^2 \\ c_2(n) = 945,1 - 2823,8n + 1415,2n^2 \end{cases} \tag{3.19}$$

As shown above, there is a relationship among PWTs, the time in minutes, and the seismic input ($\alpha$) and number ($n$) of non-functional treatment rooms. The set of the illustrated equations provides the time behavior of the ED under emergency operating conditions. By applying Equation (3.18), it is possible to calculate the length of time patients must wait to see a doctor in a treatment room at a given instant in time under emergency operating conditions.

$$WT(t,\alpha,n) = \frac{a_e(\alpha,n)}{t} * \exp\left[-0.5 * \left(\frac{\ln\left(\frac{t}{b_e(\alpha,n)}\right)}{c_e(\alpha,n)}\right)\right] \tag{3.20}$$

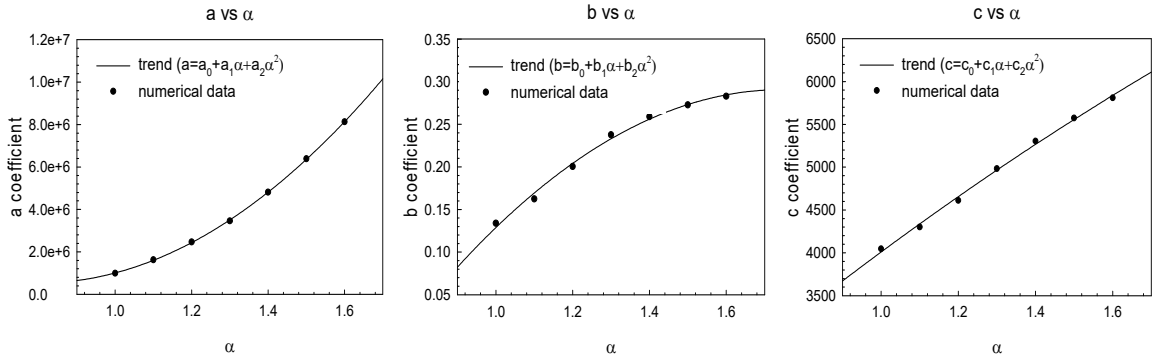

Figure 3.23  Quadratic interpolation, coefficients *a*, *b*, and *c* vs $\alpha$ for *n* = 1.

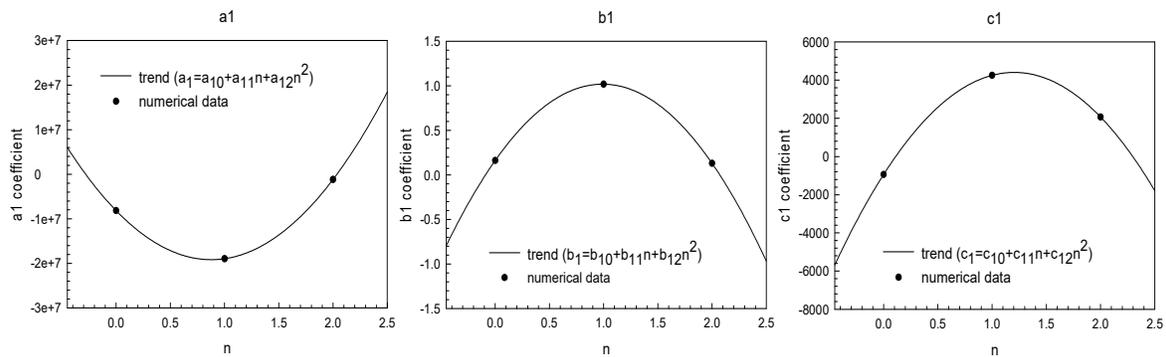

Figure 3.24  Quadratic interpolation, coefficients $a_1$, $b_1$, and $c_1$ versus *n*.

### 3.5.7 Comparison between the Meta-Model vs. DES model

The previous section presented meta-model for both normal and emergency operating conditions. This mathematical approximation avoids the problem of the computation time required for the more complex DES model. Simulation models are often too complex and require extensive computational capacity for use in conducting a sensitivity analysis based on Monte Carlo methods. Because the meta-model proposed herein is an approximation of the simulation model, its accuracy with respect to the experimental data must be determined.

To determine if the meta-model adequately represents the behavior of the output generated by a simulation model, the meta-model's output must to be compared with the results obtained from the DES. In this specific case, the PWTs trends for both the meta-model and DES model were considered. A comparison was performed between PWTs given by the DES and analytical model for a different number of non-functional treatment rooms and for different intensity of seismic input.

First, normal operating conditions were analyzed. The results for two different $\alpha$ and *n* combinations are illustrated in Figure 3.25. Next, emergency operating conditions were taken into account. A comparison between the proposed meta-model and the simulations results are shown in Figure 3.26 considering two different $\alpha$ and *n* combinations.

Considering the peak PWTs values, the error between the experimental data and the analytical model were calculated for both normal and emergency operating conditions; see Table 3.3. Results show that the meta-models for normal and emergency operating conditions provide an accurate description of ED behavior. Overall, it can be seen that the relative error between the proposed meta-model and the simulations results is rather low for all the possible $\alpha$ and $n$ combinations.

Table 3.4 illustrates that per the meta-model, the maximum estimated error under normal operating conditions is equal to 10.81% when all treatment rooms are operational, and $\alpha$ is 1.2. When the emergency plan is applied, the estimated error reaches a maximum value of 11.21% if there is one non-functional treatment room, and $\alpha$ is 1.5. Consequently, it has been concluded that the proposed meta-models show good agreement with the experimental results in describing the system response. Emergency response planners can apply the meta-models with confidence, thereby replacing the more complex DES model.

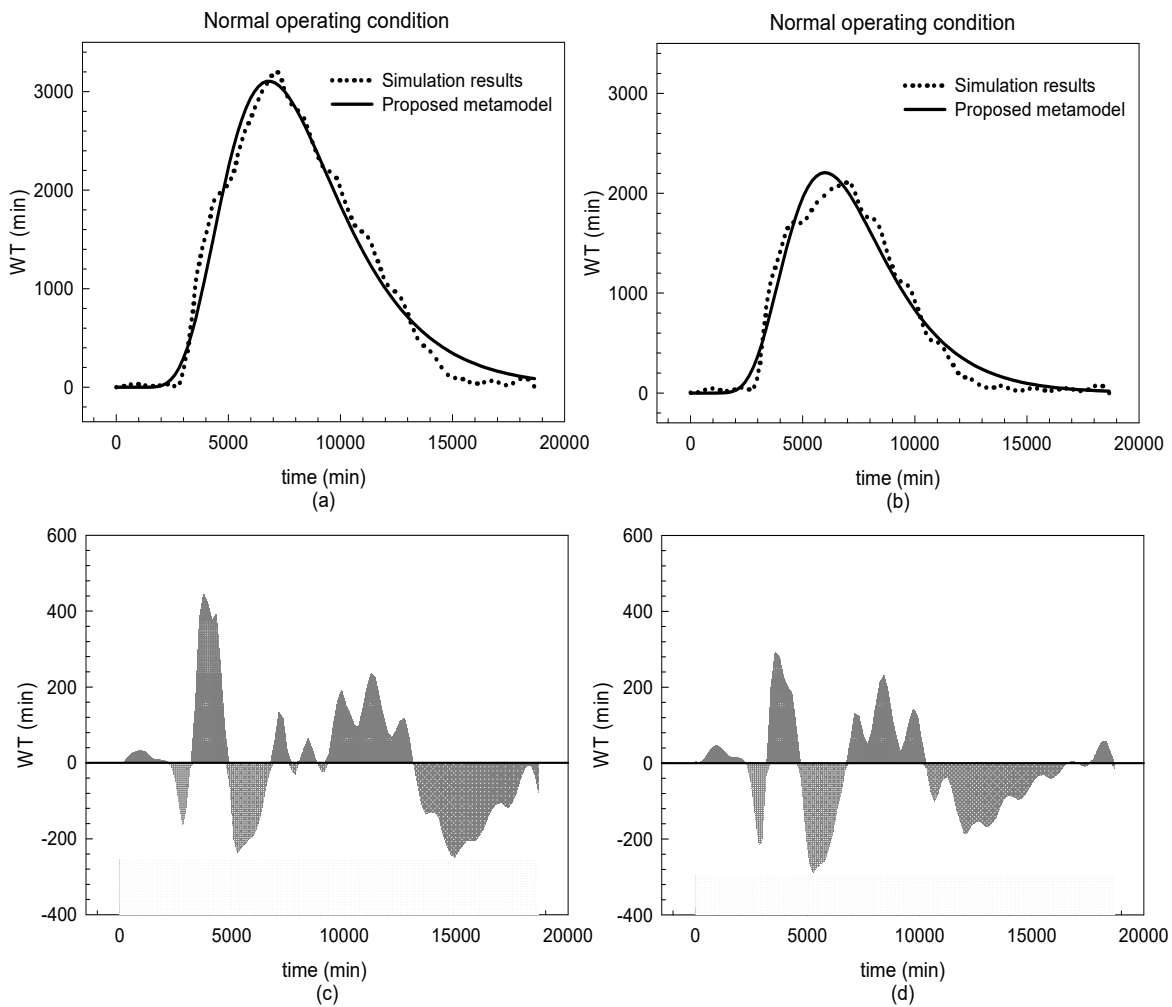

Figure 3.25  Simulation data vs. proposed meta-model under normal operating conditions for (a) $n = 0$, $\alpha = 1.6$, (b) $n = 1$, $\alpha = 1.2$. and(c), (d) error bars.

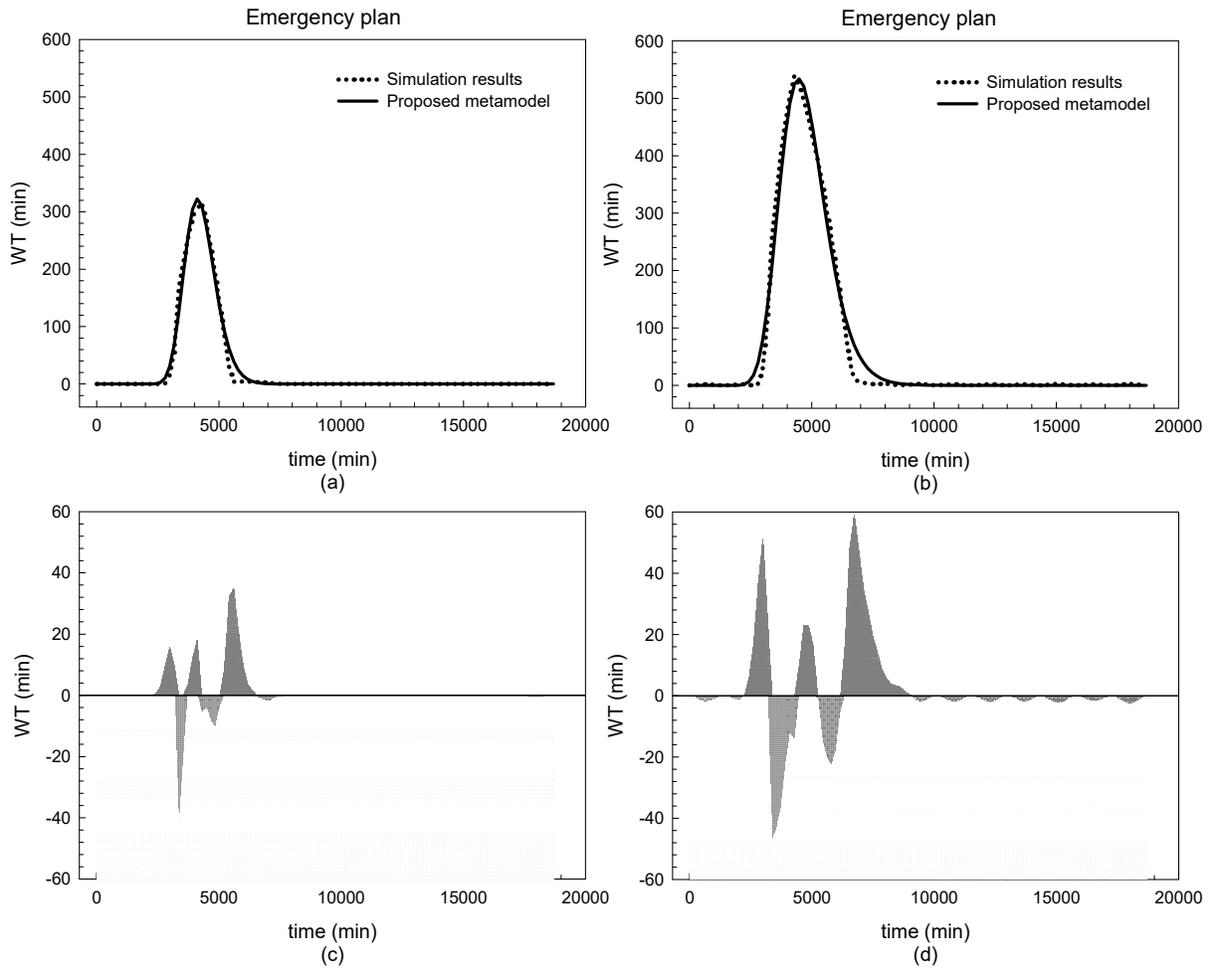

**Figure 3.26** Simulation data vs proposed meta-model under emergency plan condition for (a) $n = 0$, $\alpha = 1.6$; (b) $n = 1$, $\alpha = 1.2$; and (c) and (d) error bars.

**Table 3.3** Error between the proposed meta-model and the simulations results evaluated at the peak value under normal operating conditions.

| Scale factor $\alpha$ | Error (%), $n = 0$ | Error (%), $n = 1$ | Error (%), $n = 2$ |
|---|---|---|---|
| 1 | 5.43% | 2.94% | 7.53% |
| 1.1 | 3.84% | 8.96% | 5.44% |
| 1.2 | 10.81% | 4.35% | 1.03% |
| 1.3 | 2.23% | 0.37% | 1.11% |
| 1.4 | 2.6% | 2.72% | 4.4% |
| 1.5 | 3.22% | 1.35% | 3.26% |
| 1.6 | 0.32% | 1.002% | 3.92% |

**Table 3.4**  Error between the proposed meta-model and the simulations results evaluated at the peak value under emergency operating conditions.

| Scale factor $\alpha$ | Error (%), $n = 0$ | Error (%), $n = 1$ | Error (%), $n = 2$ |
|---|---|---|---|
| 1   | 8%    | 9.17%  | 5.31% |
| 1.1 | 15.2% | 1.05%  | 3.71% |
| 1.2 | 7.93% | 1.11%  | 0.93% |
| 1.3 | 8.13% | 5.24%  | 0.38% |
| 1.4 | 6.89% | 8.96%  | 1.63% |
| 1.5 | 7.33% | 11.21% | 1.92% |
| 1.6 | 1.89% | 9.82%  | 2.41% |

## 3.6 THE GENERAL META-MODEL

### 3.6.1 Problem Formulation

The meta-model presented above was developed to adequately represent, in real time, the dynamic response of the Mauriziano Hospital's ED considering whether or not the emergency plan is in place. It describes a number of hospital scenarios and takes into account variable seismic input and possible structural damage as a result of the earthquake, which could render some treatment rooms non-functional. Thus, the proposed Mauriziano Meta-model presents a sensitivity analysis under the considered variables.

This meta-model is valid only for the case-study hospital and cannot be used to estimate the behavior of other EDs under similar seismic conditions. In addition, long-running times complicate the evaluation of a hospital's performances using DES models. For these reasons, a general meta-model analyzes the capacity of healthcare facilities to cope with and respond to a catastrophic event such as an earthquake will provide hospitals with a tool to assess in advance a hospital's resilience.

This research has proposed and constructed a general meta-model for general application. As done for the Mauriziano meta-model, PWTs were chosen as the main parameter of response to describe the efficiency of healthcare facilities. It is dependent on both internal and external factors. The main goal is to create a model that includes significant parameters that can describe the trend of PWTs for any hospital in the event of a disaster. Considering that each hospital is substantially different from another, the problem is complex, and a considerable number of variables are needed to describe the behavior of any healthcare facility. Given that it is impossible to create a general model with the same amount of detail as a model developed specifically for one hospital, it is essential simplify the problem; therefore, the number of variables must be reduced.

Based on a study of the simulation results carried out for the Mauriziano Hospital and of the emergency plan's structure, three parameters have been identified as the most significant inasmuch they affect PWTs. The total number of emergency rooms, the number of doctors/nurses per color-coded area, and seismic input were considered as the most representative parameters to characterize a generic ED.

Therefore, the general meta-model was constructed using the following approach. First, a mathematical equation for the model was chosen, and a log-normal function was considered based on the meta-model developed for the studied hospital. Then, the parameters of the selected model were estimated. Finally, the general meta-model was validated against the results of the DES simulations. Thus, a hospital's dynamic behavior can be obtained virtually and instantaneously in real time.

### 3.6.2 Assumptions

Before developing the general meta-model, some assumptions were made. First, the total number of doctors in each color-coded area varies proportionately with the number of the available treatment rooms and which each treatment room is assigned one doctor and one or more nurses. Therefore, the meta-model assumed that the number of doctors for each ED is the same as the number of available treatment rooms. Therefore,

$$d = m \tag{3.21}$$

where $d$ represents the total number of doctors in a considered color-coded area, and $m$ is the total number of the available treatment rooms in the same color-coded area. Only doctors in the color-coded areas were considered, and the doctors working in others ED areas (triage, waiting rooms, etc.) were ignored. This hypothesis is reasonable because treatment rooms are equipped to administer care to only one patient at a time.

Second, the meta-model was constructed using as input data patient arrival rates per data used for the Mauriziano meta-model. An amplified seismic input was considered, as discussed above.

Then, the meta-model was developed only for those patients whose condition have been coded "yellow," those patients deemed critical. During an emergency situation, the greater the PWTs, the more likely an injury will worsen, therefore, PWTs are a critical parameter for evaluating hospital performance. This model assumed that all patients are coded "yellow." It is certainly possible to build a meta-model also for patients labeled "green," but it hasn't been developed herein.

### 3.6.3 Development of the General Meta-Model

As discussed earlier, a hospital's resilience during an emergency event is evaluated by how quickly injured patients receive treatment, which is directly correlated to the PWTs. Therefore, PWTs were chosen as the main response variable of the meta-model.

The *WT* function was determined considering three variables simultaneously: seismic input ($\alpha$); the total number of emergency rooms per color-coded area ($m$); and the time in minutes ($t$). A lognormal function was considered to describe the trend in PWTs, which describes the time that patients have to wait before receiving care at a given instant in time. Below is the general form of the lognormal equation used for the *WT*:

$$WT(t,\alpha,m) = \frac{a(\alpha,m)}{t} * \exp\left[-0.5 * \left(\frac{\ln\left(\frac{t}{b(\alpha,m)}\right)}{c(\alpha,m)}\right)\right] \qquad (3.22)$$

where $WT$ is the patient wait time in minutes, $\alpha$ is the seismic input, $m$ is the total number of emergency rooms per color-coded area, $t$ is the time in minutes and $a, b, c$ are nonlinear regression coefficients. The values of the coefficients were calculated considering their dependence from the seismic input. As in the previous case, $a, b, c$ coefficients are considered quadratic functions of $\alpha$.

$$a(\alpha) = a_0 + a_1\alpha + a_2\alpha^2 \qquad (3.23)$$

$$b(\alpha) = b_0 + b_1\alpha + b_2\alpha^2 \qquad (3.24)$$

$$c(\alpha) = c_0 + c_1\alpha + c_2\alpha^2 \qquad (3.25)$$

In contrast, the dependence of the coefficients $a_0$, $a_1$, $a_2$, $b_0$, $b_1$, $b_2$, $c_0$, $c_1$, and $c_2$ by the parameter $m$ cannot be analyzed as was done in the previous case. For the general meta-model, the total number of treatment rooms per colored areas ($m$) has been considered instead of the number of nonfunctional treatment rooms ($n$) that were included in the meta-model for the Mauriziano Hospital. Therefore, in order to study the reliance of $WT$ on the total number of treatment rooms ($m$), the following approach has been applied. The DES model developed for the Mauriziano's ED was used, and Monte Carlo simulations were run increasing the total number of treatment rooms under a constant number of incoming patients. This method was then applied to a number of different hospitals with a different number of surgeries, and the results from the simulations collected and analyzed.

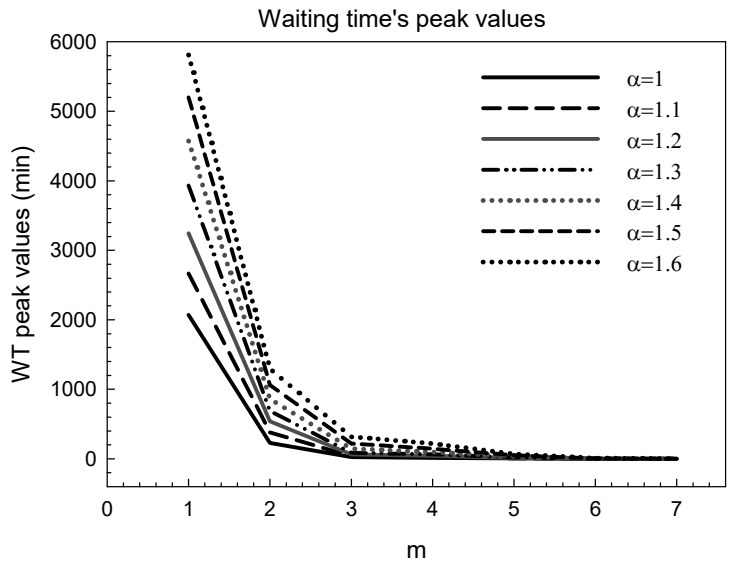

**Figure 3.27**     Variation of the peak wait-time values for a variable total number of treatment rooms.

As expected, PWTs decreased with increasing the total number of treatment rooms. In particular, by examining the output response from simulations, an exponential function is appropriate to describe the trend of the peak PWTs values when there is an increase in number of treatment rooms. The results are shown in Figure 3.28 for different $\alpha$ values.

The trend of the peak $WT$ can be approximated using a single exponential function, whose general form is:

$$WT_{max} = ae^{(-bm)} \tag{3.26}$$

where $WT_{max}$ represents the maximum waiting time corresponding to the peak value, $m$ is the total number of treatment rooms coded yellow, and $a, b$ are coefficients that depend on the seismic input. The coefficients $a, b$ are given by Equations (3.27) and (3.28):

$$a = 18534,7 + \frac{6787,1}{1 + \frac{-(\alpha - 1,28)}{0,059}} \tag{3.27}$$

$$b = 5,01 - 3,9465\alpha + 1,0981\alpha^2 \tag{3.28}$$

where $\alpha$ is the considered seismic input. In view of these considerations, coefficients $a_0$, $a_1$, $a_2$, $b_0$, $b_1$, $b_2$, $c_0$, $c_1$, and $c_2$ can be expressed in function of the total number of treatment rooms ($m$). Their mathematical expressions have been determined using curve fitting procedures. After fitting different equations, Equation (3.29) has been selected:

$$g(m) = x_1 + m^4 \left( f + \frac{d}{m} + \frac{c}{m^2} + \frac{b}{m^3} \right) \tag{3.29}$$

where $m$ is the total number of treatment rooms in the yellow-coded areas, and $x_1$, $b$, $c$, $d$, and $f$ are coefficients whose values are specified below:

$$a_0(m) = -89323896 + m^4 \left( -1106445 + \frac{14701509}{m} - \frac{69987233}{m^2} + \frac{138734467}{m^3} \right) \tag{3.30}$$

$$a_1(m) = 132611723 + m^4 \left( 2072754 - \frac{26999059}{m} + \frac{124474864}{m^2} - \frac{233300000}{m^3} \right) \tag{3.31}$$

$$a_2(m) = 16657792 + m^4 \left( -543784 + \frac{6227391}{m} - \frac{22646870}{m^2} + \frac{22339458}{m^3} \right) \tag{3.32}$$

$$b_0(m) = 5.57 + m^4 \left( 0.08 - \frac{1.04}{m} + \frac{4.89}{m^2} - \frac{9.34}{m^3} \right) \tag{3.33}$$

$$b_1(m) = -7.65 + m^4 \left( -0.12 + \frac{1.58}{m} - \frac{7.34}{m^2} + \frac{13.67}{m^3} \right) \tag{3.34}$$

$$b_2(m) = 2.79 + m^4 \left( 0.04 - \frac{0.54}{m} + \frac{2.54}{m^2} - \frac{4.78}{m^3} \right) \tag{3.35}$$

$$c_0(m) = 28475.3 + m^4 \left( 338.6 - \frac{4684.3}{m} + \frac{22726}{m^2} - \frac{43551.1}{m^3} \right) \tag{3.36}$$

$$c_1(m) = -43772 + m^4 \left( -578.5 + \frac{8013.6}{m} - \frac{38812}{m^2} + \frac{74209.6}{m^3} \right) \tag{3.37}$$

$$c_2(m) = 11604.2 + m^4 \left( 123.1 - \frac{1811}{m} + \frac{9196.2}{m^2} - \frac{18167.4}{m^3} \right) \tag{3.38}$$

Thus, a general meta-model valid for any considered hospital was developed. It provides PWTs considering representative variables. At a given instant in time, the *WT* is given by a lognormal function in which the considered parameters are the total number of treatment rooms in the yellow-coded area and the seismic input.

### 3.6.4 Meta-Model Validation

To ensure that the proposed general meta-model is generally applicable to any hospital considering all the specified assumptions, a validation process has been implemented. To apply this model generally, the meta-model must approximate real hospital behavior; this might introduce errors in the model. A meta-model without sufficient accuracy may lead to an incorrect assessment of the ED's performance. Therefore, validation of the meta-model is critical.

This research compared the experimental data obtained from the DES model with analytical results. The goal of the simulation model is to capture very closely the response of the "real" system. Two healthcare facilities were considered to verify the model's efficacy. For both hospitals, Monte Carlo simulations were performed using the same patient arrival rates and the *WT*'s average curve.

First considered was the Mauriziano Hospital; see Figure 3.28. The results from the DES model were compiled and compared with the proposed meta-model. Input data for the meta-model considered a total of three rooms in the Mauriziano's ED ($m\ \alpha = 3$). In addition, a variable seismic input was taken into account, varying the multiplicative scale-factor values from 1 to 1.6. The results are illustrated for $\alpha = 1$ and $\alpha = 1.6$.

The San Sepolcro Hospital, located in Arezzo (Italy), was also considered, and a similar comparison between the DES and the meta-model results performed. As done before, a comparison between simulation and analytical model was conducted. Input data for the meta-model considered four rooms in the San Sepolcro's ED ($m = 4$). The results are shown in Figure 3.29 for a variable seismic input.

By analyzing results shown in Figure 3.28 and Figure 3.29, note that for both hospitals the proposed meta-model was able to describe accurately the behavior of the EDs. The *WT* curves obtained from the analytical meta-model show good agreement with the experimental data.

The error between the simulation data and the analytical meta-model at the peak value is given in Table 3.5 for the different $\alpha$ values and for both Mauriziano and San Sepolcro Hospitals. The highest error at the peak value is 19.6% for Mauriziano Hospital obtained for $\alpha = 1$, whereas the maximum error for San Sepolcro Hospital evaluated at the peak value is 25.4% when $\alpha = 1.1$. Therefore, the proposed general meta-model shows a good match to the simulation results.

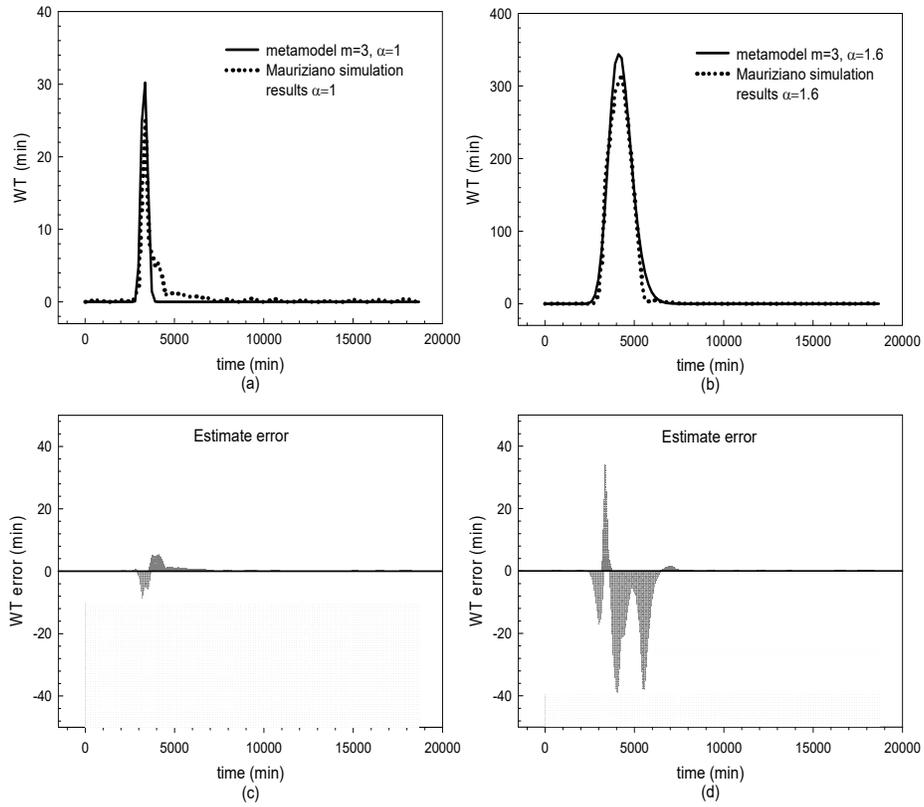

**Figure 3.28**    **Comparison between analytical meta-model and Mauriziano Hospital's experimental data for (a) $\alpha = 1$, (b) $\alpha = 1.6$, and (c) and (d) error bars.**

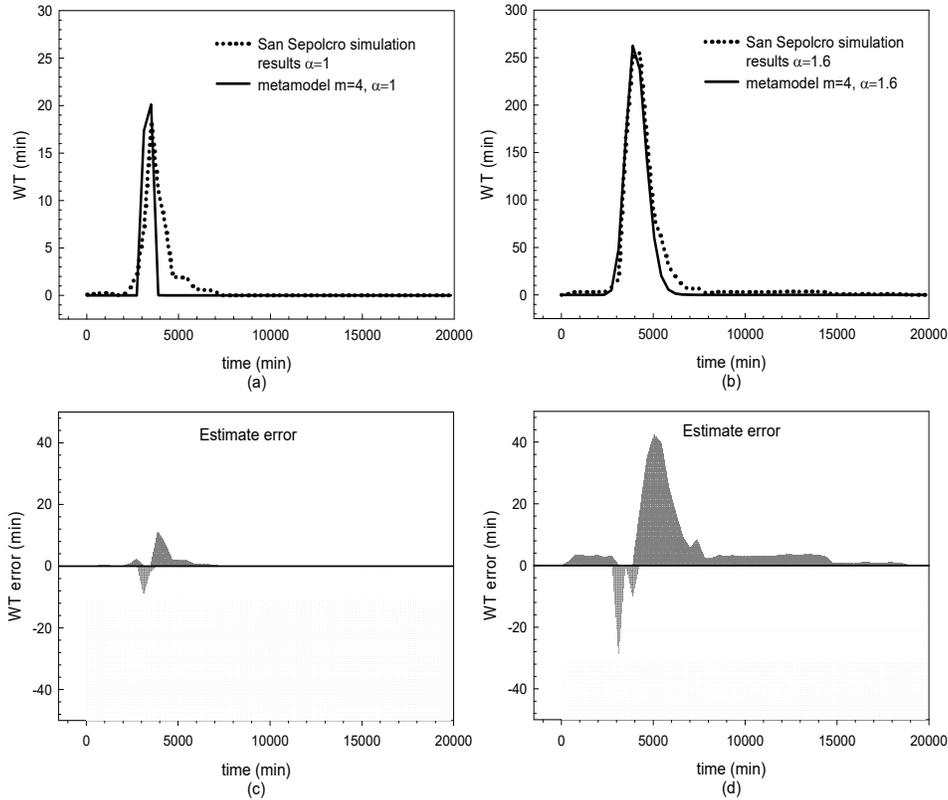

**Figure 3.29**  Comparison between analytical meta-model and San Sepolcro Hospital's experimental data for (a) $\alpha = 1$, (b) $\alpha = 1,6$, and (c) and (d) error bars.

**Table 3.5**  Error between the experimental and analytical model evaluated at the peak value for Mauriziano and San Sepolcro Hospitals.

| Scale factor $\alpha$ | Error (%) Mauriziano ED | Error (%) San Sepolcro ED |
|---|---|---|
| 1 | 19.6% | 10.7% |
| 1.1 | 16.9% | 25.4% |
| 1.2 | 13.8% | 24.3% |
| 1.3 | 9.3% | 21.2% |
| 1.4 | 17.2% | 15.3% |
| 1.5 | 13.1% | 5.1% |
| 1.6 | 5.9% | 1.7% |

## 3.7  REMARKS AND CONCLUSIONS

Healthcare facilities play a key role in our society, especially during and immediately following a disaster. A geographic region many be vulnerable to several potential hazards, and it is critical that hospitals continue to operate in the event of an earthquake. A healthcare facility must remain accessible and able to function at maximum capacity, providing its services when they are most needed. Discrete event simulation (DES) is a powerful tool to represent complex systems such as hospitals. Since the 1980s, it has been used widely in the medical industry to study hospital response. Different parameters have been considered over the years where patient wait times (PWTs) have been identified as the main parameter in order to evaluate the ED resilience. Presented herein, a hospital in Italy has been considered as a case study. A DES model has been built for the hospital's ED, taking into account the existing emergency response plan (ERP). The trends in PWTs obtained when the ERP is applied were compared with normal operating trends. Considering the same input data for both models, PWTs under emergency operating conditions are lower than those in normal operating conditions.

    Because building a DES model is time consuming; a simplified model called meta-model was developed, which is defined as an analytical framework representing the relationship between the system response and some selected variables without the need to run the model several times. In order to build the meta-model, different scenarios were considered, taking in account the amplitude of the seismic input and possible structural damage to the ED. After verification of the meta-model's accuracy of the case study, a general meta-model was developed to provide hospitals with a useful tool able to evaluate in real time the ED behavior.

# 4 Hospital Emergency Network: Earthquake Impact in San Francisco

## 4.1 INTRODUCTION

California is one of the most seismically active parts of the world. About once every 150 years there is a major seismic event on the San Andreas fault, which runs the length of state. This fault has reached a sufficient stress level that there is concern that the next "big one," a hypothetical earthquake of very high magnitude, is imminent. Such an earthquake could cause enormous damage to the major cities, especially in the San Francisco Bay Area [Poland 2008].

Even though it is impossible to predict exactly where or when the next major earthquake will occur in California, it is not impossible to predict its effects. Indeed, the ability to be prepared for large extreme events is critical for saving lives and reducing earthquake's heavy consequences. Key to the recovery efforts in the event of such a disaster is the structural integrity of healthcare facilities. Hospital management during emergencies requires coordinated processes and resources, including situational awareness [Downey et al. 2012]. That is why the concept of resilience–which may be defined as the capability of social units to adapt and maintain their functions during emergencies, organizing themselves in order to minimize the effects of disasters–has become increasingly important when planning the emergency strategies. In particular, hospitals play a key role during and after catastrophic events [Cimellaro et al 2010a]. When a disaster occur, hospitals not only have to provide care to a large number of casualties in a setting of limited resources, but in addition they have to collaborate and cooperate effectively with other healthcare facilities.

Several studies have been conducted over the years to understand how medical professionals may coordinate the rescue and relief work after a big natural disaster [Zhou et al. 2014]. Indeed, because earthquakes are one of the most catastrophic events, different approaches have been proposed to assess different disaster scenarios. For example, Hashemi and Alesheikh [2013] developed a multi-agent simulation model in order to evaluate the earthquake impact in Tehran, Iran.

In this chapter, a magnitude 7.2 earthquake on the San Andreas fault has been selected in order to analyze the types of consequences the city of San Francisco can expect following this strong event. In particular, the performance of San Francisco hospitals has been studied. The goal is to understand whether healthcare facilities are able to provide emergency care to all the injured in a timely and efficient manner. To achieve this goal, several steps have been followed in order to develop the hospitals' emergency network.

First, the number of injuries that could occur in the selected earthquake scenario has been obtained. A total number of 3650 patients were assumed, distributed by neighborhoods and taking into account the amount of damage that each neighborhood can experience.

Second, the six San Francisco's emergency departments (EDs) have been considered, and all the patients were distributed in the hospitals; it was assumed that during emergencies patients are directed to the closest hospital to them. Thus, the number of patients arriving in each hospital was obtained for the considered earthquake scenario.

Then, patient wait times (PWTs) were chosen as the most significant parameter in order to measure hospital resilience during emergencies. A maximum waiting time of 3 hours (180 min) was considered based on interviews with hospital staff; this the metric used to evaluate the performance of the San Francisco's EDs.

Finally, a hospital network was created using PWTs as input data. Two options were considered in order to mitigate the risk to exceed the maximum time each patient should have to wait. Even if the results of this analysis are limited to certain boundaries and assumptions, those limitations are manageable from the point of view of operational decision makers, which can evaluate rapidly how this urban area can manage an emergency situation.

## 4.2 METHODOLOGY

The main purpose of this study is to analyze San Francisco's hospitals network during emergencies. The goal is to check whether they will be able to work together in order to coordinate and deliver a broad spectrum of services to the San Francisco community after a seismic event. Along these lines, it is possible to relate the organizational aspects of healthcare facilities with a hospital's resilience by measuring the quality of care provided during the disaster. The quality of care provided by the San Francisco's hospitals could be defined using the PWTs in the emergency department (ED) before receiving care. So, PWTs have been chosen as the main parameter to evaluate the response of hospitals during hazardous event such as an earthquake. In particular, a maximum waiting time of 3 hours (180 min) has been considered as the maximum time patients can wait before their medical conditions begin to deteriorate.

As shown in Figure 4.1, this study selected a magnitude 7.2 earthquake on the San Andreas fault. After defining the earthquake scenario, the number of injuries that could occur in the city was obtained. All patients were distributed among the six San Francisco's EDs by considering the distance from each of them and the number of patients arriving at each hospital. Two solutions have been proposed to avoid exceeding the acceptable maximum waiting time. The first one considers the possibility of redistributing the injured through a Control Center.

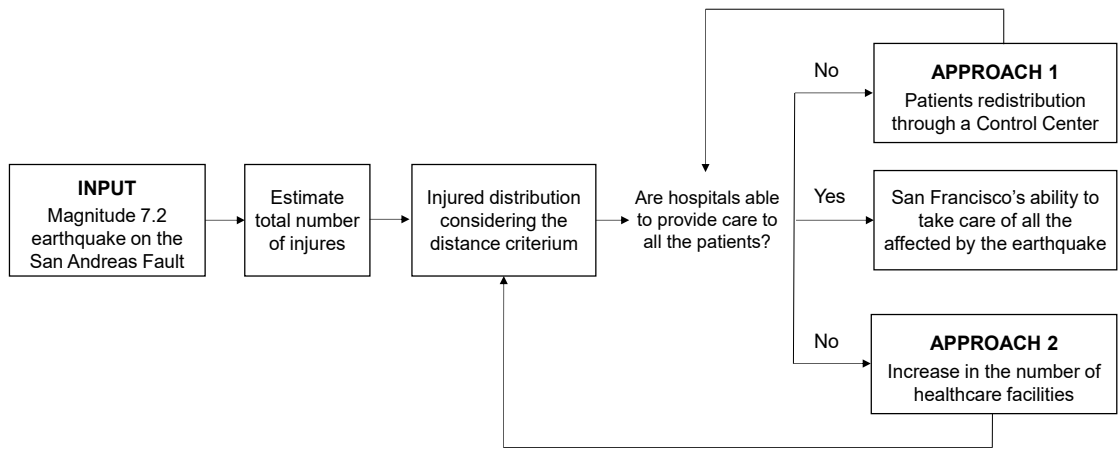

Figure 4.1     Research logic framework.

## 4.3   ASSUMPTIONS

Assumptions made before developing hospital networks analysis include:

1. The magnitude on which fault the earthquake will occur must be defined. For the purpose of defining hospital disaster resilience, a magnitude 7.2 earthquake located on the Peninsula segment of the San Andreas fault was selected. It is also used by the San Francisco's CAPSS (Community Action Plan for Seismic Safety).

2. Damages sustained by infrastructure have not been taken into account in this work. Although it is anticipated that there will be damage to freeways and overpasses in the event of a 7.2 magnitude earthquake, this is considered negligible in terms of surface street access.

3. Number of injured is distributed in proportion to the damage of residential buildings. This hypothesis is reasonable because the main reason why people are injured during a seismic event is because of damaged or collapsed buildings; thus a high relationship between these two variables could be supposed. In this study, possible casualties from other causes have not been considered.

4. The patient arrival rate used as input for the meta-model is the same collected in a California hospital during 1994 Northridge earthquake; see Cimellaro et al. [2011] for the pattern of the Northridge arrival rate, which was used for this model.

5. Only six San Francisco's hospitals were considered in this study. There are more than six hospitals in San Francisco, but some of them do not have an ED and others are specialized hospitals (e.g. children's specialized hospital, geriatric psychiatry hospital, etc.). For the purpose of this chapter, only general hospitals provided with a functioning ED were analyzed.

## 4.4 ESTIMATING DISTRIBUTION OF INJURED

A large-scale disaster like an earthquake affects a large number of people, with those sustaining injuries ranging from minor to severe. Depending on the size and time of the earthquake, the number and the severity of the injured could vary considerably. This research considers a magnitude 7.2 earthquake on the San Andreas fault occurring at rush-hour. According to Tobin and Samant [2009], four severity levels have been considered. Severity 1 refers to minor injuries, meaning that patients need basic medical care that can be administered by nurses or paraprofessionals. These kinds of patients could be associated with white or green triage codes. Severity 2 represents serious injuries that require a greater degree of medical care. These patients present a partial impairment of vital functions so they could be coupled with yellow triage codes. Severity 3 are patients with severe injuries, which means compromised vital functions. They need an immediate medical care and could be labeled with either yellow or red triage codes. Finally, Severity 4 refers to critical injuries level. These patients are mortally injured or their lives are at risk; they are associated with red triage codes. The number of injured for each severity level is shown in Table 4.1. This research considered patients with yellow and green codes only in order to study San Francisco's hospitals network; severity 4 injuries were not taken into account in this chapter.

A total of 3650 injured with a yellow triage code was considered. This number has been obtained by summing the highest injured number considering both severity 2 and severity 3 and the expected number of patients with yellow and green codes included in the severity 2 level. This means that San Francisco's hospitals have to provide care to 3650 patients with earthquake-related injuries.

Earthquakes affect different parts of a city in different ways due to each location's proximity to faults, underlying soil, and types of buildings. Therefore, for San Francisco the distribution of the injured is not homogeneous over the city. For this reason, patients have been distributed by neighborhood. The city has been divided into sixteen neighborhoods according to the Department of Public Works, including Bayview, Downtown, Excelsior, Ingleside, Marina, Merced, Mission, Mission Bay, North Beach, Pacific Heights, Richmond, Sunset, and Twin Peaks [Tobin and Samant 2009]. San Francisco's six EDs have been considered as shown in Figure 4.2.

**Table 4.1**   Estimated injuries in the considered earthquake scenario.

| Levels of severity | Casualties |
|---|---|
| Severity 1 | 3200 to 5600 |
| Severity 2 | 760 to 1300 |
| Severity 3 | 90 to 150 |
| Severity 4 | 170 to 300 |

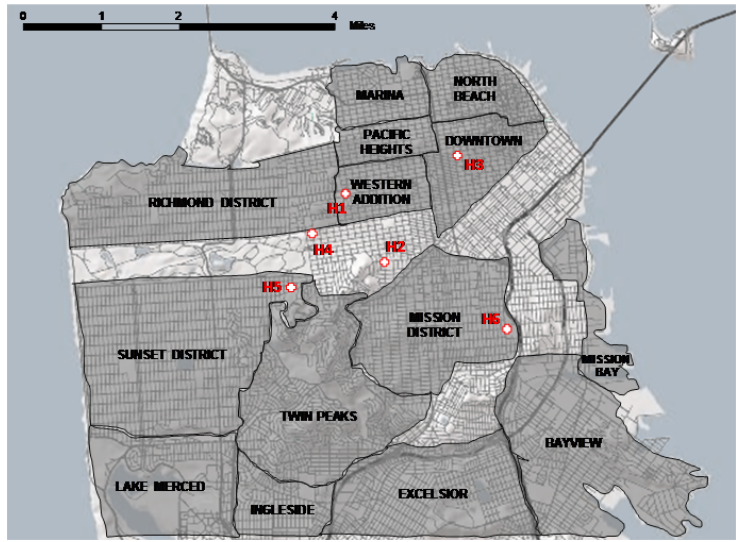

**Figure 4.2    San Francisco's neighborhoods and hospitals.**

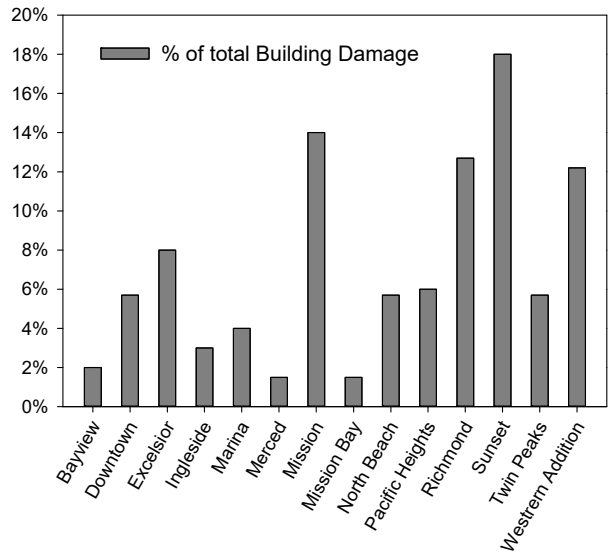

**Figure 4.3    Percentage of damaged buildings for each San Francisco's neighborhood for M7.2 scenario.**

In order to estimate injured distribution by neighborhoods, in this study the total number of injured has been distributed based on the percentage of damaged buildings for each neighborhood. According to Tobin and Samant [2009], each neighborhood's share of the total building damage in the city is shown in Figure 4.3. This figure illustrates that the level of projected damage in Mission, Sunset, and Western Addition neighborhoods represents the greatest share of the city's building damage. The Bayview, Merced and Mission Bay neighborhoods will suffer the lowest level of damage. According to the distribution of damage to buildings, the number of injuries for each neighborhood is shown in Figure 4.4. The estimated injured in Figure 4.4 include only those injuries sustained in privately-owned buildings.

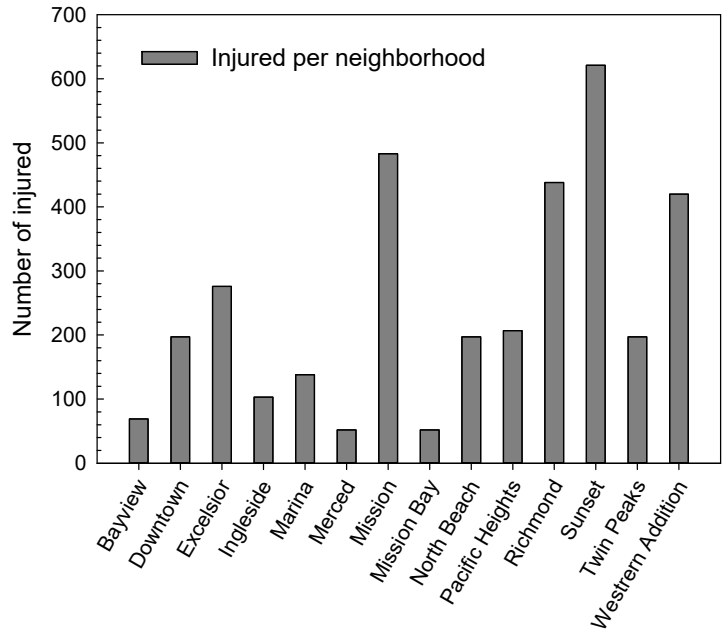

**Figure 4.4**     Number of injured per neighborhood for M7.2 scenario.

After catastrophic events, such as an earthquake, all injured must be transported to one of city's EDs. This travel has to be as short as possible because long average drive times to a trauma center could worsen a patient's conditions worse. For this reason, all the estimated injured for the magnitude 7.2 scenario have been distributed considering their proximity to the nearest hospital. Within a single neighborhood, a homogeneous distribution of the patients has been considered in order to send each of them to the nearest ED. The percentage of patients for each hospital is shown in Table 4.2.

Patient wait times has been chosen as the main response parameter used to examine hospital disaster resilience during earthquakes: a wait time of 180 minutes (3 hours) has been selected as a limit value above which the hospital is considered not resilient. Interviews with medical staff of several hospitals have shown that patients' conditions may worsen irreversibly after waiting more than 3 hours. Equation (3.23) has been used to determine PWTs. It describes the time that patients have to wait before receiving care at a given instant in time considering simultaneously two different variables: the total number of treatment rooms per color area ($m$) and the seismic input ($\alpha$).

A survey for each hospital included in this study has been conducted a specially designed questionnaire and relevant data including the number of emergency rooms ($m$) have been collected for each hospital. By knowing the number of patients arriving at each ED, the seismic input ($\alpha$) can be obtained using Equation (4.1):

$$\alpha \propto NP \tag{4.1}$$

where $NP$ is the number of patients for the considered $\alpha$ value. It has been assumed that 559 is the number of patients arriving at the ED corresponding to $\alpha = 1$. The $\alpha$ values obtained for the six analyzed hospitals are listed in Table 4.3.

By knowing the *m* and *α* values, the trend of PWTs for each considered hospital has been obtained. Then, a comparison between the estimated waiting times and the maximum acceptable PWTs value (3 hours) has been done for each hospital; see Figure 4.5. As shown in the figure, hospital 1, hospital 3, and hospital 4 are unable to provide care to all the patients arriving at the ED. In particular, the average PWTs for hospital 1 reaches a peak value of about 280 min while hospital 3 and hospital 4 of about 4700 min. This because the hospitals' capability is determined in terms of their size measured by the number of treatment rooms available.

Two approaches have been considered in order to ensure that all patients receive emergency care within the maximum acceptable PWT. One assumes a functionalist perspective in which the capacity of one of the others healthcare facilities (hospitals 2, 5, 6) is used to guarantee emergency care to all the patients that cannot be treated in the hospital closest to them. This implies the presence of a Control Center to redistribute the flow to other hospitals. The other considers the possibility of increasing the number of hospitals by using other healthcare facilities already existing in the San Francisco's area that are not equipped with an efficient ED. This means that some aspects of basic emergency planning need to be implemented in these structures.

Table 4.2    Estimated percentage of injured for each analyzed hospital.

| Hospital | % of Injured |
|---|---|
| H1 | 15.45% |
| H2 | 15.95% |
| H3 | 15.55% |
| H4 | 15.35% |
| H5 | 16.35% |
| H6 | 21.35% |

Table 4.3    *α* values for each considered hospital.

| Hospital | N. of Patients | *α* |
|---|---|---|
| Hospital 1 | 564 | 1.01 |
| Hospital 2 | 582 | 1.04 |
| Hospital 3 | 568 | 1.02 |
| Hospital 4 | 560 | 1 |
| Hospital 5 | 597 | 1.08 |
| Hospital 6 | 779 | 1.39 |

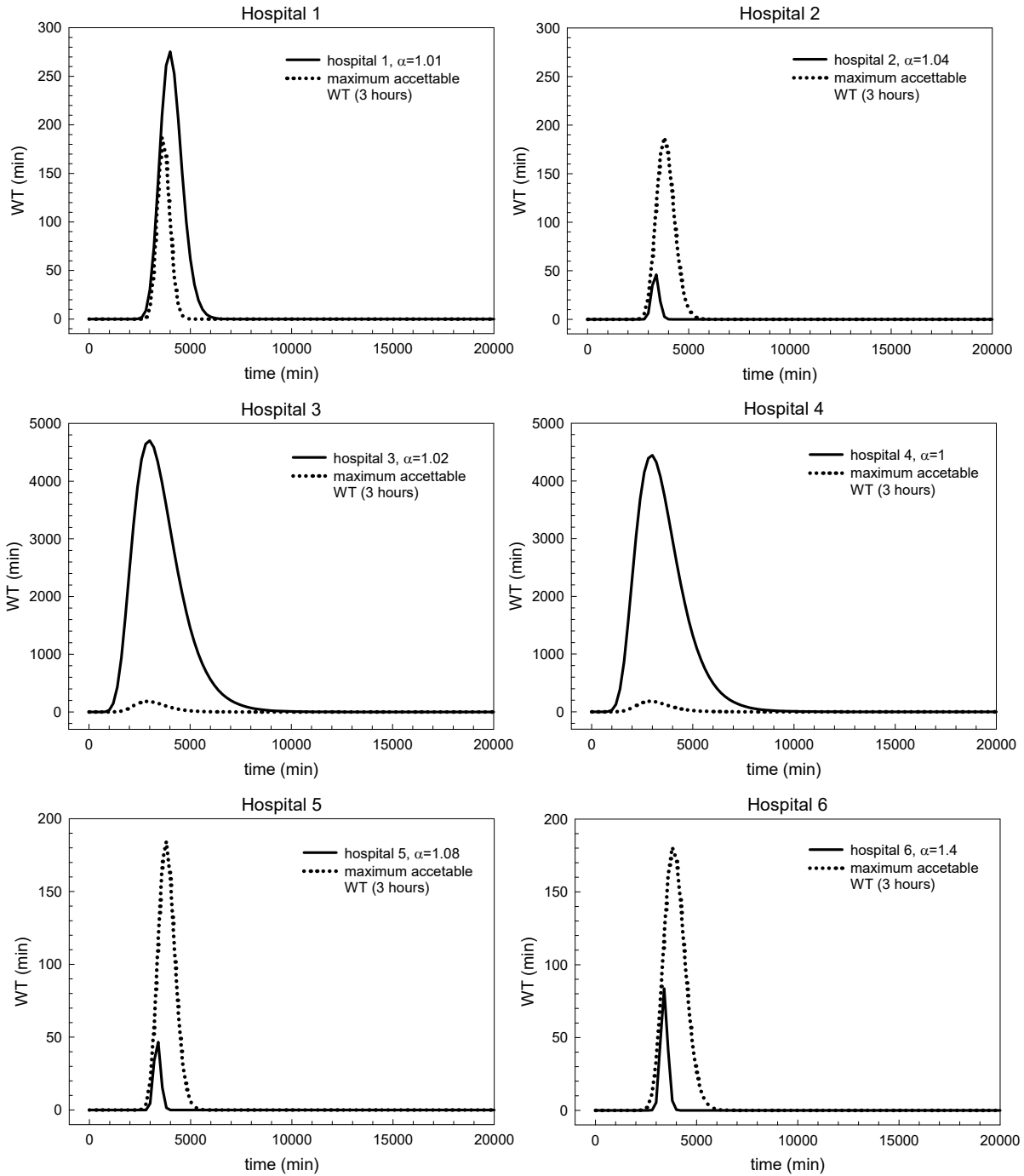

Figure 4.5  Patient's estimated wait times vs maximum acceptable wait times (3 hours).

## 4.5 APPROACH 1: REDISTRIBUTION OF PROPOSED PATIENTS

The patient redistribution approach has been developed to guarantee emergency care to all the injured in the San Francisco's area assuming a 7.2 magnitude earthquake on the San Andreas fault. Patients have been distributed in the six San Francisco's hospitals on the basis of two criteria, including each hospital's capacity and the distance from the place in which the patient is located to the ED. A hospital's capacity represents the number of patients that the hospital is able to treat so that PWTs do not exceed the maximum acceptable.

First, patients were distributed considering only the distance criterion. Then, the total numbers of patients arriving at the ED and the hospital's capacity were compared. As aforementioned, three hospitals had received too many patients and had reached capacity. Thus, a Control Center has been assumed to manage the flow of patients in the hospitals. The maximum acceptable WT has been used as the primary source for measurement of emergency care hospital capacity. Therefore, the maximum number of patients that can be treated in each of the three hospitals has been obtained, and the remainder of the patients has been distributed in the hospitals with higher capacity. New $\alpha$ values were calculated. Considering the new $\alpha$ values, the trend of PWTs for each considered hospital after the redistribution was calculated, and a comparison between the estimated PWTs and the maximum acceptable PWT value (3 hours) has been done for each hospital; see Figure 4.6. The redistribution of patients has a marked influence on the city's capacity to provide emergency care to all earthquake injuries. As shown in Figure 4.6, PWTs never exceed the maximum acceptable PWT of 3 hours.

Although redistributing the patients according to the capacity of each hospital addressed the issue of the maximum acceptable PWTs (3 hours), the travel time to reach hospitals increased. For this reason, the rate of increase in patients' travel time has been considered in order to evaluate the accuracy of the proposed solution. The maximum travel time between hospitals and their areas of expertise has been calculated considering normal San Francisco's traffic conditions in rush hour. These travel times have been amplified by 50% in order to take into account the traffic congestion caused by the emergency. Results are listed in Table 4.5, which gives the rates of increase in patient travel time for each considered hospital. The highest increase rate is about 9 min, validating the initial analysis that the proposed solution is a viable option.

That said, this approach presents some limitations. The presence of a Control Center is mandatory in order to manage patient flow to the hospitals. Installing a Control Center might be prohibitively expensive and given the issue of communication breakdowns in the event of an earthquake, directing people to the right hospital is a stumbling block. For this reason, another approach has been developed, as discussed below

Table 4.4  α values after redistribution.

| Hospital | No. of Patients | α |
|---|---|---|
| Hospital 1 | 492 | 0.88 |
| Hospital 2 | 654 | 1.17 |
| Hospital 3 | 341 | 0.611 |
| Hospital 4 | 341 | 0.611 |
| Hospital 5 | 816 | 1.46 |
| Hospital 6 | 1006 | 1.8 |

Table 4.5  Maximum travel time between hospitals and their areas of expertise calculated considering normal San Francisco traffic conditions in rush hour.

| Hospital | Travel time before redistribution | Travel time after redistribution | Increase in patients travel time |
|---|---|---|---|
| Hospital 1 | 21-29 min | 21-29 min | 0 min |
| Hospital 2 | 10-17 min | 19-26 min | 9 min |
| Hospital 3 | 20-28 min | 20-28 min | 0 min |
| Hospital 4 | 17-26 min | 17-26 min | 0 min |
| Hospital 5 | 19-25 min | 24-30 min | 5 min |
| Hospital 6 | 24-30 min | 29-38 min | 8 min |

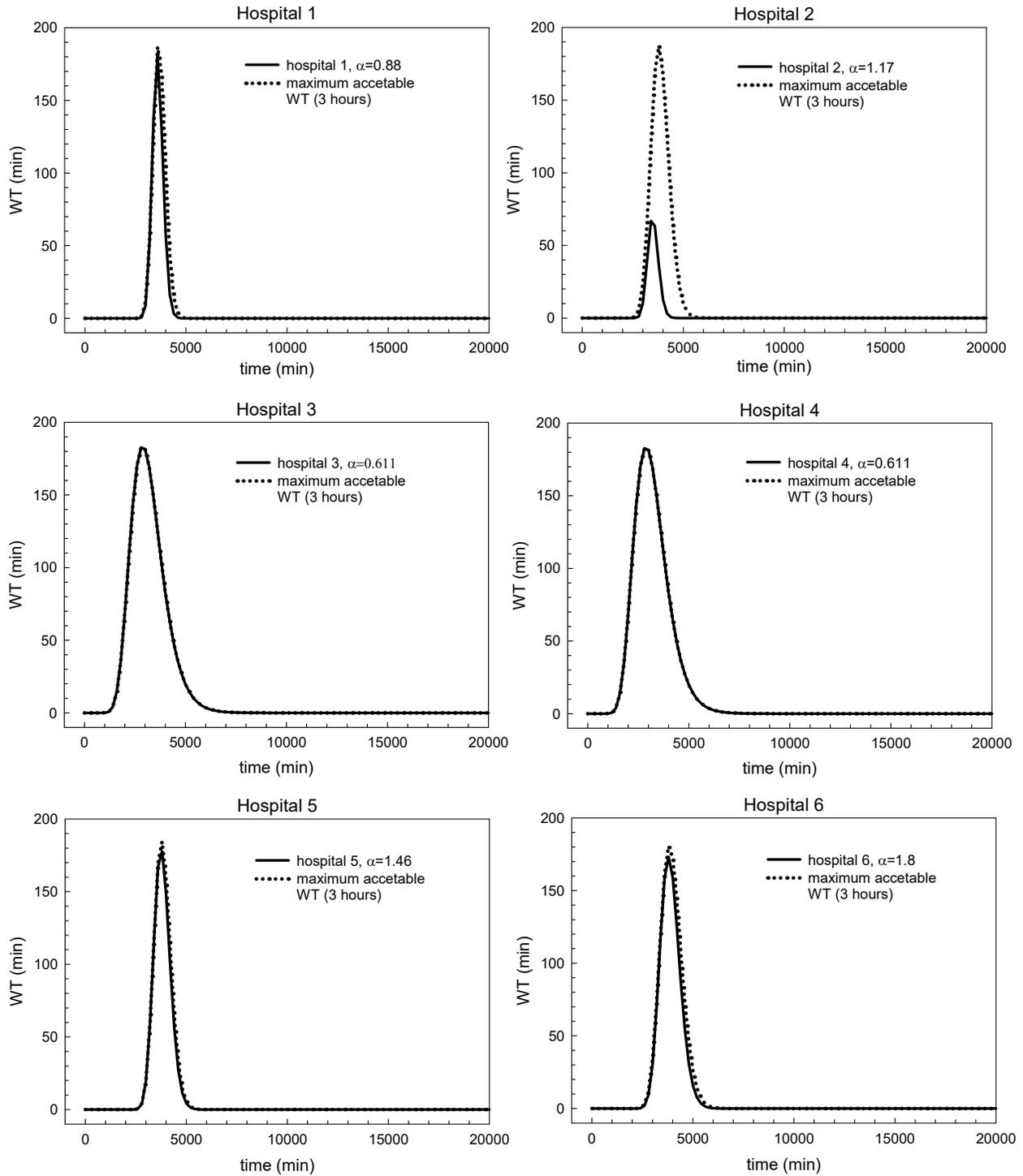

**Figure 4.6** Patient's estimated wait times after redistribution vs maximum acceptable wait times (3 hours).

## 4.6 APPROACH 2: INCREASE IN THE NUMBER OF HEALTHCARE FACILITIES

An increase in the number of healthcare facilities has been considered as a second solution to guarantee emergency care to all the injured in the San Francisco area within a defined PWT. In particular, two possibilities were analyzed. The first one includes the use of one or more already existing structures in San Francisco and adapting them in order to accommodate a flow of patients caused by an emergency such as the hypothesized earthquake. This aspect has been not examined in detail in this study. The second possibility proposes the creation of a new hospital to provide safe and timely care to all the injured. This solution will be explored below.

### 4.6.1 Area Identification

The identification of the area in which the new hospital should be located is the first important step. As illustrated in Section 4.4, three San Francisco hospitals (hospitals 1, 3, 4) are unable to provide care to all the patients arriving at their EDs after a large earthquake. This means that the new healthcare facility has to be placed in the area of expertise of those three hospitals, including Richmond district, Western Addition, Pacific Heights, Marina, North Beach, and Downtown neighborhoods. In this way, all the injured that cannot be treated in these three hospitals could be redirected to the new healthcare facility.

In order to find the most appropriate location, the density of injured in each district has been considered. In particular, the technique used in making location decisions is the "center of gravity" method. With this method, the coordinates for the optimal location are chosen as an average of the coordinates of the various neighborhoods, which are weighted according to the number of injured expected from each neighborhood. Two steps have been followed in the area identification process. First, the neighborhoods' center-of-gravity (G) has been calculated using a Cartesian reference system. The center-of-gravity represents the geometric center of the neighborhood considering a uniform density of patients in each area. In a subsequent step, the total number of patients per neighborhood has been used as weight. The center-of-mass has been calculated, and the coordinates of the location of the new hospital obtained. As shown in Figure 4.7, the structure should be placed in the northern part of the Richmond neighborhood. Figure 4.7 illustrates the position of the centers-of-gravity as well as the number of patients per neighborhood.

Locating the new hospital at exactly (24.14, 8.38) on the coordinate system may not be possible, since there may not be any available land at that location. Decision makers may choose a location that is feasible and near that location. This calculation is simply based on estimates on number of injured, which may vary in accuracy. In addition, the distances to the neighborhoods are more than a function of map coordinates since roads may be direct or may wind around; however, this method provides useful information that can help locate an appropriate site.

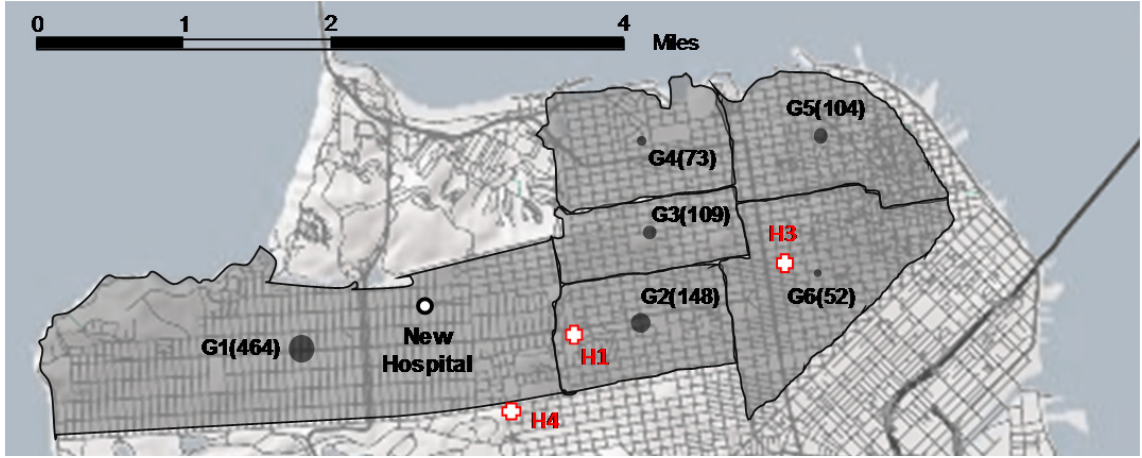

Figure 4.7    New hospital location.

### 4.6.2 Hospital Size

Once the most suitable position for the new facility has been identified, the number of patients arriving at each ED has been calculated considering the distance criterion. Results are listed in Table 4.6.

In this section, the number of patients arriving at the new healthcare facility has been used as an input parameter to define the hospital size. The capacity planning issue is primarily driven by the general meta-model described above. Using Equation (3.23), we obtain the minimum number of treatment rooms that should be allocated to the new ED to meet efficiency targets. Reversing the aforementioned equation, the minimum number of treatment rooms ($m$) can be expressed in function of patient waiting time ($WT$) and the seismic input ($\alpha$):

$$m = f(WT, \alpha) \tag{4.2}$$

A waiting time of 180 minutes (3 hours) has been considered as an input value since it represents the maximum time patients can wait before make their conditions irremediably worsen. Knowing the number of patients arriving at the EDs allows us to determine the seismic input ($\alpha$). Results are listed in Table 4.7.

By knowing the $\alpha$ and $WT$ values, the minimum number of treatment rooms has been obtained using Equation (4.1). Results indicate that a minimum three EDs is required to accommodate all the patients arriving at the new ED within a maximum waiting time of 3 hours. By knowing the number of treatment rooms required for the new facility, the peak value of PWTs for each considered hospital is shown in Table 4.8. This table compares the estimated PWTs and the maximum acceptable PWT value (3 hours) hospitals could provide timely and efficient care to all the patients arriving at their EDs, demonstrating that this is a viable option for addressing the overcrowding of current ED capacity in the event of a major earthquake.

**Table 4.6** Estimated number of injured for each analyzed hospital.

| Hospital | No. of patients |
|---|---|
| Hospital 1 | 432 |
| Hospital 2 | 582 |
| Hospital 3 | 336 |
| Hospital 4 | 319 |
| Hospital 5 | 597 |
| Hospital 6 | 779 |
| New Hospital | 605 |

**Table 4.7** $\alpha$ values for each hospital.

| Hospital | No. of patients | $\alpha$ |
|---|---|---|
| Hospital 1 | 432 | 0.77 |
| Hospital 2 | 582 | 1.04 |
| Hospital 3 | 336 | 0.61 |
| Hospital 4 | 319 | 0.57 |
| Hospital 5 | 597 | 1.07 |
| Hospital 6 | 779 | 1.39 |
| New Hospital | 605 | 1.08 |

**Table 4.8** Peak value of patient wait time for the considered hospitals.

| Hospital | Average WT peak |
|---|---|
| Hospital 1 | 127 min |
| Hospital 2 | 66 min |
| Hospital 3 | 98 min |
| Hospital 4 | 81 min |
| Hospital 5 | 167 min |
| Hospital 6 | 173 min |
| New Hospital | 53 min |

## 4.7 REMARKS AND CONCLUSIONS

As recent events have shown, even moderate damage from earthquakes can become a catastrophe. This is why the concept of resilience—*the ability of a system to mitigate and contain the effects of disasters while plan effective strategies to recover fast minimizing social disruption*—has gained more and more attention. In this work, a methodology to evaluate San Francisco's preparedness for the next big earthquake in the Bay Area has been proposed. Considering a magnitude 7.2 earthquake scenario, a hospital's ability to provide care to all the injured arriving at the EDs has been evaluated. Results show that three of the six San Francisco's hospitals considered are not able to meet the increased demand for medical care. For this reason, two solutions have been proposed to better address a post-earthquake scenario. The approaches evaluate the optimal recovery plan that maximizes the resilience index of the San Francisco's healthcare facilities so that all the EDs can guarantee timely and efficient care to all the injured from the earthquake.

Approach 1 includes a redistribution strategy whereby patients are rerouted to hospitals with higher capacity. The benefits gleaned in preserving and using already existing healthcare facilities avoid the extensive planning phase for a new healthcare facility, the time necessary for the construction, and the extremely high costs. The installation of a Control Center could reroute the injured in real time to hospitals with greater capacity, thus ensuring that PWTs never exceed the maximum value of three hours. Coordinating such a relocation might prove difficult, assuming there might be some breakdown in the communication system as cell phone users overload the system. In addition, damage to local streets from broken sewers line, debris from collapsed building fronts, etc., might make local infrastructure difficult to navigate, increasing travel times to unacceptable levels.

Approach 2 identifies behaviors and resources that contribute to a system's ability to respond to the unexpected. By building a new hospital, the resources that are needed for resilient adaptation have been qualified, and patients have been sent to the nearest facility to them. This approach negates the need for a Control Center. However, building a new healthcare facility requires a major financial investment. Although this addresses the issue of ED capacity during a major earthquake, it might be too costly to implement.

Future research plans to extend this methodology to include the insertion of components of the physical infrastructure dimension (e.g., transport network, hospitals structures, etc.) as well as the financial aspects.

# 5 The Application of Factor Analysis to Evaluate Hospital Resilience

## 5.1 INTRODUCTION

Natural and manmade disasters worldwide have steadily increased, becoming more frequent and more intense due to increased rates of urbanization, environmental degradation, and to changes in climate variables (such as higher temperatures or extreme precipitation). Healthcare facilities and emergency services are at the front lines of this increased demand. Hospitals differentiate themselves from other healthcare organizations as they play an important role during the immediate aftermath of an emergency by providing continued access to care, serving as a safety net, and preparing and responding to disasters. In order to ensure that hospitals respond appropriately, it is necessary to define and analyze complex scenarios so that the system does not collapse in the event of a dire emergency.

Hospital disaster resilience focuses on a system's overall ability to prepare and plan for, and recover, from catastrophic events as well as sustain required operations under both expected and unexpected conditions. However, a hospital's adaptive behaviors depend on several variables related to the complexity of the system. Hospital disaster resilience must be measured separately, using multiple concepts such as hospital safety, cooperation, recovery, emergency plans, business continuity, critical care capacity, and other important dynamics. Thus, the overall resilience level can be obtained by combining the resilience of each individual variable in order to take into account hospital's response ability at all levels of the system [Zhong et al. 2014]. Several methods have been proposed that measure a hospital's ability to provide emergency care to all the injured in an extreme situation. In this study, a framework for hospital resilience has been developed, using empirical data from hospitals in the San Francisco Bay Area (California). To achieve this goal, data from a survey questionnaire were analyzed in order to extract key factors for hospital resilience measurement.

Factor analysis is a means of describing a characteristic that is not directly observable based on a set of observable variables. It has been used largely to analyze and measure latent factor in several different fields [Li et al. 2013]. Herein, it has been conducted using principal components analysis. In this specific case, eight variables were considered as those most representative of a hospital's performance during an emergency. Three factors explaining over 80% of variance were found, including cooperation and training management, resources and equipment capability, and structural and organizational operating procedures. Each of these factors can be analyzed separately as a means of determining which part of the hospital's internal system needs modification. Next, a score model was established to measure the level of hospital disaster resilience. The model provides an analytical expression for hospital resilience ($R$), combining linearly the three extracted factors. The weight for each factor was obtained, and the overall resilience for the considered hospitals was calculated.

## 5.2 METHODOLOGY

A study was conducted on Tertiary Hospitals in the San Francisco Bay Area, California. A questionnaire was developed to gather relevant data for the hospitals' resilience analysis. The survey was conducted between April 2014 and July 2014. A total of 16 hospitals completed the questionnaire, representing a 69% response rate. The 16 hospitals are shown in Figure 5.1. The survey was conducted by in-person interviews of the hospital's emergency staff or by sending the questionnaire by e-mail. Each hospital selected a responder who was familiar with emergency planning to respond to the questionnaire. All of the hospitals were informed of the research objectives. Before starting the factor analysis, the collected questionnaires were reviewed to check their completeness and consistency. The questionnaire consists of 33 questions and 8 sections in total. All the questions were multiple choices, where the only two possible answers were "yes" or "no". Option "yes" was assigned a score of "1"; option "no" received a score of "0". "Yes" answers represents the hospital's ability to resist and absorb the shock of disasters; "no" answers deemed the hospital's behavior as "not resilient." The total score of each section was obtained by summing the score of each question. The higher the total score, the better the hospital's disaster resilience.

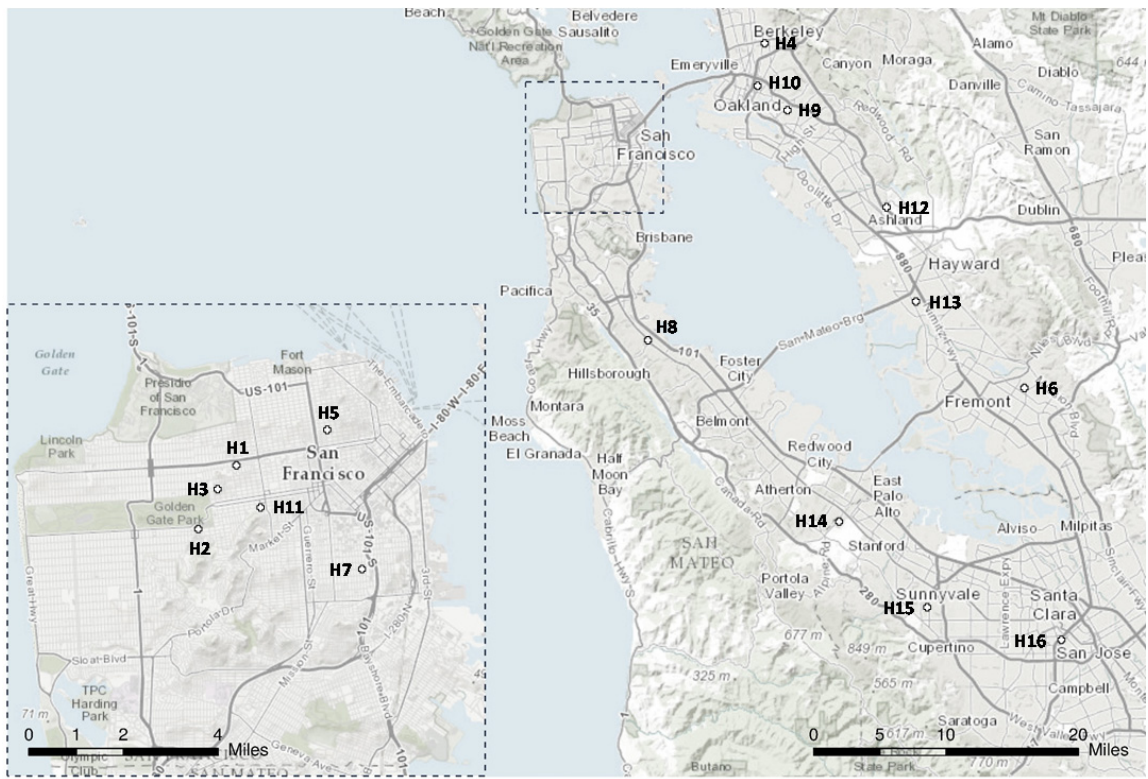

Figure 5.1    Tertiary hospitals in the San Francisco's Bay Area.

Eight items were selected to reflect hospital's behavior during an emergency. In order to simplify the analysis, each item has been replaced with a letter, as following:

(a) Hospital safety

(b) Hospital disaster leadership and cooperation

(c) Hospital disaster plan

(d) Emergency stockpiles and logistics management

(e) Emergency staff

(f) Emergency critical care capability

(g) Emergency training and drills

(h) Recovery and reconstruction

All the collected data from the survey were analyzed to determine a lower number of unobserved variables that reflected a hospital's disaster resilience. A Microsoft Excel 2013 database was created, and the answers were transferred to a spreadsheet to be checked and analyzed. Factor analysis was performed using IBM SPSS Statistic version 21, downloaded on May 15, 2014.

Several steps were considered in building a valid framework to measure a hospital's disaster resilience. The basic idea was to represent all the variables included in the hospital's resilience analysis with a smaller number of variables because some of them could be linearly related to each other. Thus, first the presence of significant correlations between the items was checked. Second, initial factor loadings were calculated using the principal component method. Once the initial factor loadings were calculated, the factors were rotated in order to find those factors that were easy to interpret. The goal of rotating factors was to ensure that all variables had high loadings only on one factor. Varimax rotation was used to rotate the extracted principal components. Then, factors scores were obtained, and the numbers of factors were chosen, selecting those eigenvalues that were greater than 1. Finally, a framework for hospital disaster resilience was obtained as linear combination of the extracted factors taking into account the calculated weights.

## 5.3   FACTOR ANALYSIS

Factor analysis is a statistical method used to investigate whether a number of variables of interest, say, $Y_1$, $Y_2$, ... , and $Y_n$, are related to a smaller number of unobserved variables $F_1$, $F_2$... $F_m$, which are called factors. A factor is a hypothetical variable that influences scores on one or more observed variables. The factor analysis's first goal is to determine how many factors are sufficient to include all the information contained in the original set of statements.

Different methods exist for estimating the parameters of a factor model. The research reported herein used the Principal Component method. It consists in an orthogonal transformation that converts a number of possibly correlated variables into a set of factors that are linearly uncorrelated and of high variance. These factors are called principal components. Therefore, each of the considered items can be expressed as a linear combination of a number of common factors:

$$z_j = k_{j1}F_1 + k_{j2}F_2 + \ldots + k_{jh}F_m \tag{5.1}$$

where $z_j$ is the *j*th standardized variable, $F_1$, $F_2$, ..., $F_m$ are common factors independent and orthogonal to each other (with $m < n$), and $k_{jh}$ are the calculated coefficients.

By applying the inverse factor model, it is possible to obtain factor equations as a linear combination of the original variables:

$$\begin{aligned} F_1 &= c_{11}z_1 + c_{12}z_2 + \ldots + c_{1n}z_n \\ F_2 &= c_{21}z_1 + c_{22}z_2 + \ldots + c_{2n}z_n \\ &\ldots \\ F_m &= c_{m1}z_1 + c_{m2}z_2 + \ldots + c_{mn}z_n \end{aligned} \tag{5.2}$$

In order to extract key component factors, three steps were considered. First, the relationships among variables were studied; second, the factors were extracted; and third, an analytical solution for hospital resilience was proposed.

### 5.3.1 Correlation Analysis

As discussed above, the goal of a factor analysis is to obtain factors that represent the correlations between variables; therefore, these variables have to be connected to each other somehow. If the relationships between variables are weak, it is unlikely that common factors exist. Two tests have been used in order to verify the presence of significant correlations between the items. The Kaiser-Meyer-Olkin test (KMO) is used to check whether the sample size is big enough. The sample is adequate when the KMO value is greater than 0.5. Bartlett's test of sphericity compares the observed correlation matrix to the identity matrix. In particular, it checks if the correlation matrix is an identify matrix, thus implying that all of the variables are uncorrelated. In this study, the KMO value is greater than 0.5, and the Bartlett's test indicates that some variables are not independent, suggesting that the data are suitable for a factor analysis as shown in the correlation matrix in Table 5.1.

By visualizing the correlation matrix, it has been demonstrated that some variables are correlated. In fact, the absolute values outside the main diagonal are often high (e.g., b and d: 0.813; g and b: 0.764). This means that these variables are useful for a factor analysis. In addition, the table of communalities has been examined to test the goodness of fit. Indeed, this table shows how much of the variance in each of the original variables is explained by the extracted factors. If the communality for a variable is less than 50%, it is a candidate for exclusion from the analysis because the factor solution contains less than half of the variance in the original variable. For this reason, higher communalities are desirable. In this case, the extracted communalities for all the testing variables are greater than 70%, which indicates that the extracted components represent the variables well.

Table 5.1    Correlation matrix.

|  |  | a | b | c | d | e | f | g | h |
|---|---|---|---|---|---|---|---|---|---|
| Correlation | a | 1.000 | 0.494 | -0.374 | 0.321 | 0.550 | -0.040 | 0.377 | 0.549 |
|  | b | 0.494 | 1.000 | 0.356 | 0.813 | 0.471 | 0.000 | 0.764 | 0.301 |
|  | c | -0.374 | 0.356 | 1,000 | 0.304 | 0.120 | 0.000 | 0.102 | -0.346 |
|  | d | 0.321 | 0.813 | 0.304 | 1,000 | 0.000 | -0.111 | 0.745 | -0.079 |
|  | e | 0.550 | 0.471 | 0.120 | 0.000 | 1,000 | 0.292 | 0.441 | 0.686 |
|  | f | -0.040 | 0.000 | 0.000 | -0.111 | 0.292 | 1,000 | 0.186 | 0.553 |
|  | g | 0.377 | 0.764 | 0.102 | 0.745 | 0.441 | 0.186 | 1,000 | 0.318 |
|  | h | 0.549 | 0.301 | -0.346 | -0.079 | 0.686 | 0.553 | 0.418 | 1,000 |

Table 5.2    Table of communalities.

|  | Initial | Extraction |
|---|---|---|
| a | 1.000 | 0.869 |
| b | 1.000 | 0.929 |
| c | 1.000 | 0.834 |
| d | 1.000 | 0.891 |
| e | 1.000 | 0.711 |
| f | 1.000 | 0.752 |
| g | 1.000 | 0.775 |
| h | 1.000 | 0.910 |

### 5.3.2 Factor Extraction

As mentioned above, the principal component method was used to extract the independent factors and eigenvalues were obtained. They indicate the amount of variance explained by each principal component or each factor so that the sum of the eigenvalues is equal to the number of variables. The number of factors determined considering the number of eigenvalues exceeds 1.0, according to the method proposed by Kaiser [1960]. Lower values account for less variability than does a single variable. In this study, three factors have an eigenvalue greater than 1 as shown in Table 5.3. By analyzing Table 5.3, the three extracted factors appear to be representative of all domains. The cumulative variance of the three-factor solution is up to 83%; therefore, the factors are adequate to describe the hospital's performance.

Table 5.3    Total variance explained.

| Component | Initial eigenvalues | | | Un-rotated factors | | | Rotated factors | | |
|---|---|---|---|---|---|---|---|---|---|
| | Total | % of variance | Cumulative % | Total | % of variance | Cumulative % | Total | % of variance | Cumulative % |
| 1 | 3.356 | 41.950 | 41.950 | 3.356 | 41.950 | 41.950 | 2.951 | 36.891 | 36.891 |
| 2 | 2.075 | 25.932 | 67.882 | 2.075 | 25.932 | 67.882 | 2.101 | 26.257 | 63.148 |
| 3 | 1.241 | 15.507 | 83.388 | 1.241 | 15.507 | 83.388 | 1.619 | 20.241 | 83.388 |
| 4 | 0.779 | 9.735 | 93.124 | | | | | | |
| 5 | 0.308 | 3.851 | 96.975 | | | | | | |
| 6 | 0.191 | 2.382 | 99.357 | | | | | | |
| 7 | 0.046 | 0.569 | 99.926 | | | | | | |
| 8 | 0.006 | 0.074 | 100.000 | | | | | | |

### *5.3.2.1 Rotation*

The rotation phase of factor analysis attempts to transform the initial matrix into one that is easier to interpret. Varimax rotation has been used in order to improve the interpretability and utility of the results. Indeed, the relationship between the initial items and the extracted factors is not clear after the extraction of factors. For this reason, rotation has been used in an effort to find another set of loadings that fit the observations equally well but can be more easily interpreted. After a Varimax rotation, each original domain tends to be associated with one of the three extracted factors, and each factor represents only a small number of items. Orthogonal rotations keep factors uncorrelated while increasing the meaning of the factors. A rotated component matrix was obtained that helped to determine what each factor represents. The total amount of variation explained by the three factors remains the same, and the total amount of the variation explained by both models is identical.

Table 5.4    Rotated component matrix.

| | Component | | |
|---|---|---|---|
| | 1 | 2 | 3 |
| a | 0.479 | 0.216 | **0.770** |
| b | **0.947** | 0.178 | 0.006 |
| c | 0.353 | 0.014 | **-0.842** |
| d | **0.919** | -0.202 | -0.083 |
| e | 0.359 | **0.733** | 0.211 |
| f | -0.120 | **0.834** | -0.208 |
| g | **0.836** | 0.270 | 0.053 |
| h | 0.118 | **0.822** | 0.469 |

Table 5.4 shows a new set of values for each of the three extracted factors. The bolded values represent the percentage of the extracted information from each item. The first factor is strictly connected to three items, including hospital disaster leadership and cooperation (0.947), emergency stockpiles and logistics management (0.919), and emergency training and drills (0.836). As each of these variables increase, so do the other three. Three variables are also included in the second factor that is primarily a measure of emergency staff (0.733), emergency critical care capability (0.834), and recovery and reconstruction (0.822). The third factor contains information mainly from two items: hospital safety (0.770) and hospital disaster plan (0.842). Note that all items have high-factor loadings on only one factor. The first factor includes all the items related with the hospital management mechanisms during emergencies. The second factor is representative of the emergency department's (EDs) capability in terms of human and financial resources as well as hospital's facilities (number of beds, rooms, etc.). The third factor focuses on all the hospital's prevention strategies (structural and organizational).

In this way, the three extracted factors have been identified and named:

- ($F_1$) Cooperation and Training Management
- ($F_2$) Resources and Equipment Capability
- ($F_3$) Structural and Organizational Operating Procedures

A linear combination of these three factors represents the hospital's resilience.

### 5.3.3 Factor Solutions

After interpreting the factors, an analytical expression for all the items has been obtained. Each of the variables has been expressed as a linear combination of the extracted factors taking into account the weight of the factors obtained by use in the component matrix; see Table 5.5:. Using the coefficients provided by the component matrix, the eight initial variables could be determined as:

$$a = 0.7F_1 - 0.292F_2 - 0,541F_3 \tag{5.3}$$

$$b = 0.877F_1 + 0.4F_2 + 0,002F_3 \tag{5.4}$$

$$c = 0.086F_1 + 0.637F_2 + 0,648F_3 \tag{5.5}$$

$$d = 0.643F_1 + 0.676F_2 - 0,146F_3 \tag{5.6}$$

$$e = 0.715F_1 - 0.388F_2 + 0,223F_3 \tag{5.7}$$

$$f = 0.259F_1 - 0.489F_2 + 0,668F_3 \tag{5.8}$$

$$g = 0.842F_1 + 0.255F_2 + 0,029F_3 \tag{5.9}$$

$$h = 0.624F_1 - 0.716F_2 + 0,094F_3 \tag{5.10}$$

In order to build an evaluation framework to determine the resilience of a hospitals, an analytical expression of the factors obtained considering that different factors have different contributions to overall resilience. The score for each factor was obtained using a regression analysis based on factor score coefficient matrix shown in Table 5.6.

Table 5.5   Component matrix.

|   | Component | | |
|---|---|---|---|
|   | 1 | 2 | 3 |
| a | 0.700 | -0.292 | -0.541 |
| b | 0.877 | 0.400 | 0.002 |
| c | 0.086 | 0.637 | 0.648 |
| d | 0.643 | 0.676 | -0.146 |
| e | 0.715 | -0.388 | 0.223 |
| f | 0.259 | -0.489 | 0.668 |
| g | 0.842 | 0.255 | 0.029 |
| h | 0.624 | -0.716 | 0.094 |

Table 5.6   Factor score coefficient matrix.

|   | Component | | |
|---|---|---|---|
|   | 1 | 2 | 3 |
| a | 0.142 | -0.063 | 0.479 |
| b | 0.322 | 0.006 | -0.038 |
| c | 0.135 | 0.122 | -0.579 |
| d | 0.349 | -0.184 | -0.032 |
| e | 0.056 | 0.331 | 0.010 |
| f | -0.120 | 0.505 | -0.286 |
| g | 0.273 | 0.059 | -0.021 |
| h | -0.042 | 0.358 | 0.172 |

Thus, expressions of factors can be determined using the coefficients provided by the matrix above:

$$F_1 = 0.142a + 0.322b + 0.135c + 0.349d + 0.056e - 0.120f + 0.273g - 0.042h \quad (5.11)$$

$$F_2 = -0.063a + 0.006b + 0.122c - 0.184d + 0.331e + 0.505f + 0.059g + 0.358h \quad (5.12)$$

$$F_3 = 0.479a - 0.038b - 0.579c - 0.032d + 0.010e - 0.286f - 0.021g + 0.12h \quad (5.13)$$

These factors represent the basic structure for a hospital's resilience. The overall level of hospital disaster resilience ($R$) can be calculated using the three extracted factors combined linearly.

$$R = \alpha F_1 + \beta F_2 + \chi F_3 \tag{5.14}$$

where $F_1$, $F_2$, and $F_3$ are the extracted factors calculated using Equations (5.11), (5.12), and (5.13), respectively; $\alpha$, $\beta$, and $\chi$ are the weights of each factor.

The coefficients $\alpha$, $\beta$, and $\chi$ were calculated by the proportion between the percentage of variance explained by each factor and the cumulative variance contribution of the three primary factors. Then, the expression for hospital disaster resilience ($R$) has been obtained:

$$R = \frac{(0.503 F_1 + 0.311 F_2 + 0.186 F_3)}{4} \tag{5.15}$$

The weight for hospital cooperation and training management is 0.503, 0.311 for hospital resource and equipment, and 0.186 for hospital structural and organizational operating procedures. This means that the first factor is more relevant at the time to assess the resilience of healthcare facilities, representing about half of the hospital emergency preparedness and response.

Three levels for hospital disaster resilience have been identified. Upon completion of the questionnaire, each of the eight items will have an overall score obtained by summing the scores of each question (with the score "0" or "1"). Using these scores, the three extracted factors can be calculated. Knowing the value of the factors, hospital disaster resilience ($R$) can be obtained using Equation (5.15). The $R$ value is in the range:

$$0 < R < 1 \tag{5.16}$$

where "0" represents "no resilience," "and "1" means "maximum level of resilience" corresponding to the ability to absorb any damage without suffering complete failure. If the resilience value is over 0.75, the hospital is very resilient to emergencies. If the resilience value is below 0.25, the hospital is not able to absorb adequately disastrous impacts. The score of each factor provides information about a hospital's behavior. It is possible to evaluate the differences between the different factors (cooperation and training management, resources and equipment capability, and structural and organizational operating procedures) to determine which part of the hospital's system has the lowest level of resilience.

Table 5.7  Levels of hospital disaster resilience.

| Low level of resilience | Moderate level of resilience | High level of resilience |
|---|---|---|
| $R \leq 0.25$ (25%) | $0.25 < R < 0.75$ (25%-75%) | $R \geq 0.75$ (75%) |

## 5.4 GENERAL DISCUSSION

The framework developed above can be applied to each hospital in this study to determine the level of resilience in the considered geographical area. As shown in Table 5.8, the disaster resilience score has been calculated for each hospital. For confidentiality reasons, all hospitals are represented by a number. According to the table, 10 hospitals, which account for about 62.5% of the sample, have an high level of resilience ($R \geq 0.75$) while six hospitals, representing the remaining 37.5%, are in the moderate resilience zone ($0.25 < R < 0.75$). There are no hospitals whose resilience score is under 0.25, which means that there are no healthcare facilities with an insufficient level of resilience. These results indicate that the San Francisco Bay Area's hospitals have generally a high level of resilience.

Furthermore, the scores of the three extracted factors were calculated to identify areas with a lower resilience level. The results have been standardized so the scores range from 0 to 1. The values of factors values are listed in Table 5.9, which shows that the three extracted factors have a generally acceptable level of resilience. Note that among them, factor $F2$ (resources and equipment capability) has the lowest resilience level. This means that if the hospital wants to increase the overall level of resistance, the resources and equipment capability area should be studied.

Table 5.8    Disaster resilience scores for the considered hospitals.

| Hospital | R | Hospital | R |
|---|---|---|---|
| 1 | 0.836 | 9 | 0.871 |
| 2 | 0.813 | 10 | 0.681 |
| 3 | 0.771 | 11 | 0.787 |
| 4 | 0.772 | 12 | 0.607 |
| 5 | 0.391 | 13 | 0.739 |
| 6 | 0.831 | 14 | 0.663 |
| 7 | 0.904 | 15 | 0.892 |
| 8 | 0.818 | 16 | 0.581 |

Table 5.9    Extracted factors scores for the considered hospitals.

| Hospital | F1 | F2 | F3 | Hospital | F1 | F2 | F3 |
|---|---|---|---|---|---|---|---|
| 1 | 1 | 0.48 | 1 | 9 | 1 | 0.57 | 0.93 |
| 2 | 0.88 | 0.71 | 0.78 | 10 | 0.78 | 0.56 | 0.71 |
| 3 | 1 | 0.52 | 0.85 | 11 | 0.89 | 0.62 | 0.86 |
| 4 | 0.55 | 0.62 | 0.93 | 12 | 0.77 | 0.48 | 0.64 |
| 5 | 0.67 | 0.29 | 0.57 | 13 | 0.77 | 0.58 | 0.92 |
| 6 | 1 | 0.47 | 0.93 | 14 | 0.66 | 0.53 | 0.72 |
| 7 | 0.98 | 0.76 | 0.92 | 15 | 0.89 | 0.71 | 1 |
| 8 | 0.99 | 0.48 | 0.94 | 16 | 0.56 | 0.57 | 0.65 |

## 5.5 REMARKS AND CONCLUSIONS

The last few years has witnessed a number of devastating events throughout the world. Hospitals and other health facilities are vital assets for communities when disaster strikes. The overall impact of a disaster is strongly influenced by how long hospitals take to recover and provide lifesaving medical care. The recovery time is strictly related to the level of hospital resilience. Resilient facilities have the ability to govern, resist, and recover after a disaster has struck.

This chapter developed a framework for measuring hospital disaster resilience. Factor analysis was used to analyze multivariate empirical data, and a three-factor solution was obtained. Three main objectives were archived in this study.

First, hospital disaster resilience has been conceptualized as a multidimensional construct with three subordinate dimensions corresponding to the three extracted factors: cooperation and training management, resources and equipment capability, and structural and organizational operating procedures. Each of the questionnaire's items represented a particular aspect of hospital's performance during emergencies. Combining all these domains linearly gives the multidimensional measure of a hospital's overall ability to cope with disasters. Results from factor analysis provided not only a measurement of a hospital's preparedness before catastrophic events, but also shows all the aspects covered by the definition of hospital disaster resilience, including hospital disaster leadership and cooperation, emergency plans, emergency stockpiles and logistics management, emergency training and drills, and critical care capability.

Second, a framework to measure hospital disaster resilience was developed. It provides an analytical expression in which the three extracted factors are linearly combined. The weight of each factor has been assigned by the proportion of variance. The results from factor analysis demonstrate that cooperation and training management ($F1$) is the most highly weighted factor as it includes hospital coordination meeting with other EDs, emergency drugs, or emergency materials as well as emergency training programs. Three levels for hospital disaster resilience were identified, including low, moderate, and high levels of resilience

Third, the proposed framework was used to assess the level of resilience of San Francisco Bay Area hospitals. Using the considered questionnaire, relevant data was collected and analyzed, and results revealed a generally high level of resilience of hospitals in the considered geographic area. Moreover, scores on particular areas of resilience were calculated in order to identify areas for improvement. This particular study was conducted on a small sample size based on the number of hospitals in the San Francisco Bay Area. It is recommended that the sample size be increased and a similar study conducted to validate the results presented herein.

# 6 Agent-Based Models to Study the Resilience of Socio-Technical Networks during Emergencies

## 6.1 INTRODUCTION

Infrastructure has become the central nervous system of modern society, and they are organized in a complex network. Infrastructure has always had a certain degree of interdependency. In fact, once a community is subjected to a shock (earthquake, terrorism, hurricanes, floods, explosions, etc.) it is more vulnerable when this degree of interdependency is greater [Cimellaro et al. 2013]. Naturally, when contemplating the loss of life due to catastrophic events, highly accurate predictions of how such infrastructure will respond is critical. Simulation is an ideal technique as it can accommodate randomness and details required in such models, and it enables a form of experimentation not possible with the actual incidents [Shendarkar and Lee 2008].

This chapter aims at developing an agent-based model (ABM) of an infrastructure in order to study the dynamics of evacuation considering human behavior. This model guides designers and legislators on how determine if a building is safe and if the occupants will be able to evacuate in an emergency situation, or highlights the need to improve the infrastructure response during an emergency situation.

Two different ABM evacuation models of a museum and a metro/train station were developed. Each of these models is divided in two phases: the no-emergency dynamic and the evacuation process after the occurrence of a blast. Although agents are able to make informed decisions and act on them, both models contain specific irrational behaviors. The evacuation time is the main parameter of response, and it is used to evaluate the efficiency and safety of the infrastructure.

## 6.2 STATE-OF-THE-ART ON AGENT-BASED MODELS AND HUMAN BEHAVIOR

This section describes the state-of-art of human-behavior modeling, which is grouped in three parts: (i) agent-based modeling (ABM), (ii) human-behavior modeling, and (iii) belief-desire-intention paradigm (BDI).

## 6.2.1 Agent–Based Modeling

Agent-based modeling (ABM) is a key approach to model complex and heterogeneous systems. It can be applied in a vast range of fields, including biology, business, ecology, social science, technology, earth science, and network theory. To develop an ABM, many agent-based platforms can be used. Among these, the most widely and common used are Netlogo, Repast, Anylogic, Mason. Each software has different programming language and features; therefore, selection of one type of software over another is dependent more on the available computing power and time. The ABM simulates actions and interactions of autonomous individuals in an environment, who are called "agents," with a view to assessing their effects on the system as a whole. An "agent" is a discrete individual with a set of characteristics and rules that govern his/her behaviors and decision-making capability; who interacts with other agents; who is flexible, and who has the ability to learn and to adapt behaviors based on experience.

Depending on the environment, an "agent" may represent individuals, groups, companies, infrastructures, etc. Modeling agent behaviors and the reciprocal interactions are possible using rules or logical operations that can be formalized by equations. It is possible to consider individual variations in the behavioral rules ("heterogeneity") and random influences or variations ("stochasticity"); see Helbing and Balietti [2012]. Furthermore, ABM can be combined with other simulation methods used in the natural and engineering sciences, including statistical physics, biology, and cybernetics.

### 6.2.1.1 Modeling Evacuation using Agent-Based Models

Emergency evacuation is the movement of people from a potentially dangerous place to a safe refuge due to threat or occurrence of a disastrous event. The possible causes for evacuation include earthquakes, building fires, military attacks, etc. [Yuan and Tan 2007]. Evacuation models, with the aim of quantifying and modeling human movement and behavior, have been underway for at least 30 years. In the light of tightened homeland security, research on evacuation is gaining impetus and attracting more attention. It guides designers and legislators on how to determine if a building is safe and if the occupants will be able to evacuate in an emergency situation, or highlights the need to improve the infrastructure response during an emergency situation. The research into this field has progressed along two routes. The first focuses on the behavior of people under normal condition without specific threats. The latter focuses on the movement of people in critical situations and evaluation of hazardous scenarios.

Gwynne et al. [1999] state that the problem of evacuation can be analyzed in terms of: *optimization, simulation,* and *risk assessment*. The first approach is used when a large number of people are considered, especially if it is assumed that they behave in a homogenous fashion. The second approach is used when it is important to analyze decisions made and paths followed during an evacuation. The third approach identifies hazards arising from each kind of incident and evaluates variations in the results that can be associated with changes in infrastructure design.

Effective evacuation strategies require accurate prediction of the environment impact on the agents and crowd behavior. Therefore, a valid ABM incorporates human behavior into the emergency dynamic. Consequently, for cases where loss of life is a potential outcome, highly accurate predictions are mandatory. Therefore, simulation is an ideal technique as it can accommodate randomness and detail needed in such models, especially in cases where it is impossible to validate in the laboratory [Shendarkar and Lee 2008].

### *6.2.1.2 Human Behavior in ABM*

In the past, the majority of ABMs focused on disaster planning emphasized either evacuation plans or human behavior. In an evacuation-centric ABM, the agent behavior is static and follows predetermined rules. In general, human-behavior-centric ABMs simulate individual or crowd behavior.

**Behavior of an individual:** Xi et al. [2011] developed a behavioral ABM that simulates the decision-making process of in motion pedestrians based on the extended decision field theory. A multi-scale modeling framework was proposed and developed using AnyLogic® software to mimic crowd-crossing behaviors under the considered situations.

Miyoshi et al. [2011] developed a multi-agent model of evacuation behavior of passengers in the event of an aircraft evacuation. In this model, agents mimic the behavior and mental state of passengers in the cabin, including mental stress, strong fear, or anxiety, which generate time delays in the evacuation. Factors that influenced panic include the remaining time, the frequency of waiting, and the difficulty of finding an exit. The results of the run simulation collected the dependencies of these factors on the time needed to complete an evacuation.

**Behavior of a crowd:** Shendarkar and Lee [2008] conducted a simulation of crowd evacuation management under terrorist bomb attacks in public areas. The limitation of this ABM is that although they simulated the behavior accurately, the agents moved according to predetermined paths.

**Complete behavioral ABM:** Luo [2008] has developed the most complete model in this field, which considers both individual and crowd behaviors. Thanks to a three-layered framework that reflects the natural pattern of human-like decision-making process of an individual in a crowd, which generally involves a person's awareness of the situation and the consequent changes on the internal attributes than emotional and social group attributes. The case study used a public transportation system and human behaviors and modeled behaviors in both normal and emergency situations.

## 6.2.2  Classification of Human Behaviors

Because human behavior is a complex mechanism influenced by culture, attitudes, emotions, values, ages, perception, etc., it is important to classify the context and the situation in order to analyze it. Analysis of such a complex system requires dividing it into its simple parts. In an evacuation scenario, the main factors that influence the behavior of an individual are:

6. *THE STATE,* which includes the role performed in the evacuation and the age of the individuals. These characteristics involve different static and basic behavior;

7. *INDIVIDUAL BEHAVIOR,* which considers the emotional aspects of a person. This is the most variable and unpredictable aspect; and

8. *CROWD BEHAVIOR*, where it has been widely demonstrated in literature that individuals in a crowd behave in certain ways. These behaviors are mostly influenced by kinship, aggregation phenomena, or collision events.

## 6.3 AN AGENT-BASED MODEL FOR EVACUATION OF PEDESTRIANS

This ABM presented below can be used to simulate an evacuation process of pedestrians due to a catastrophic event; this specific case defines the event as an explosion. It can be applied to every geometric configuration, and it facilitates making important considerations about the evacuation process. First, the plan metric configuration of the rooms as well as the location of the emergency exits should be determined. In consideration of the geometry of the rooms and the general architecture of a building, determination of locations where an explosion would inflict maximum damage and mitigation of that damage should be considered. Identification of the optimal distribution of users and performing various simulations with different distributions of rescuers or users should be considered. Determination of these factors aids in reducing the risk in emergency conditions, thus ensuring greater safety. The core of this ABM is the consideration of human behavior, To obtain realistic assessments of the evacuation process, agents are able to make decisions and to feel emotions during the emergency evacuation. Two case studies are considered herein: the Ursino Castle, a museum in Catania, Italy (Case Study 1) and a metro station of the Paris train station, the Gare De Lyon (Case Study 2).

### 6.3.1 ABM Phases

In both ABMs, the dynamic simulation includes two phases: (i) *Normal phase,* and (ii) *Emergency phase.*

#### 6.3.1.1 Normal Phase

The *normal phase* is defined as agents performing typical actions germane to that environment. In Case Study 1, agents move in sequence from one artwork to another. In Case Study 2, agents wait for trains to arrive inside the station, the train arrives and stops, and they enter the train car. At the same time, agents who are inside the trains disembark and move towards the exits.

#### 6.3.1.2 Emergency Phase

The emergency phase starts with an explosion. The severity of the explosion determines the number of deaths and injured at different entity levels. After the explosion, the evacuation process starts. First, the evacuation dynamics include agents moving towards the emergency exits depending on their mobility and proximity to those exits. Otherwise, they run through the path made previously, following a path already known. Some simulations include agents labeled "security," whose duty is to rescue the severely injured and support their evacuation.

### 6.3.2 Human Behavior in the ABM

Challenger et al. [2009] studied crowd behavior, basing his conclusions on the evidence collected by real evacuation processes. This study showed that a frequent phenomenon during the evacuation process is the leader–follower relationship. More importantly, a controlled evacuation run by a leader could lead to a success. Given that phenomenon of altruism is widespread, it should be included in any evacuation simulation. During the evacuation process the agent's behavior is modeled using BDI, which simulates the decision-making process performed by humans in the real world. Furthermore, the behavior of an individual in a crowd in emergency conditions is generally driven by the avoidance of social norms.

#### 6.3.2.1  Static Attributes

The static attributes called "agent's state" are a series of default features that remain unchanged throughout the simulation:

- Role played in the simulation: there are agents who have a specific role, e.g., officers inside the museum or the train station whose job duties include rescuing the injured. An added factor is that they know the environment and placement of the emergency exits.
- Health conditions: After the explosion, agents can be categorized as injured or not-injured. Each of these agents, depending on the degree of severity and type of damage reported, have a different behavior.
- Leaders–Followers: In emergency situations there are two main psychological profiles. Some people tend to play the role of leaders, i.e., these agents make decisions independently and proceed to evacuate on their own. These agents are often followed by other individuals—the followers—who tend to imitate actions performed by the leaders.

#### 6.3.2.2  Individual Behavior

First of all, the behavior of agents is driven by static rules. If an emergency exit is in their visual field in the event of a blast, they will move toward the emergency exit. Otherwise, the agents will retrace their previous path rather than running an unknown path. Unpredictable behaviors, which are driven by emotions, are categorized using a BDI paradigm defined earlier that simulates the human decision-making process. This mathematical theory is able to simulate the decision making-process of leaders and followers, as well as to simulate altruism, i.e., the instinct to help others who need help.

### 6.3.3  Comparison between Agent–Based Models against an Agent-Based Model Categorizing Human Behavior

To quantify the significance of incorporating human behavior, an ABM was compared to three other models with different characteristic. Three models were developed: (i) a base model; (ii) a deterministic model; and (iii) a probabilistic model.

#### 6.3.3.1  Base Model

The basic model—a very simplified version of ABM—can be compared to a hyper-rational model because the human behavior is characterized by a set of static rules. Evacuation is governed by one simple rule: each agent evacuates through the emergency exit placed in his her radius of vision. If agents cannot see any emergency exit, they retrace their previous path to escape through the exit where they entered. Obviously, this model reflects an ideal situation that does not reflect reality.

### *6.3.3.2 Deterministic Model*

The deterministic model incorporates leader–follower dynamics; however, the rules are static. This means that if an agent sees a leader in his radius of vision, automatically he becomes a follower. Agents cannot decide how to behave, which makes it different from the probabilistic model where the leader–follower behavior is not automatic, but agents have a choice to follow or not follow a leader.

### *6.3.3.3 Probabilistic Model*

The probabilistic model is integrated with the BDI paradigm whereby leader–follower behaviors and altruism are simulated. First, psychological profiles emerge after the explosion, and a variable percentage of agents become leaders. The other agents as followers are classified, but this does not mean that they follow a leader. In fact, they can decide if they choose to follow or not follow a leader depending on rational or irrational factors. The rational factors include the health status of the agent and the location of the emergency exits. Irrational factors include emotional responses; i.e., responses that are dependent on the emotional state of the agent. Altruism is also a factor, whereby an agent sees another injured agent in his/her radius of vision and stops to help evacuate the injured agent.

### 6.3.4 Agents

The categories of agents are the same for all models (base, deterministic, and probabilistic) and the two case studies. However, a category may behave differently, depending on the complexity of the model. Agents are categorized at the beginning of the simulation. During the evacuation process, their categorization may change permanently or temporarily as a function of the human behavior. Each agent has its own intrinsic characteristics. In order to understand visually evacuation dynamics, icons for each agent were developed:

VISITORS: These agents in good health after the explosion. In this agent category, various psycological figures are defined. A percentage of "visitor" agents will temporarily change their status and become another agent category. For example, a percentage of them will assume the role of leader and a percentage the role of follower. Some may decide to become altruist and help another agent.

SECURITY: Security staff of the infrastructure are also called guards. During the evacuation process, their duty is to help those agents who are not severely injured but still require assistance to evacuate, i.e., those agents with broken bones or burns. Security's behavior is completely static. Their role after the blast is to rescue the injured continued until all the injured have been evacuated.

DEAD: Dead agents after the explosion are generated according to a survivability contour. During the running simulation they do not move.

SEVERELY INJURED: These agents have experienced serious injuries and cannot evacuate on their own; they need specialist rescue because their injuries are life-threatening.

INJURED BROKEN BONES / BURNS: These agents are injured but not severely and are able to evacuate but the speed at which they can evacuate is compromised by their injuries. This category of agents may be assisted in their evacuation path by altruists or rescurers.

INJURED EARDRUM RUPTURE: Because of the high-pressure shock waves from the blast, these agents have suffered eardrum rupture, and as a result experience a sense of disorientation and partial loss of balance. These agents walk at a reduced speed due to their loss of balance. If they become followers, then it is possible that they may increase their walking speed due to the leader guide, which may compensate for their disorientation and loss of balance.

LEADER: Leaders are randomly generated according to fixed percentages. Thanks to their automatic assumption of a leadership role, they may be followed by followers. They also might transform into Altruists if they stop and assist the injured in evacuating.

FOLLOWER: Followers are visitors or those whose eardrums have ruptured who decide to follow a leader. This role may be temporary if they lose sight of the leaders they were following; at that point they transform into Stop Following agents.

STOP FOLLOWING: These agents are followers who have lost the leaders they were following. Within a time interval of a few seconds, they may transform into followers if they are able to reconnect with their leader. Otherwise, they return to being visitors.

ALTRUISTS: Altruists are leaders or visitors that see an injured agent and stop to assist the injured agent in evacuating the site; therefore, the speed at which they evacuate is compromised.

RESCUERS: Rescuers are security agents that assist other agents; therefore, the speed at which they evacuate is compromised.

### 6.3.5 Measured Parameters and Variables

The parameters used in order to measure the effectiveness or otherwise of the BDI paradigm are as follows:

- *Numerical:* number of agents for each category

- *Geometric:* placement of the bomb rather the distance between bomb and emergency exit
- *Time:* evacuation times (partial and total)

*6.3.5.1  Definition of Evacuation Time*

Evacuation time is an important time concept in the evacuation process. Evacuation time is defined as the elapsed time between the instant that all the occupants who are able to evacuate receive an alarm and their arrival at a nominally safe destination inside or outside the infrastructure. Occupants are those agents who are not severely injured and are able to able to evacuate without assistance. Partial evacuation time is defined as the elapsed time between the instant that all the agents belonging to specific category are able to evacuate to a safe area.

## 6.4   BELIEF-DESIRE-INTENTION PARADIGM (BDI)

The human behavior that we want to simulate is the decision-making process of an individual in a crowd by using the BDI paradigm, which can be employed in order to imitate the human reasoning and the decision-making process. As shown in Figure 6.1, the BDI paradigm can be divided into sub-modules as follows: Belief-sub-module, Desire-submodule, Intention-sub-module, and Decision-making sub-module. Beliefs are information that an individual possesses regarding a situation. They may be incomplete or incorrect due to the nature of human perception. Desires are the states of affairs that a human would wish to see manifested. Intentions are desires that a human is committed to achieve. Zhao and Son [2008] extended the intentions sub-module to include other sub-modules. As shown in Figure 6.2, Zoumpoulaki et al. [2010] extended the classic BDI framework by incorporating Personality and Emotions.

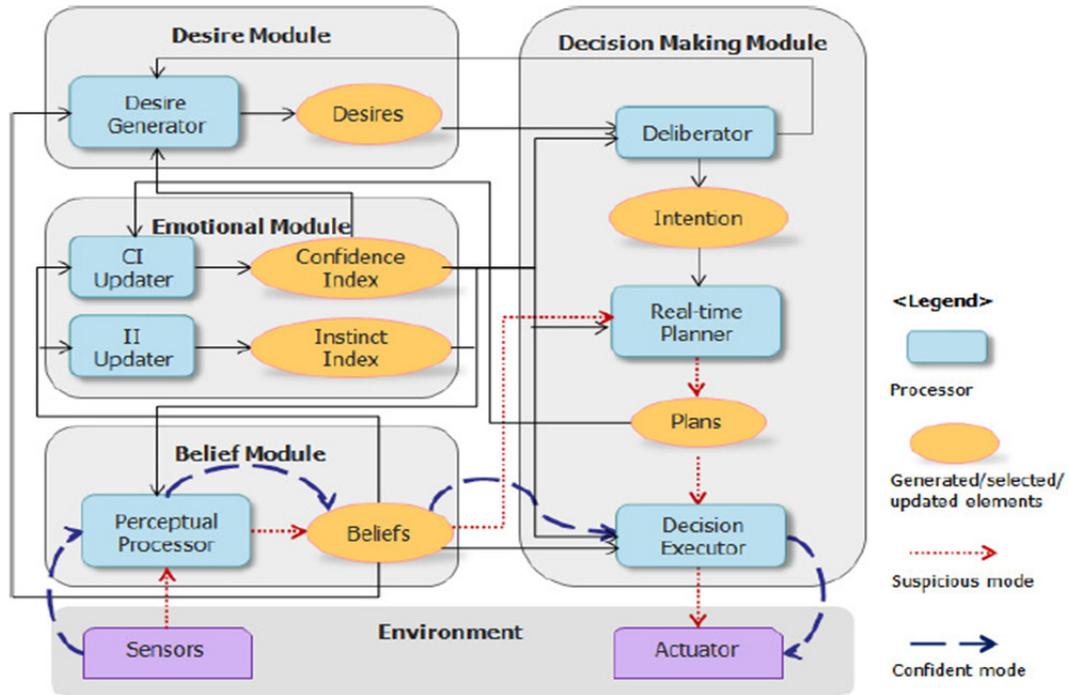

Figure 6.1     Depiction of BDI paradigm.

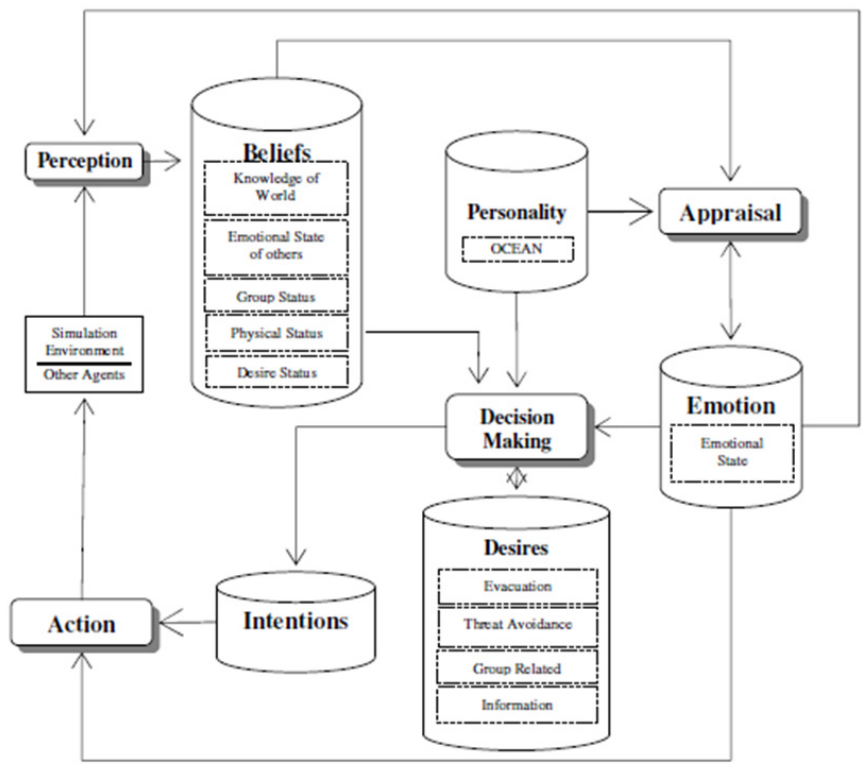

Figure 6.2     Depiction of Zoumpoulaki BDI framework.

### 6.4.1 The Decision Making Module and the Extended Decision Field Theory

The decision making module is the core of a BDI paradigm. This module is implemented through the decision field theory developed by Busemeyer [2002] and subsequently extended by Lee [2009]. The extension allows for updates in the model under a dynamic environment, using the *subjective evaluation* and the *attention weights* for alternatives. It provides a dynamic and probabilistic mathematical approach in order to simulate human deliberation process in making decision under uncertainty. Dynamics due to the variable "time" is a factor that affects decisions as well as the changing of the environment.

The EDFT calculates the dynamic evolution of preferences among *n* options expressed by an agent over time using the linear system formulation expressed in Equation (6.1) The parameters that could change the preference choice are the *m* attributes

$$P(t+h) = SP(t) + CM(t+h) \cdot W(t+h) \tag{6.1}$$

where $P(t)_T = [P_1(t), P_2(t), ..., P_n(t)]$ represents the preference state, and $P_i(t)$ is the strength of preference corresponding to option *i* at time *t*. The preference state is updated at every time step *h*.

The first member of the equation represents the memory effect. It is the product of the preference chosen at the previous state $P(t)$ and the *stability matrix*, *S*. The diagonal elements of *S* are the memory for the previous state preferences and all are assumed to have the same value. The off-diagonal elements are the inhibitory interactions among competing options. These assumptions ensure that each option has the same amount of memory and interaction effects.

*C* is the *contrast matrix* comparing the weighted evaluations of each option. If each option is evaluated independently, then *C* will be *I* (identity matrix). In this case, the preference of each option may increase simultaneously; see Equation (6.1). Alternately, the elements of the matrix *C* may be defined as $c_{ii} = 1$ and $c_{ij} = -1/(n-1)$ for $i \neq j$, where *n* is the number of options.

*M* is the value matrix, (*n* × *m*) vector, where *n* is the number of options, and *m* is the number of attributes; this represents the subjective evaluations (perceptions) of a decision-maker for each option of each attribute. If the evaluation value changes according to the environment, the matrix *M* is constituted with multiple states.

Finally *W* is the *weight vector,* (*m* × 1) vector, where *m* is the number of attributes. It allocates the weights of attention corresponding to each column (attribute) of *M*. In the case where *M* is constituted with multiple states, each weight $W_j(t)$ corresponds to the joint effect of the importance of an attribute, and the probability of a state. *W(t)* changes over time according to a stationary stochastic process.

## 6.5 QUESTIONNAIRE TO CALIBRATE THE HUMAN BEHAVIOR UNDER EMERGENCY CONDITIONS

### 6.5.1 Objective

In order to simulate the leader–follower dynamics and altruism among agents, it is necessary to define the most common behaviors of leaders and followers under emergency conditions. Therefore, a survey was created to assess the percentage of leaders and followers based on a sampling of individuals. In order to evaluate the behavior as a function of different boundary conditions, it examines how agents (leaders or followers) would behave during the evacuation if they saw an injured or another agent who needed help. A statistical distribution of the answers was calculated using the vector *W* and preferences *P*.

### 6.5.2 Structure of the Questionnaire

To compose the questionnaire, surveys used in other scientific research were taken as examples. The guidelines provided by the theory of planned behavior (TPB) [Ajzen 1991] were used. The process for constructing a TPB survey is described next.

First, the behavior being investigated must be clearly defined. The categories of agents are:

*LEADER:* The leader tends to go towards the emergency exit, he/she accepts a higher risk and he is determined.

*FOLLOWER*: The follower tends to follow a group of people or a leader. The boundary conditions that can affect an agent's decision are as follows:

- The agent LOCALIZES an EMERGENCY EXIT
- The agent is INJURED.
- The agent MEETS another agent on his/her way who NEEDS HELP.

Therefore, the behavior of two types of agents—leaders / followers—inside an environment after an explosion must be simulated.

#### 6.5.2.1 Specifying the Research Population

The population of interest to the investigators must be clearly defined. The sample size consisted of a minimum of 100 individuals, ranging in ages between 15 and 75, and males and females of all educational levels.

#### 6.5.2.2 Question Types

The response options were mutually exclusive, close-ended questions:

Ordinal-polytomous, where the respondent has more than two ordered options.

| Very unlikely | Unlikely | Neutral | Likely | Very Likely |
|---|---|---|---|---|
| ○ | ○ | ○ | ○ | ○ |

| Definitely | Probably | Not sure | Probably not | Definitely not |
|---|---|---|---|---|
| ○ | ○ | ○ | ○ | ○ |

(Bounded) Continuous, where the respondent is presented with a continuous scale.

| 1 | 2 | 3 | 4 | 5 |
|---|---|---|---|---|
| ○ | ○ | ○ | ○ | ○ |

The boundary conditions or possible attributes were as follows:

- S = Health status (injured-not injured)
- E = Emergency exit location (see the emergency exit-do not see the emergency exit)
- I = Presence of injured (encounter an injured- do not encounter an injured)
- For both cases possible conditions are $2^3= 8$

Table 6.1 below lists all possible combinations among the three options and their applications that were used in order to assess agent behavior. The possible agent behaviors are:

- the agent evacuates by his/her own
- the agent follows a leader
- the agent stops to help another agent

**Table 6.1    Options combinations and relative survey questions.**

| NOT INJURED<br>DO NOT SEE EMERGENCY EXIT<br>DO NOT SEE INJURED | NOT INJURED<br>DO NOT SEE EMERGENCY EXIT<br>SEE INJURED | NOT INJURED<br>SEE THE EMERGENCY EXIT<br>DO NOT SEE INJURED | NOT INJURED<br>SEE EMERGENCY EXIT<br>DO NOT SEE INJURED |
|---|---|---|---|
| $\begin{pmatrix} QUESTION\ 1 \\ QUESTION\ 2 \\ NONE \end{pmatrix}$ | $\begin{pmatrix} QUESTION\ 1 \\ QUESTION\ 2 \\ QUESTION\ 3 \end{pmatrix}$ | $\begin{pmatrix} QUESTION\ 11 \\ QUESTION\ 9 \\ NONE \end{pmatrix}$ | $\begin{pmatrix} QUESTION\ 11 \\ QUESTION\ 9 \\ QUESTION\ 13 \end{pmatrix}$ |
| INJURED<br>DO NOT SEE EMERGENCY EXIT<br>DO NOT SEE INJURED | INJURED<br>DO NOT SEE EMERGENCY EXIT<br>SEE INJURED | INJURED<br>SEE EMERGENCY EXIT<br>DO NOT SEE INJURED | INJURED<br>SEE EMERGENCY EXIT<br>SEE INJURED |
| $\begin{pmatrix} QUESTION\ 5 \\ QUESTION\ 6 \\ NONE \end{pmatrix}$ | $\begin{pmatrix} QUESTION\ 5 \\ QUESTION\ 6 \\ QUESTION\ 7 \end{pmatrix}$ | $\begin{pmatrix} QUESTION\ 12 \\ QUESTION\ 10 \\ NONE \end{pmatrix}$ | $\begin{pmatrix} QUESTION\ 12 \\ QUESTION\ 10 \\ QUESTION\ 7 \end{pmatrix}$ |

### 6.5.3 Survey

The questionnaire using the online form-building tool Adobe Forms central ® was developed as a submission-enabled PDF form or a web form was distributed. The subjects filling out the survey were asked to imagine being in an environment with certain boundary conditions (injured or not injured / see an emergency exit or do not see an emergency exit / there is an injured on his path- there is not an injured on his path). They were asked to choose the following options: to follow another individual, to evacuate on their own, or to stop to assist a person who needed help. To make the answers truthful as possible and provide psychological context, before filling out the survey subjects were asked to watch a video of some real explosions.

### 6.5.4 Introduction to the Survey

Before filling in the questionnaire, subjects were informed about the reason and purpose of the research study. A privacy statement and contacts for requests for information were provided.

**Informed Content:** You are being invited to take part in a research study. The information in this form is provided to help you decide whether or not you want to take part.

**What is the purpose of this research study?** The objective of this project is to use an emergency scenario survey to analyze how people evaluate emergencies situations and to develop an accurate simulation model of an emergency evacuation, such as from an explosion, taking into account human behavior.

**Will the information that is obtained from me be kept confidential?** No personal information of yours will be collected. The only persons who will know that you participated in this study will be the research team members; specifically, the Principal Investigators and the advisor. Your responses will be confidential. You will not be identified in any reports or publications resulting from the study.

**May I change my mind about participating?** Your participation in this study is voluntary. You may stop the study at any time. Also any new information discovered about the research will be provided to you.

**Whom can I contact for additional information?** You can obtain further information about the research or to voice concerns or complaints about the research by calling the Principal Investigators.

Sex, age, and education questions.

| Gender | Age | Education |
|---|---|---|
| ○ male | ○ 15-30 | ○ elementary school |
| ○ female | ○ 31-45 | ○ middle school |
| | ○ 46-60 | ○ high school |
| | ○ 61-80 | ○ bachelor degree |
| | | ○ master degree or higher |

### 6.5.5 Leader–Follower Categorization

The first step was to create questions that effectively self-selected Leaders or Followers. Three specific questions were prepared. Question 1 asked the following:

Overall imagine the scenario described at the beginning. So, you're in a museum and a blast occurs, regardless of whether you are injured or not:

| | Definitely it's safer evacuate by my own | Probably it's safer evacuate by my own | I don't know | Probably it's safer to join a group | Definitely it's safer to join a group |
|---|---|---|---|---|---|
| Do you think it is safer find an emergency exit by yourself or to join to a group that is evacuating? | ○ | ○ | ○ | ○ | ○ |

Potential leaders were individuals who responded to the question above with the following answers:

- Probably it's safer evacuate by my own
- Definitely it's safer evacuate by my own

Questions 1 and 2 were designed to verify the responder as a leader.

| | Very unlikely | Unlikely | Neutral | Likely | Very Likely |
|---|---|---|---|---|---|
| How likely is that you would decide to evacuate on your own? | ○ | ○ | ○ | ○ | ○ |

| | Definitely | Probably | Not sure | Probably not | Definitely not |
|---|---|---|---|---|---|
| Would you decide to join a group of people who are running away from the blast? | ○ | ○ | ○ | ○ | ○ |

A leader must show a strong propensity to evacuate by his/her own and a low propensity to follow a group of people that is evacuating. This means that the subject has to respond to Questions 1, "How likely is that you would decide to evacuate on your own?" as follows:

- Likely
- Very likely

In addition, answers to the Question 2, "Would you decide to join a group of people who are running away from the blast?" should be as follows:

- Probably not
- Definitely not

The survey results formed the basis of a numerical analysis using the mathematical operations presented earlier: the "weight vectors" $W$, the "preference vector" $P$, and through an inverse formula the $M$ matrices. To obtain the most accurate results from the survey, those responders who answered NEUTRAL to Question 2 the answer NEUTRAL were categorized as leaders were also considered agents. This less restricted categorization was only used for calculations.

### 6.5.6 Survey Scenarios

The scenarios were presented in the second person singular, with each followed by a series of questions asking the participant to imagine themselves in the scenario. Several survey studies have demonstrated the feasibility of this design in measuring respondent behaviors from different aspects proposed by TPB, including intention, attitude, subjective norms, and perceived behavior control.

### 6.6    IMPLEMENTATION OF THE HUMAN BEHAVIOR MODULE

In this section, implementation of the human behavior module will be discussed. The goal of this research is to study evacuation dynamics under emergency conditions and to see how they change depending on human behavior. The evacuation time is the main parameter of response and is used to evaluate the efficiency and safety of an infrastructure.

The first step was to define a block framework in order to model the human decision-making process represented in Figure 6.3. The agent is immersed in an environment. It is able to perceive within a given radius of vision the environment's characteristics. The perception by the beliefs possessed by each agent is filtered. Beliefs are information that a human possesses about a situation; they may be incomplete or incorrect due to the nature of human perception. This module together with the influence of personality and the emotional evaluation of the agent are incorporated into the decision-making module. Mathematically, all these blocks modifies the matrices and vectors that make up the DFT mathematical model. The decision-making process allows an agent to generate intentions. Intentions are the completion of the agent's desires. Finally, the intentions were carried out as real actions in the environment. Therefore, the detailed blocks of the decision-making process are discussed below.

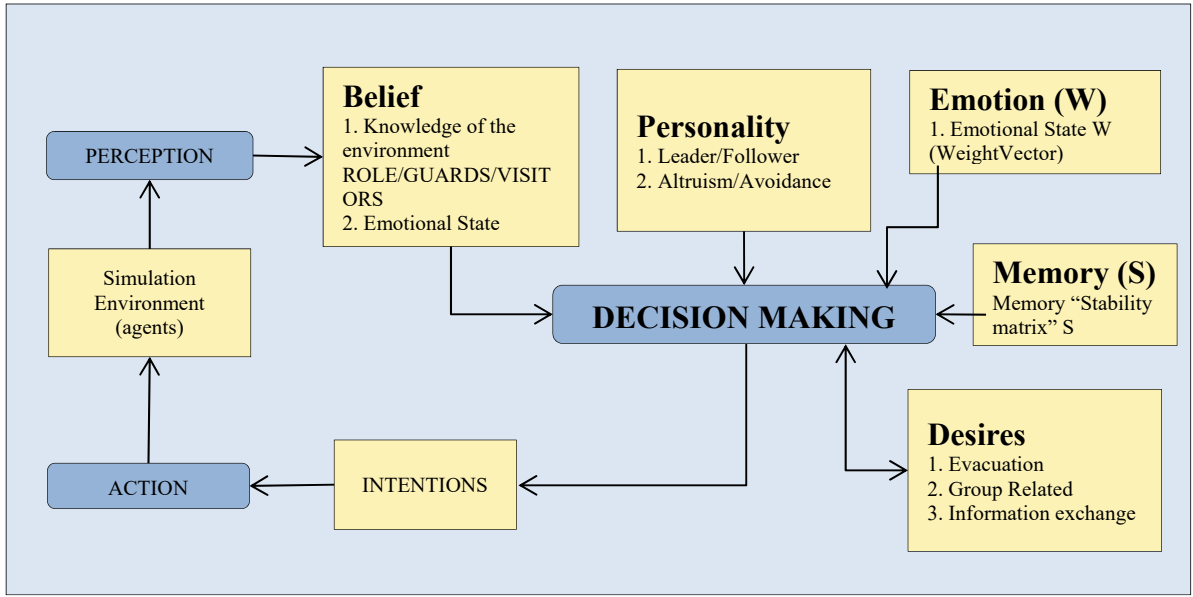

**Figure 6.3    Extended BDI architecture.**

### 6.6.1  Simulation Environment

The simulation environment represents a virtual discretization of the real world. In Case Study 1 it is a museum; in Case Study 2 it is a metro station. Constructing the environment is one of the main steps of the simulation because it provides the stimulus to the agents. These stimuli are absorbed by the agents as perceptions.

### 6.6.2  Perception

All the stimuli and information from the virtual world by the agent are acquired. This phase is called Perception. In the simulations, agents are able to perceive injured persons, emergency exits, or other persons as leaders or guards.

### 6.6.3  Belief

The agent's perception of the environment is filtered through their belief system.

1. First, the agent is influenced by his/her role and knowledge of environment, e.g., a guard knows the location of emergency exits while the visitor may not;
2. The emotional state may influence the danger perception of an agent;
3. The physical condition of an agent (injured / not injured) may affect human behavior; and
4. Status of the desire corresponds to memory. The agent's preference is influenced by the preference expressed at the previous instant: the preference at the previous instant is the desire to perform an action.

### 6.6.4 Personality

Personality is the most important aspect of the decision-making process. Each agent has a specific psychological profile that leads him/her to perform actions. In terms of characterizing agents, the first parameter to determine is whether an agent is a leader or a follower. Leaders tend to go forward to the emergency exit; they tend to accept higher risks and are more determined. Therefore, people tend to gather around and follow his/her lead. In contrast, followers tend to follow a group of people or a leader and it is unlikely that they will evacuate on their own. Agents are grouped into three categories: absent, avoidant, or altruistic towards other agents who needs help; they decide on an individual basis whether or not to stop and help other agents.

### 6.6.5 Emotions

Emotions are the irrational component of the human psyche and is simulated by vector $W$. First of all, the weight vector is calculated based on the responses to the survey and will have a form like this:

$$W(t+h) = \begin{pmatrix} a_1 - a_2 \\ b_1 - b_2 \\ c_1 - c_2 \end{pmatrix} \tag{6.1}$$

where the weight vector $W$ in the form of intervals is defined. Through a stochastic process, the preference expressed by agent is calculated and a numeric value is extracted and obtained. This extraction is conducted for each agent in order to simulate the "personality" of each agent and its unique decision.

### 6.6.6 Memory

Even the short-term memory may have an influence on the decision-making process. The equation's first member takes into account the memory effect, which is the product of the previous chosen preference $P(t)$ and of the stability matrix $S$. The $S$-diagonal elements are the memory for the previous preferences. The off-diagonal elements are the inhibitory interactions among competing options.

### 6.6.7 Desires

Desire also influences the decision making process and depends on the situation in which agents are involved. For example, it may be the need to help an injured person or a family member, or even try to find an emergency exit.

### 6.6.8 Decision Making

Decision-making is the block where all the previous modules converge. All these modules (in the form of boundary conditions or numbers) constitute the inputs of this module. Thanks to the mathematical theory of the DFT, a preference of an agent to perform a specific action is calculated. Therefore, the output of this module is a preference in the performance of an action.

### 6.6.9 Intention

The preference expressed in the previous block by the human brain in an action is processed.

### 6.6.10 Action

In this module, the intention into an effective action is realized. The action undertaken within the environment will have an impact. This means that the cycle begins again until the performance of a new action.

## 6.7 CASE STUDIES

The evacuation model implemented NetLogo 5.0.5 (December 2013), which is a multi-agent based programmable simulation software developed by Wilensky [1999]. Two different ABM evacuation models of a museum and of a metro/train station were developed. Each of these models was divided in two phases: the normal dynamic and the evacuation process of the agents after the occurrence of a blast. In both of them, the agents were able to perform decisions, but some specific irrational behaviors by agents were also included. The evacuation time is the main parameter of response, and it is used to evaluate the efficiency and safety of the infrastructure.

### 6.7.1 Case Study 1: Ursino Castle Museum

First, an agent-based evacuation model of a museum was developed. The model is divided into two phases of fruition: the normal dynamic and the evacuation process of agents after the occurrence of a blast. The evacuation time is the main parameter of response and is used to evaluate the efficiency and safety of the infrastructure. Different scenarios in order to assess what happens in different situations were analyzed.

#### 6.7.1.1 Dimensional Data

This model was developed based on the geometry of the Ursino Castle Museum, which is located in Catania (Sicily, Italy). Figure 6.4 shows the plan of the museum. The castle consists of 10 rooms located around a central courtyard. The main entrance, represented with a red arrow, is the only emergency exit that allows exit from the castle. The other two emergency exits, depicted in yellow, have access to the courtyard. To simulate a real situation, it was hypothesize that three possible scenarios of explosion occurred and for each scenario the probability of death, injuries, and ruptured eardrums were calculated.

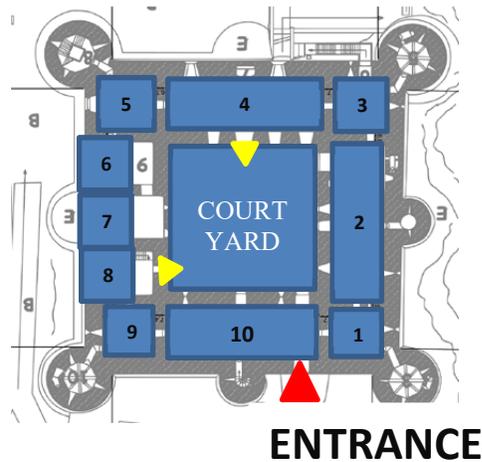

**ENTRANCE**

Figure 6.4    Depiction of the museum plan.

### Case Study 1 scenarios

Different scenarios to determine what happens in different situations were analyzed.

### Scenario 1 – Room 4

In the first scenario, an explosion in the middle of the fourth room was simulated.

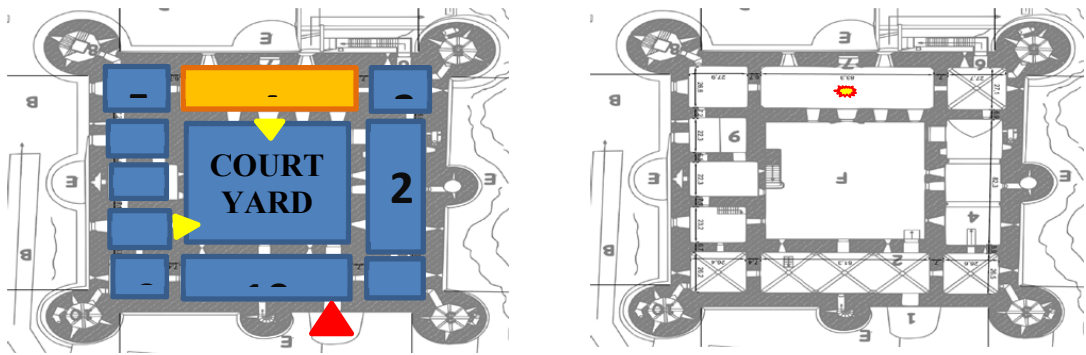

Figure 6.5    Scenario 1: Museum's plan with location of explosion.

### Scenario 2 – Room 10, near the emergency exit

In the second scenario, the blast was located in the 10th room near the emergency exit.

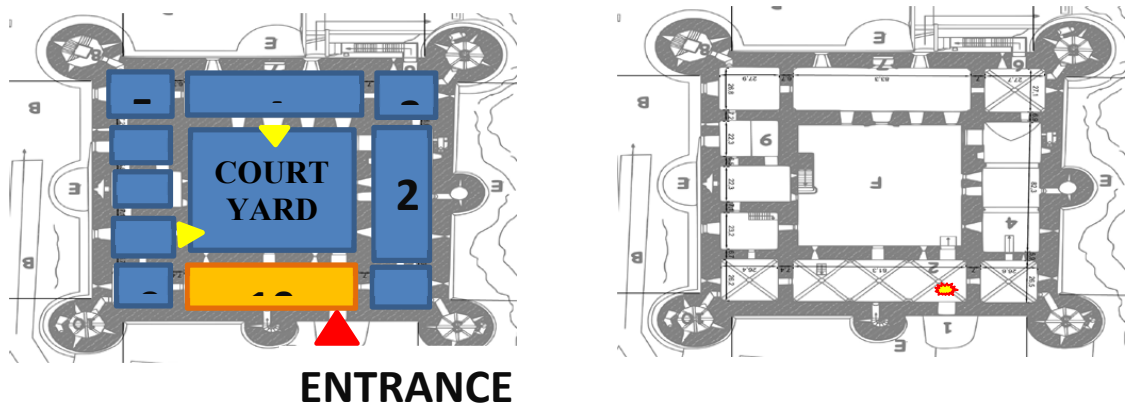

Figure 6.6　　Scenario 2: Museum's plan with location of explosion.

### Scenario 3 – Rooms 2 and 4

The third scenario simulated the worst configuration in terms of dead and injured because two consequential explosions occurred. As shown in Figure 6.7, the first blast occurred in Room 2 and the second explosion in Room 4.

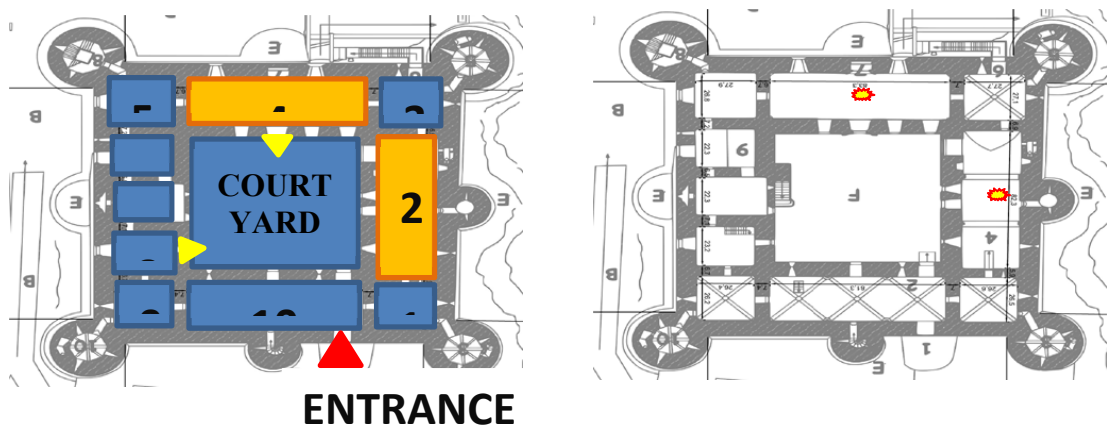

Figure 6.7　　Scenario 3: Museum's plan with location of explosions.

**Generating Injuries**

Based on a study conducted by the Army, Navy, and Airforce that determined blast effects, the distribution of the dead and injured was calculated. For each scenario, a survivability contour was calculated and the range of probability of survival established; see Figure 6.8. Specifically, a random number for each agent was extracted. In the red contour, if this number is in the range between 1 and 99, the agent dies. In the yellow contour, if the number is between 1 and 50, the agent will be severely injured to the point that their survival is questionable. In blue contour, survival rates vary between 50 and 99% of the cases, and the agents who do survive will suffer from broken bones and burns. The last contour ensures a survival rate of 99%, who will report eardrum rupture, resulting in loss of orientation and partial loss of balance. Each scenario holds its own survivability contours with relative diameters.

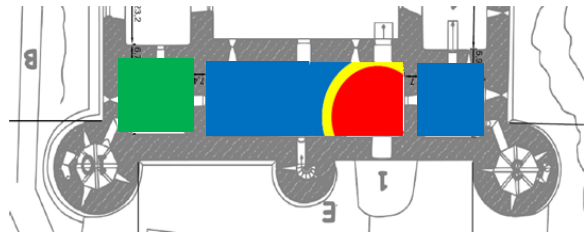

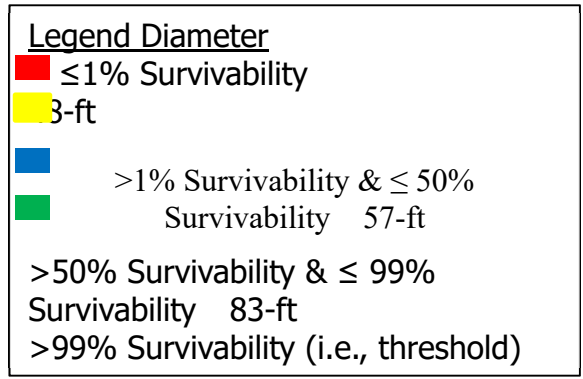

Figure 6.8    Depiction of an example of survivability contours.

### 6.7.2  Case Study 2: Gare de Lyon Station

The evacuation model developed was applied to the Gare De Lyon metro station, which is a node of the Paris Metrò infrastructure encompassing lines 1 and 14. The model is divided into two phases: the normal dynamic and the evacuation process of the agents after the occurrence of a blast. The evacuation time is the main parameter of response, which is used in order to evaluate the efficiency and safety of the infrastructure.

#### 6.7.2.1  Dimensional Data

The Metro 14 line station is rectangular in shape with sides 142 m to 18 m (465 ft 10 in. to 59 ft). Access stairs to the platform are located at the two ends of the rectangle. The station has two tracks on either side of a large central platform shielded by automatic platform screen doors. To access the trains, there are 14 double doors of 1.5 m (4 ft, 11 in.) on both sides.

#### 6.7.2.2  Model Implementation

The agent-based evacuation model was implemented based on Larcher et al.'s [2011] investigation of an explosion inside a rail system .Fluid–structure interaction calculations (FSI) calculated the displacements and probable structural failure, and the probabilities of death and eardrum rupture. The outputs of these simulations as input of the model of evacuation were used. As shown in Figure 6.9, several considerations were made in generating the geometry of the model using NetLogo, which is based on a two-dimensional grid of patches; see Section 7.2.1. Patches size is measured in terms of pixels and turtles, and movements are measured in patches. Each patch measures 0.5 m × 0.5 m, which fits the needs of this simulation perfectly. A person fills approximately 0.5 m². The turtles have been set for a size equal to a patch (i.e., 0.5 m²). By using a "while" loop, geometries were defined. The image b shows the metro station geometry developed in NetLogo where:

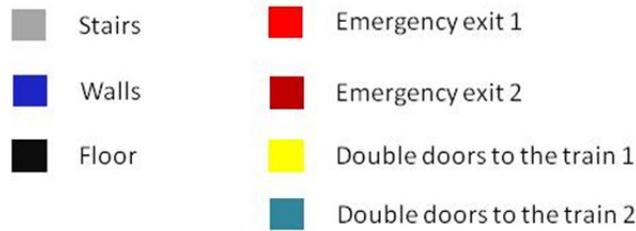

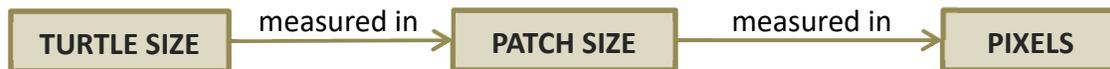

Figure 6.9    Depiction of metro station in NetLogo with legend.

#### 6.7.2.2.1  Normal Dynamic

The simulation starts when the *GO button* is pressed. The simulation begins from the two entrances whereby a flow of travelers enters the metro station and will focus on the track depending on their destination and nearest free door. In the normal fruition, trains arrive every 3 min. at 100% punctuality. To simulate a real situation, a slider button that allows for selecting the interval of arrival time between the two trains in the opposite direction was created. The cursor, which can be adjusted, has an arrival interval between 1 sec and 90 sec.

#### 6.7.2.2.2  Emergency Dynamic

Three models (base, deterministic, and probabilistic) were developed. The emergency phase of each model is identical; however, each of these models includes their own unique human behavior simulation dynamics. In the Metro Station, the simulation considered 180 agents subjected to a blast from 50 lbs of explosive. At that point, the dynamics of evacuation begin.

### 6.7.2.3  Agent's Categories and Evacuation Speed

According to the simulation the explosion caused both death and injuries. Those agents who were severely injured are unable to move. Those agents with broken bones and burns are in the category of seriously injured that are able to evacuate but at a slower pace than others. If helped by other agents, these agents may improve their evacuation speed. Those subjected to eardrum rupture are considered visitors who suffer from with consequent loss of sense of direction, which may impact their ability to evacuate. The evacuation speed of each agent is categorized as follows:

- Uninjured visitors: 1.5 m/sec
- Those agent injured due to eardrum rupture suffer from disorientation and problems with balance and therefore walk unsteadily at a speed of 1 m/sec
- Those agents helping other injured agents: 0.6 m/sec
- Those agents injured due to broken bones and/or burns injured: 0.3 m/sec

Speeds are scaled according to the size of the patches (0.5 m x 0.5 m). The NetLogo command for setting the speed of agents is $f_{ad}$. To avoid possible program bugs, $f_{ad}$ should be set < 1. A higher value would result in the agents moving beyond one patch, thus being able to penetrate walls or other obstacles. Therefore, the maximum speed at which agents are categorized (1.5 m/sec) is designed by the command $f_d = 1$, i.e., agents move one patch (0.5 m) per second; thus, agents move with a speed of 0.5 m/sec. Therefore, evacuation times must be scaled by factor $f_s$ $\left( f_s = \dfrac{1.5 \text{ m/sec}}{0.5 \text{ m/sec}} = 3 \right)$.

Table 6.2    Metro station agents speed.

| Agents | Real speed | Scaled speed | $f_d$ |
|---|---|---|---|
| Visitor (emergency) | 1.5 m/sec | (1.5 m/sec)/3 =0.5 m/sec | $f_d$=1 |
| Visitor (normal phase) | 1 m/sec | (1 m/sec)/3 =0.33 m/sec | 1 : 0.5(m/sec) = x : 0.33 (m/s) → $f_d$ = 0.66 |
| Eardrum rupture | 1 m/sec | (1 m/sec)/3 = 0.33 m/sec | 1 : 0.5(m/sec) = x : 0.33 (m/sec) → $f_d$ = 0.66 |
| Helped injured | 0.6 m/sec | (0.6 m/sec)/3 = 0.2 m/sec | 1 : 0.5(m/sec) = x : 0.2 (m/sec) → $f_d$ = 0.4 |
| Injured | 0.3 m/sec | (0.3 m/sec)/3 = 0.1 m/sec | 1 : 0.5(m/sec) = x : 0.1 (m/sec) → $f_d$ = 0.2 |

## 6.8    NUMERICAL RESULTS

The analysis described below focused on how the human behavior (BDI) affects the evacuation times.

### 6.8.1  Museum Ursino Castle

The Ursino Castle model, unlike the Metro Station model, was not based on an existing model. This affected the level of detail available for the Museum Ursino analysis, and fewer simulations were conducted. The base model simulation results for the three scenarios show the dependence of evacuation time depending on bomb placement. Each agent according to their respective health status had different walking speeds. Therefore, the greater the distance of the bomb from the emergency exit more, the greater the speed of the evacuees. For example, Scenarios 1 and 3 show comparable results because bomb placement was identical. In Scenario 2, the evacuation time is lower compared to Scenarios 1 and 3 because the bomb was placed near the emergency exit.

In the Ursino Castle model, the leader–follower dynamics are valuable because the infrastructure environment is complex. This means that the emergency exits are not always visible; therefore, agents follow a leader. The results of the Ursino Castle deterministic model show that the follower evacuation time is lower than the eardrum rupture and visitors partial evacuation time. Therefore, a leader–follower dynamic can have a significant effect. Finally, combining altruism and leader–follower dynamic, the evacuation time decreases. The same simulation with 20 Leaders had a positive impact only in Scenario 2.

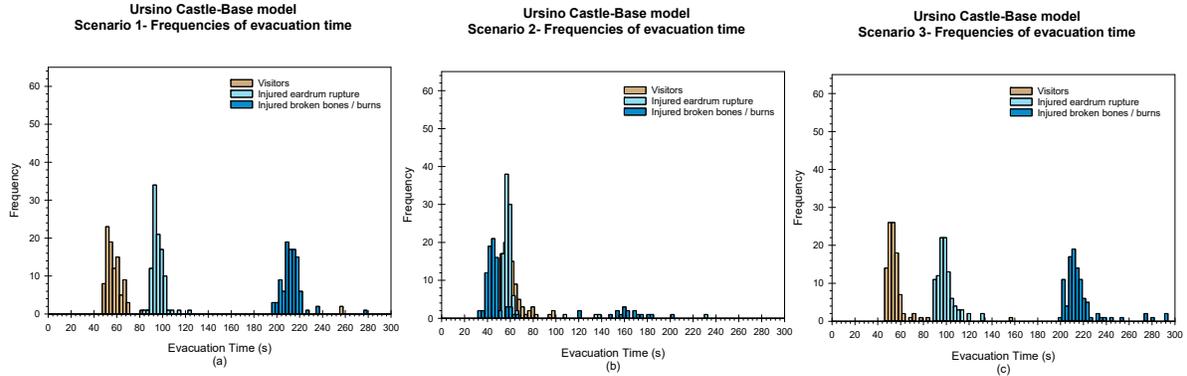

**Figure 6.10** Base models evacuation time. (a) Scenario 1, bomb in the Room 4; (b) Scenario 2, bomb in Room 2; and (c) Scenario 3, bombs in Room s and 4.

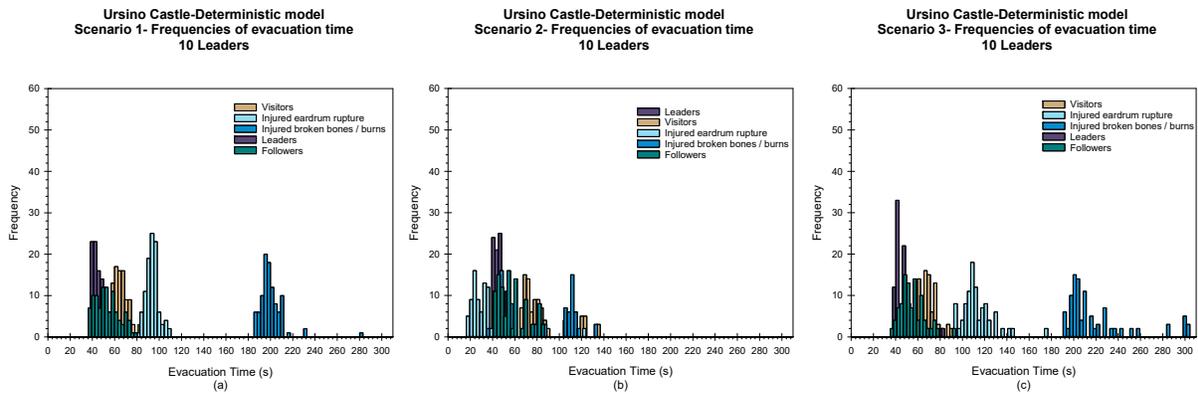

**Figure 6.11** Ursino Castle deterministic model considered 10 leaders: (a) Scenario 1; (b) Scenario 2; and (c) Scenario 3.

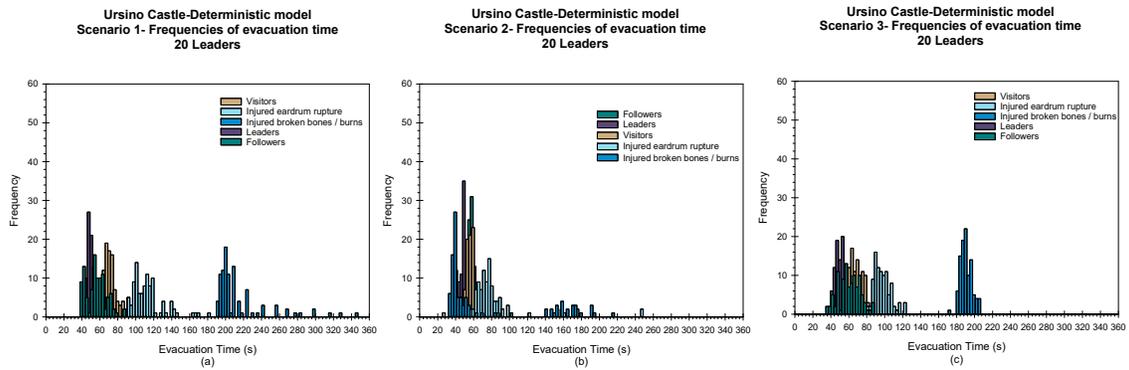

**Figure 6.12** Ursino Castle deterministic model considering 20 leaders.

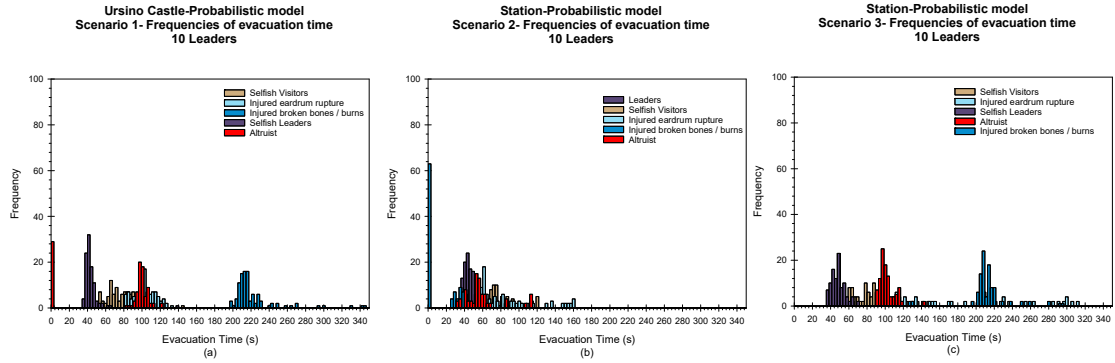

Figure 6.13    Probabilistic model considering 10 leaders.

### 6.8.2  Metro Station Gare de Lyon

For the metro station simulations, the evacuation time is also dependent on the distance of the bomb blast from the emergency exit. In the base model, factors that significantly influence the behavior of the agent—and thus the evacuation time—are injuries sustained from broken bones, burns, etc. Given their injuries, their speed is compromised, and they are usually the last agents to exit from the station. Because the distance between the bomb blast and the emergency exit decreases as we move from Scenario 1 to Scenario 3, naturally, the evacuation time reflects the same trend.

In the deterministic model, the variation due to the human behavior is considered, with the introduction of "leaders" and "followers." This variation between the base model and the deterministic model does not alter the outcome by a significant factor. Regardless of leader–follower dynamics, the injured remain the last agent category to evacuate the area. Therefore, an evolved probabilistic model was developed, and security, a new agent category, was incorporated, whose function is to help injured to evacuate and accelerate the evacuation process. Despite followers having a lower evacuation time compared to leaders, at this state the leader–follower dynamic is not yet significant. Increasing the number of leaders to 20 results in a marked improvement in the evacuation dynamics. Therefore, at this point in the simulation, all human behavior dynamics are inserted. The probabilistic model includes both altruistic behavior and leader–follower behavior.

Note that in these simulations that although there are leaders in the Metro Station, the leader–follower dynamic are totally absent; therefore, followers and their evacuation time were not recorded. This phenomenon is not isolated, but it occurs throughout the station probabilistic model. Because the metro station is a simple rectangle where the emergency exits are clearly visible, all agents move towards the emergency exits without following any other agent; see Figure 6.19. Also worth noting that in this figure, altruism does not occur in all scenarios; see Scenario 1. This phenomenon occurs only if the injured agent is placed between leaders and/or visitors and the emergency exit along the evacuation route. If the injured are in the central part of the station (Scenario 1) nobody stops to help other agents because nobody will retrace steps already taken. Hypothetically, if all agents were altruistic and assisted all injured agents, the evacuation time for all scenarios would be halved.

Doubling the number of leaders has no effect on the evacuation times, with the only change registered is that frequency of the injured, i.e., a more effective evacuation would occur if all the broken-bones/burns injured were rescued. For this reason, a variable number of security agents were included in the simulations. Security is a class of agents whose duty is to rescue any agents in trouble, and their role in the simulation does not end until all people have been rescued.

The results show a reduction in evacuation time. In those cases where the number of security guardians is minimal, the evacuation time increases. This is because several trips are necessary to rescue the injured, with a consequent increase in the evacuation time. As the number of leaders increase, the phenomena of altruism also increases. Therefore, the intervention of the security is very effective. In scenarios where the emergency exits are clearly visible, the leader–followers dynamic is more likely to occur, with an accompanying increase in altruism. According to a preliminary analysis, combining altruism and the intervention of the security can lead to a 50% reduction in evacuation times.

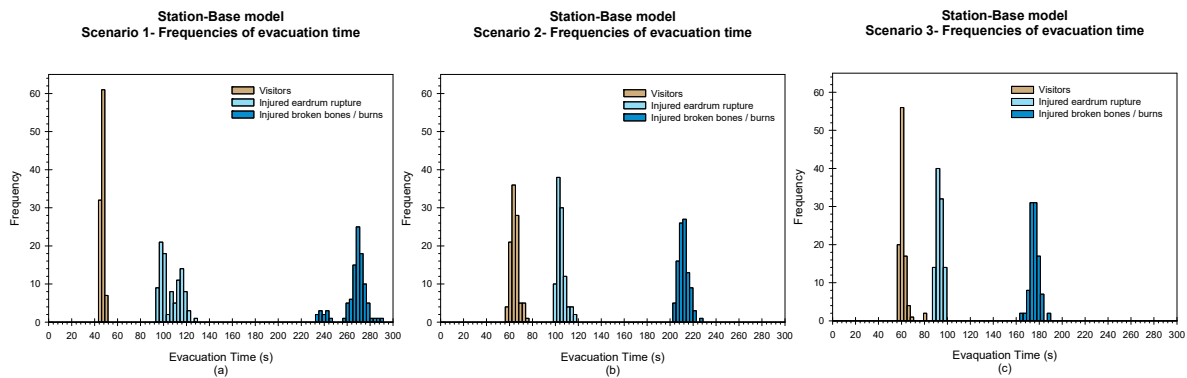

**Figure 6.14** Base models evacuation time. (a) Scenario 1 blasts 60 m away from the emergency exits; (b) Scenario 2: blasts 40 m away from the emergency exits; and (c) Scenario 3: blasts 30 m away from the emergency exits.

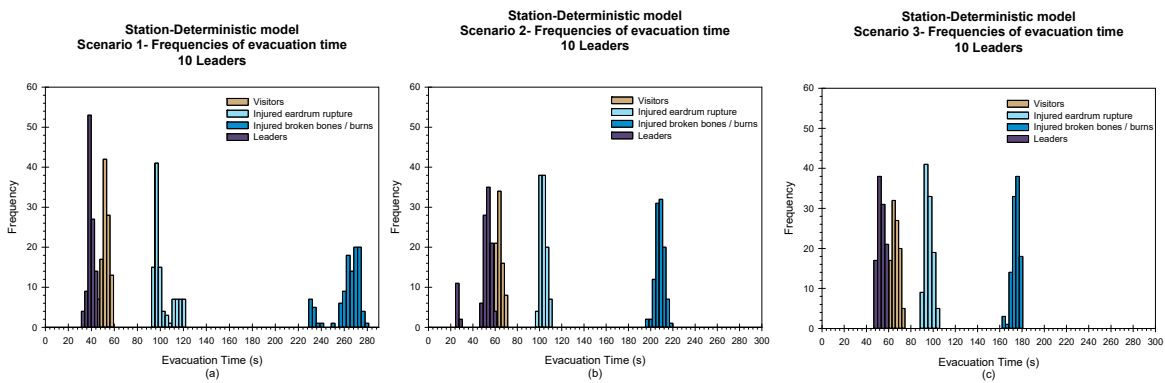

**Figure 6.15** Deterministic models evacuation time for 10 Leaders without security: (a) Scenario 1 blasts 60 m away from the emergency exits; (b) Scenario 2 blasts 40 m away from the emergency exits; and (c) Scenario 3 blasts 30 m away from the emergency exits.

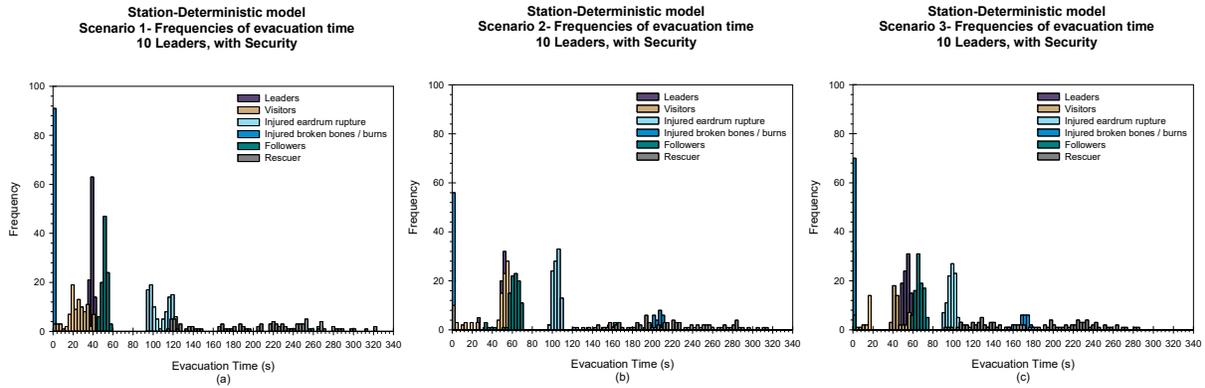

Figure 6.16  Deterministic models evacuation time for 10 Leaders with security: (a) Scenario 1 blasts 60 m away from the emergency exits; (b) Scenario 2 blasts 40 m away from the emergency exits; and (c) Scenario 3 blasts 30 m away from the emergency exits.

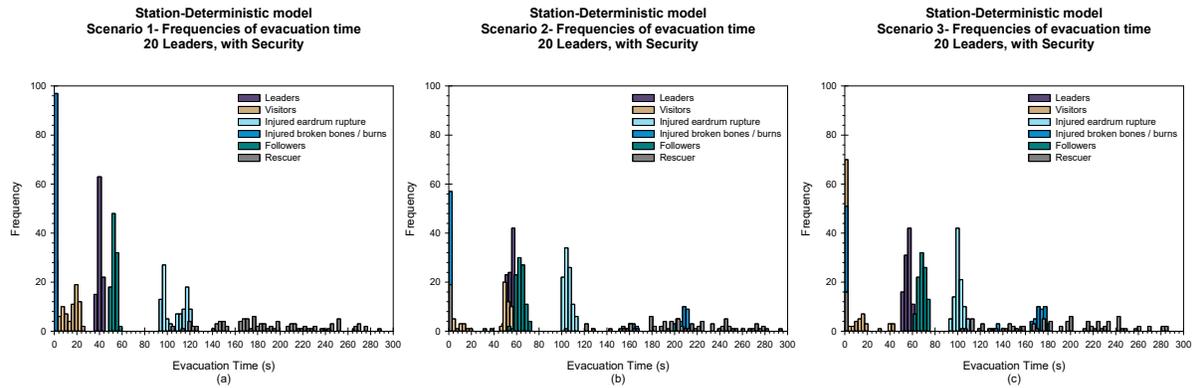

Figure 6.17  Deterministic models evacuation time for 20 Leaders with security, (a) Scenario 1 blasts 60 m away from the emergency exits; (b) Scenario 2 blasts 40 m away from the emergency exits; and (c) Scenario 3: blasts 30 m away from the emergency exits.

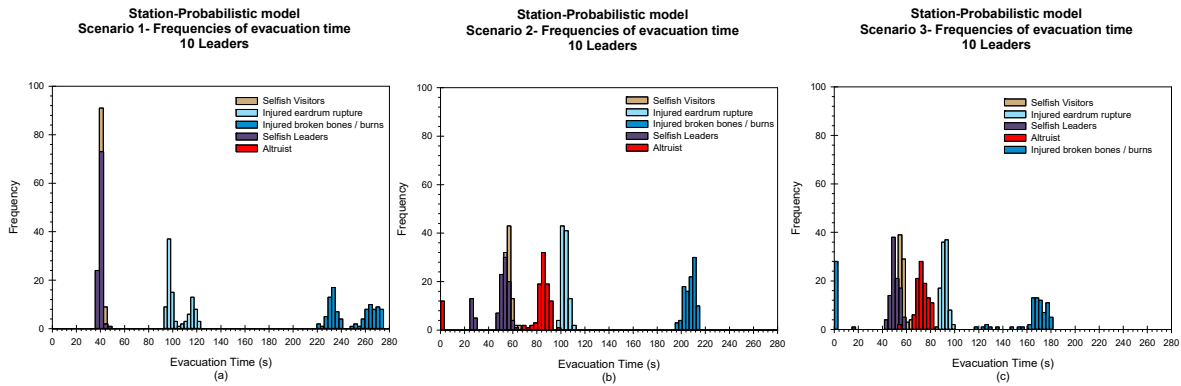

Figure 6.18  Probabilistic models evacuation time for 10 Leaders without security, (a) Scenario 1 blasts 60 m away from the emergency exits; (b) Scenario 2 blasts 40 m away from the emergency exits; and (c) Scenario 3 blasts 30 m away from the emergency exits.

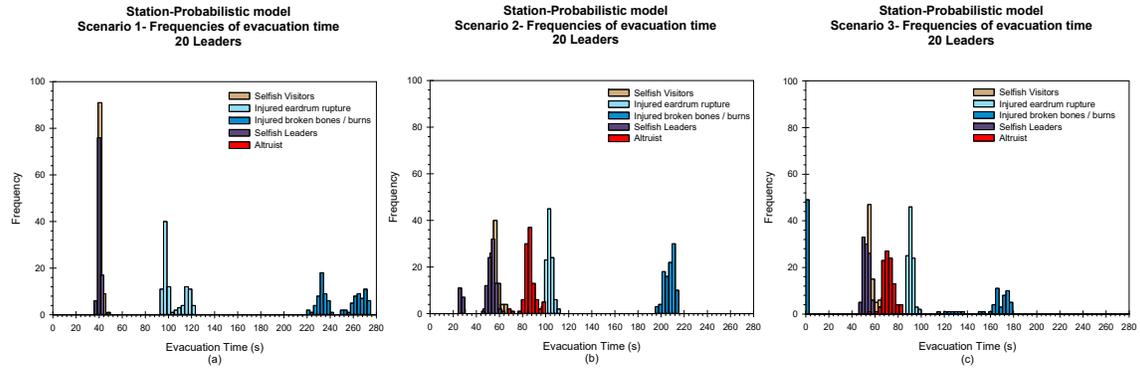

Figure 6.19    Probabilistic models evacuation time for 20 Leaders without security, (a) Scenario 1 blasts 60 m away from the emergency exits; (b) Scenario 2 blasts 40 m away from the emergency exits; and (c) Scenario 3 blasts 30m away from the emergency exits.

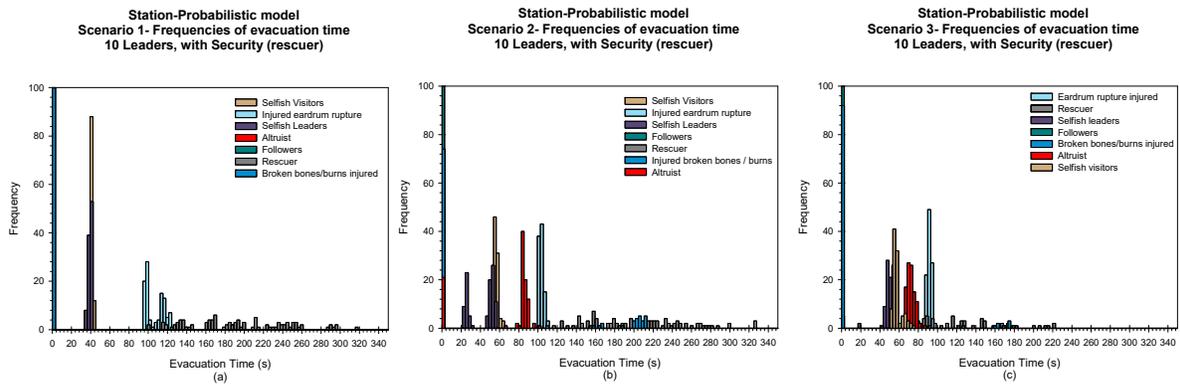

Figure 6.20    Probabilistic models evacuation time for 10 Leaders with security; (a) Scenario 1 blasts 60 m away from the emergency exits; (b) Scenario 2 blasts 40 m away from the emergency exits; and (c) Scenario 3 blasts 30 m away from the emergency exits.

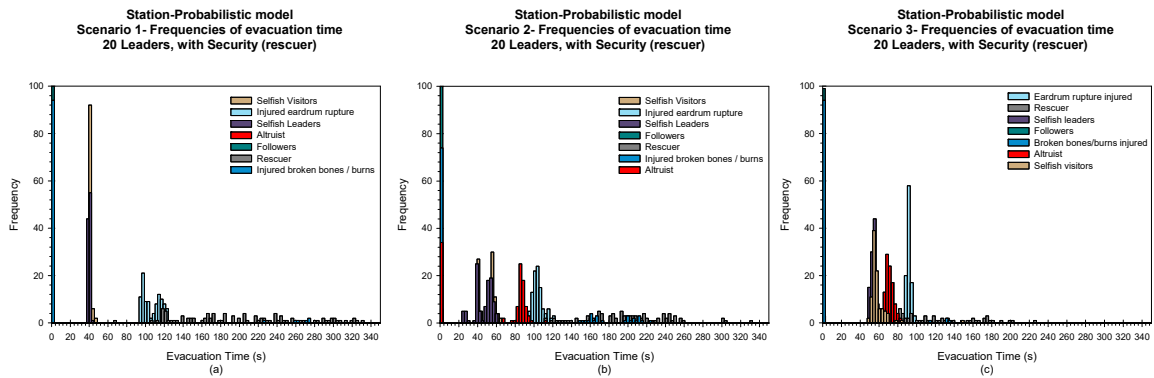

Figure 6.21    Probabilistic models evacuation time for 20 Leaders without security: (a) Scenario 1: blasts 60 m away from the emergency exits; (b) Scenario 2: blasts 40 m away from the emergency exits; and (c) Scenario 3: blasts 30 m away from the emergency exits.

## 6.9 REMARKS AND CONCLUSIONS

In this chapter, an Agent-based Model (ABM) of an infrastructure was studied to determine evacuation dynamics that considers human behavior in the case of a serious event. The two locations chosen for application of this ABM were a museum and a metro train station. These two facilities were specifically chosen because of their geometries. The Ursino Castle Museum consists of a series of rooms located successively around a squared courtyard. Therefore, the access from one room to another is through small doors, and the emergency exits are not always visible. In contrast, the geometry of the metro station is simple. It is rectangular in shape, with the easily visible emergency exits located at the two ends of the rectangle.

We began with a basic model and made it more complex by including the influence of human behavior in an ABM. The results are as follows: although the geometry may be the same, it cannot be assumed that all agents will react in the same manner. This was demonstrated in the study of the train station where the emergency exits are clearly marked and always visible, and the geometry is fairly simple. Therefore, leader–follower dynamics do not occur and agents tend to evacuate on their own. In the museum study, which has a complex geometry, the leader–follower dynamics occurs and agents follow a leader.

In general, regardless of the geometry between the two behavioral dynamics implemented (*altruism & leader–follower*), the less significant is the *leader–follower* dynamic. Instead, *altruism* is a key determinant in the accurate estimation of the evacuation time. Security guards have been also modeled in the analysis. When considering both behavioral dynamics models in ABM, regardless of the geometry considered, incorporating the human-behavior model affects the evacuation process.

For example the ABM of the Ursino Castle Museum without considering the human behavior estimate results in an evacuation time of 1 minute and 40 seconds (100 sec). If the human-behavior factor is included, the evacuation time is 4 minutes and 10 seconds (250 sec). The same trend can be observed in the Gare De Lyon metro station where the evacuation time without considering human behavior is approximately 1 minute (60 sec). If the human-behavior factor is included, the evacuation time increases to 2 minutes and 50 seconds. In both cases not including the human-behavior model underestimated the evacuation time.

Sensitivity analyses can also be performed by varying the number of agents, emergency exits, security guards, and the location of the blast. Further research can be analyzed by including other aspects of the human behavior such as kingship, age, panic, etc.

# 7 Restoration Fragility Functions of an Emergency Department

## 7.1 INTRODUCTION

Hospitals are critical facilities that affect the emergency response after a catastrophic event such as a strong earthquake. The non-functionality of an emergency department (ED) during a crisis might significantly impact healthcare services and affect the recovery process. A hospital's capability to remain accessible and ability to function at maximum capacity, providing its services to the community when they are needed most, can be evaluated using resilience indicators. A possible resilience indicator for healthcare facilities is the waiting time, which is the time the patient waits from the moment he/she walks in the ED until he/she first receives care from medical personnel [Cimellaro et al. 2010; Cimellaro et al. 2011]. A key role in the evaluation of the resilience indicator is dependent on the recovery time and the shape of the restoration curve because they are both uncertain quantities.

Presented herein is a procedure for building fragility curves of restoration processes called restoration fragility functions (RFFs), which can be adopted for resilience analysis. Restoration fragility functions take into account the uncertainties of the restoration process. In detail, RFFs are defined as the probability of exceedance of a given restoration process when a certain damage state occurs. To calculate the RFFs, it is necessary to define the functionality ($Q$) of the system considered and the recovery time. The ED of the Umberto I Mauriziano Hospital in Italy is presented as a case study. After building and calibrating a discrete event simulation (DES) model of the ED using real data collected on site, different case studies have been tested by modifying patient arrival rates and changing the number of available emergency rooms. In this research, two scenarios were considered: the ED with emergency plan applied and the ED in under normal conditions. The RFFs of both scenarios are compared.

## 7.2 STATE-OF-THE ART

### 7.2.1 Resilience of Hospital Facilities

A system is usually designed to behave in a certain way under normal circumstances. When disturbed from equilibrium by a disruptive event, the performance of the system will deviate from its design level. The resilience of the system is its ability to reduce both the magnitude and duration of the deviation as efficiently as possible to its usual targeted system performance levels [Proag 2014]. Bruneau et al. [2003] define resilience as the ability of a system to reduce the chances of a shock, to absorb a shock if it occurs (abrupt reduction of performance), and to recover quickly after it. By improving seismic resilience, we can minimize loss of life, injuries, economic losses, and any reduction in quality of life due to earthquakes.

They state that the strategy for measuring community resilience is to quantify the difference between the ability of a community's infrastructure to provide community services prior to the occurrence of an earthquake and the expected ability of that infrastructure to perform after an earthquake. More specifically, a resilient system should show reduced failure probabilities, reduced consequences from failures, and reduced time to recover to its normal level of performance [Bruneau et al. 2003]. However, according to Stevenson et al. [2014], few studies consider the ways social and economic connectivity of organizations of all kinds (private, public, for profit, and nonprofits) shape organizational and community recovery after a disaster. Also Yavari et al. [2010] stated that typically past earthquakes have focused on the performance of structural systems and not on the ability to provide necessary services after an earthquake.

Rose [2005] and Rose et al. [2007] distinguished post-disaster organizational capacities as static resilience: (1) the ability to absorb impacts and maintain function when shocked by making use of the resources available at a given time; and (2) dynamic resilience: the speed at which an entity or system recovers from a shock to attain a desired state [Park et al. 2011; Wein and Rose 2011]. It is logical to begin analyzing resilience by focusing on organizations whose functions are essential for community well-being in the aftermath of earthquake disasters. [Bruneau et al. 2003].

As part of the National Science Foundation's Earthquake Engineering Research Centers program, the Multidisciplinary Center for Earthquake Engineering Research (MCEER) has conducted a project to increase resilience by developing seismic evaluation and rehabilitation strategies for the post-disaster facilities and systems that society expects to be operational following an earthquake. The program was divided into several major research thrust areas, and one was the seismic retrofit of acute-care facilities. The research tasks focused on the integrated issues of the performance of hospital buildings, including both structural and nonstructural systems, and components and their functionality. (http://mceer.buffalo.edu/research/hospitals/default.asp) in the event of a disaster.

Hospitals have been recognized as critical buildings in hazardous events; hospitals must continue to function when an emergency occurs and must supply essential health services to the community at a time of disaster. The deaths and injuries after earthquakes are due to a variety of factors, but there are two main issues: the vulnerability of buildings and how quickly the emergency responders get to victims. Within a short time, hospitals have to provide care to a large number of injured whose lives are at risk, and they must have the ability to expand their services quickly beyond normal operating conditions to meet an increased demand for medical care.

Healthcare key factors are often classified into two categories: physical and social. The physical category includes structural and nonstructural parts, while the social category encompasses staff and administrative parts. A typical healthcare facility depends on the state of its building; the continuity of its utility supplies; availability and sufficiency of staff, equipment, and medical supplies; and easy accessibility for its daily operation. The failure of any of these components affects the continuity of medical care.

The lack of easy access to information and reports stating the experience of previous hospitals means that many hospital buildings are still very vulnerable to earthquakes [Achour et al. 2011]. Indeed keeping a hospital safe from natural or human-made disasters goes beyond protection of its physical structure. It requires preservation of its infrastructure as well as emergency staff trained to keep the facilities operational, collaboration with the network of health facilities, emergency plans, business continuity, and others specific abilities.

## 7.2.2 Functionality and Performance Levels of Infrastructures

Most research in this rapidly evolving field has focused on the evaluation side of resilience (i.e., defining and measuring resilience) [Mieler 2015]. Less attention has been paid to design side issues: for example, if a community wants to improve its resilience to earthquakes or other hazards, exactly what levels of performance are required from its buildings and lifelines? In fact, Mieler explains that performance levels established by modern building codes reflect choices that balance the desire to minimize initial construction costs with the need to ensure adequate levels of safety for the building's occupants [BSSC 2009]. From a resilience perspective, this performance objective, when aggregated across a community's entire building stock, can impede recovery after an earthquake. He developed a set of performance targets for important systems such as hospitals, schools, etc.

Bruneau et al. [2003] stated that the actual or potential performance of any system can be measured as a point in a multidimensional space of performance measures. Over time, performance can change, sometimes gradually, sometimes abruptly. Abrupt changes in performance occur in the case of disastrous events like a major earthquake. The performance of a system over time can be characterized as a path through the multidimensional space of performance measures. They defined a measure, which varies with time, for the quality of the infrastructure of a community. Specifically, performance can range from 0% to 100%, where 100% means no degradation in service and 0% means no service is available. If an earthquake occurs at time $t_0$, it could cause sufficient damage to the infrastructure such that the quality is immediately reduced. Restoration of the infrastructure is expected to occur over time until time $t_1$, when it is completely repaired.

Zhu and Frangopol [2014] defined the probability of a bridge experiencing different performance and functionality levels (e.g., one lane closed, all lanes closed). Inspired by the Federal Highway Administration [FHWA 2010] and ATC-13 [1999], they modeled the restoration process of bridge functionality by a normal cumulative distribution function corresponding to each bridge damage state considered. Recovery functions are highly dependent on their associated damage states. For example, a bridge categorized in a severe damage state may need more time to be restored to its full functionality compared to a bridge slightly damaged.

Padgett and DesRoches [2007] conducted analogous research to assess the probability of meeting various damage states expressed in terms of restoration of functionality and, subsequently, facilitate the refinement of component limit-state capacities for analytical fragility curve development of bridges. Padgett's research was limited to gathering information relating bridge damage to functionality by soliciting expert opinion. The FEMA-funded ATC-13 project recognized the need for this type of data in California and data on loss of function and restoration time for lifeline facilities was collected [ATC 1985]. One of the results of the damage-functionality survey was to quantify the probability of having a given restoration function or capacity over time. Padgett then developed fragility curves indicating the probability of the bridge being damaged beyond a given state for various levels of ground-motion intensity.

Yavari et al. [2010] introduced a methodology for anticipating the post-earthquake functionality of hospitals in a region. He defined performance levels for interacting systems (structural, nonstructural, lifeline, and personnel) in a hospital, and then probabilistically modeled them using damage data from past earthquakes. He proposed four potential functionality classes: fully functional (FF), functional (F), affected functionality (AF), and not functional (NF). Next, he proposed an overall measure of hospital functionality and related this overall performance to that of the interacting systems.

Cimellaro et al. [2010a] defined functionality of a hospital as the combination of a qualitative functionality related to the quality of service (QS) and a quantitative functionality that is related to losses in the healthy population. The qualitative functionality is related to the QS and can be defined using the waiting time (WT) spent by patients in the emergency room (ER) before receiving care. The WT is the main parameter used to evaluate the response of the hospital during normal and hazardous event operating conditions.

The quantitative functionality is considered when the maximum capacity of the hospital is reached. In this condition, the hospital is not able to guarantee a normal level of QS because the main goal now is to provide treatment to the most number of patients. In this case, the number of patients treated is a good indicator of functionality. Previously, Holmes and Burkett [2006] suggested classifying structural and nonstructural damage into different levels: None, Minor, Affecting Hospital Operations, and Temporary Closure. They used historical data on seismic vulnerability of hospitals to define performance levels, but they didn't take into account the role of personnel.

Formerly Nuti and Vanzi [1998] had synthetically defined the performance of hospital facilities as the time elapsed before a casualty is treated. The efficiency of the system is measured in terms of the mean distance for persons injured by the earthquake to be treated and by damage to the system. Their model aimed at minimizing the distance covered by each earthquake victim to reach a hospital. This calculation was then used to determine possible retrofitting strategies and to evaluate the effect of different post-earthquake emergency measures.

Jacques et al. [2013] presented a standardized methodology to analyze the impact of disasters on the functionality of healthcare systems. They developed a survey tool that collects field data on the performance of critical building systems and infrastructure, assesses the impact of system and infrastructure failure on the ability of hospitals to keep functioning, and provides data that can enhance existing and future tools to assess the performance of healthcare facilities. The survey assumed that hospital functionality is dependent on the physical infrastructure (e.g., continued functionality of electricity and water) and human infrastructure (e.g., reporting of healthcare providers as well as support staff). They created an event tree of hospital services based on lessons from seismic events (Bío-Bío, Baja California, and Christchurch earthquakes of 2010 and 2011); the event tree traces the failure or reduction of different hospital services examined.

### 7.2.3 Fragility Functions

Shinozuka et al. [2000] define fragility curves as functions that represent the probability that a given structure's response to various seismic excitations exceeds performance limit states. As such, fragility curves are a measure of performance in probabilistic terms. Fragility curves can be generated using actual damage data collected from existing structures that have already been subjected to earthquake loads. As these are scarce and rarely available for the areas of interest, the data is usually obtained from computational simulations [Koutsourelakis 2010].

In the current state-of-art, fragility functions describe the conditional probability that a structure, a nonstructural element or in general a system, will exceed a certain damage state, assuming a certain demand parameter (e.g., story drift, floor acceleration, etc.) or earthquake intensity level [e.g., peak ground acceleration (PGA), peak ground velocity (PGV) or spectral acceleration ($S_a$)] is reached. Usually, fragility functions take the form of lognormal cumulative distribution functions, having a median value $\mu$ and logarithmic standard deviation, $\beta$ [Porter et al. 2007]. According to Baker et al. [2011], the probability of collapse at a given $S_a$ level, $x$, can be estimated as the fraction of records for which collapse occurs at a level lower than $x$. A lognormal cumulative distribution function is often fit to this data, to provide a continuous estimate of the probability of collapse as a function of $S_a$.

Kafali and Grigoriu [2005] measured seismic performance using fragility surfaces instead of fragility curves. A fragility surface is the probability of system failure as a function of moment magnitude and site-to-source distance, consequences of system damage and failure, and system recovery time following seismic events. They used a Monte Carlo simulation and crossing theory of stochastic processes to calculate fragility surfaces for different limit states [Kafali and Grigoriu 2005].

Up until now, most of the studies on fragility curves focused on building and developing fragility functions using data from nonlinear dynamic structural analysis. There are a number of procedures for performing nonlinear dynamic structural analyses to collect the data for estimating a fragility function. One common approach is incremental dynamic analysis (IDA), where a suite of ground motions are repeatedly scaled in order to find the intensity measure (IM) level at which each ground motion causes collapse [Vamvatsikos and Cornell 2002; FEMA 2009]. A second common approach is multiple stripes analysis, where analysis is performed at a specified set of IM levels, each of which has a unique ground-motion set [Jalayer 2003].

In this work, the procedure to calculate the probability of exceedance is different from the usual method: the data taken into account to build RFFs are related to the performance of the ED during an extreme event when a certain damage state occurs. The procedure to estimate the parameters and the use of a lognormal cumulative distribution function to fit the data are similar to the methodology presented by Baker [2014]. The recovery functions are computed for three different damage states (DS), *no damage*, *moderate damage*, and *complete damage*. For each DS, a characteristic restoration curve is defined.

### 7.3    DEFINITION OF RESTORATION FRAGILITY FUNCTION

Presented herein is a procedure for building fragility curves of restoration processes that can be adopted for resilience analysis. The RFF is the probability of exceedance of a given restoration curve (*rf*) when a certain damage state (*DS*) occurs for a given earthquake intensity measure *I*. The general definition of RFF based on earthquake intensity *I* is given by:

$$RFF(i) = P\left(RF_j \geq rf_{DS1} \middle| DS = DS1, I = i\right) \quad (7.1)$$

where the $RF_j = j$th restoration function; $rf_j$ = restoration function associated to a given damage state *DS* (1,2,…n); *I* is an earthquake intensity measure, which can be represented by PGA, PGV, the pseudo velocity spectrum (PVS), the modified Mercalli scale (MMI), and *I*, which is a given earthquake intensity value. The main difference between RFFs and standard fragility functions is that the RFF is correlated to a given *DS*. In other words, the RFF is conditional on *DS* and *I*, while standard fragility curves are only conditional on the intensity measure *I*.

## 7.4 METHODOLOGY

The RFFs are evaluated using the experimental data of the restoration curves collected by the numerical analyses of the model considered. Different outputs can be considered, but in this specific case, the waiting time (*WT*) spent by patients in the ER before receiving care is considered as an indicator of functionality [Cimellaro et al. 2010]. In particular, the following relationship has been used to define its functionality *Q*:

$$Q = \frac{WT_0}{WT} \quad (7.2)$$

where $WT_0$ is the acceptable waiting time in regular conditions when the hospital is not affected by a catastrophic event, and *WT* is the waiting time collected during the simulation process. When the *WT* is less or equal to $WT_0$, the value of *Q* is equal to 1, meaning that the hospital's functionality is at its maximum.

Different restoration functions (*rf*) associated at different damage states have been chosen. Then, for each simulation, the probability of exceedance of a given restoration curve (*rf*) has been calculated. The frequency of exceedance at a given instant is defined as

$$f = N/N_{tot} \quad (7.3)$$

where *N* is the number of times when the restoration curves exceed the restoration curve associated at a given damage state; $N_{tot}$ is the number of simulations.

Finally, the probability of exceedance of a given restoration state is calculated by

$$P_{ex} = \frac{\sum f_i}{T} \quad (7.4)$$

where $\sum f_i$ is the sum of the frequencies at each time instant, and *T* is the length of the simulation (e.g., *T* =13 days in the case study).

Finally, different methods to fit fragility curves are compared such as:
- -MLE method: maximum likelihood method
- -SSE method: sum of squared errors

## 7.5 CASE STUDY: THE MAURIZIANO HOSPITAL

The Umberto I Mauriziano Hospital shown in Figure 7.1 is used as a case study to show the applicability of the methodology. Located in Turin, Italy, the hospital stands out in the landscape of healthcare facilities in the Piedmonte region as a medical facility; it provides both basic care and several areas of specialization. The hospital is located in the southeast part of the city, almost 3 km from downtown. It was built in 1881 and was bombed several times during World War II. This explains why several buildings have been rebuilt or added. Presently the hospital includes 17 units corresponding to different departments, covering an overall surface of 52827 m$_2$. Only the ED (building 17) has been modeled.

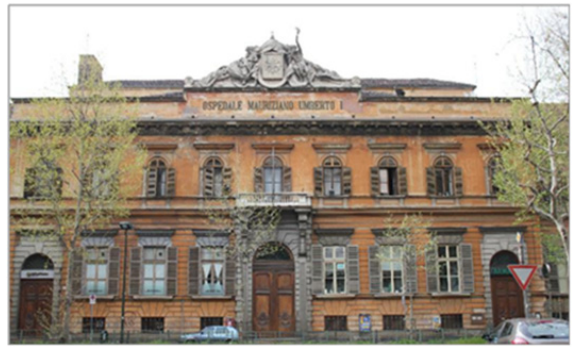

Figure 7.1     Umberto I, Mauriziano hospital in Turin.

### 7.5.1 The Emergency Department

A hospital's ED is the most affected area of a medical facility in the event of a disaster. Emergency departments play a pivotal role in the delivery of acute ambulatory and inpatient care, providing immediate assistance request during 24-hour period post-event [Morganti et al. 2013]. Experience has shown that the effectiveness of rescue operations in the first 24 to 48 hours, and especially the capacity of the medical system, can considerably reduce the number of deaths [ATC 2002].

The ED is a complex and dynamic setting, with multiple interactions between patients, doctors, nurses, technicians, and different departments, and multiple patient paths to categorize [Anders et al. 2006]. Moreover, hospital EDs are also critical pressure points during disasters since their non-functionality might significantly impact the healthcare services and affect the recovery process.

The ED consists of an entrance area in which a procedure called "triage" is carried out, and four macro areas corresponding to the four different color codes that represent the severity of injury. In particular, these four color codes are red, yellow, green, and white. Red codes (emergency) identify patients with compromised vital functions, already altered or unstable whose lives are at risk. Yellow codes (urgency) are patients who are not in immediate danger of life but present a partial impairment of vital functions. Green codes (minor urgency) are patients that are not in critical condition; their lives are not at risk, and their injuries do not affect vital functions. White codes (no urgency) include all patients who do not have neither serious nor urgent problems and who do not really need to be in the ED, but their need of care could be provided by a general doctor. This research used only those patients designated with yellow codes to develop RFFs.

Considering this classification, the ED is normally divided in four main areas but when the emergency response plan (ERP) is applied, the number of areas is reduced to three. This is because under emergency conditions white codes are sent to another facility outside the ED. Furthermore, under normal operating conditions, yellow- and green-coded patients share the same area and, consequently, treatment rooms. Under emergency conditions, the red-code area is located immediately in front of the ambulance entrance and contains two rooms where patients receive preliminary treatment. Parallel to this area is the yellow-coded area, which is composed of three treatment rooms. Separate from this zone, the green-coded area is situated perpendicular to yellow- and red-coded areas, and includes two treatment rooms. Each area is provided with waiting rooms where patients can stay before being treated. There are also a number of recovery rooms located inside the ED where patients can stay before being discharged or transported to another part of the hospital for extended recovery time; see Figure 7.2. To assess the impacts on the length of stay in the ED, a computer simulation model was developed that varies the number of the available treatment rooms relative to the seismic input.

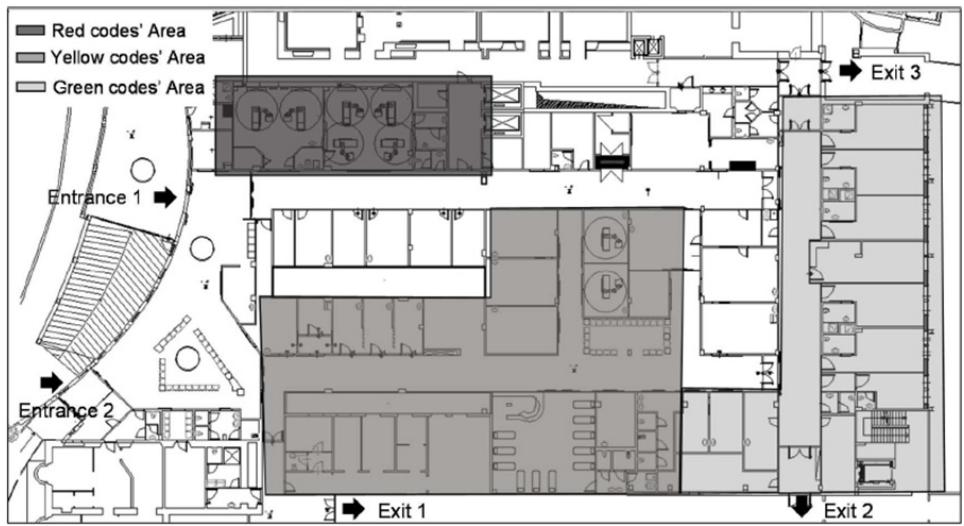

Figure 7.2    Emergency Department color-coded areas.

### 7.5.2 Discrete Event Simulation Model

A discrete event simulation model (DES) of the ED was developed (Figure 7.3) using ProModel version 7.0, downloaded on February 15, 2014. In particular, different scenarios have been analyzed considering some structural damage to specific parts of the building due to a catastrophic event and variable patient arrival rates depending on the seismic intensity. Discrete event simulation models represent useful tools to test emergency-response plans under a rapid increase in the volume of incoming patients. Using discrete-event Monte Carlo computer simulations, hospital administrators can model different scenarios of the hospital to see how they compare to the desired performance [Morales 2011]. Moreover, DES model allows investigation and planning of the use of hospital resources [Šteins 2010].

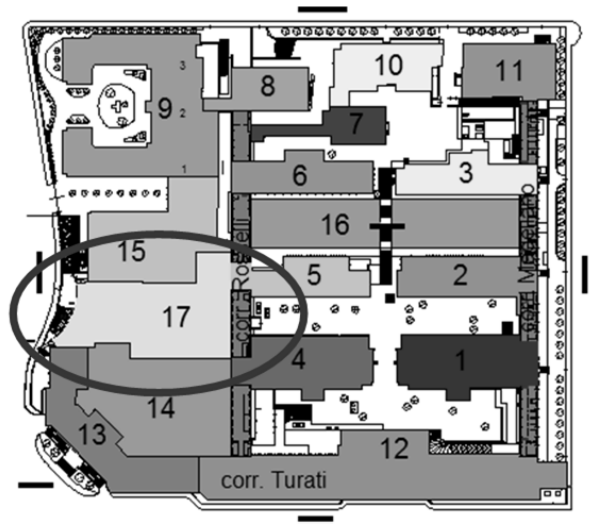

Figure 7.3    Hospital's units: Emergency Department building.

Healthcare systems show multiple interactions between patients, doctors, nurses, technicians, different departments, and circulation patterns. Therefore, difficulty often exists in evaluating how each of these components can affect the whole system while taking into account all the multiple interactions. A healthcare system could be represented as a chronological sequence of events that occur at a definite instant in time and mark changes in the system. In this way, the end of each event marks the start of the next event. Testing the ERP requires identifying those factors that represent the quality of healthcare services and can best describe ED performance during a dramatic event. Different parameters can be used to evaluate the efficacy of emergency plans and, among these parameters, the most representative one is patient wait times (PWTs).

Several steps were taken to build a valid analytical framework. First, a generic simulation model of the ERP was developed taking into account the hospital resources, emergency rooms, circulation patterns, and patient codes. One hundred simulations were carried out using as input data the seismic arrival cycle.

### 7.5.3 Comparison between Emergency Response Plan and Normal Operating Conditions

An emergency plan consists of a number of procedures designed to respond to those situations where the standard operation procedures are not be able to provide essential health services. It was developed to assure adequate medical resources during an emergency for the maintenance of patient care and equipment, availability of treatment supplies, and appropriate interaction with others critical infrastructures. Generally, an ERP is activated when the number of ill or injured exceeds the normal capacity of the ED or the normal operations of multiple departments to provide the quality of care required. Its very nature means that it is impossible to test the effectiveness of an emergency plan before a disaster occurs.

This research analyzed two case scenarios. First, 100 simulations of the ED under normal operating conditions were carried out; next, another 100 simulations were conducted with the ERP applied. The output of the model is a record of PWTs that can be used to develop functionality curves for each case scenario. The ERP is considered effective if the PWTs obtained when the emergency plan is applied is significantly lower than the PWTs obtained under emergency conditions when the emergency plan is not active.

### 7.5.4 Input Data of the Model and Assumptions

The data input of the model are patient arrival rates under normal operating conditions; this data has been extracted by the hospital's register statistics. A 13-day simulation was conducted that considered the occurrence of a seismic event after two days of simulation. In order to simplify the model, it was assumed that the ERP was applied also in the first 2 days of simulation, even though the minimum required conditions for application were not satisfied. Considering that there are no experimental data available, the probability values entered for the construction of the model were obtained from interviews with the hospital's medical staff. If the response of the medical staff were found not satisfactory, a probability of 50% was considered.

The ED can be characterized by the number of operating rooms, the number of resources (doctors, nurses, and healthcare operators), and the procedures available inside the different rooms, as well as the circulation patterns and the patient arrival rates. Patient arrivals in the ED vary from hour to hour and, in order to determine the patient arrival distributions, an arrival cycle was defined using data from the hospital's register.

In order to take into account the increase of the patients flow due to a catastrophic event and the consequent crowding of the ED, a seismic input was considered. The data collected from a California hospital during 1994 Northridge, California, earthquake were used in the model to simulate the seismic event. Northridge's arrival rate was selected because it is the only documented event to date [Stratton et al. 1996; Peek-Asa et al. 1998; and McArthur et al. 2000] where patient arrival rates have been collected. The pattern of the Northridge patient arrival rates is given in the work of Cimellaro et al. [2011]. Then the patient arrival rates was scaled to the seismic hazard in Turin using a procedure based on the MMI scale. An earthquake with a return period of 2500 years was considered, assuming a nominal lifespan for a building of strategic importance that is 100 years old according to the Italian seismic standards [NTC-08 2008].

The purpose of the research was to build RFFs of the ED. Therefore, the data collected was related to increasing seismic intensities. To study the effect of the seismic arrival rate on PWTs, an amplified seismic input was considered. The seismic arrival rate was amplified in order to analyze the sensibility of the ED towards the amplitude of the earthquake. Multiplicative scale factors ranging from 1.1 to 1.6 were used to amplify the input data. The factors ($\alpha$) used for the analysis are shown in Figure 7.4.

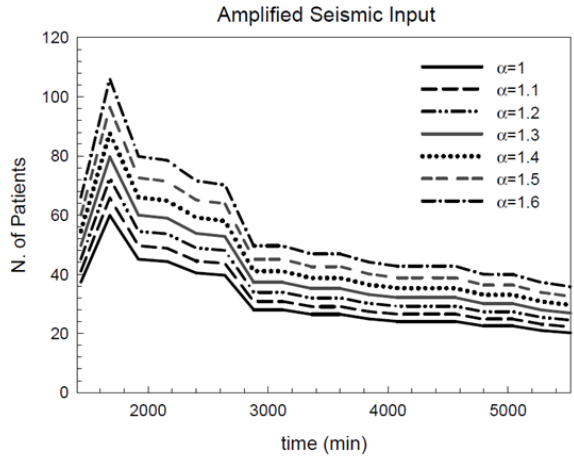

**Figure 7.4** Amplified seismic input for different scale factors $\alpha$.

### 7.5.5 Simulation and Output of the Model

The simulation lasted twelve days and assumed that the first two days proceeded under normal operating conditions. These two days are needed in order to make the system stable and avoid any influence by initial conditions. A three-day period followed the normal operating conditions, which assumed that the ED was operating under emergency conditions determined using the scaled arrival rate calculated in the "seismic input" paragraph. After the seismic input, another eight-day period was run assuming normal operating conditions; this is because the system needs some time to return to its previous steady state before the earthquake occurred. The twelve-day-long simulation was run 100 times for each different scenario. Table 7.2 presents each different case scenario and the number of simulations that were conducted. Three different damage states (*DS*) were considered:

- DS=Fully operational/No Damage ($n = 0$);
- DS=Moderate Damage ($n = 1$);
- DS=Severe Damage ($n = 2$);

where $n$ is the number of treatment rooms not functioning because they were damaged by the earthquake. The outputs of the model are the waiting time of each patient and the time instant in which the patient walks into the ED.

Table 7.1    Case scenarios and the corresponding simulation number.

| ED with Emergency Plan | | | ED in normal operating conditions | | |
| --- | --- | --- | --- | --- | --- |
| | α | number of simulations | | α | number of simulations |
| Fully operational/No Damage (n=0) | 1 | 100 simulations | Fully operational/No Damage (n=0) | 1 | 100 simulations |
| | 1.1 | 100 simulations | | 1.1 | 100 simulations |
| | 1.2 | 100 simulations | | 1.2 | 100 simulations |
| | 1.3 | 100 simulations | | 1.3 | 100 simulations |
| | 1.4 | 100 simulations | | 1.4 | 100 simulations |
| | 1.5 | 100 simulations | | 1.5 | 100 simulations |
| | 1.6 | 100 simulations | | 1.6 | 100 simulations |
| Moderate Damage (n=1) | α | number of simulations | Moderate Damage (n=1) | α | number of simulations |
| | 1 | 100 simulations | | 1 | 100 simulations |
| | 1.1 | 100 simulations | | 1.1 | 100 simulations |
| | 1.2 | 100 simulations | | 1.2 | 100 simulations |
| | 1.3 | 100 simulations | | 1.3 | 100 simulations |
| | 1.4 | 100 simulations | | 1.4 | 100 simulations |
| | 1.5 | 100 simulations | | 1.5 | 100 simulations |
| | 1.6 | 100 simulations | | 1.6 | 100 simulations |
| Severe Damage (n=2) | α | number of simulations | Severe Damage (n=2) | α | number of simulations |
| | 1 | 100 simulations | | 1 | 100 simulations |
| | 1.1 | 100 simulations | | 1.1 | 100 simulations |
| | 1.2 | 100 simulations | | 1.2 | 100 simulations |
| | 1.3 | 100 simulations | | 1.3 | 100 simulations |
| | 1.4 | 100 simulations | | 1.4 | 100 simulations |
| | 1.5 | 100 simulations | | 1.5 | 100 simulations |
| | 1.6 | 100 simulations | | 1.6 | 100 simulations |

## 7.6    HOSPITAL PERFORMANCE AND RESTORATION FUNCTIONS

Generally, the performance of a hospital under seismic hazard is quantified considering all its possible damage states. The performance of the ED is quantified within this work by mapping the current damage state to a value between 0 and 1.0. Assuming a certain damage state occurs in the hospital, different restoration functions (*rfs*) can be applied to the damaged structure to restore its functionality. However, the restoration functions (*rfs*) of the ED are highly dependent on their associated damage states. For example, an ED categorized in a severe damage state may need more time to be restored to its full functionality compared to an ED that is only slightly damaged; therefore some *rf*s have a higher probability of occurring with respect to others.

## 7.7 METHODOLOGY

### 7.7.1 Calculation and Definition of Functionality

The RFFs were evaluated using the experimental data of the restoration curves collected by the numerical analyses of the model considered. Different output can be considered, but in this specific case, the waiting time (*WT*) spent by patients in the ED before receiving care is considered as an indicator of functionality [Cimellaro et al. 2010]. In particular, the following relationship has been used to define its functionality *Q*:

$$Q = \frac{WT_0}{WT} \tag{7.5}$$

where $WT_0$ is the acceptable waiting time under regular conditions when the hospital is not affected by a catastrophic event, and *WT* is the waiting time collected during the simulation process. When the *WT* is less or equal to $WT_0$, the value of *Q* is equal to 1, meaning that the hospital's functionality is at its maximum. Other procedures to define functionality can be used.

### 7.7.2 Smoothing Procedure

The functionality values obtained directly from the model's outputs have been plotted with the program SigmaPlot 12.3. The data obtained from the model were subjected to a smoothing procedure via SigmaPlot 12.3. Smoothing is used to elicit trends from noisy data. This procedure is used when it is necessary to smooth data to remove high-frequency component or to resample observations. In this case, a local smoothing technique that computes the median of the values at neighboring points with a bandwidth method of the nearest neighbors was adopted. Figure 7.5 shows the functionality over time of 100 simulations considering a damage state of no damage (*DS*=0) subjected to an earthquake of magnitude VIII-IX.

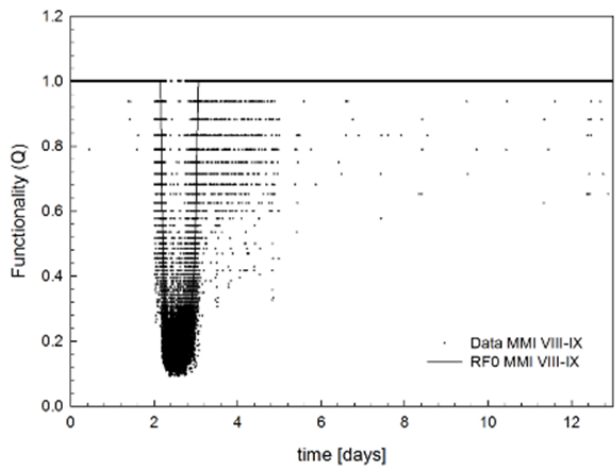

**Figure 7.5** Restoration function (*DS* = 0) assuming an earthquake of magnitude VIII-IX obtained by a smoothing procedure of the output data from the model.

### 7.7.3 Restoration Functions and Damage State Definitions

Different restoration functions (*rfs*) associated with different damage states were calculated. Three different damage states (DS) were considered:

- DS=Fully operational/No Damage ($n = 0$)
- DS=Moderate Damage ($n = 1$)
- DS=Severe Damage ($n = 2$)

where *n* is the number of treatment rooms not functioning because they were damaged by the earthquake.

Function $Q(t)$ has been plotted for each damage state at different earthquake intensity measures. Figure 7.6 shows the *rf* for the three different damage states, assuming the occurrence of an earthquake of magnitude VIII–IX in the MMI. Two case studies were considered: the ED with and without an ERP.

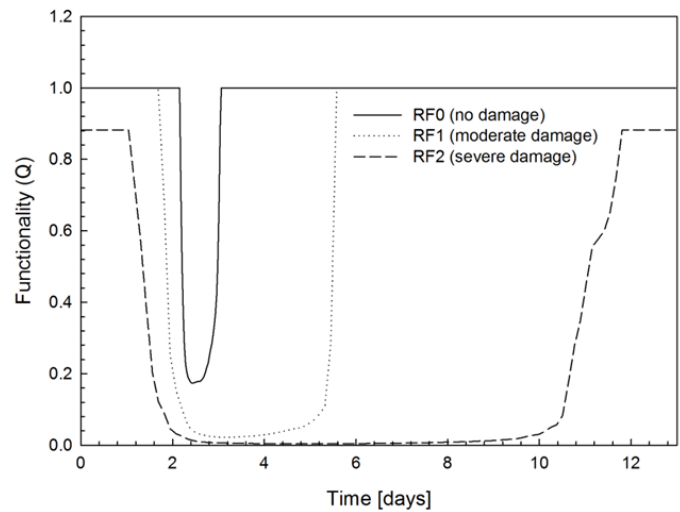

**Figure 7.6**   Restoration functions assuming earthquake of magnitude VIII–IX (the Emergency Department with the emergency plan applied).

### 7.7.4 Reference Restoration Functions

To calculate the probability of exceedance of a given *rf*, it is necessary for comparison's sake to define reference *rf*s. Three *rf*s associated with specific damage states have been chosen to calculate the fragility restoration curve. The *rf*s chosen in this study refer to the functionality curve assuming no damage (RF0), moderate damage (RF1), and complete damage (RF2) for an earthquake of magnitude VI in the MMI. As shown in Figure 7.7, RF0 has a restoration time of one day, while RF1 and RF2 have restoration times of two days and six days, respectively. The restoration time $t_r$ specifies how long the ED takes to recover from a disaster.

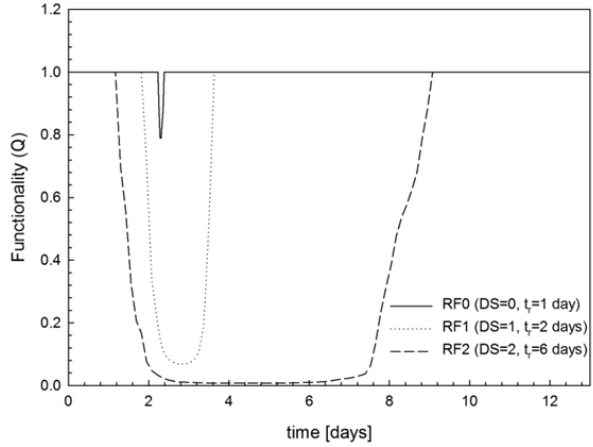

Figure 7.7    Restoration functions assuming earthquake of magnitude VI.

### 7.7.5 Frequency of Exceedance

Three different *rf*s of each damage state were compared for increasing values of seismic intensity measures with the ERP applied or under normal operating conditions,. For each simulation, the probability of exceedance of a given restoration curve was calculated. The frequency of exceedance at a given instant is defined as

$$f = N/N_{tot} \tag{7.6}$$

where $N$ is the number of times when the restoration curve exceeded the restoration curve associated at a given damage state; $N_{tot}$ is the number of simulations. For each case scenario, damage state, and seismic intensity measure, the frequency over time was plotted. In Figure 7.8, the frequency of exceedance of the RRF RF0, for the damage state *no damage*, was plotted, and numerical results are presented for additional case studies.

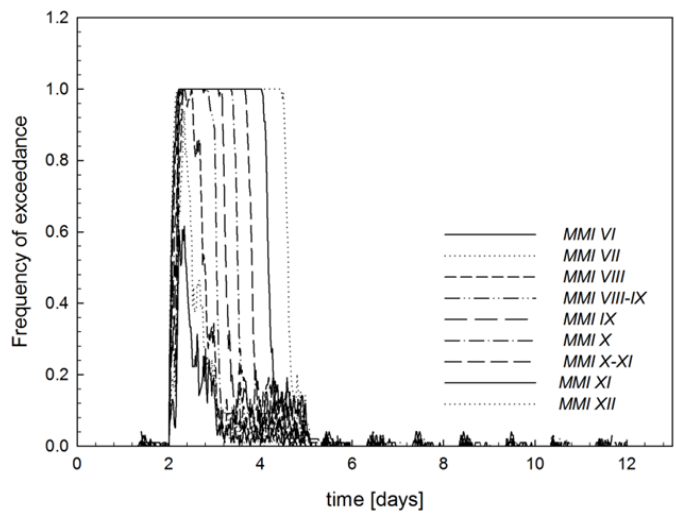

Figure 7.8    Frequency of exceedance of reference restoration function RF0, no damage (*n* = 0).

### 7.7.6 Data Interpolation

One hundred simulations for each case study were conducted. The data obtained from the model were the time *t* at which a single patient arrived at the ED associated with PWTs before receiving care. In order to calculate the frequency of exceedance defined above, it was necessary to have at a given instant *t*, 100 values of functionality *Q(t)* that could be compared to the value of the *rf* at the same given instant *t*. A data interpolation function was used. According to Meijering [2002], interpolation is a procedure where an approximating function is constructed in such a way as to agree perfectly with the usually unknown original function at the given measurement points. It is a method of constructing new data points within the range of a discrete set of known data points. The interpolation function *interp1* of Matlab-r 2011b was used to interpolate the data.

$$Q_{new} = \mathrm{interp1}(t, Q_{old}, t_1, 'linear') \tag{7.7}$$

where $Q_{new}$ is the new value of functionality that has been interpolated; *t* is the time instant associated with the value of functionality; and $t_1$ is a new vector of time instants at which $Q_{new}$ is evaluated.

### 7.7.7 Probability of Exceedance of a Given Restoration State

The probability of exceedance of a given restoration state has been calculated by

$$P_{ex} = \frac{\sum f_i}{T} \tag{7.8}$$

where $\sum f_i$ is the sum of the frequencies at each time instant, and *T* is the length of the simulation (e.g., *T* = 13 days in the case study). The probability of exceedance has been calculated for increasing seismic intensities. Two different methods to fit fragility curves are compared:

- -MLE method: maximum likelihood method
- -SSE method: sum of squared errors

## 7.8 NUMERICAL RESULTS

As outputs of the model, the PWTs of the ED when the ERP was activated have been collected for different scenarios. Three different damage states (DS) were considered:

- DS=Fully operational/No Damage (*n* = 0);
- DS=Moderate Damage (*n* = 1);
- DS=Severe Damage (*n* = 2);

where *n* is the number of treatment rooms not functioning because they were damaged by the earthquake.

For each DS, several simulations were conducted by changing the intensity of the seismic event using the methodology described above. The intensity was increased by means of scale factors that multiplied the patient arrival rates. Three different restoration functions (RFs) were chosen as comparison. The functionality $Q$ of the ED was evaluated for increasing seismic intensities based on the MMI scale. In this case study, the MMI scale was adopted, but other parameters such as PGA, PGV, or $S_a$ can be used as well. Each graph shows different damage states:

- Emergency plan fully operational with $n = 0$, where $n$ is the number of treatment rooms not available because they have been damaged by the earthquake (Figure 7.9)
- ERP affected by moderate damage ($n = 1$) (Figure 7.10)
- ERP affected by severe damage ($n = 2$) (Figure 7.11)

As shown in the graphs, the functionality is reduced, and the recovery time increased when two treatment rooms are not operative. The functionality is also dependent on the seismic intensity. As the seismic intensity increased, the restoration curves reflect that the recovery time to return to their initial functionality increased. As shown in Figure 7.11 for higher seismic intensities, the functionality at the end of the simulation doesn't reach the ideal value, showing that the ED has not totally recovered from the seismic event.

Three reference *rfs* associated with specific damage states were chosen to calculate the fragility restoration curve; in this study, they refer to the functionality curve assuming no damage (RF0), moderate damage (RF1), and complete damage (RF2). As shown in Figure 7.12, RF0 has a restoration time of 1 day, while RF1 and RF2 have restoration times of 2 days and 6 days respectively. The restoration time $t_r$ specifies how long the ED takes to recover from a disaster. The frequency of exceedance of the RFFs has been calculated for each case study. The frequency is plotted for each time instant. The results are shown in Figure 7.13 to Figure 7.21.

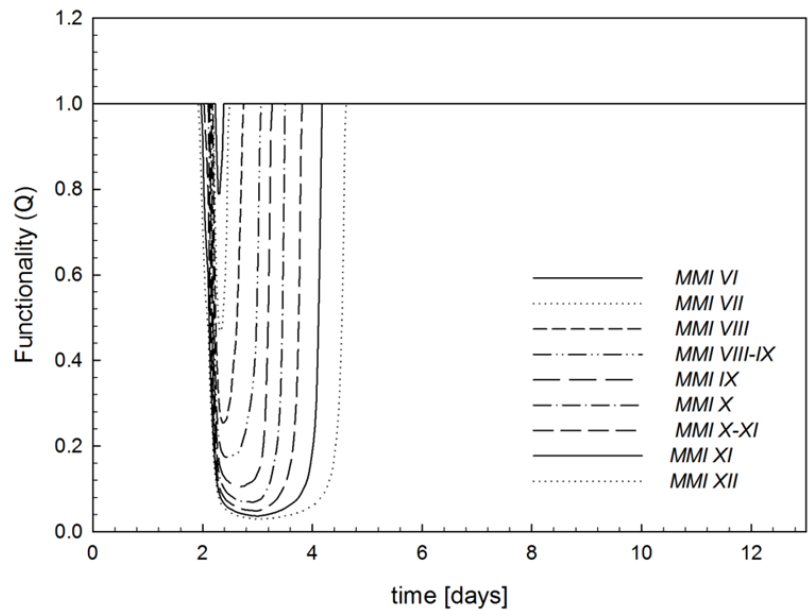

**Figure 7.9** Functionality curves as a function of seismic intensity, no damage ($n = 0$).

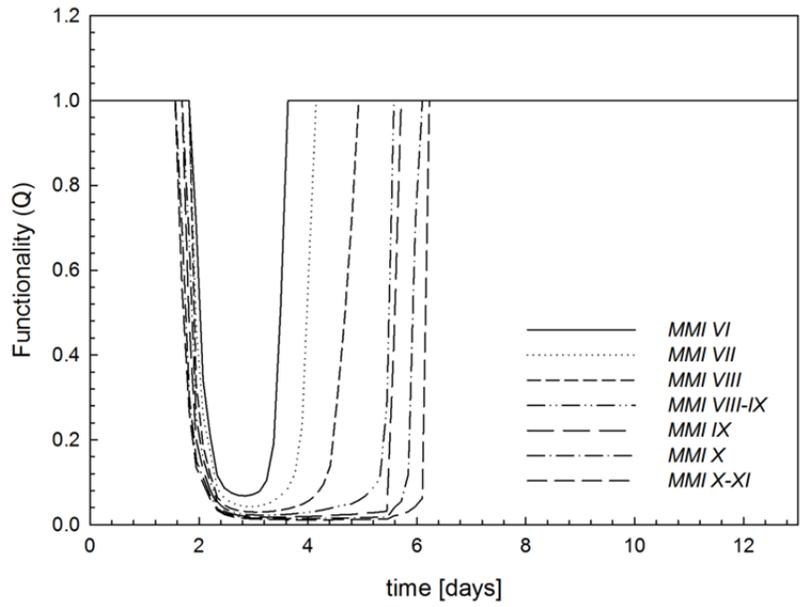

**Figure 7.10** Functionality curves as a function of seismic intensity, moderate damage ($n = 1$).

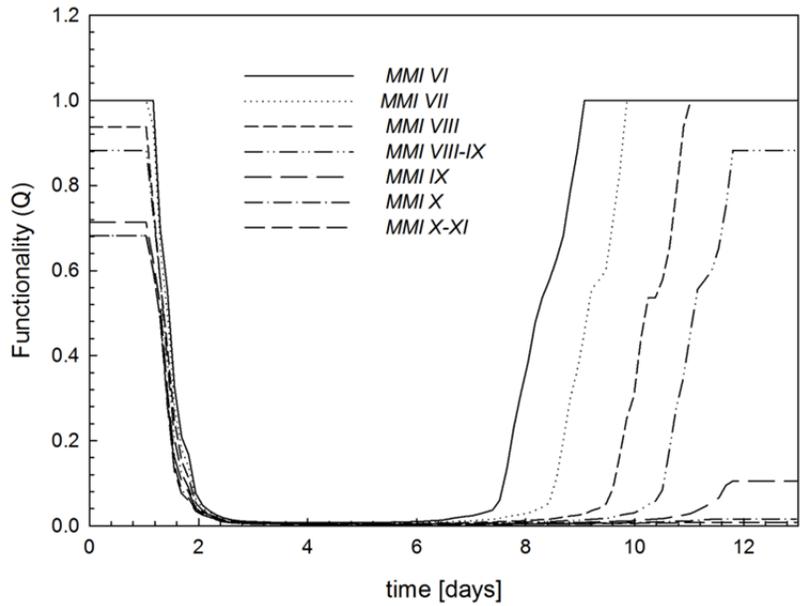

**Figure 7.11** Functionality curves as a function of seismic intensity, severe damage ($n = 2$).

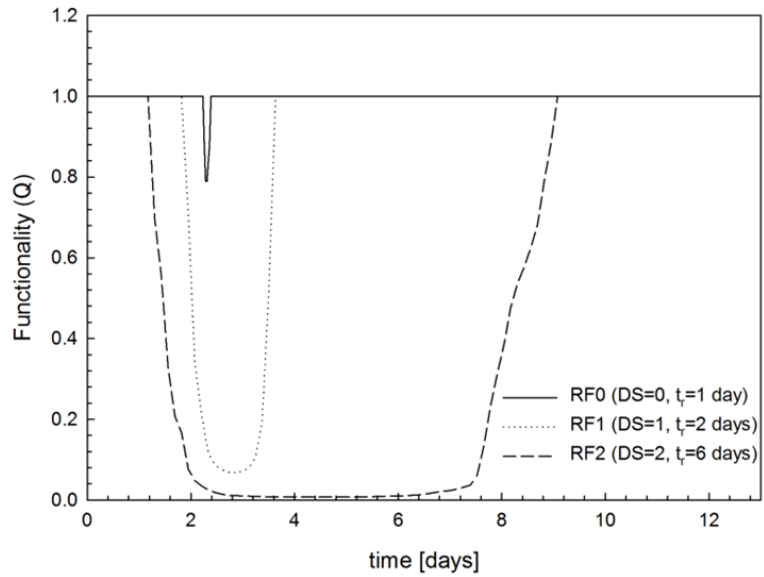

Figure 7.12    Restoration functions assuming earthquake of magnitude VI.

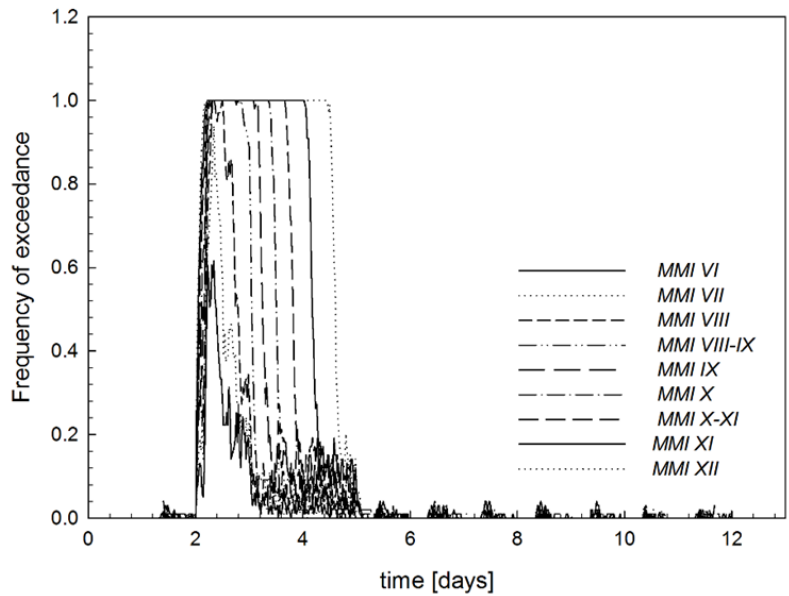

Figure 7.13    Frequency of exceedance of reference restoration function RF0, no damage ($n = 0$).

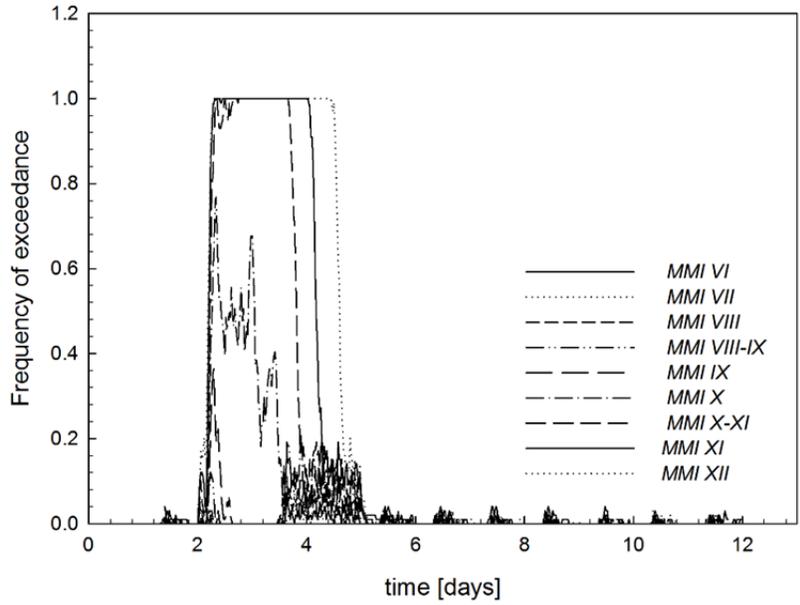

**Figure 7.14**   Frequency of exceedance of reference restoration function RF1, no damage ($n = 0$).

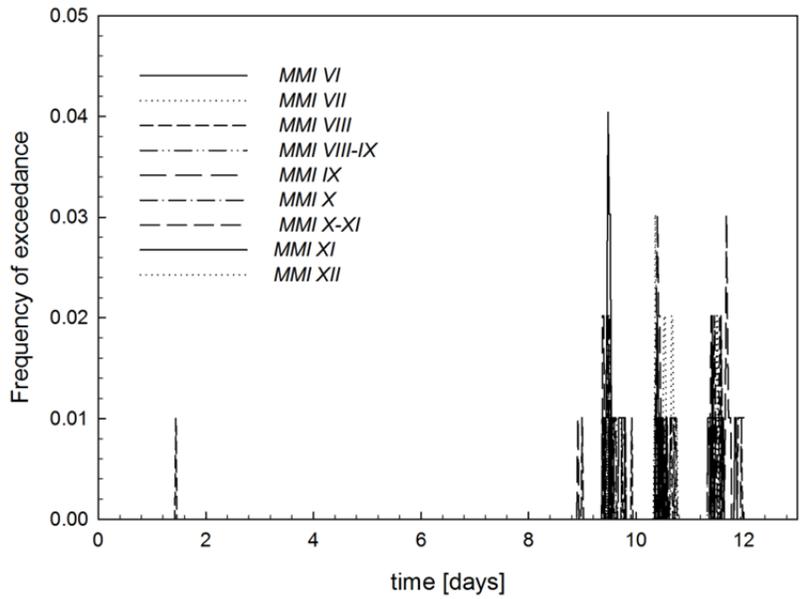

**Figure 7.15**   Frequency of exceedance of reference restoration function RF2, no damage ($n = 0$).

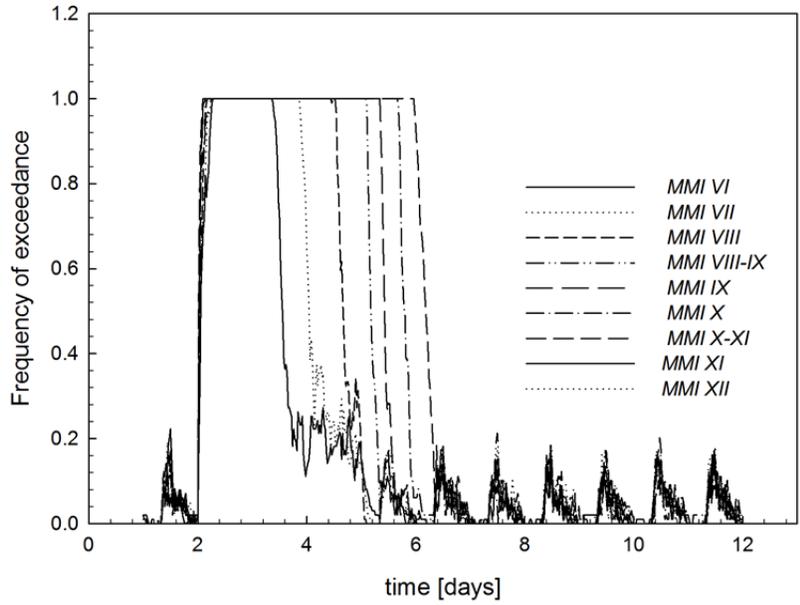

**Figure 7.16** Frequency of exceedance of reference restoration function RF0, moderate damage (*n* = 1).

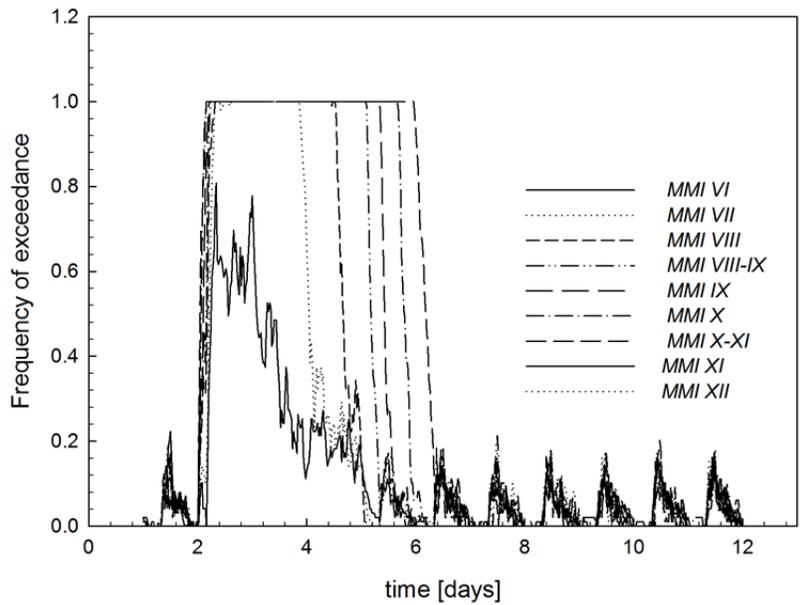

**Figure 7.17** Frequency of exceedance of reference restoration function RF1, moderate damage (*n* = 1).

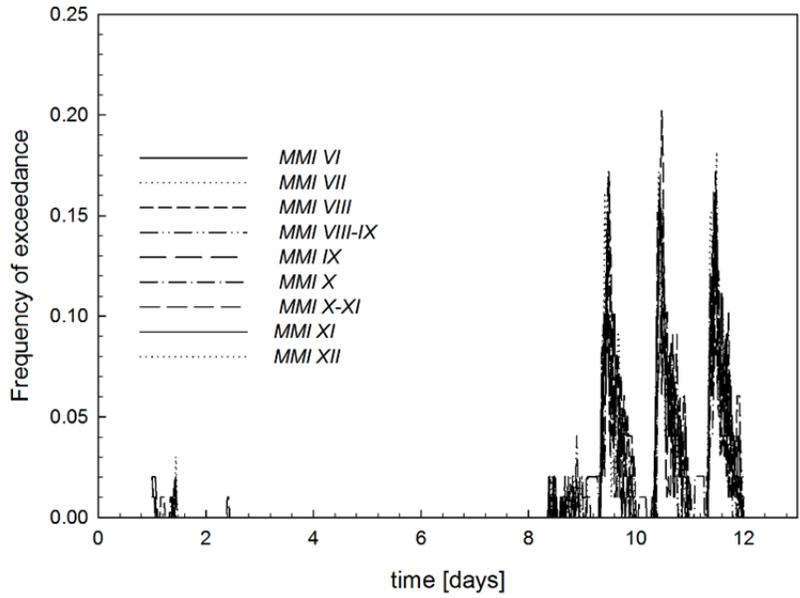

**Figure 7.18** Frequency of exceedance of reference restoration function RF2, moderate damage ($n = 1$).

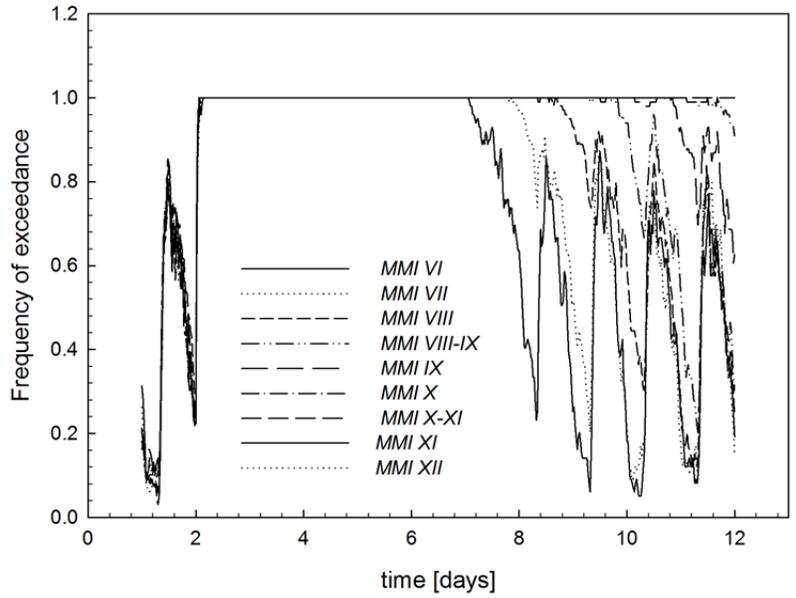

**Figure 7.19** Frequency of exceedance of reference restoration function RF0, severe damage ($n = 2$).

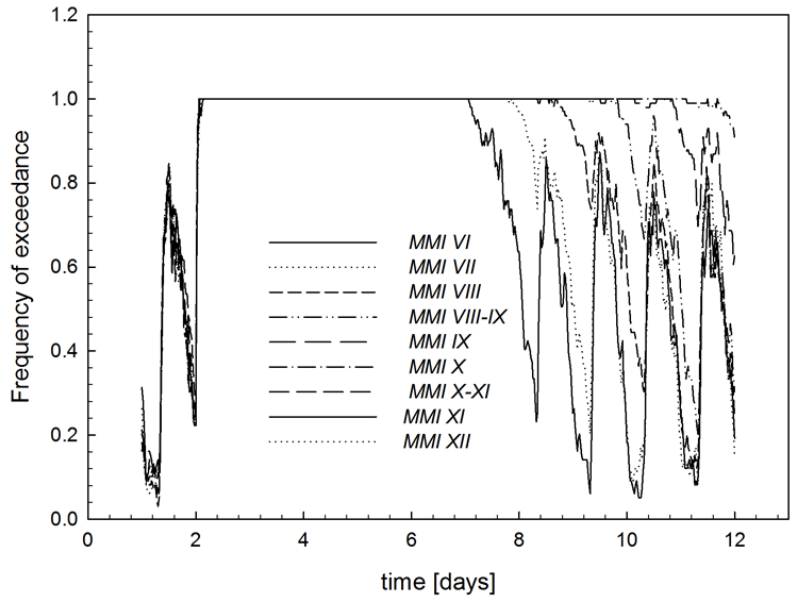

**Figure 7.20** Frequency of exceedance of reference restoration function RF1, severe damage ($n = 2$).

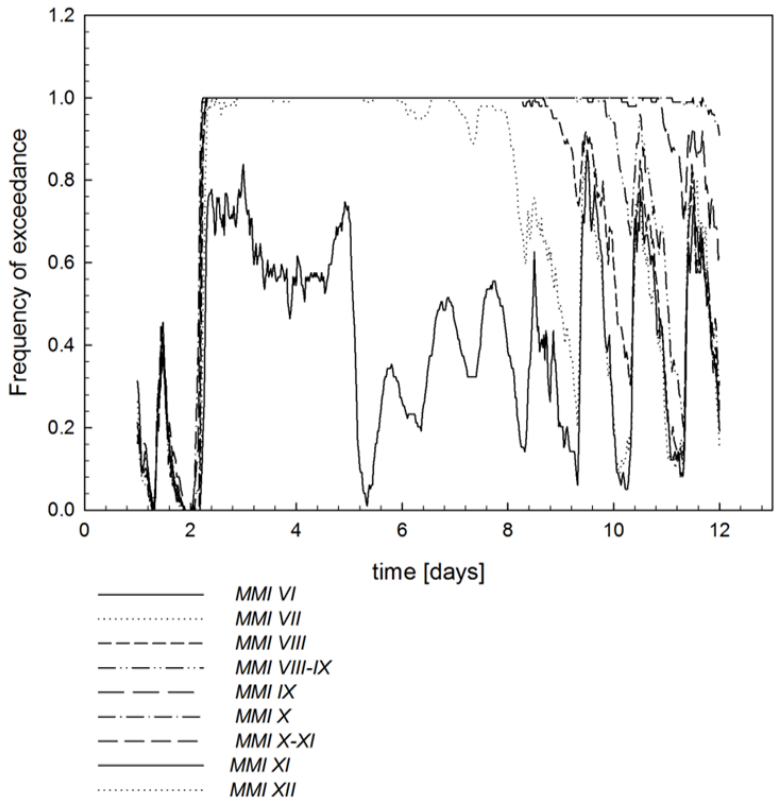

**Figure 7.21** Frequency of exceedance of reference restoration function RF2, severe damage ($n = 2$).

The RFFs for each damage-state scenario were calculated. The RFF is the probability that a given restoration function *rf* (Figure 7.12) is reached when a certain damage state occurs for a given earthquake intensity measure *I*. In Figure 7.22–Figure 7.27 the probability of restoration is plotted. The lognormal cumulative distribution function is used to fit the data to provide a continuous estimate of the probability of restoration as a function of MMI. Described in Baker [2013], two different methods to fit fragility curves are compared:

- -MLE method: maximum likelihood method
- -SSE method: sum of squared errors

As shown in Figure 7.22–Figure 7.27, the two fitting methods produced similar results. In Figure 7.22 and Figure 7.23, the probability of exceedance of the *rf* RF0 increases with the increment of the MMI. For higher MMI, the probability of exceedance of RF1 reaches the probability of exceedance of RF0. In Figure 7.22 and Figure 7.23, the same behavior can be observed. Note that the probabilities of exceedance of RF0 and RF1 overlap. The RRF related to RF2 increases considerably with respect to the previous damage states. Figure 7.28–Figure 7.30 show the RFFs related to the ED without the ERP.

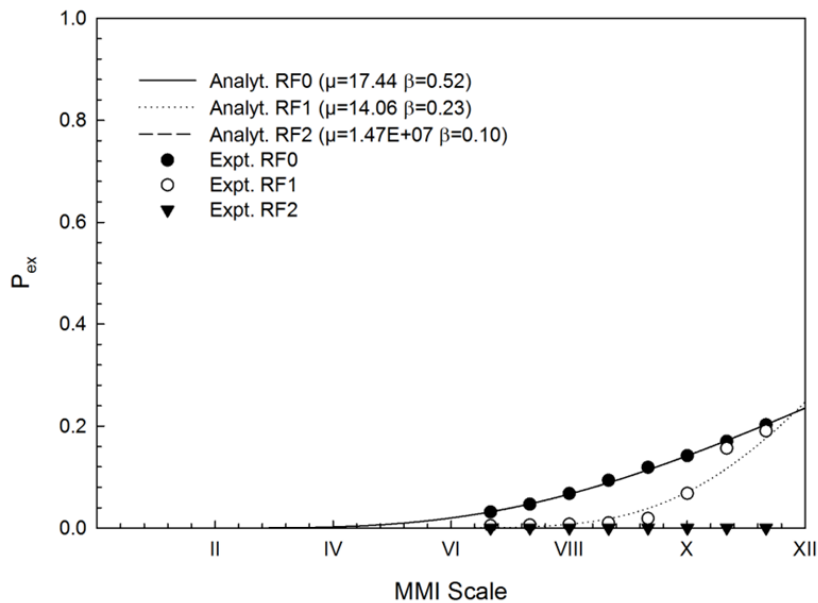

**Figure 7.22** Reference restoration function given DS = 0(no damage) using MLE method (the Emergency Department with emergency response plan applied).

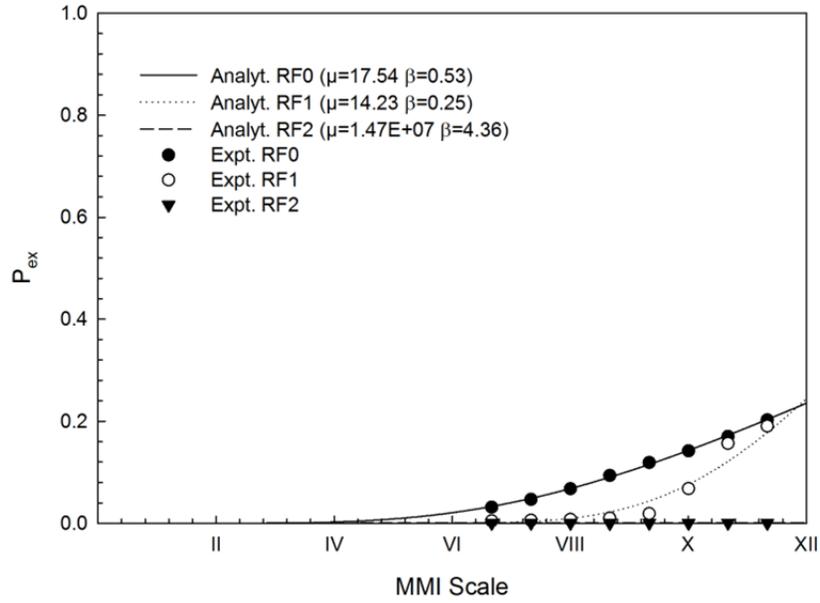

**Figure 7.23**   Reference restoration function given a DS=0(no damage) using SSE method (the Emergency Department with emergency response plan applied).

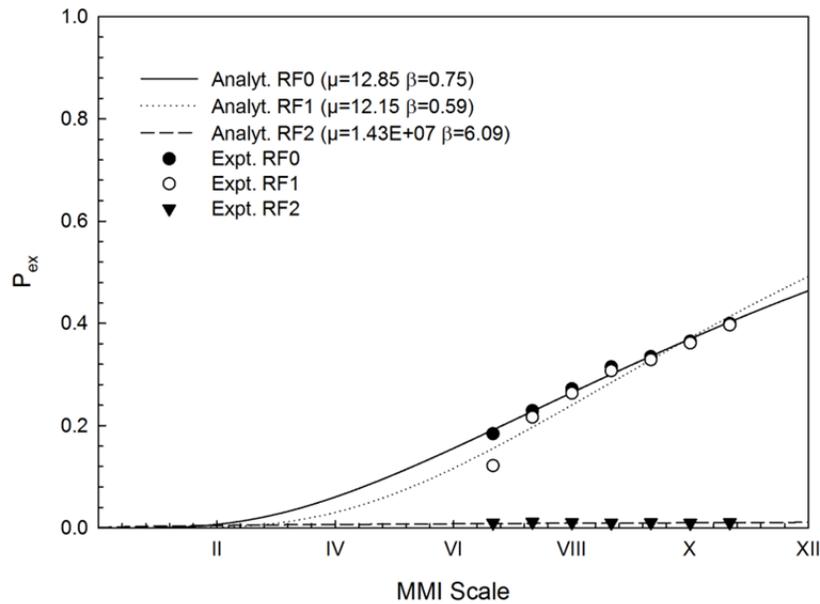

**Figure 7.24**   Reference restoration function given a DS = 1(moderate damage) using MLE method the (Emergency Department with emergency response plan applied).

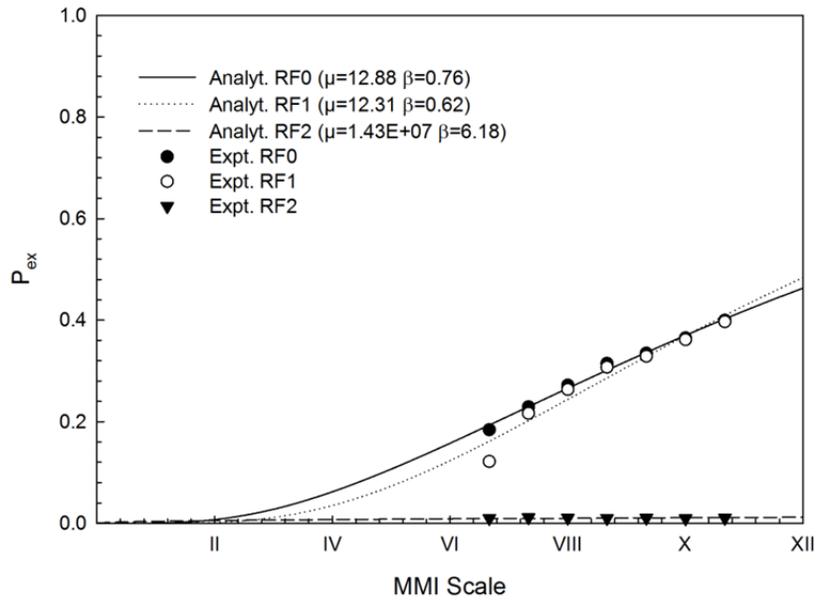

**Figure 7.25**  Reference restoration function given a DS = 1(moderate damage) using SSE method (the Emergency Department with emergency response plan applied).

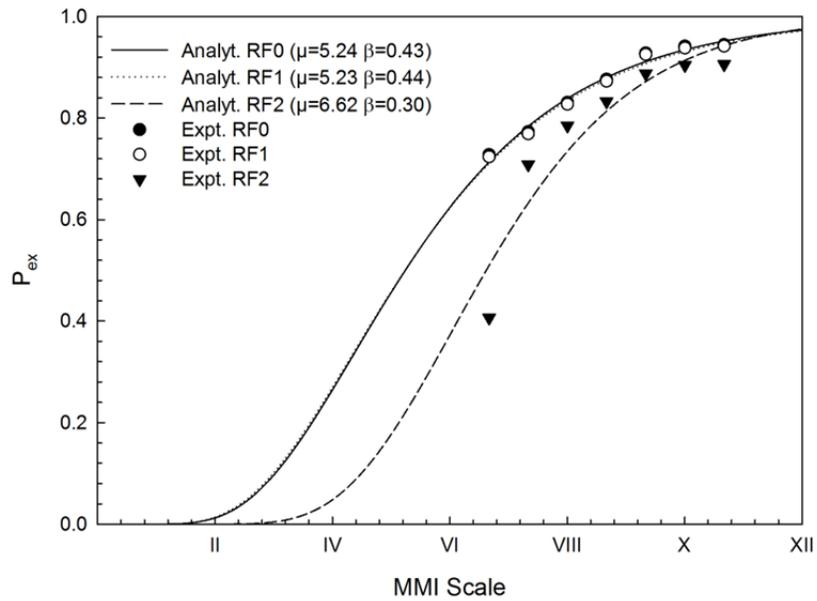

**Figure 7.26**  Reference restoration function given a DS = 2 (severe damage) using MLE method the (Emergency Department with emergency response plan applied).

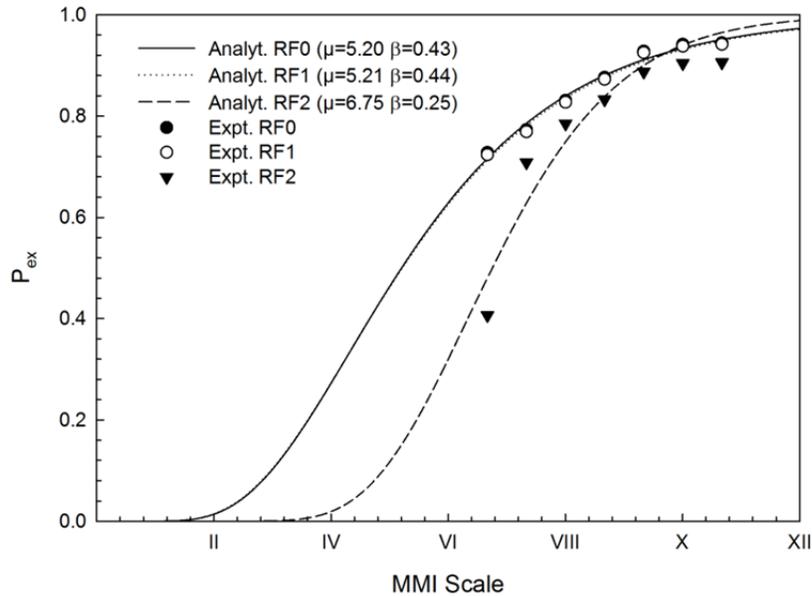

**Figure 7.27** Reference restoration function given a DS = 2 (severe damage) using SSE method (the Emergency Department with emergency response plan applied).

## 7.9 COMPARISON BETWEEN EMERGENCY DEPARTMENT WITH AND WITHOUT EMERGENCY PLAN APPLIED

The probability of exceedance of a given restoration curve is higher without the ERP than when the ERP is applied. Therefore, the ERP can be considered effective since PWTs when the ERP is applied is significantly lower than PWTs without the ERP. However, the only exception is when the damage state is severe (DS = 2). In that case, the RRFs of both case scenarios mainly overlap.

## 7.10 COMPARISON BETWEEN SSE AND MLE FITTING METHODS

As can be seen in Figure 7.31–Figure 7.33, the two fitting methods produced very similar results. Table 7.3 compares the parameters estimated with MLE method and SSE method. Note that the values slightly differ from one another.

**Table 7.2** Number of simulations conducted: no damage.

| | RFF DS no damage | | | |
|---|---|---|---|---|
| | MLE | | SSE | |
| | $\vartheta$ | $\beta$ | $\vartheta$ | $\beta$ |
| RF0 | 17.44 | 0.52 | 17.54 | 0.53 |
| RF1 | 14.06 | 0.23 | 14.23 | 0.25 |
| RF2 | 1.47E+07 | 0.10 | 1.47E+07 | 4.36 |

Table 7.3        Number of simulations conducted: moderate damage.

| | RFF DS moderate damage | | | |
|---|---|---|---|---|
| | MLE | | SSE | |
| | $\vartheta$ | $\beta$ | $\vartheta$ | $\beta$ |
| RF0 | 12.85 | 0.75 | 12.88 | 0.76 |
| RF1 | 12.15 | 0.59 | 12.31 | 0.62 |
| RF2 | 1.43E+07 | 6.09 | 1.43E+07 | 6.18 |

Table 7.4        Number of simulations conducted: severe damage.

| | RFF DS severe damage | | | |
|---|---|---|---|---|
| | MLE | | SSE | |
| | $\vartheta$ | $\beta$ | $\vartheta$ | $\beta$ |
| RF0 | 5.24 | 0.43 | 5.20 | 0.43 |
| RF1 | 5.23 | 0.44 | 5.21 | 0.44 |
| RF2 | 6.62 | 0.30 | 6.75 | 0.25 |

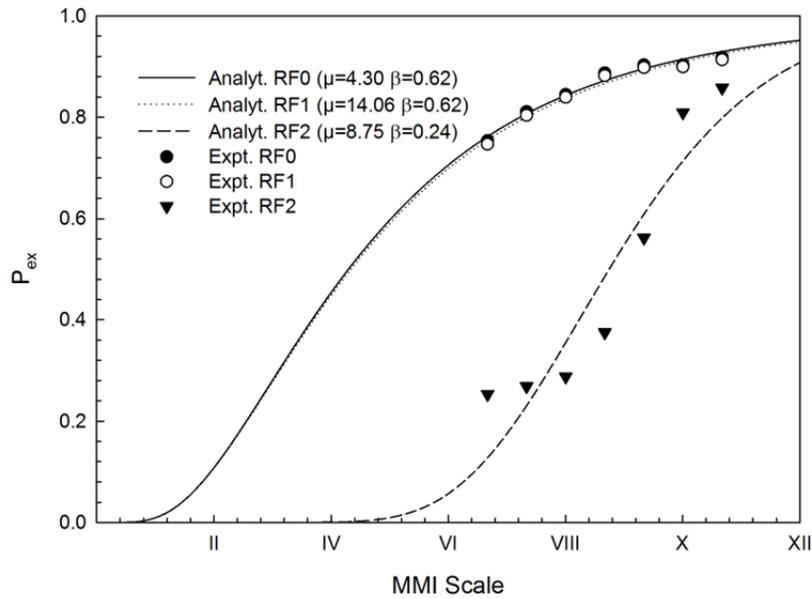

Figure 7.28    Reference restoration function given a DS = 0 (no damage) using MLE method (the Emergency Department without emergency response plan).

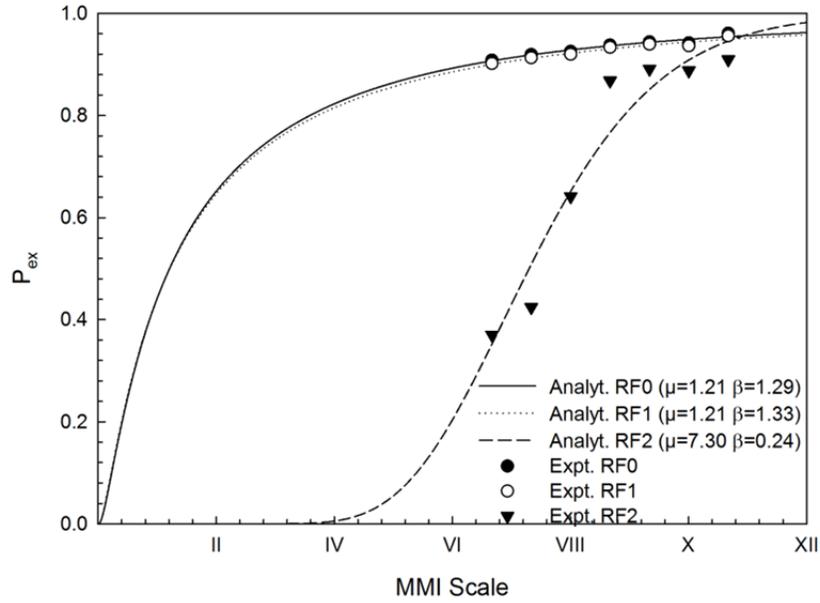

Figure 7.29  Reference restoration function given a DS = 1 (moderate damage) using MLE method (the Emergency Department without emergency response plan).

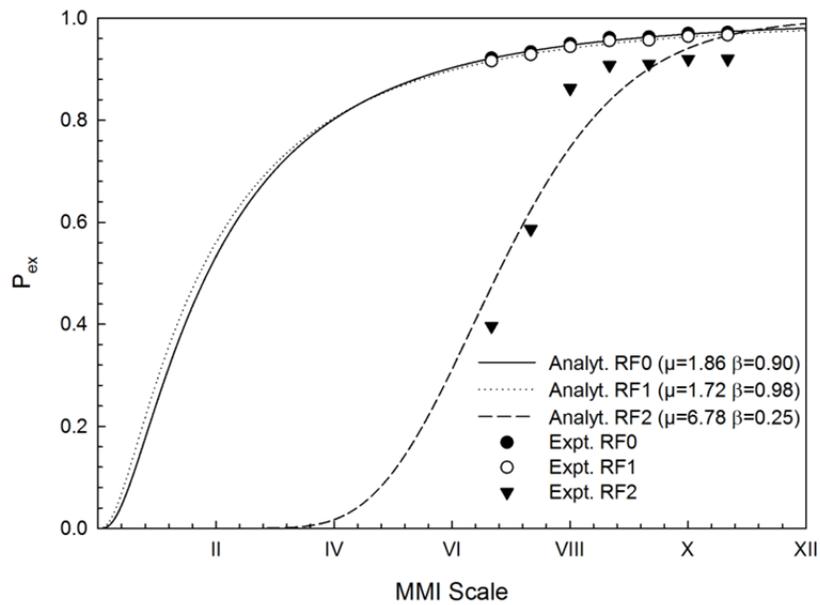

Figure 7.30  Reference restoration function given a DS = 2 (severe damage) using MLE method (the Emergency Department without emergency response plan).

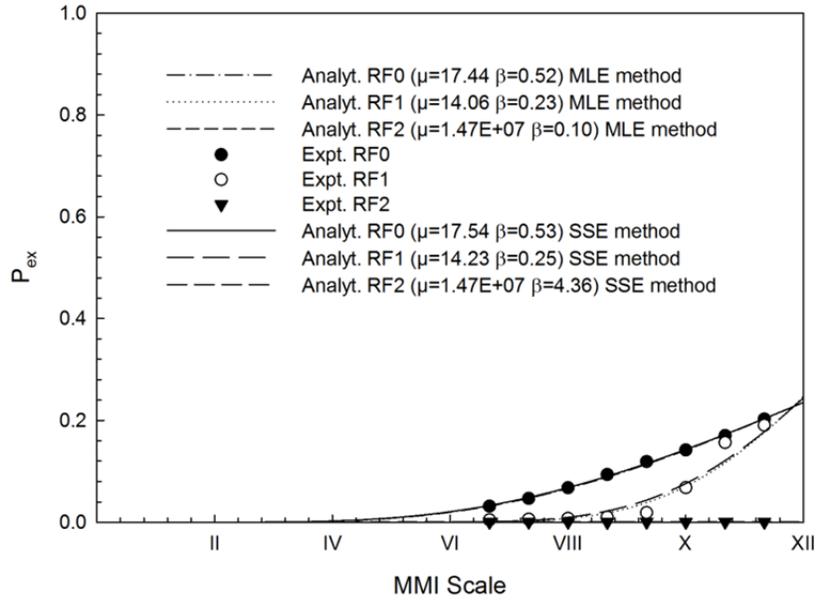

**Figure 7.31** Reference restoration function given DS = 0(no damage) using MLE and SSE methods (the Emergency Department with the emergency response plan applied).

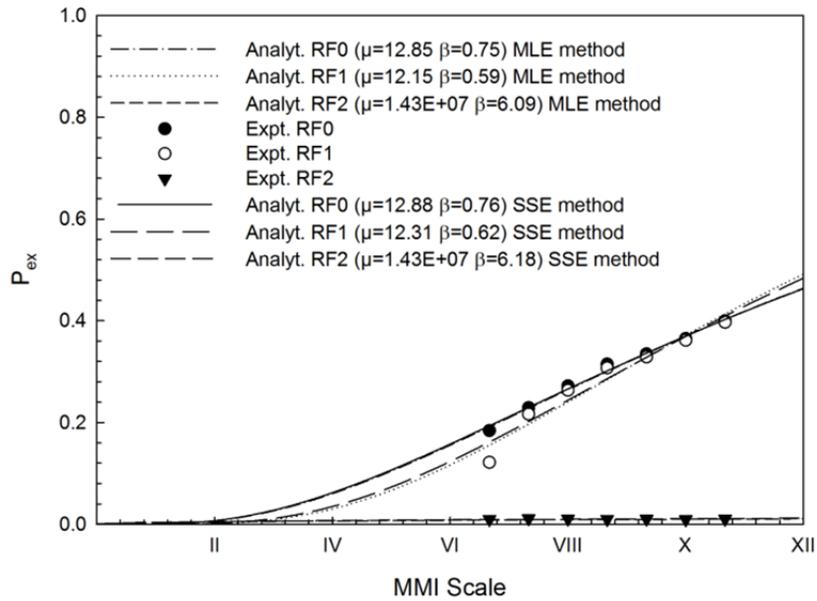

**Figure 7.32** Reference restoration function given a DS = 1(moderate damage) using MLE and SSE methods (the Emergency Department with emergency response plan applied).

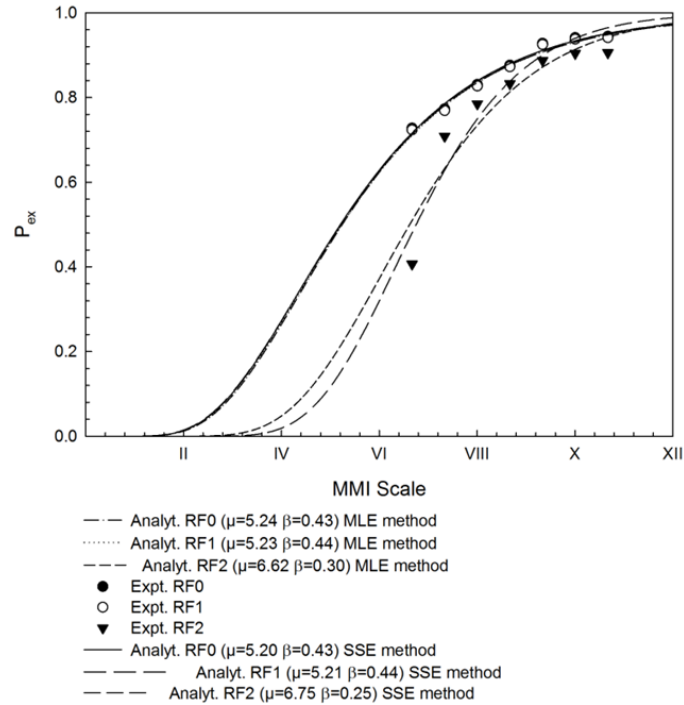

Figure 7.33  Reference restoration function given a DS = 2 (severe damage) using MLE and SSE methods (the Emergency Department with emergency response plan).

## 7.11 REMARKS AND CONCLUSIONS

Presented herein a methodology for building restoration fragility functions (RFFs) that describe the probability of exceedance a given restoration curve associated with a given damage state. The restoration process is one of the most uncertain variables in the resilience analysis, therefore, it is necessary to consider it in probabilistic terms and introduce RFFs.

Restoration fragility functions can be a useful tool in defining resilience of a hospital network. For example, they can be used to estimate the restoration process of an ED as a function of the seismic intensity. The main difference between RFFs and a standard fragility function is that the RFF is correlated to a given damage state. In other words, RFFs are conditional on DS and $I$, while standard fragility curves are only conditional on the intensity measure $I$. The method has been applied to the model of an ED of an existing hospital during a crisis with and without an emergency response plan (ERP).

A discrete event simulation model (DES) of the ED was developed and patient wait times (PWTs) collected as output. The RFFs for three different damage states were developed. The MLE and SSE methods of estimation of the lognormal cumulative function parameters were compared, which showed similar results. The ERP was considered effective since the RFFs are significantly lower when the ERP is applied. However, the only exception is when the damage state is severe (DS = 2), because in that case the RRFs of both cases mainly overlap. Further studies could lead to the development of fragility restoration functions based on different resilience indicators beside PWTs, e.g., structural damage and economic losses.

# 8 Modeling the Interdependence of Lifelines at the Regional and Local Level with Temporal Networks

## 8.1 INTRODUCTION

Recent interest from both scientists and policy makers has focused on increasing community resilience and identifies lifelines as one of the most important areas to consider. Lifelines can be defined as critical infrastructure systems—network-structured and interdependent—that provide a reliable flow of services and goods essential to the economic, social, and political security of a community. Disruptive events like the 9/11 terrorist attack and the Fukushima Dai-chi nuclear disaster, along with high costs related to recovery and reconstruction phases, have highlighted the importance of this issue. Some examples are the creation of the U.S. President's Commission on Critical Infrastructure Protection and the European Program on Critical Infrastructure Protection. This increased government attention urges universities and research centers to increase efforts in the study and modeling of lifelines behavior to improve the resilience of communities.

Two of the most studied and analyzed features of lifelines are interdependencies and temporal effects. Critical infrastructure systems are not isolated but highly interdependent and mutually interconnected. The links among different networks increase the potential of cascading failures, which can bring to catastrophic amplification of the impact. In the case of large-scale failures, they destabilize the system's environment, making it non-stationary and subject to temporal effects. Time-dependent analysis is required when temporal inhomogeneity matters and the sequence of event is important, which is usually the case in an emergency situation. If these aspects are ignored, system performance may be greatly overestimated. In the emergency response phase, these considerations become even more relevant because of the reduced period of time considered and of the density of events that populate it.

The aim of this work is to implement a temporal network approach in the civil engineering and the emergency management sector. First, we introduce the state-of-the-art and the features of lifeline systems. Second, we analyze and describe the methodology adopted for modeling a lifeline system in an effective and appropriate manner. Finally, this work is validated using a real case study.

## 8.2 MODELING TEMPORAL NETWORKS

Like many other networks, lifelines can change their characteristics and topology over time. Being critical systems, they are designed to be reliable even under stress, and in the best case scenario there are usually backup systems in place. To take into account the effect of time on networks, it is necessary to have a model capable of representing at every time step the current condition of the system.

The problem will be analyzed studying the connectivity of the system rather than the physical phenomena involved. First, a brief presentation of existing methods for evaluating the connectivity features of networks is introduced that focusses on the input-output inoperability method (IIM), an important methodology for evaluating cascading effects in a system. This system has some limitations; therefore, application of some of its methodology is utilized in conjunction with the development of a new methodology to model certain features, like the temporal variability of the topology. The suggested method is then compared to the probability risk assessment (PRA) method used for the analysis of critical sites. The results suggest running analyses at different scales of detail, from the regional level to the local level.

### 8.2.1 Existing Interdependence Models

This section groups and reviews the existing modeling and simulation approaches used for interdependence analysis. They are broadly categorized into six types: system dynamics-based models, network-based approaches, empirical approaches, agent-based approaches, and economic theory-based approaches. An analysis of each approach is conducted, highlighting the strengths and weaknesses of each approach [Ouyang 2014].

### 8.2.2 System Dynamics-Based Models

System dynamics-based approaches model the dynamic and evolutionary behavior of the interdependent lifelines by capturing important causes and effects under disruptive scenarios. System dynamics-based approaches use a top-down method to manage and analyze complex adaptive systems involving interdependencies. Feedback, stock, and flow are the basic concepts in this type of approach. Feedback loops indicate connection and direction of effects between the infrastructure's system components. Stocks represent quantities or states of the system, the levels of which are controlled over time by flow rates between stocks. System dynamics-based approaches model the interdependent infrastructures using two diagrams: (1) a causal-loop diagram capturing the causal influence among different variables; and (2) a stock-and-flow diagram describing the flow of information and products through the system.

This type of approach has some weaknesses. First, as the causal loop diagram is established based on the knowledge of a subject-matter expert, it is also a semi-quantitative method. Thus, many parameters and functions in the models require calibration, which need a huge amount of data that is not easily accessed. Lastly, due to the difficulty to obtain relevant data, validation efforts usually consist of conceptual validation so there is relatively limited validation of the model. These weaknesses call for integrating other modeling approaches in a uniform analysis framework for overall decision support [Bush et al. 2005].

### 8.2.3 Network-Based Models

As already affirmed before, infrastructure systems can be described by networks where nodes represent different system components and links mimic the physical and relational connections among them. Network-based approaches model single examples of infrastructure by networks and describe the interdependencies by interlinks, providing descriptions of their topologies and flow patterns. Performance response of lifeline systems to hazards can be analyzed by first modeling the component failures from hazards at component level and then simulating the cascading failures at the system level [Patterson et al. 2007]. Depending on whether or not the particle flow is modeled, network-based studies are broadly grouped into topological models and physics-based models.

### 8.2.4 Empirical Approach

The empirical approach analyzes lifelines interdependencies according to historical events or disaster data and expert experience. Studies using this type of approaches can identify frequent and significant failure patterns, and quantify interdependence strength metrics to inform decision making. Historical interdependence incidents can be used to uncover the interdependence structures or relationships between critical infrastructures under extreme events, such as the 2011 Tohoku earthquakes in Japan. Establishing special databases from the incident reports and then analyzing the data can help identify the frequent and significant failure patterns. Usually, interdependency incident records are collected from newspapers, media reports, internet news outlets, official ex-post assessments, and utility owners and operators.

This type of approach has several weaknesses. First, due to the bias of reporting, there may exist underreporting of some frequent interdependence failures that may have significant impact. Second, scholars use different databases to collect failure data without a standardized data-collection methodology for interdependent critical infrastructure performance. Third, the reliance of the empirical approaches on previous failure records may not provide good predictions for new disasters. These weaknesses call for other modeling and simulation approaches for additional decision support [McDaniels et al. 2007].

### 8.2.5 Agent-Based Models

Agent-based approaches are an effective way to model critical infrastructure systems and the related decision-making process that characterize them during an emergency. These approaches adopt a bottom-up method and assume the complex behavior or phenomenon emerges from many individual and relatively simple interactions of autonomous agents. Each agent interacts with others, and its environment based on a set of rules that mimic the way a real counterpart of the same type would react. Most critical infrastructure components can be viewed as agents.

Agent-based approaches model the behaviors of decision-makers and the main system participants in the interdependent lifelines in order to capture all types of the interdependencies among lifelines by discrete-event simulations, provide scenario-based what-if analysis and the effectiveness assessment of different control strategies, and can be also integrated with other modeling techniques to provide more comprehensive analysis. However, this type of method has some weaknesses. First of all, the quality of simulation is highly dependent on the assumptions made by the modeler regarding agent behaviors, and such assumptions may be difficult to justify theoretically or statistically. Then, calibrating the simulation parameters is a challenge due to lack of relevant data and the difficulty in modeling participant behaviors; detailed information about each critical infrastructure system is considered highly desirable by utilities managers [Bonabeau 2002].

#### *8.2.5.1 Economic Theory-Based Models*

Lifeline systems interdependence can be analyze through models of economic interdependencies. In the existing literature, two types of economic theories are employed to model lifelines interdependencies: input–output and computable general equilibrium [Rose 2005]. Inoperability input-output models can easily analyze how perturbations propagate among interconnected infrastructures and how to implement effective mitigation efforts. This model will be presented in detail in the next section. Computable General Equilibrium-based methods extend the capacities of the input-output methods, capture the nonlinear interactions among infrastructure systems, provide resilience or substitution analysis of single infrastructure and the whole economy, and are able to capture different types of interdependencies in a single framework. The weaknesses of these types of models are related to calibration and data acquisition [Partridge and Rickman 1998].

### 8.2.6 Input-Output Inoperability Method: Analysis

Developed by Haimes and Jiang [2001], the IIM model is an adaptation of the Leontief's input-output (I-O) analysis of economic interdependencies [1986]. However, instead of focusing on the economic impact of a perturbation, the IIM proposed in this work is intended to simulate the propagation of risk of inoperability in the infrastructure sector.

Inoperability is defined by the authors as "the inability for a system to perform its intended function." It is quantified by a value between 0 and 1, determined from considerations on the likelihood and the level of failure. When the inoperability of an element is 0 it means that it is working at the top of its potentialities; when it is 1, it is completely inoperative. These risks of inoperability are propagated between different networks following interdependency patterns. The equation describing the IIM is as follows:

$$q = [I - A]^{-1} \cdot c \tag{8.1}$$

where $q$ is the damage vector that contains the inoperability values for the $n$ infrastructures considered; $A$ is a matrix that depicts the extent of interdependence between infrastructures and is the transpose of the adjacency matrix that describes the topology of the system; $I$ is an identity matrix, and $c$ is the scenario vector that includes the effects of the perturbation (e.g., natural disasters, man-made attacks, intrinsic failures, etc.) on each infrastructure.

The damage vector $q$ is the output of the model and quantifies the level of inoperability of the infrastructures composing the system following a perturbation that propagates according to the topology described by the interdependency matrix, $A$. Each element of this matrix quantifies the level of influence of the $j$th infrastructure on the $i$th infrastructure. They can be a value between 0 and 1, which is complete propagation of the scenario from $j$ to $i$, and 0, whereby there is no propagation from *j* to *i*. Thus, the *A* -matrix represents the probability of transferring inoperability across different infrastructures.

To give an example of how the IIM works, the case of a six-node network developed by Valencia [2013] is used, which will be referred to as Example 1; see Figure 8.1. This example is concerned with two networks: an electric network and a water network, which serve three buildings.

The hazard considered is infrastructure aging. To measure the impact of individual node decay across the network, the column summation of the damage vector *q* of each node *i* at each time *t* is computed; see Figure 8.2. To calculate the decay score, we used the following equation:

$$dc\_s_i(t) = \sum_{i=1}^{j} q_i(t) \qquad (8.2)$$

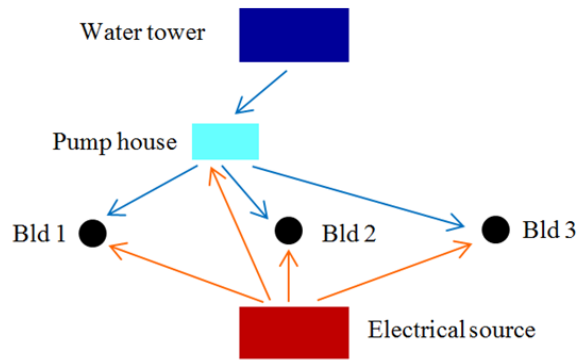

**Figure 8.1**    **Graph representing Example 1 topology.**

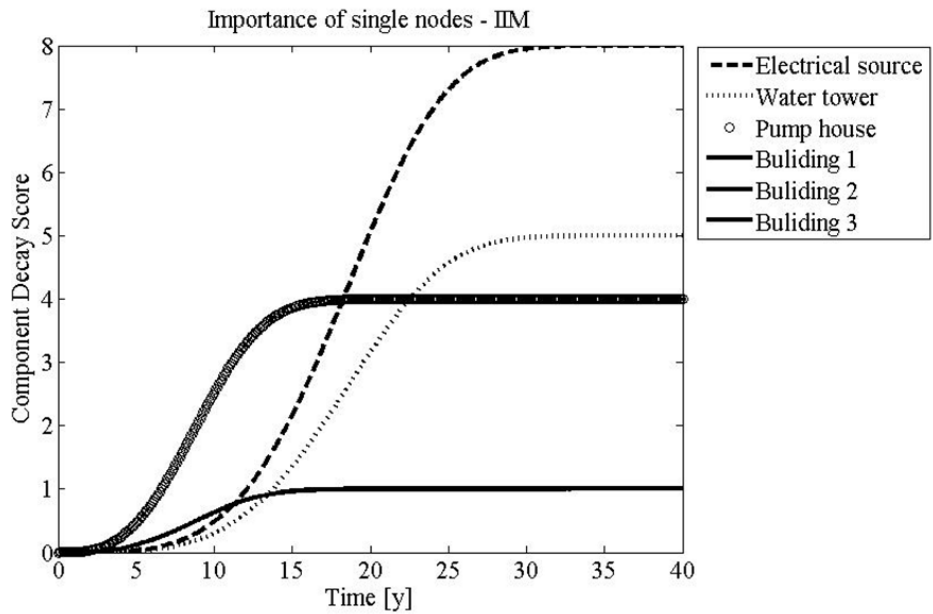

**Figure 8.2**    **Component decay score for Example 1.**

This approach, which was applied to a complex infrastructure network, presents three severe limitations: (1) it does not take into account the redundancies of the system; (2) because it is a static model, it does not consider the temporal evolution of the system and does not account for temporal effects that can disrupt the system; and (3) its inputs and outputs are not user-friendly.

To address the limitation relative to redundancies, a simple solution has been implemented. Regard the topology of the system in Figure 8.1. If a new pump house is added in parallel to the first one, the network presents a redundancy, which will be referred to as Example 2. Figure 8.3 shows the new topology of Example 2. It is clear that the performance of the system is improved with respect to the previous case because both the pump houses can perform the same work, and their simultaneous failure is more unlikely than the failure of just one.

We expect that the impact of the water tower and the electrical source remains the same, while the impact of each of the pump houses decreases. If we compare the results obtained in Figure 8.4 in the case where there is only one pump for the system, we can see how the expected trends are not present. Electrical source and water tower $dc\_s$ increase and the pumps decay score doesn't change.

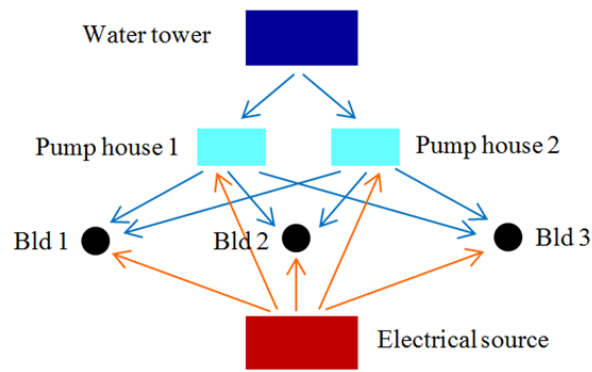

**Figure 8.3**  Graph representing Example 2 topology.

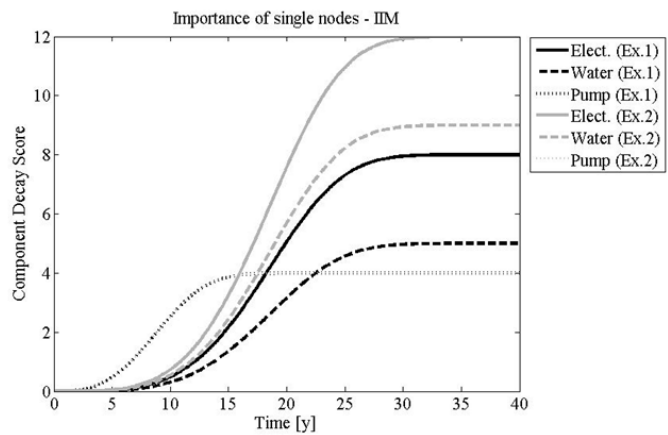

**Figure 8.4**  Comparison of decay scores of Example 1 and Example 2, using the traditional IIM.

To solve the problems related to redundancies, probabilities of nodes in parallel can be combined properly. In the previous case, the $dc\_s$ of electrical source and water tower are increasing because the algorithm sees another node (the new pump) that needs to be operated as well. To avoid this, we introduce the Series-Parallel Vector:

$$SP = \begin{Bmatrix} 1/n_1 \\ 1/n_2 \\ \vdots \\ 1 \quad (for\ BLD) \end{Bmatrix} \quad (8.3)$$

where $n_i$ is the number of nodes redundant of node $i$. After having expanded it to the $n$-dimension, it is possible to add it to the damage vector equation:

$$\boldsymbol{SP}^* = \boldsymbol{SP} \times \{1 \quad 1 \cdots 1\}_{1 \times n} \quad (8.4)$$

$$q_i(t) = \left[\boldsymbol{I} - \boldsymbol{A} \cdot \boldsymbol{SP}^*\right]^{-1} \cdot c_i \quad (8.5)$$

After this operation, the values of the electrical source and water tower return to the right values. In regards to the pumps, if we assume that failure of the two pumps is stochastically independent, the probability that both fail at the same time is given by:

$$P(A \cap B) = P(A) \cdot P(B) \quad (8.6)$$

Result of this implementation reflects the initial expectations about the effect of redundancy and is shown in Figure 8.5.

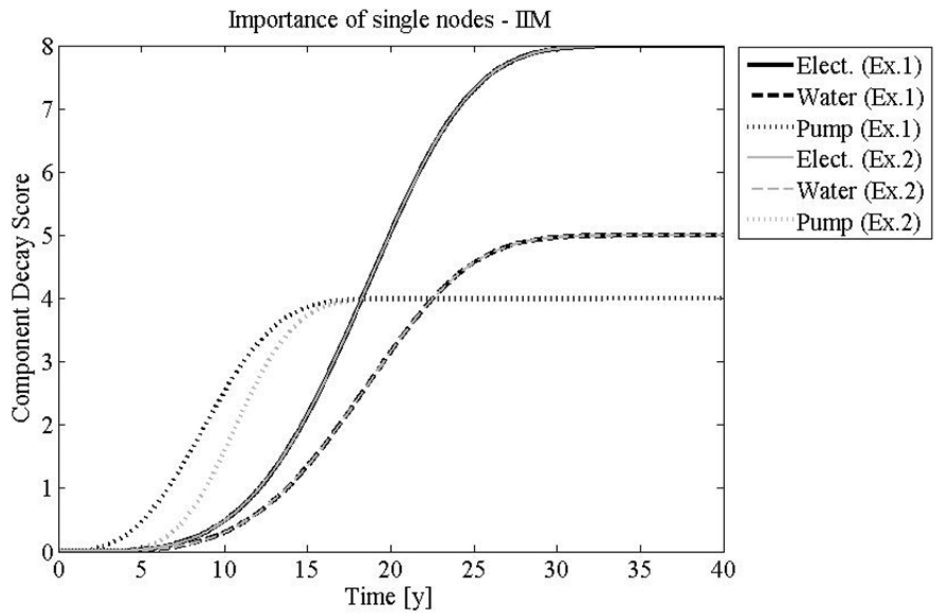

**Figure 8.5**   Comparison of decay scores of Example 1 and Example 2 using the implemented IIM.

Now that it is clear how the introduction of an additional pump can improve the performance of a system, it is necessary to introduce an index able to represent the performance of the entire system and not just single nodes. Thus, the system score is introduced. It is a dimensionless risk index that varies in the range $0 \div \infty$ and expresses the rating of a system of infrastructures at a specific time $t$ it is defined as follows:

$$sys\_s(t) = \sum_k \frac{\sum dc\_s_{k,i}(t)}{n°_k} \tag{8.7}$$

where $k$ is the type of node (i.e., electrical sources, water towers, and pump houses). The final targets (i.e., buildings) are not considered when calculating the $sys\_s$. A low value of $sys\_s$ indicates that the system of infrastructure has low risk of failure at the target nodes, while a high value indicates high risk. A threshold separating the low-risk region from the high-risk region needs to be calibrated on the basis of the importance of the system and minimum acceptable performance.

Through this new index, a sensitivity analysis can be performed to establish which intervention better improves the performance of the system. A ten-building system is considered. Analyzed improvements are: (1) modification of the system into two smaller systems; (2) the introduction of redundant nodes; and (3) a plan of maintenance interventions.

As shown in Figure 8.6, the positive effects of the intervention is represented by the drastic lowering of the plateau of $sys\_s$ in Figure 8.7. Studying the effect of applying the redundancy intervention of Figure 8.8 instead demonstrates how the plateau does not vary, but the $sys\_s$ in Figure 8.9 decreases in the short-term segment of the function.

Comparing these two scenarios with the one obtained from planned maintenance interventions, it is clear that modification of the system is the best solution for the long term. In the short term, the addition of redundancies is very effective. Maintenance interventions have a relevant positive effect in both the short and the long term, and are probably the most feasible solution from an economic point of view. Results of the sensitivity analysis are shown in Figure 8.10.

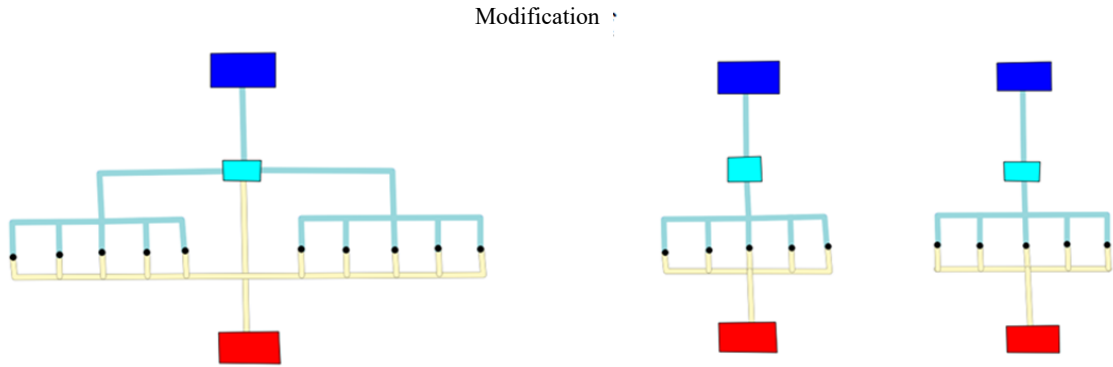

**Figure 8.6**     **Modification of Example 1 system.**

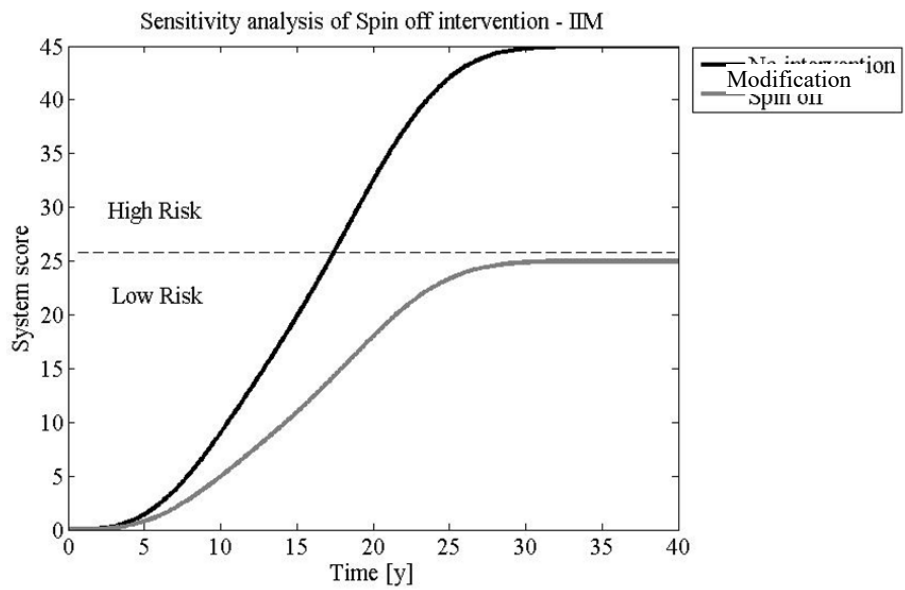

**Figure 8.7** System score before and after the modification intervention applied to Example 1.

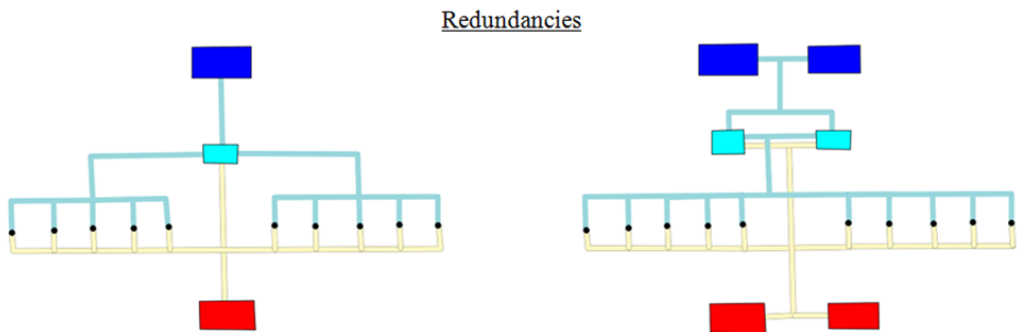

**Figure 8.8** Redundancy intervention applied to the Example 1 system.

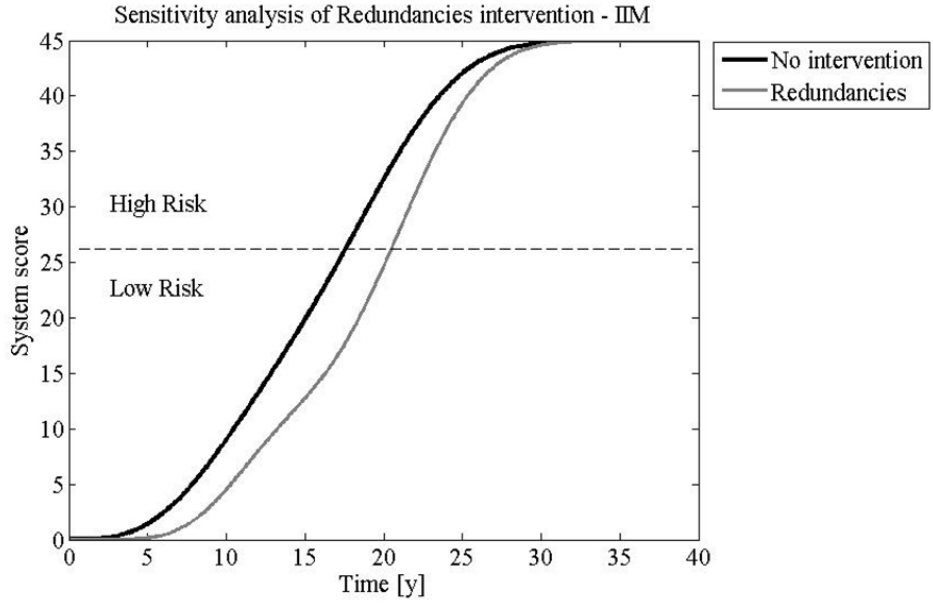

**Figure 8.9** System score before and after the redundancy intervention applied to Example 1.

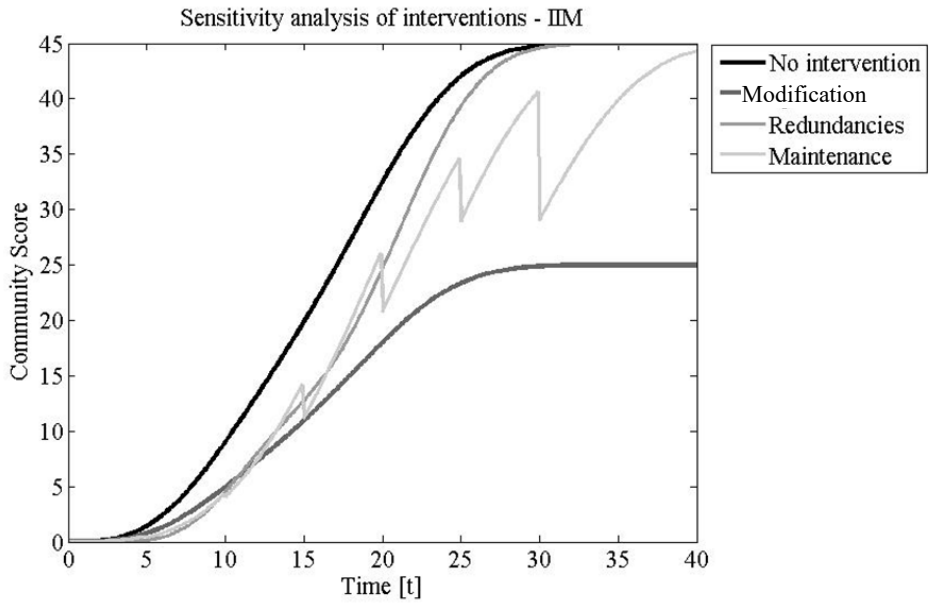

**Figure 8.10** Sensitivity analysis of interventions for risk mitigation on Example 1 system.

### 8.2.7 Modified IIM for Temporal Networks

Many extensions of the model have been proposed such as the Dynamic IIM (DIIM) and the Multi-Regional IIM (MR-IIM). However, the modified IIM model presented hereafter tries to overcome some of the limitations of the methodology proposed by Haimes et al. [2005], while capturing the key aspects of infrastructure behavior during an emergency. A presentation of the theoretical framework of the model and then three different implementations of the IIM are presented below.

#### 8.2.7.1 Topology formalization using graph theory

Graph theory has been used to model infrastructure networks. The geographical, topological, and flow information of a network can be represented with a graph $G(V, E)$, which is formed by a set $V$ of vertices, herein called nodes, and a set $E$ of edges. The definition of the nodes is dependent on the spatial scale of the problem considered, which might be an entire infrastructure (e.g., electric network, water network, or gas network) [Cimellaro et al. 2014], a sub-system (e.g., wind turbine), or even a unit (e.g., a gearbox of the wind turbine). Each node can be attributed with specific features such as hierarchy, resistance, and autonomy; edges do not have any features assigned in the proposed model, but they are oriented. The edges can link nodes to the intra-network (i.e., within a specific infrastructure) or to inter-networks (i.e., across different infrastructures). An inter-network represents the interdependencies described in the $A$-matrix. Instead of attempting to specify the likelihood and the degree of interdependency in the $A$-matrix, this model defines an inter-network link as Boolean, either 0 or 1. Thus, $a_{x_i y_j}$ values will be 0 if the $x$th node belonging to $i$th infrastructure is dependent on the $y$th node belonging to the $j$th infrastructure.

With respect to existing formulations, the concept of chains is introduced in the model. A chain is a sequence of nodes from one vertex to another using the edges. The chains of interest are those that connect a source (i.e., a node without inflows) to a sink (i.e., a node without outflows). The task of every source is to feed all the sinks of the network if the topology allows doing so. If it doesn't, the source is called partial. An example of partial source is a photovoltaic plant on the roof of a building. This plant belongs to the general electric network of a city or a block, but it only feeds the building where it is located and not any other structure. It is assumed that every node of a chain must have at most one inflow edge, but it can have multiple outflow edges. This means that different supply lines exist in a critical infrastructure system. For example, besides the main supply line, backup lines exists that can substitute in case of failure or malfunction. Each of these chains can guarantee the operability of the network even though they are mutually exclusive.

The hierarchy of their operation is defined by the design of the infrastructure. There are two types of hierarchy. The source hierarchy corresponds to the rank of priorities for the entry into operation of the sources. The path hierarchy corresponds to the rank of priorities for the activation of different possible paths. It is assumed that source hierarchy is stronger than path hierarchy. This means that if the first chain fails, the network tries to maintain operation starting from the previous source that then inquires if new paths are available (if possible). If no other path is available for that source, then it skips to the next path. This theoretical framework and notation will be adopted while discussing methodology implementation to the IIM model.

### 8.2.7.2 Probabilistic Formulation of Inputs and Outputs

The proposed methodology modifies the IIM deterministic formulation in probabilistic terms because while the damage score just gives a snapshot of the cascading propagation of inoperability, it does not say anything about the final state of the network. The probability of failure of a single node is obtained by combining the natural hazard with the infrastructure vulnerability and refers to the status (fully operative or failed) of the node itself after the perturbation. Hereinafter, it will be called self-failure probability ($P_{sf}$) and will substitute for the scenario vector, $c$.

The hazard component is represented by an event vector $E n \times 1$, where $n$ is the number of nodes in the system. At a given time $\bar{t}$, every node will be disrupted by a natural event (e.g., earthquakes, tsunamis, fires, sabotage, etc.). The elements of the $E$-vector can be physical quantities such as the, *pga*, *pgv*, and *pgd*, earthquake magnitudes, height *hw* of a tsunami wave, the megatons *Mt* of an explosion, etc. These quantities can be different from node-to-node because infrastructures usually have a large spatial extension; see Figure 8.11(a).

By performing different simulations, using different $E$-vectors, it is possible to approach the problem in probabilistic terms. Each simulation has a weight, which corresponds to the probability of occurrence of the event of a certain magnitude that is directly taken from the hazard curves. The vulnerability component of each node is represented by the fragility curves, which define the probability of failure of each node depending on the type of hazard considered; see Figure 8.11(b). Thus, there are as many fragility curves for each node as the type of acting hazard. Only complete-failure fragility curves are used; intermediate damage levels are not considered at this stage.

The probability of failure $P_{sf}$ of a node under a specific event $E$ is obtained by inserting the value of the $E$-vector in the node fragility curve. Proposed by Valencia [2013], this approach of summing up the elements of the $q$-vector to obtain a final score to evaluate the interdependency performances has obvious limitations because they are not normalized to the dimension of the system (e.g., the longer will be the chain, the higher the score will be). Moreover, as pointed out previously, the index proposed by Valencia does not take into account the benefits provided by redundancies present in the infrastructure. In the modified IIM proposed, the probability of failure $P_f$ of every node is obtained by combining the $P_f$ with the cascading failure probability $P_{cf}$, which is transmitted by the upstream nodes and is calculated using a step-by-step approach that takes into account the ramifications of the system; see Figure 8.11(c). In other words, $P_f$ is the probability of failure of each node that is obtained as a result of all the disrupting events and the cascading propagation effects.

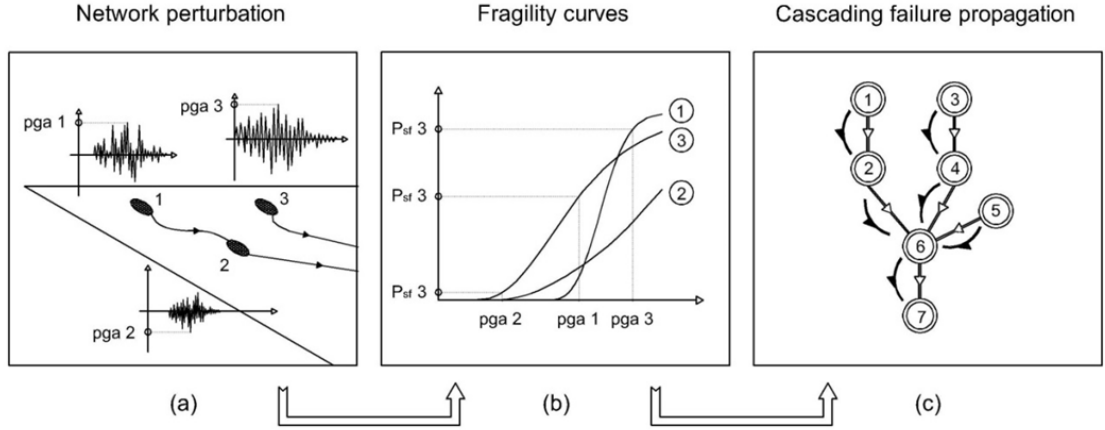

**Figure 8.11** Flowchart of the probabilistic approach: (a) a network subjected to perturbation: (b) probabilities of failure of nodes are computed on fragility curves; and (c) propagation according to the topology of the network.

### 8.2.7.3 Multilayer Approach for Spatial Interdependency

The more intuitive approach for analyzing a system of infrastructure is solving each network separately and then considering their interaction. Infrastructure networks can be seen as layers that overlap each other and share some nodes presented in both networks, and are virtually connected by inter-infrastructure edges. It is said "virtually" because it is not a real physical connection. Let's consider the network of Example 1. Operation of the element pump needs both electricity and water; therefore, it belongs to both the electric and water networks. By visualizing a layer, a single node will be projected in the two layers, and a virtual edge will link the two projections; see Figure 8.12.

The IIM model is incapable of dealing with different layers and adjacency matrixes. In fact, it needs to store the topology in one general matrix and considers the entire system as single network. This is because the IIM can only use square matrices, while the inter-networks matrices are usually rectangular. To overcome this limitation, Valencia [2013] suggested introducing the *I*-matrix. These are $n \times m$ matrices, where $n$ is the number of nodes of the $j$th infrastructure and $m$ the number of nodes of the $i$th infrastructure, which are dependent on the $j$th. The idea is to increase the values of *c*-vector of infrastructure $i$, by adding the *q*-vector computed for the $j$th network, as expressed in Equation (8.8):

$$c_{j \to i} = I_{j \to i}^T \cdot q_j + c_i \tag{8.8}$$

Inserting the output of the first network into the input of the second dependent network is the correct approach for evaluating cascading effects. However, because this formulation starts from the same deterministic values as before, it cannot be considered satisfactory. The current method involves the combination of $P_{cf}$, for upstream and downstream networks:

$$P_{cf_i}^* = \left(I_{j \to i}^T \cdot P_{cf_j}\right) \cup P_{cf_i} \tag{8.9}$$

where $P_{cf_i}^*$ can be considered cascading-failure probability that incorporates in the node all the information coming from upstream networks and upstream nodes of the current network that converge at a certain point.

This multilayer approach brings many benefits:

(i) It acknowledges the analysis and results of layers and interdependencies, and helps understand where criticalities are located and which are the tighter and more stressed inter-links. While the evaluation of single infrastructure is mature, the interdependency studies are still at a development stage and determining their connectivity is the real issue;

(ii) Moreover, including possibilities to each infrastructure manager running the model of a given layer and then controlling the interaction between the different layers at a higher level is closer to the professional practice adopted during an emergency response phase. A model that considers all the elements of the system simultaneously won't be used in the real practice because none has the authority or expertise to manage all the data; and

(iii) In the end, the diffusion of informatics tools, like Geographic Information Systems (GIS), in both the emergency response and the risk planning sector suggests adoption of a unified methodology. The GIS platform has great potentialities, and it can be effectively used to organize input data and visualize outputs. Their relational databases are shaped with a layer structure that is in accordance with the one proposed above.

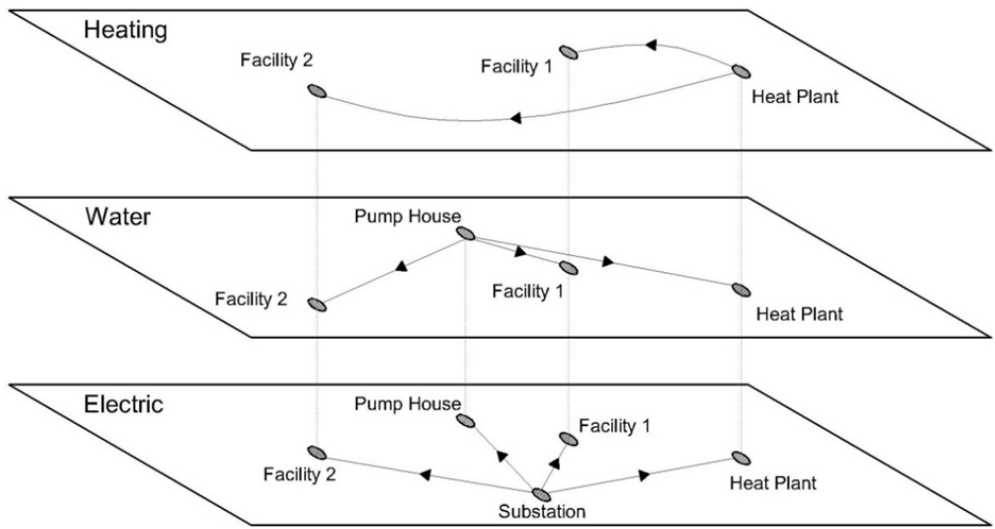

**Figure 8.12** Example of layer subdivision for interdependent networks.

### 8.2.7.4 Tensor Notation for Accounting Temporal Effects

What has not been addressed yet is temporal dimension. Compared to the traditional static IIM the first add-on is the introduction of a timeline $\tau = \{t_0, t_1, t_2, \cdots, T\}$, where the range $t_0 \div T$ must be extended enough to include all the events and their effects. Time step $\Delta t$ of the elements of the $\tau$-vector represents the time necessary for the propagation of the events across the entire system. This means that if at time $\bar{t}$ a landslide overwhelms an electrical pylon before time $\bar{t} + \Delta t$, the pylon will fail, and the effect of this failure must propagate throughout all the system. Therefore, the transmission of information in the system is immediate. After having solved the system at the time $\bar{t}$, the final situation will be the initial condition at time $\bar{t} + \Delta t$. Given this timeline, it is clear that each event must be associated with a time of occurrence, and that the model must run at every time step.

Now the model is not stationary but is composed of temporal networks, denoted by $G(t) = G[V, E(t)]$. The $P_f$ of nodes changes over time in accordance with the sequence of events, including the existence of edges. Changes in the status of the nodes may result in changes in the topology of the system. For example, let us consider the node $V_a$, which is a water purification plant; the node $V_b$ is a water collection pit, and the node $V_c$ is an aqueduct. $V_c$ is usually fed by $V_a$ through the edge $E_{ac}$, but if $V_a$ fails, the edge $E_{ac}$ disappears and the edge $E_{bc}$ is activated. The active chain has shifted from $a \rightarrow c$ to $b \rightarrow c$.

From this example, it can be inferred that different chains of a network are not only spatial layers but also temporal layers. Although the multilayer approach was effective in modeling interdependencies among different networks, here the networks are mutually exclusive and not linked. The solution adopted is to pass from bi-dimensional matrixes to a tri-dimensional tensor notation. The topology of every network is now described by an adjacency tensor whose elements are $a_{x_i y_j}(t)$. Each different temporal layer of the $A$-tensor represents a possible chain. The first in hierarchy is the ordinary supply line, while the others are backup lines. Figure 8.13 shows the three different possible functional configurations that the seven-node network examined can assume; this case will be referred to as Example 3.

To better understand which of the chains is active at time $\bar{t}$, the probability of occurrence of a specific configuration $P_{occ}$ is assigned to every layer; see Figure 8.14. This value determines if the layer is "on" ($P_{occ} = 1$), or if it is "off" ($P_{occ} = 0$), at the considered time step. In the current configuration, the condition for being "on" is that target nodes of the network do not fail and that configurations with higher degree of hierarchy are "off." Transferring this concept in the probabilistic model means that values of $P_{occ}$ become probabilities of being active. The sum of the probability of occurrence of a network is $0 \leq \sum P_{occ} \leq 1$, and the value $1 - \sum P_{occ}$ represents the percentage of lost capacity of the network ($LoC$).

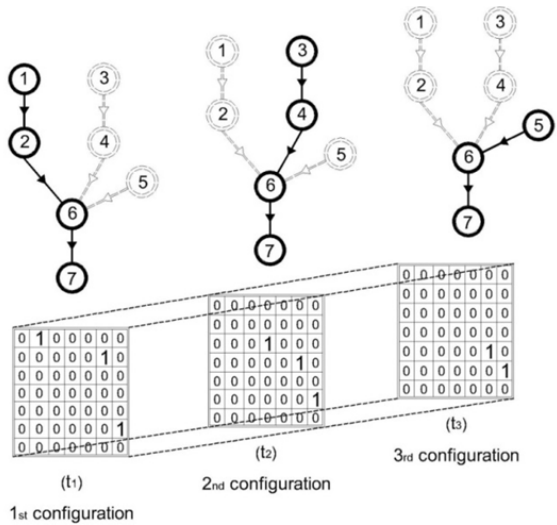

**Figure 8.13** Tensor notation for a network. An adjacency matrix is associated with each of the possible and mutually exclusive configurations of the network.

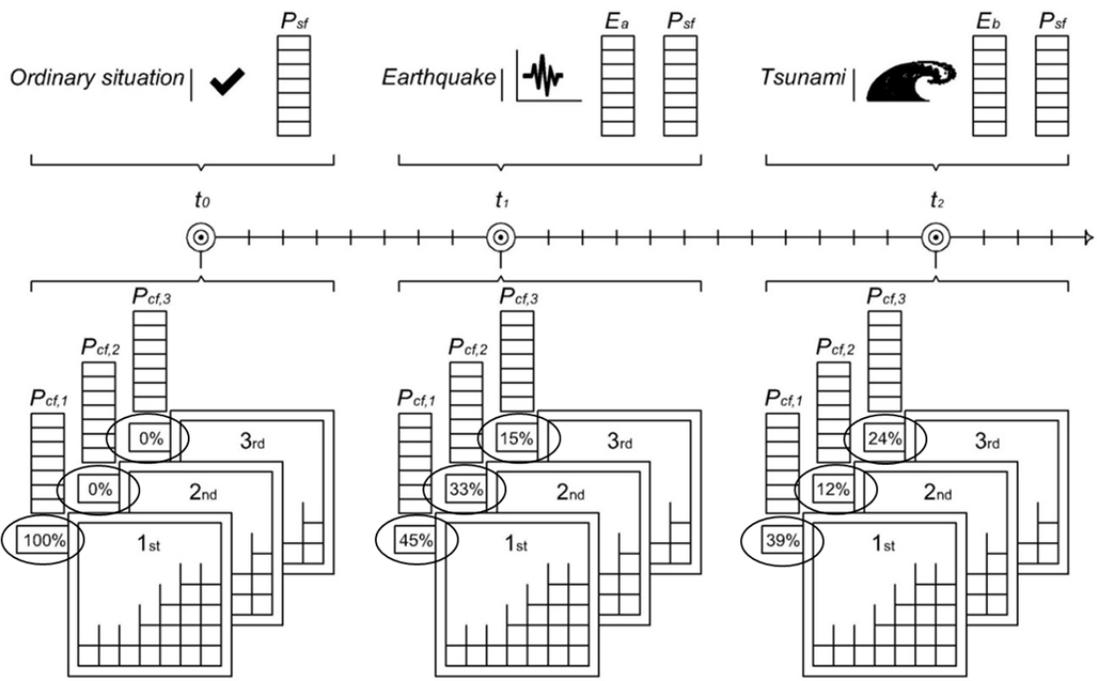

**Figure 8.14** Temporal variance of Operability Labels. Disruptive events tend to change their value and decrease the overall probability of operability $\left( \sum P_{occ} \right)$.

Once it is determined how probable the activation of a chain is, it is useful to determine additional time effects. For example, let's consider that the primary power source of a hospital is not working due to a blackout. The emergency power generator is activated to maintain the operability of the hospital, which is the target node. This UPS is fed by the fuel contained in a tank. Unless the time considered in the model is much smaller than the runtime of the tank, the unloading of the tank must be taken into account. The run out of autonomy of a node cannot be classified as an event, but its effects are well documented by various famous disasters. This work increases the importance of backup systems in lifelines because most backup sources have a capacity that can be considered limited in time. What emerges is that nodes can have temporal features—like the autonomy—that influence their status.

This temporal tensor notation has many advantages compared to static bi-dimensional of before:

(i) It is able to describe changes in the topology of the system that usually occurs after individual node failures;

(ii) The separation of chains in different layers allows computation of cascading failure responsibility of each node without considering the presence of parallel branches. The propagation of cascading effects is linear, and the results of each layer should be weighted with respect to the value of their operability label;

(iii) The $P_{occ}$ furnishes direct information on the activity of each chain and allows the evaluation of time-related effects, like autonomy;

(iv) It is possible to use the value $1-\sum P_{occ}$ as an index for quantifying the loss of capacity of the network; and

(v) The possibility of varying the topology of the system provides an opportunity to add new layers to existing networks. For example, in the case of recovering the operation of a network, rescue teams can modify its path or add a new provisional source, adding a new layer in the tensor.

#### 8.2.7.5 Comparison with Traditional IIM

Comparing the IIM with the modified IMM provides additional insight into the value of the suggested modifications. To evaluate the performance of the overall system, the system score $sys\_s$ was introduced. It is the sum of the terms of the ***q***-vectors. Its weakness is that it doesn't represent the circumstances at the target nodes, and that it needs a threshold to evaluate the level of risk; see Figure 8.15. The calibration of this threshold is problematic.

The modified IIM assesses the probability of failure, which don't need to be interpreted, and furnish a precise and mathematical measure of the risk. In an evaluation of the performance of the system, it is possible to observe the $P_f$ functions of target nodes. Figure 8.16 shows the $P_f$ for the electric supply and the water supply of each building of Example 2. Note how easy it is to identify that their functionality is limited ($P_f = 1$) after 20 years; in Figure 8.15 this is not clear.

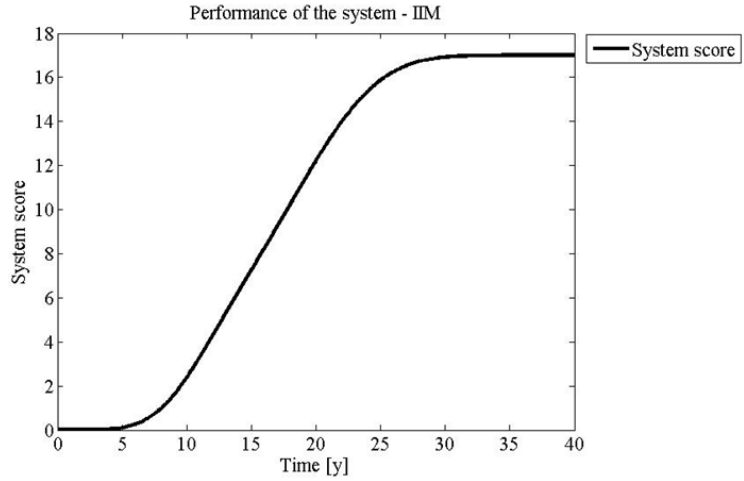

**Figure 8.15** Performance of the Example 2 network, according to the IIM.

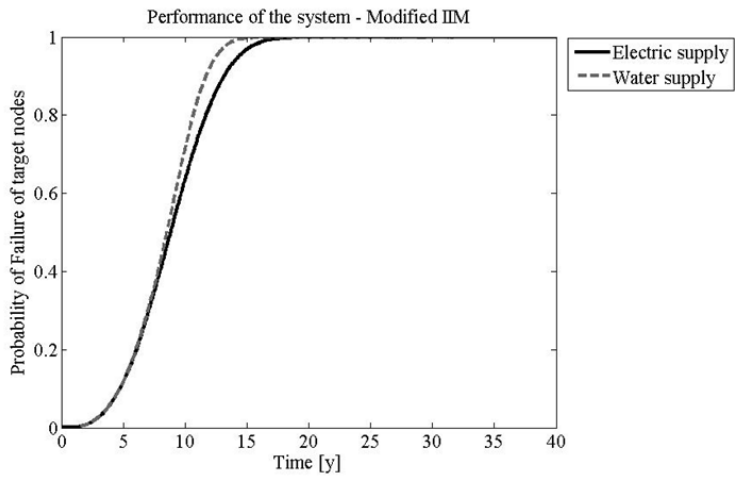

**Figure 8.16** Performance of the Example 2 network according to the modified IIM.

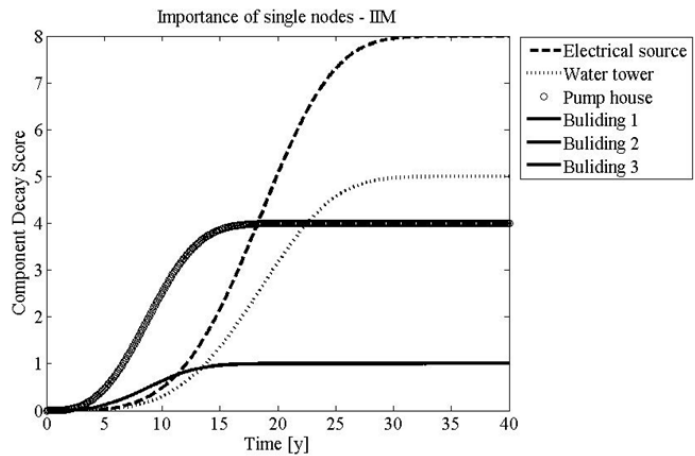

**Figure 8.17** Importance of nodes for the Example 2 network according to the IIM.

The component decay score, $dc\_s$ is a measure of the importance of single nodes in the network. This index does not consider the mutual effects of nodes in the network but multiplies the probability of self-failure $P_{sf}$ for the number of nodes that are topologically located downstream; see Figure 8.17. Using the modified IIM instead, it is possible to determine the influence each node has on the final failure of target nodes. If we compute the probability of failure of target nodes of a system in both cases and that the node of interest is subjected to an event ($P_{sf} \geq 0$) and not ($P_{sf} = 0$), it is possible to obtain the curves shown in Figure 8.18. The difference between the two functions at each time step represents the effect that the node has on the entire system. The higher the difference, the most relevant is the damaging effect that the considered node propagates to the target nodes. If we plot this difference, we obtain the curves shown in Figure 8.19 and Figure 8.20.

Considering Example 2, results show that damage to a building's electrical supply is much more relevant than damage to the electrical source. This is because fragility curves of buildings go to 0 more rapidly than the curves for the electrical source once buildings are down as no electric furniture will be available to users. If we conduct the same analysis for the water network, the results are the opposite. Here, the buildings are less important compared to the electrical network because the combination of probability of failure of electrical source, water tower, and pumps is higher than the probability of failure of the electrical source itself. This result would be totally unpredictable if the traditional IIM had been applied.

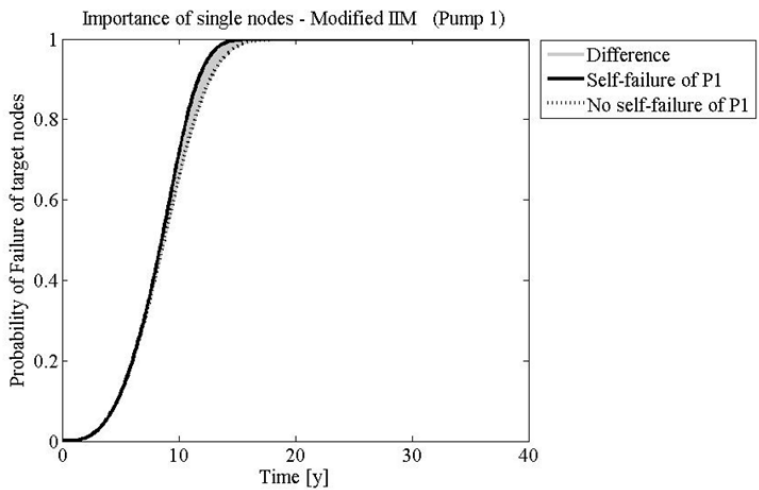

**Figure 8.18**     Differences induced in the system by Pump 1 for Example 2 network.

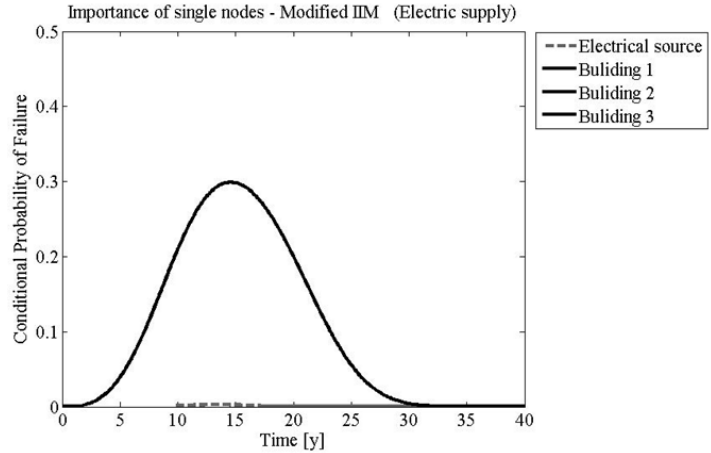

**Figure 8.19** Importance of nodes for the Example 2 network electric supply according to the modified IIM.

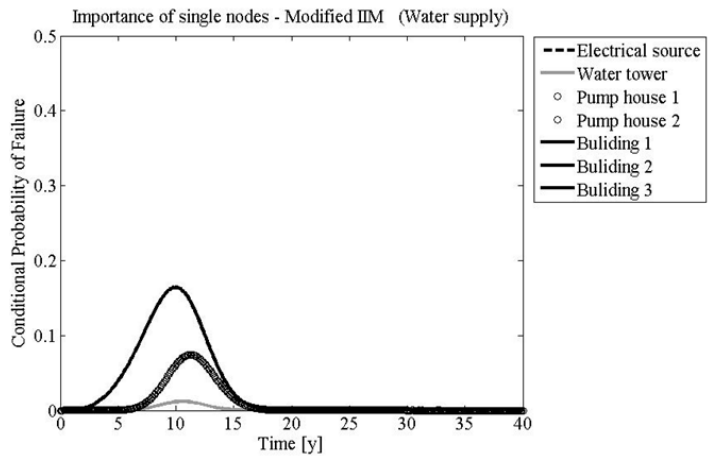

**Figure 8.20** Importance of nodes for the Example 2 network water supply according to the modified IIM.

### 8.2.8 Probability Risk Assessment

After having introduced the modified IIM for lifeline networks, next we compare it to the PRA method, which regards specific critical sites. The aim is to determine if the modified IIM can model both regional-scale and local-scale networks. The PRA is a systematic and comprehensive methodology to evaluate risks associated with every life-cycle aspect of a complex engineered technological entity (e.g., power plant, facility, or spacecraft) from concept definition, through design, construction and operation, and up to removal from service. In a PRA, risk is characterized by two quantities: (1) the magnitude or severity of the adverse consequences that can potentially result from the given activity or action; and (2) by the likelihood of occurrence of the given adverse consequences. If the measure of the severity of the consequences includes the potential for the number of people injured or killed, risk assessment becomes a powerful analytic tool to assess safety performance. A PRA usually answers three basic questions:

1. What can go wrong with the studied technological entity, or what are the initiators or undesirable initiating events that lead to adverse consequences?

2. What and how severe are the potential detriments, or the adverse consequences that the technological entity may be eventually subjected to as a result of the occurrence of the initiator?

3. How likely to occur are these undesirable consequences, or what are their probabilities or frequencies?

The answer to the first question requires technical knowledge of the possible causes leading to detrimental outcomes of a given activity or action. Probabilistic risk assessment studies can be performed for internal initiating events as well as for external initiating events. Here, internal initiating events are defined to be hardware or system failures or operator errors in situations arising from the normal mode of operation of the facility. External initiating events are those encountered outside the domain of the normal operation of a facility. Initiating events associated with the occurrence of natural phenomena (e.g. earthquake, storms, etc.) are typical examples of external initiators.

The answers to the second and third questions are obtained by developing and quantifying accident scenarios, which are chains of events that link the initiator to the end-point detrimental consequences. Focusing on the third question, the answer is obtained by using Boolean logic methods for model development and by probabilistic or statistical methods for quantification of the model. Boolean logic tools include inductive logic methods like event tree analysis (ETA) and deductive methods like fault tree analysis (FTA). It is easy to confuse these two techniques. In fact, the two are complimentary and are often used together, but each technique focuses on opposite sides of an undesired event. Figure 8.21 shows how they fit together. A more comprehensive description of these methods is discussed below.

Probabilistic risk assessment studies require special but often very important analysis tools like human reliability analysis (HRA) and dependent-failure or common-cause analysis (CCF). Human reliability analysis models human error while CCF evaluates the effect of inter-system and inter-component dependencies that tend to cause significant increases in overall system or facility risk. The final result of a PRA is given in the form of a risk curve and the associated uncertainties, which is generally the plot of the frequency of exceeding a consequence value as a function of the consequence values [ICAO 2014].

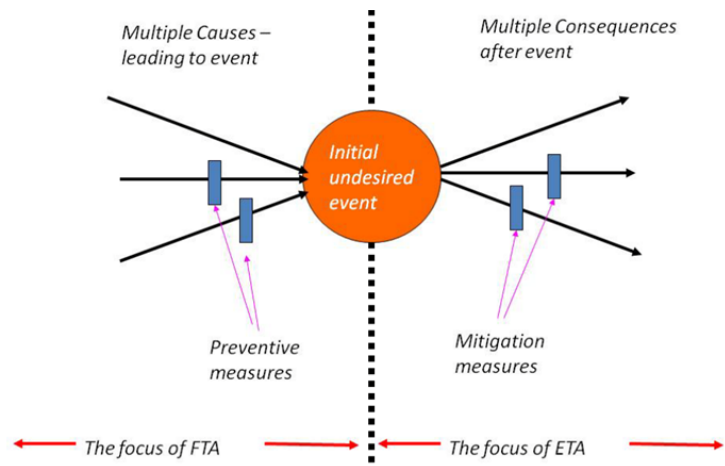

**Figure 8.21** Looking at undesired event using failure tracing methods (source: ICAO).

### *8.2.8.1 Fault-Tree Analysis*

The fault tree is a logic diagram based on the principle of multi-causality, which traces all branches of events that could contribute to an accident or failure. A fault-tree diagram is drawn from the top down. The starting point is the undesired event of interest, which is referred to as the "top event." The process consists in determining in sequence the immediate contributory fault conditions leading to that event. These may each in turn be caused by other faults and so on. The difficulty in constructing a fault-tree diagram is determining the correct sequence of failure dependencies [NASA 2014].

### *8.2.8.2 Event-Tree Analysis*

This is a complimentary technique to fault tree analysis (FTA) but defines the consequential events that flow from the primary 'initiating' event. Event trees are used to investigate the consequences of loss-making events in order to investigate ways of mitigating rather than preventing losses. The process for constructing an event tree analysis (ETA) is as follows:

1. Identify the primary event of concern.
2. Identify the controls that are assigned to deal with the primary event such as automatic safety systems, alarms on operator actions.
3. Construct the event tree beginning with the initiating event and proceeding through failures of the safety functions.
4. Establish the resulting accident sequences.
5. Identify the critical failures that need to be addressed.

There are a number of ways to construct an event tree. They typically use Boolean logic gates (i.e., a gate that has only two options, such as success/failure, yes/no, on/off). They tend to start on the left with the initiating event and progress to the right, branching progressively. Each branching point is called a node. Simple event trees tend to be presented at a system level and are not detailed [NASA 2014].

### 8.2.8.3 *An Example of a Probabilistic Risk Assessment*

To clarify what is a PRA and how FTA and ETA work, a simple numerical example is presented, referred to as Example 4. The aim is to determine the frequency over the course of a year of being late at work because of oversleeping. It is possible to construct a simple event tree model by defining an initiating event (i.e., it is a work day) and mitigating systems (i.e., an alarm clock and a backup person); see Figure 8.22. After having defined the initiating event frequency, the model is solved by the determination of branch probabilities, which may require constructing a fault tree; see Figure 8.23. As a convention for the ETA, the upper branches are considered that success has been achieved (green probabilities), while lower branches (the red probabilities) indicate an unsuccessful outcome.

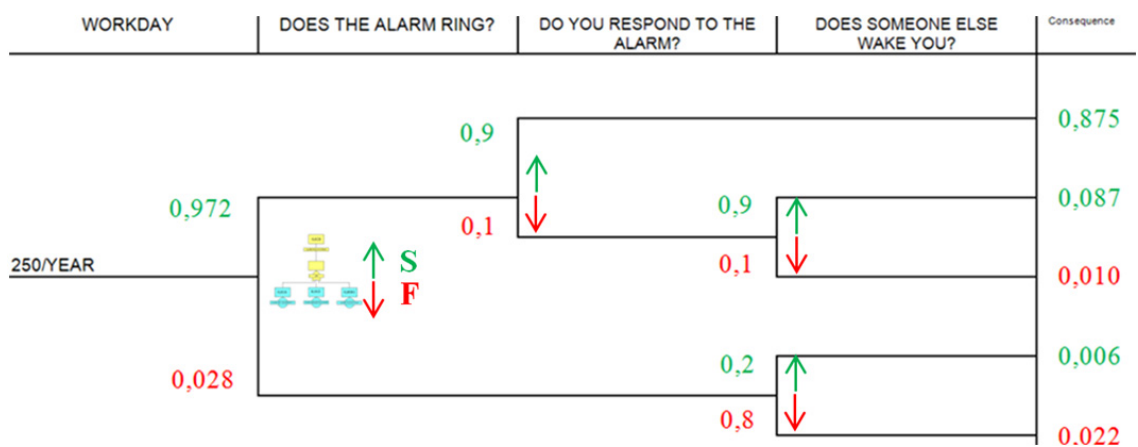

**Figure 8.22** Event tree analysis of Example 4 to determine the probability of being late at work because of oversleeping over a year-long period (source: U.S. NRC).

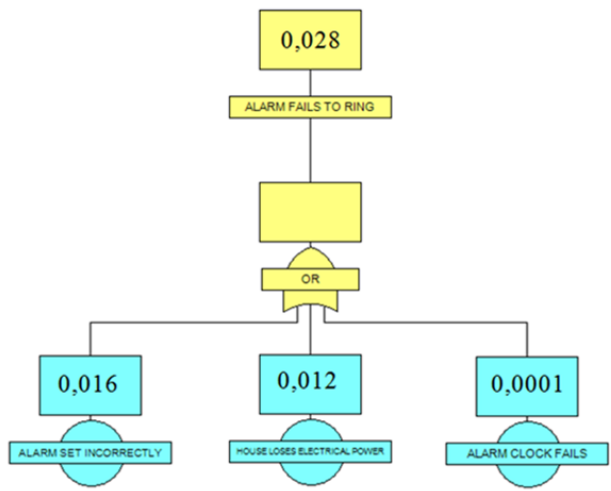

**Figure 8.23** Fault tree analysis of the probability of failure of the alarm clock over a year-long period. The results of this analysis are used by the ETA of Figure 8.22 (source: U.S. NRC).

### 8.2.8.4 Limitations of PRA

In general, PRAs suffers from analytical limitations. The three main problems are:

1. It cannot account for the indirect, nonlinear, and feedback relationships that characterize many accidents in complex systems;
2. It does a poor job of modeling human behavior and their impact on known, let alone unknown, failure modes;
3. It is conceptually impossible to be complete in a mathematical sense in the construction of event trees and fault trees.

The modified IIM suggested in the previous subsection overcomes problems related to feedback relationships and rigorous mathematical formulation because it allows loops of interdependence and uses analytical relations described by the IIM constitutive equation. In the case of modeling human behavior, agent-based models can furnish more reliable simulations. In conclusion, a modified IIM for temporal networks that incorporates agent-based models to simulate human behavior is believed to be superior in determining the risk related to a system. Next, a comparison between a PRA and the modified IIM is presented.

### 8.2.8.5 Comparison with the Modified IIM

Even if inoperability and risk are similar measures to assess the performance of a system, no correlation between the IIM and PRA has ever been done. As already shown, the modified IIM for temporal networks can model infrastructure networks at a regional scale. It is also possible to apply it to analyze infrastructure networks at the local level.

One of the main advantages of the modified IIM is its ability to take into account the positive effect of redundancies in the system. Redundancies are computed in both the PRA and the modified IIM through the logic operator "OR." Figure 8.24 shows how the tensor notation of the modified IIM is useful to this application.

Figure 8.25 compares the numerical results regarding Example 4. The example of the clock failure presented in Figure 8.23 is replicated with a simple 3-node network. A probability of self-failure, $P_{sf}$, is applied to each node equal to the frequency of failure of the oversleeping example. With the IIM algorithm, inoperability is propagated in the network, and the result is given by the probability of failure $P_f$. The result obtained with the two methods is identical. In conclusion, the combination of events and cascading effects done by the modified IIM follows the same logic approach of FTA and give the same results.

Finding a similar correlation between the modified IIM and the ETA is more difficult. The concept of sequences of events was not addressed by the traditional IIM and was introduced in the proposed implementation with the tensor formalism. The sequences referred to in the modified IIM are the occurrences of different configurations of the network. If an event tree refers to sequences of event that cannot be identified by the occurrence of a configuration, the modified IIM cannot obtain the same results if it starts with the same data. In conclusion, there are some event sequences that can be simulated through the suggested model and others that cannot.

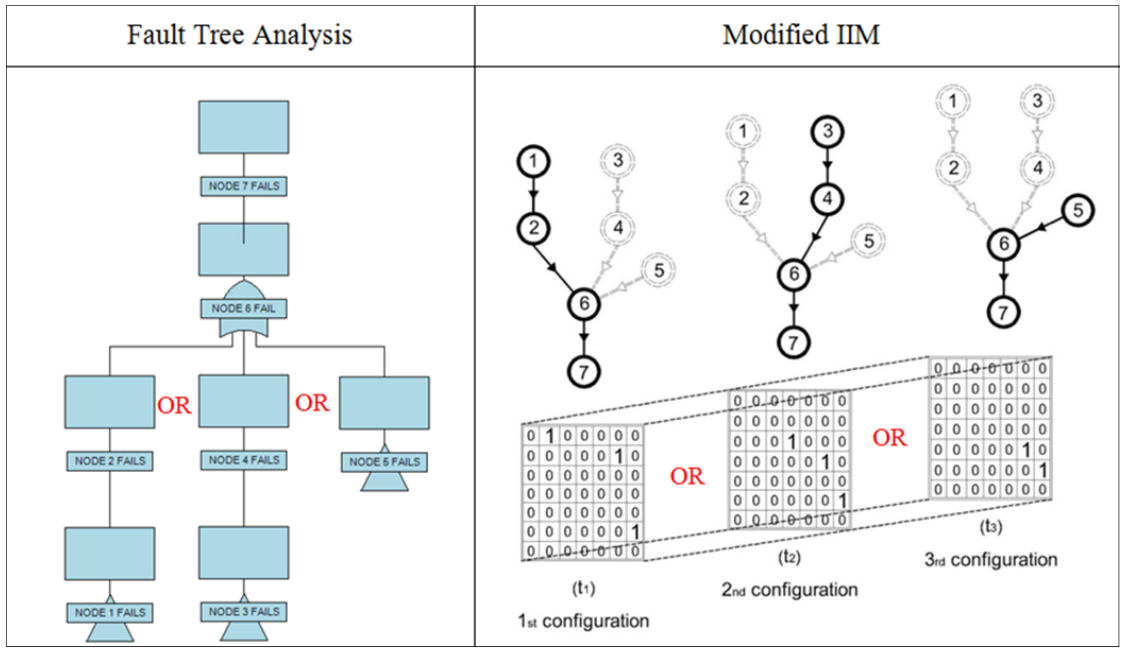

**Figure 8.24**   Structure to obtain an "OR" operator with FTA and modified IIM.

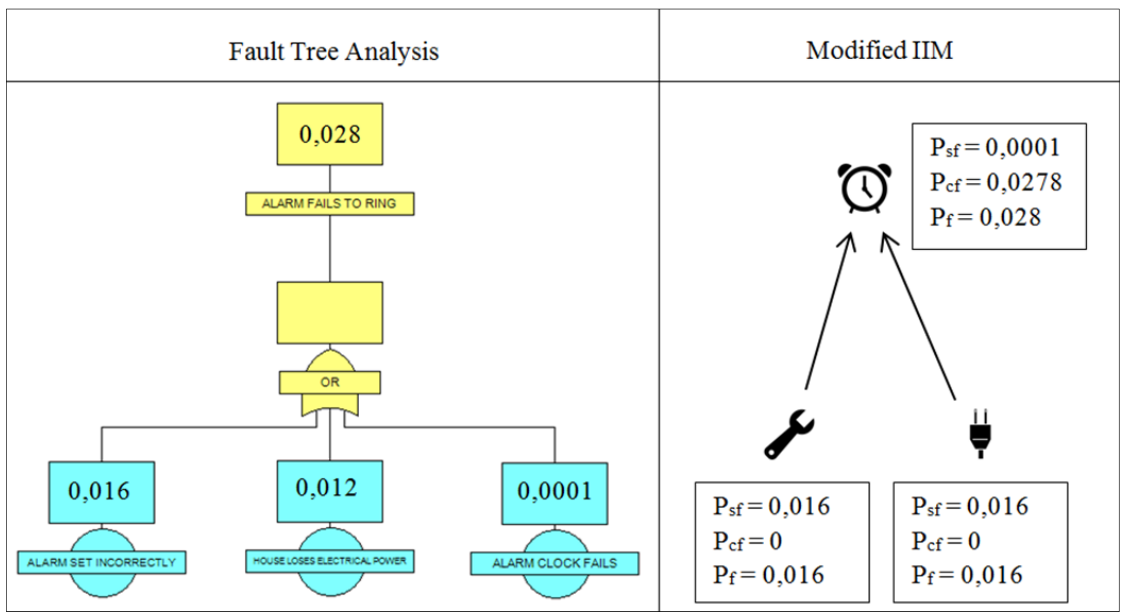

**Figure 8.25**   Comparison between the FTA and the modified IIM for Example 4.

If we look at the example of the oversleeping risk assessment and model it using the network of Figure 8.26, the events (see Figure 8.22) are not the success/failure of the three configurations of Figure 8.26, but the success/failure of single nodes. A modified IIM can compute conditioned probability of occurrence of every configuration but not the conditioned probability of single nodes.

If the event tree is structured in a way that every temporal sequence refers to the success/failure of a configuration of the system, the results obtained with ETA and modified IIM are the same. For Example 3, Figure 8.27 shows how the $P_{occ}$ and the $LoC$ are the same for both methods. An event sequence is acceptable if it is mutually exclusive from the others and if it represents the complete flow form the source to the sink, i.e., per the modified IIM. Because the probabilities present in the event tree of Example 4 are conditioned, they cannot be used as input in the modified IIM.

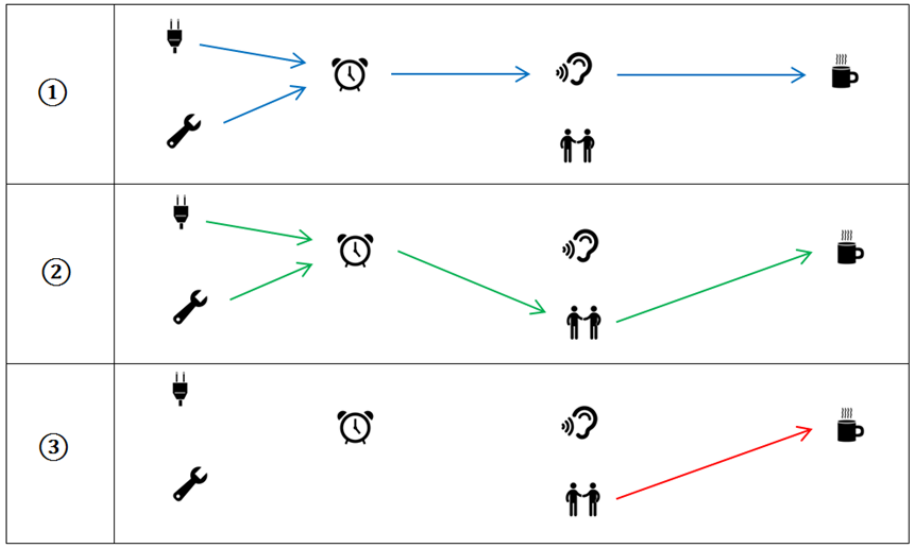

Figure 8.26    Different configurations of a network simulating Example 4.

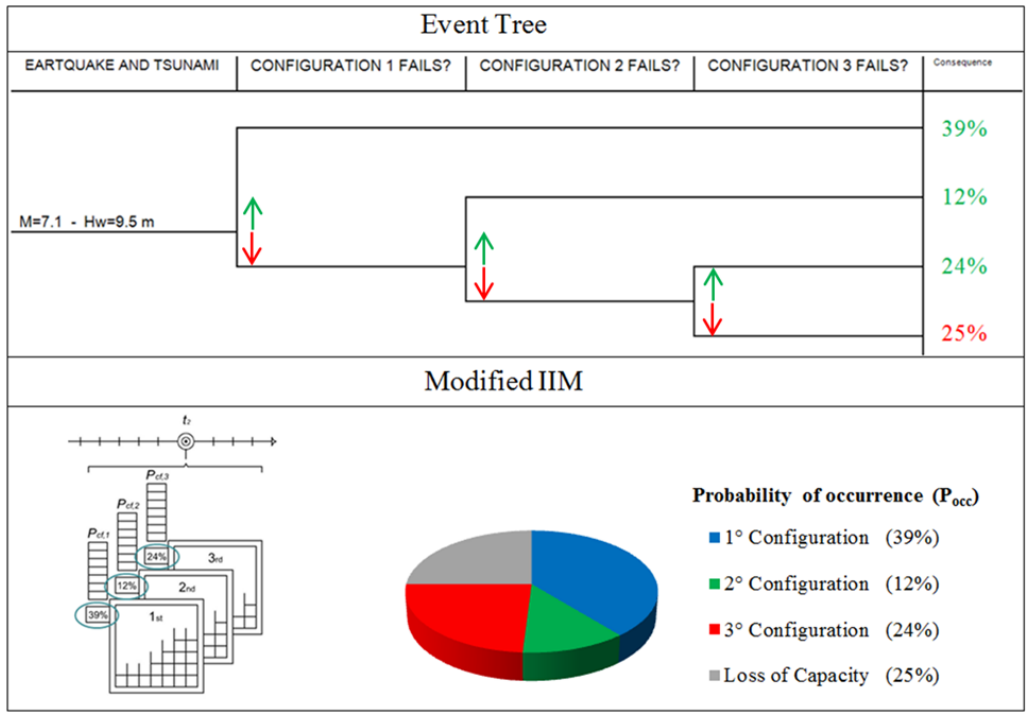

Figure 8.27    Comparison between the ETA and the modified IIM for Example 3.

## 8.3 CASE STUDY: LIFELINES SERVING A NUCLEAR POWER PLANT

Next we focus on lifeline systems serving critical sites by applying the proposed methodology to a nuclear power plant (NPP). Nuclear power plants are dependent on extended regional scale infrastructures, which can be analyzed using an IIM, but at the same time are strategic sites and have service plants at the local scale, which can be analyzed with a PRA. It will be shown that the modified IIM model can perform both a regional- and local-level risk assessment and that a separate PRA is not necessary.

This case study will use the 2011 Fukushima NNP disaster to illustrate the method. This disaster is one of the most complete examples of failure due to interdependence and temporal effects. After a brief presentation of the real event, a description and calibration of models adopted follows. Once the results of this analysis are presented, the proposed models are compared and followed by some general remarks.

### 8.3.1 The 2011 Fukushima Nuclear Power Plant Disaster

The 2011 Fukushima Daiichi disaster that followed the Tohoku, Japan, earthquake and tsunami was a serious indictment of safety engineering up to that point. The complexity of the events and of the system were not modeled accurately by risk planners; the cascading effects resulted in major damage to the NNP.

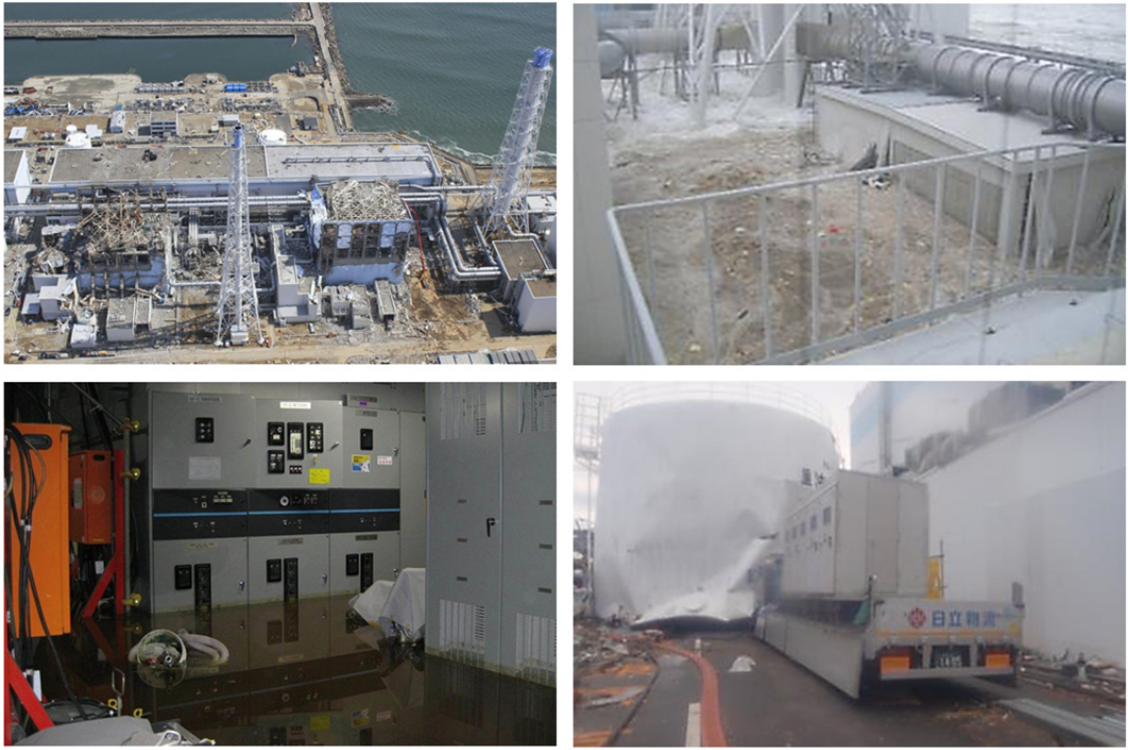

**Figure 8.28**   Photographs of the Fukushima nuclear power plant after the 2011 Tōhoku earthquake and tsunami.

*8.3.1.1 Sequence of the Accident*

The huge earthquake (Figure 8.29) and tsunami (Figure 8.30) that struck Japan's Fukushima Daiichi nuclear power station on March 11, 2011, knocked out backup power systems that cooled the reactors at the plant, which caused three of them to undergo fuel melting, hydrogen explosions, and radioactive releases. Although radioactive contamination from the Fukushima plant forced the evacuation of communities up to 25 miles away and affected up to 100,000 residents, it did not cause any immediate deaths.

The Tokyo Electric Power Company (TEPCO) operates the Fukushima nuclear power complex in the Futaba district of Fukushima prefecture in northern Japan, which consists of six nuclear units at the Fukushima Daiichi station and four nuclear units at the Fukushima Daini station. All the units at the Fukushima complex are boiling water reactors, with reactors 1 to 5 at the Fukushima Daiichi site being of General Electric Mark I design; this design has been used for NNPs located in the U.S. The Fukushima Daiichi reactors entered commercial operation from 1971 (reactor 1) to 1979 (reactor 6). When the earthquake struck, Fukushima Daiichi units 1, 2, and 3 were generating electricity and shut down automatically. The earthquake caused offsite power supplies to be lost, and backup diesel generators started up as designed to supply backup power. However, the subsequent tsunami flooded the electrical switch-gear for the diesel generators, causing most AC power in units 1 to 4 to be lost. Because Unit 4 was undergoing a maintenance shutdown, all of its nuclear fuel had been removed and placed in the unit's spent fuel storage pool. One generator continued operating to cool units 5 and 6.

The loss of all AC power in Units 1 to 3 prevented valves and pumps from removing heat and pressure being generated by the radioactive decay of the nuclear fuel in the reactor cores. As the fuel rods in the reactor cores overheated, they reacted with steam to produce large amounts of hydrogen, which escaped into Units 1, 3, and 4 reactor buildings and exploded (the hydrogen that exploded in Unit 4 is believed to have come from Unit 3). The explosions interfered with efforts by plant workers to restore cooling and helped spread radioactivity. Cooling was also lost in the reactors' spent fuel pools, although recent analysis has found that no significant overheating took place.

Radioactive material released into the atmosphere produced extremely high doses of radiation near the plant and left large areas of land uninhabitable, especially to the northwest of the plant. Contaminated water from the plant was discharged into the sea, creating international controversy [CRS 2012]. A complete timeline of all the events occurring at Unit 1 is presented in Figure 8.31.

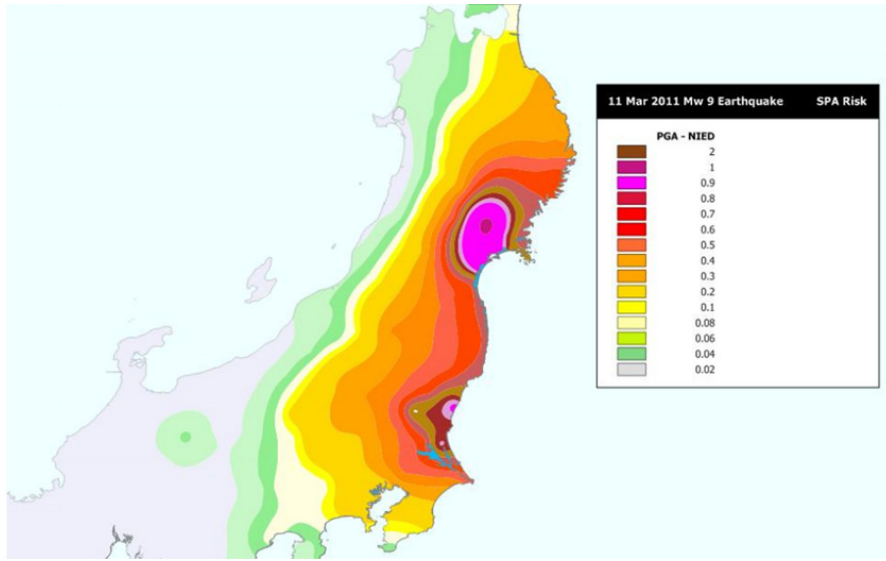

**Figure 8.29** Shake map of the Eastern Japan Coast after Tohoku earthquake (source: Scawthorn [2011]).

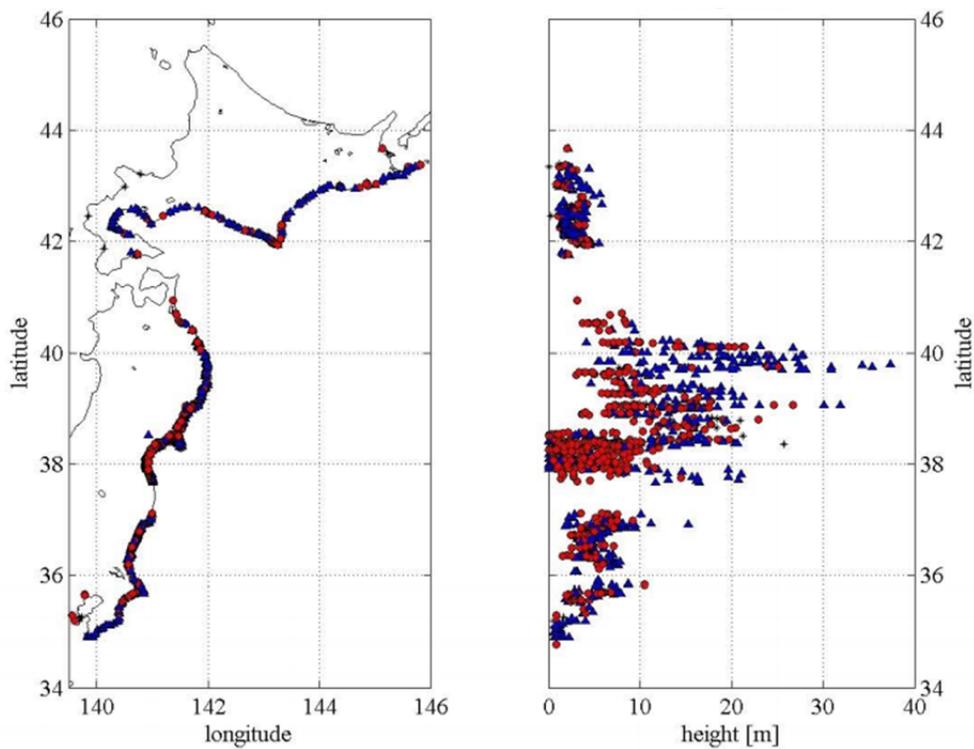

**Figure 8.30** Inundation map of the Eastern Japan Coast after Tohoku tsunami (source: Scawthorn [2011]).

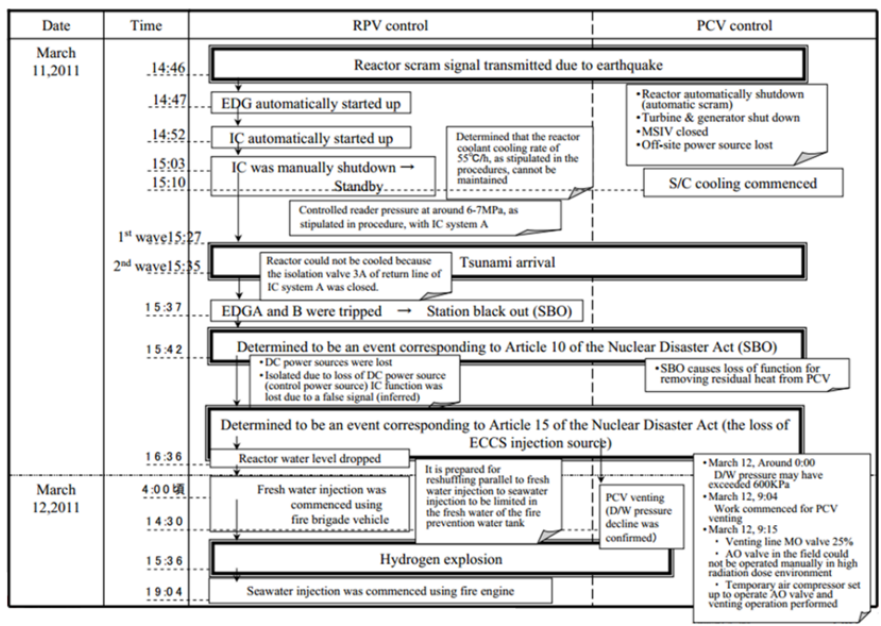

**Figure 8.31** Timeline of events occurring at Unit 1 of Fukushima Daiichi nuclear power plant (source: TEPCO [2012]).

### 8.3.2 Modeling the Nuclear Power Plant

Replicating what happened at Fukushima requires constructing models of the lifeline networks serving a generic NNP. This work does not specifically model the Unit 1 of Fukushima Daiichi NPP because data regarding this case study is not available. The aim is to model a NPP equipped with a Boiling Water Reactor (BWR) loosely based on Unit 1 at Fukushima. The topology and data regarding disrupting events affecting the system are like those of the Fukushima disaster, but parameters of the component of the system are generic and taken from literature.

#### 8.3.2.1 Topology

The proper functioning of a NNP requires substantial infrastructure whose complexity is difficult to convey to people not in the industry. Focusing on the connectivity of these networks, it is possible to construct logic schemes while setting the topology of a model. The plant scheme furnished by TEPCO was used as a reference for building this model; see Figure 8.32.

This logic scheme considers the electric network, the water network, and the steam network. For the purposes of this analysis, the steam network has not been considered by itself, but it has been included with the water network. Cooling circuits are considered closed to prohibit the spread of radioactive substances. These loops have been modeled with one-direction links from the source to the reactor core. Apart from the water network, which is denoted as being part of the local/building scale, there is also the electric network, which expands from the regional scale to the local one. Although this lifeline is more important than the cooling plant, this case study made no distinction. Because the task is to run a performance analysis of all systems serving the reactor core, all components considered important for the success or failure of the reactor cooling have been modeled by nodes and connected following logical assumptions. Two different models will be presented. The first is a simplified version of the scheme of Figure 8.32; the second model is more detailed and is integrate to the physical infrastructures and emergency responders' networks.

The simplified model shown in Figure 8.33 is composed of an electric and a water network. In order of priority, the sources of the electric network are the external electric network, diesel generators, and DC batteries. All these possible configurations converge into a power panel that then feeds the pumps of the ordinary water network. The source of this water network is the ocean, which is considered to have unlimited autonomy as does the external electric network. The first emergency cooling system is located in the Isolation Condenser (IC), which cools the steam coming from the reactor in a pool that doesn't need electricity because the flow is gravity-driven. The high-pressure coolant injection (HPCI) system can cool the core in an emergency condition. It draws water from the condensate storage tank (CST) or from the suppression pool (SP), and injects it into the core form to reduce the internal pressure. Because the pump used by this system is steam-driven, it feeds automatically once the plant is started. Thus, we have three possible configurations for the electric network and four configurations for the water-cooling network.

The model above is too limited to apply to Fukushima Unit 1. Both the electric and the water networks connections are more complex than the one presented above. Starting from electric sources, the self-sustainment guaranteed by the NPP power plant is introduced, where all sources and relative paths feed particular target nodes but not all of them. For example, the ordinary cooling line is only fed by the NPP turbine and the off-site AC power, while diesel generators feed the residual heat removal (RHR) cooling system. The IC and HPCI systems, which were considered not dependent on electricity in the simplified model, are now indirectly dependent on it because their activation is performed by valves that can be remotely controlled only with a functioning electrical supply. The DC battery is responsible for functionality of these valves.

To better model possible human interventions on the system, three additional networks have been added: the telecommunication network, the transportation network, and the emergency service network. Some nodes and edges of these networks may not be active for certain periods of time. All the layers of these new models are interdependent as shown in Figure 8.34. The connections among nodes can be present at the regional scale (Figure 8.35), at the power plant scale (Figure 8.36), and at the reactor scale (Figure 8.37).

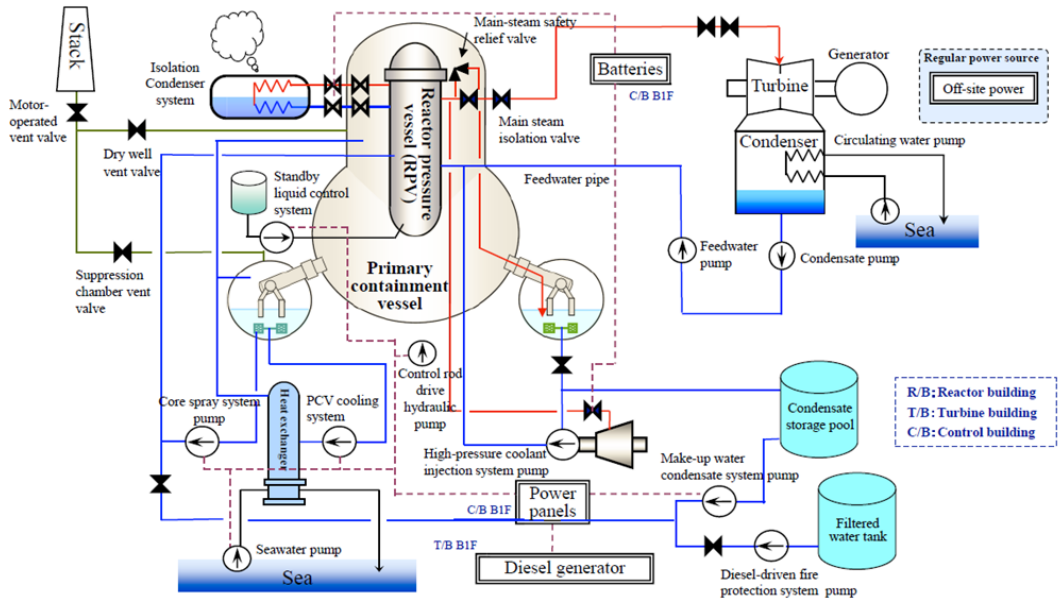

**Figure 8.32** Scheme of the Unit 1 Fukushima Daiichi reactor (source: TEPCO [2012]).

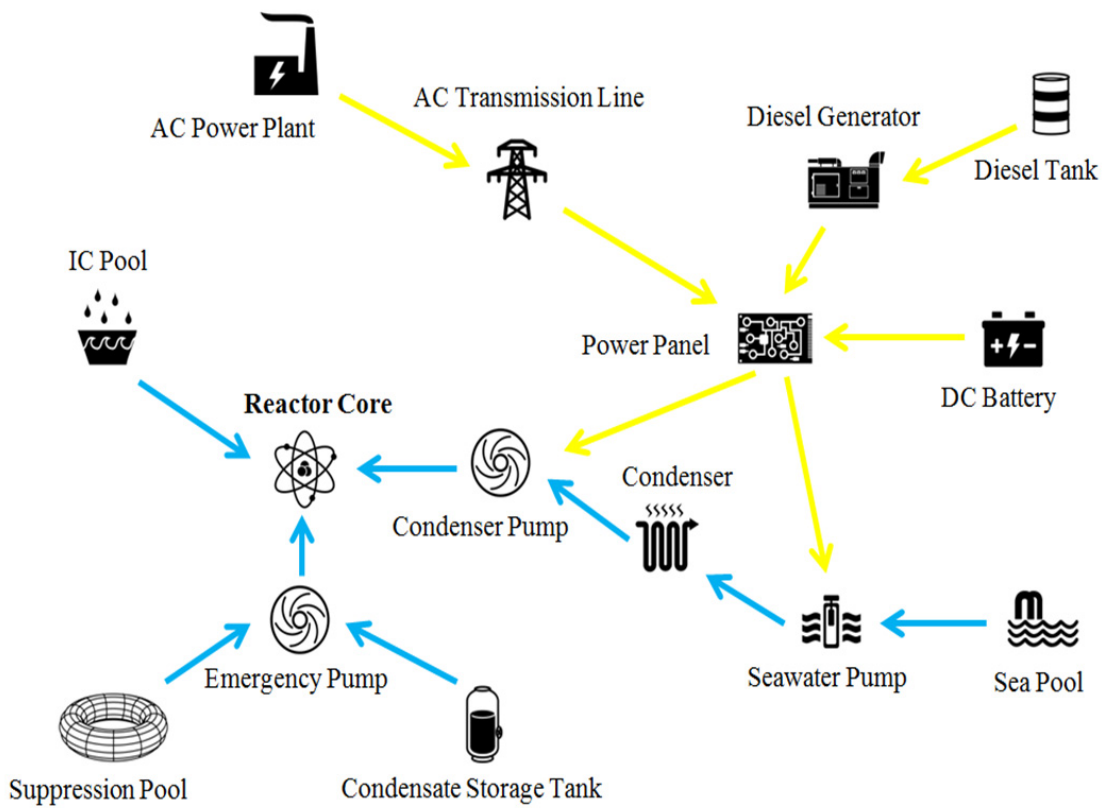

**Figure 8.33** Simplified model for lifelines serving a nuclear power plant.

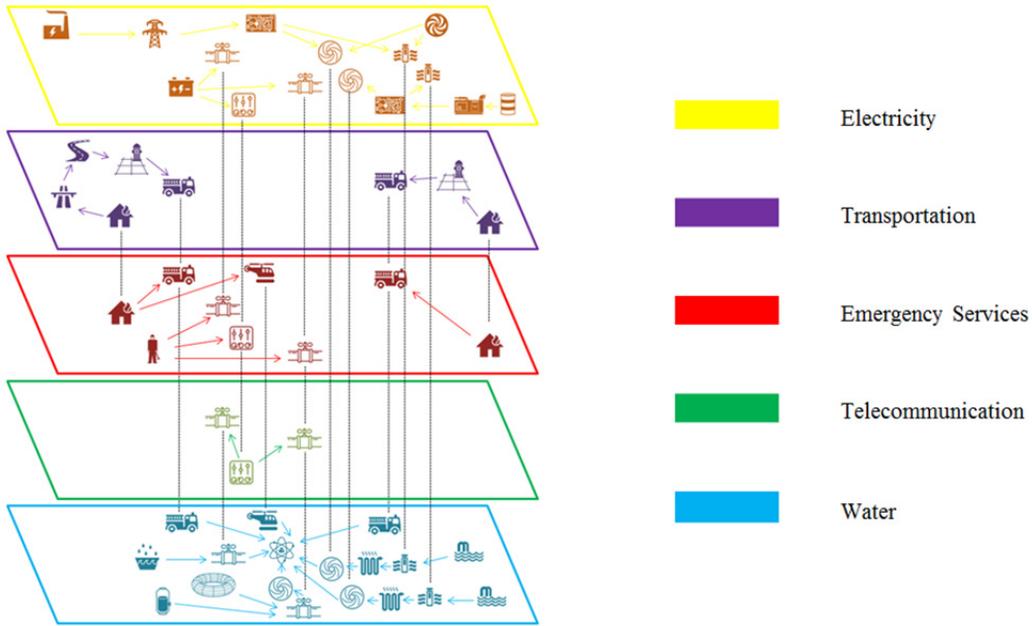

**Figure 8.34** Interdependent layers of the detailed model for lifelines serving a nuclear power plant.

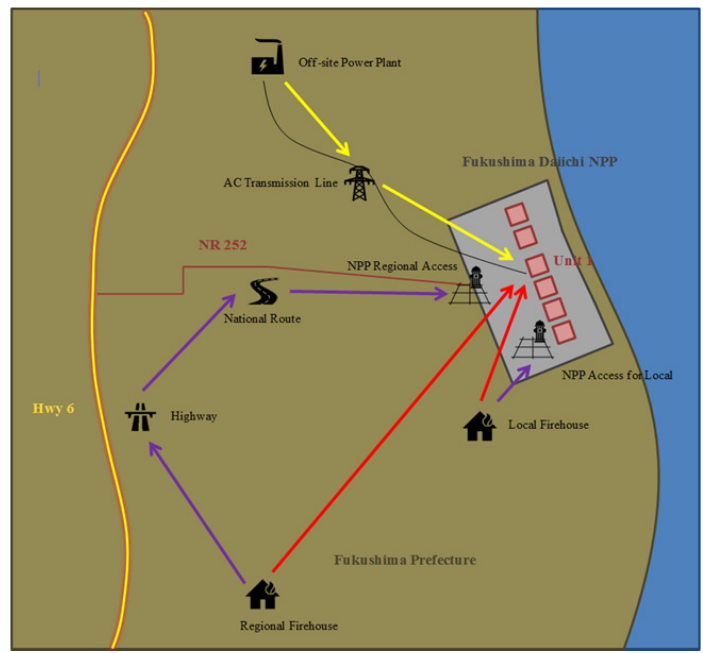

**Figure 8.35** Detailed model for lifelines serving a nuclear power plant at the regional scale.

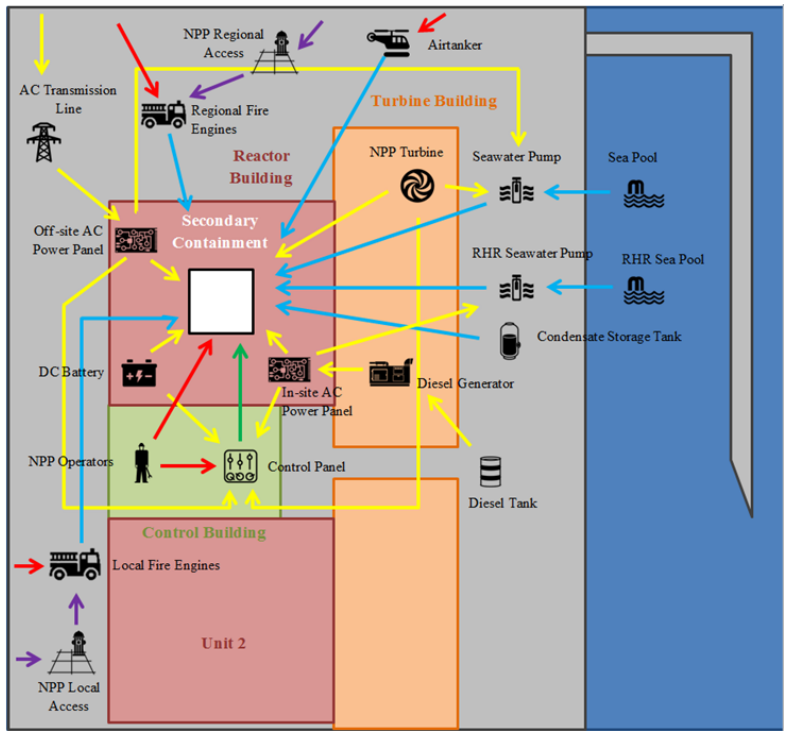

**Figure 8.36** Detailed model for lifelines serving a nuclear power plant at the power plant scale.

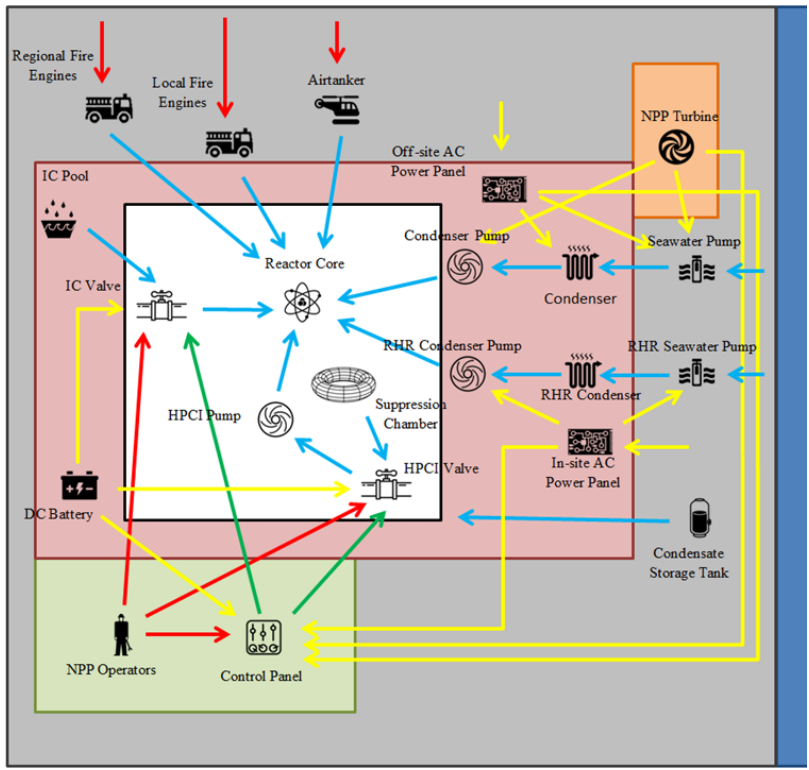

**Figure 8.37** Detailed model for lifelines serving a nuclear power plant at the reactor scale.

*8.3.2.2 Hazards*

The hazards considered for the analysis are an earthquake and the tsunami. Intensities of these hazards are considered deterministic and taken from the real event. To fill in Table 8.1, nodes were classified according to their location and altitude. Figure 8.38 and Figure 8.39 show a plan and cross-section of the NNP used for this analysis. It is then possible to estimate through data obtained from shake and flow maps what is the intensity of the event for the considered node. To assemble the event matrix, *E*-vectors should be positioned at the correct time step. Time steps are defined in accordance with the timeline present in Figure 8.31.

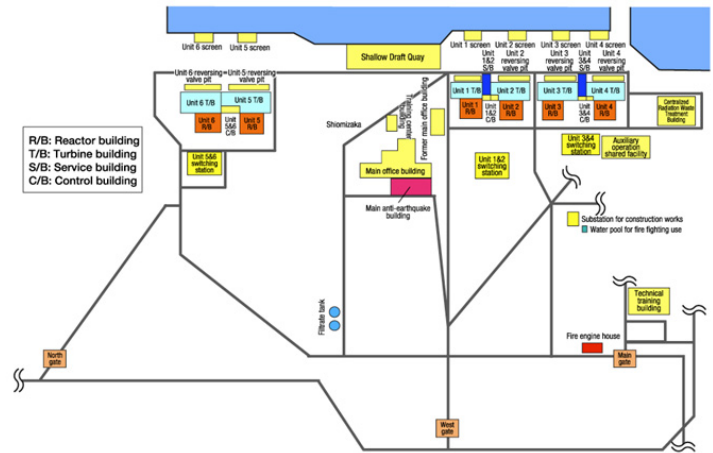

**Figure 8.38**  Location of facilities at the Fukushima Nuclear Power Plant (source: TEPCO [2012]).

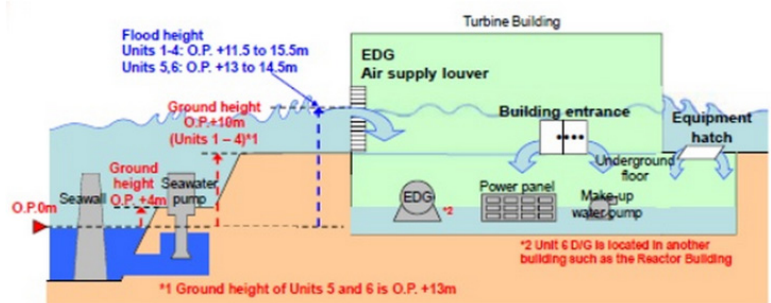

**Figure 8.39**  Path of inundation for the Fukushima Nuclear Power Plant (source: Ansaldo Nucleare [2013]).

**Table 8.1**     Intensity of hazards for the nuclear power plant models.

| Node | Location | Altitude (m) | Earthquake PGA (g) | Tsunami depth of water (m) |
|---|---|---|---|---|
| NPP Turbine | Turbine Building | 10 | 0,469 | 6 |
| AC Power Plant | Hinterland | >50 | 0,415 | - |
| AC Line | Hinterland | >50 | 0,415 | - |
| Off-site AC Power Panel | Turbine Building | 10 | 0,469 | 6 |
| Diesel Tank | NPP Apron | 10 | 0,469 | 3 |
| Diesel Generator | Reactor Building | 10 | 0,469 | 9 |
| In-site AC Power Panel | Reactor Building | 10 | 0,469 | 9 |
| DC Battery | Reactor Building | 10 | 0,469 | 9 |
| Highway | Hinterland | >50 | 0,415 | - |
| Road | Hinterland | >50 | 0,415 | - |
| Local Firehouse | Hinterland | >50 | 0,469 | - |
| NPP Local Access | NPP Apron | 10 | 0,469 | 3 |
| Local Fire Engines | - | - | 0,469 | - |
| Regional Firehouses | Hinterland | >50 | 0,415 | - |
| NPP Regional Access | NPP Apron | 10 | 0,469 | 3 |
| Regional Fire Engines | - | - | 0,415 | - |
| Airtanker | - | - | 0,415 | - |
| NPP Operators | - | - | 0,469 | - |
| Control Panel | Control Building | 10 | 0,469 | 6 |
| Sea Pool | Wharf | 4 | 0,469 | 5 |
| Seawater Pump | Wharf | 4 | 1 | 3 |
| Condenser | Reactor Building | 10 | 0,469 | 9 |
| Condenser Pump | Reactor Building | 10 | 1 | 9 |
| RHR Sea Pool | Wharf | 4 | 0,469 | 5 |
| RHR Seawater Pump | Wharf | 4 | 0,469 | 3 |
| RHR Condenser | Reactor Building | 10 | 0,469 | 9 |
| RHR Condenser Pump | Reactor Building | 10 | 0,469 | 9 |
| IC Pool | Reactor Building | 10 | 0,469 | 9 |
| IC Valve | Reactor Building | 10 | 0,469 | 9 |
| Condensate Storage Tank | NPP Apron | 10 | 0,469 | 3 |
| Suppression Pool | Reactor Building | 10 | 0,469 | 9 |
| HPCI Valve | Reactor Building | 10 | 0,469 | 9 |
| HPCI Pump | Reactor Building | 10 | 0,469 | 9 |
| PCV | PCV | 20 | 0,469 | 16 |

*8.3.2.3 Parameters*

In general, is difficult to find reliable data of sensitive structures; in the case of NPPs, it has proved even more difficult. The lack of data available for Fukushima necessitated analyzing a generic NNP instead of a specific one. That said, by evaluating different sources, credible parameters were assigned to each component. Most of the earthquake and tsunami fragility functions were taken from ATC-13 [1985]. These data are old and generic but are still broadly employed in the absence of more reliable and updated sources, and provide damage probability matrices for structural/non-structural components based on expert opinion. Given that the goal of the analysis is to evaluate the performance of the systems, some assumptions about the level of damage of the component that caused its inoperability were made. Other earthquake fragility curves were taken from the ALA report [2001] and from the HAZUS database [FEMA, 2013]. Tsunami fragility curves are considered linear functions between two values obtained from the ATC-13 recommendations and consideration about the robustness of buildings. Autonomy curves were estimated to be step functions, where the step is located in correspondence with the nominal value indicated by Hitachi-GE [2011]. The figures below illustrate the fragility and temporal effects curves for the nodes of the water network. Similar functions were adopted for all the other networks.

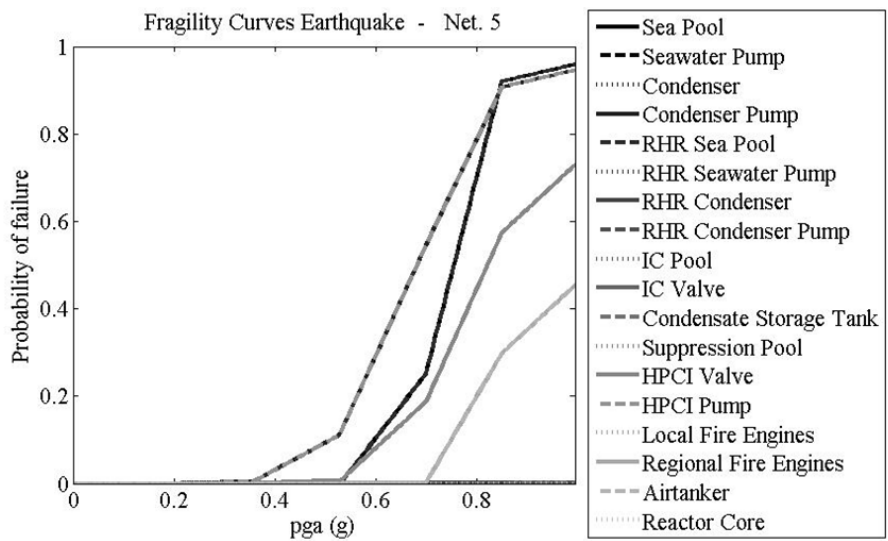

**Figure 8.40**     **Earthquake fragility curves for the nodes of the water cooling network.**

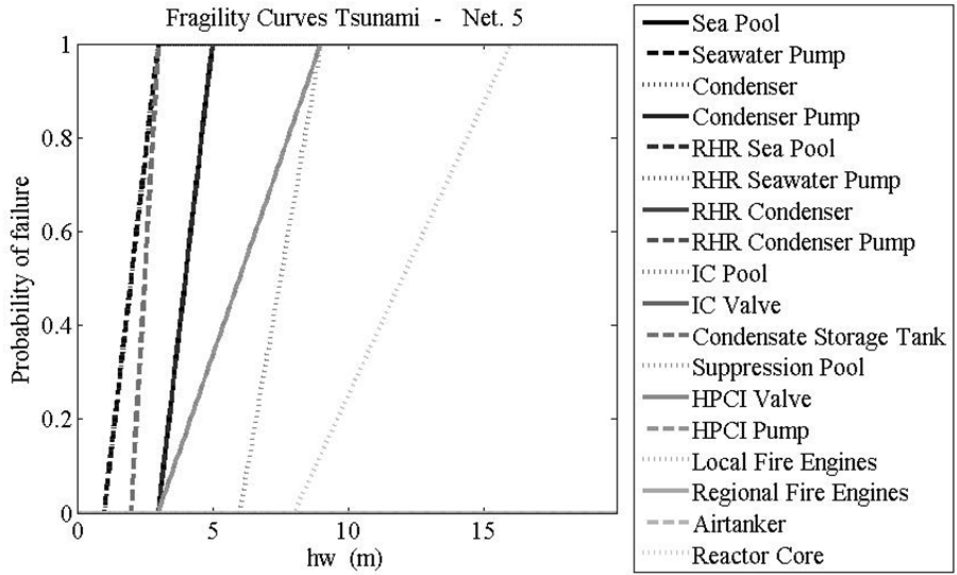

Figure 8.41  Tsunami fragility curves for the node of the water cooling network.

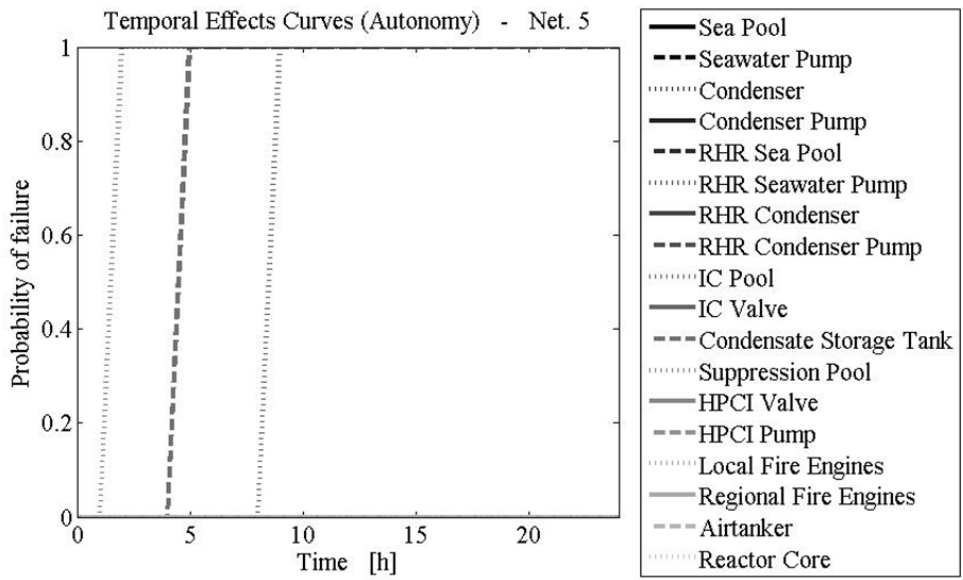

Figure 8.42  Autonomy curves for the nodes of water cooling network.

### 8.3.3  Analysis of the System

Having defined all the necessary input, it is possible to compute the various probability of failure of the systems. This section shows the results relative to the electric network and the water network of both the simplified model and the detailed model of the NPP. Analysis of other networks composing the detailed model are not reported because they are not directly comparable with the simplified model; however, they have been computed and influence the results of the water network presented. The results of the results of the simulations are presented next and compared with the real event that occurred in Japan.

*8.3.3.1 Simplified Model*

In the simplified model, the earthquake is responsible for the shutdown of other power plants and the collapse of the AC transmission line; this implies the loss of off-site AC power, which represents the first configuration of the electric network. The electric network changes configuration, and the power supply is guaranteed to the water network. Other components suffer little damage.

The arrival of the tsunami wave drastically changes the situation. Tanks placed in the NPP apron are swept away as well as the sea pumps. Diesel generators and the ordinary cooling line are out of order. Batteries are damaged too, but there is no need of them anymore since pumps they were feeding have failed. The IC backup cooling system, which doesn't need electricity because it is gravity-driven, is initiated and the probability of failure of the reactor core cooling is still close to 0.

Ten hours after the earthquake, the autonomy of the IC begins to decrease, and it is substituted by the HPCI system, which doesn't need electricity because it is equipped with a steam-driven pump. As the probability of failure of the IC increases because of the failure of autonomy run increases, the probability of failure of the reactor core increases as well because it now relies solely on the HPCI system, which was potentially damaged by earthquake and the tsunami; see Figure 8.46, Figure 8.51, and Figure 8.53.

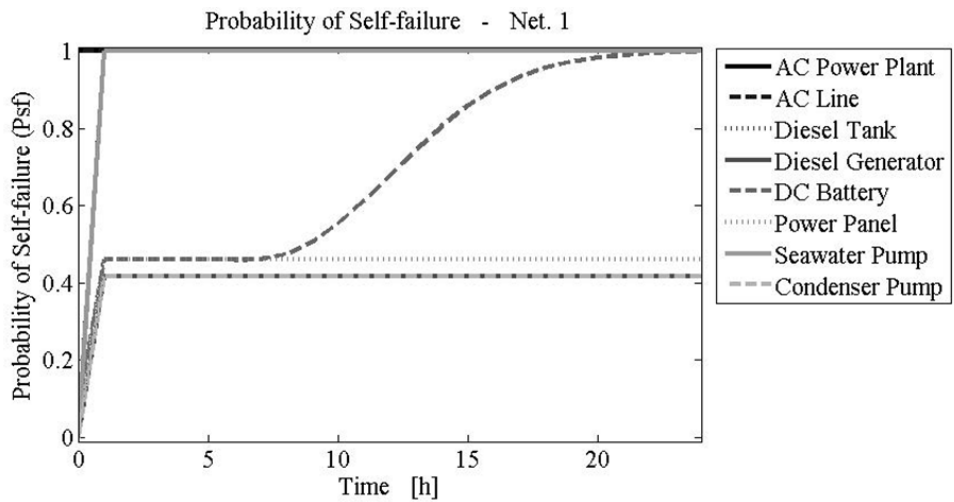

**Figure 8.43** Probability of self-failure of nodes of the electric network for the simplified model.

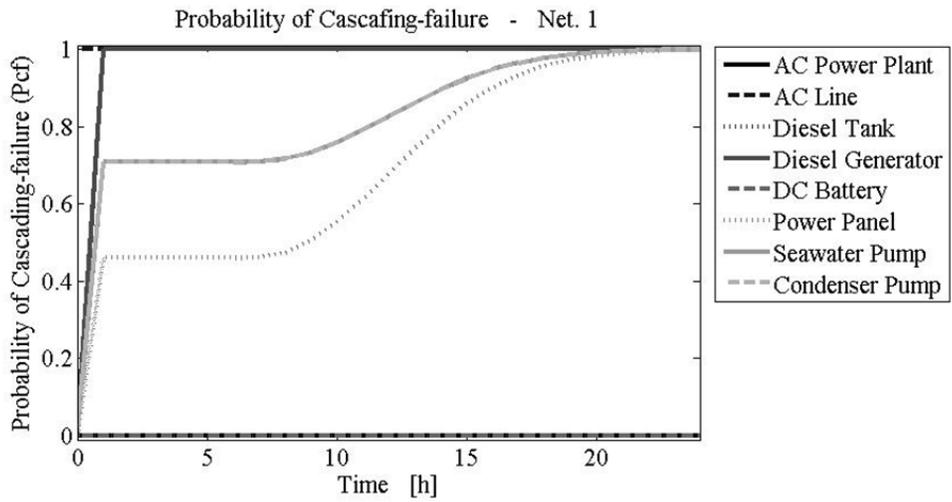

**Figure 8.44** Probability of cascading-failure of nodes of the electric network for the simplified model.

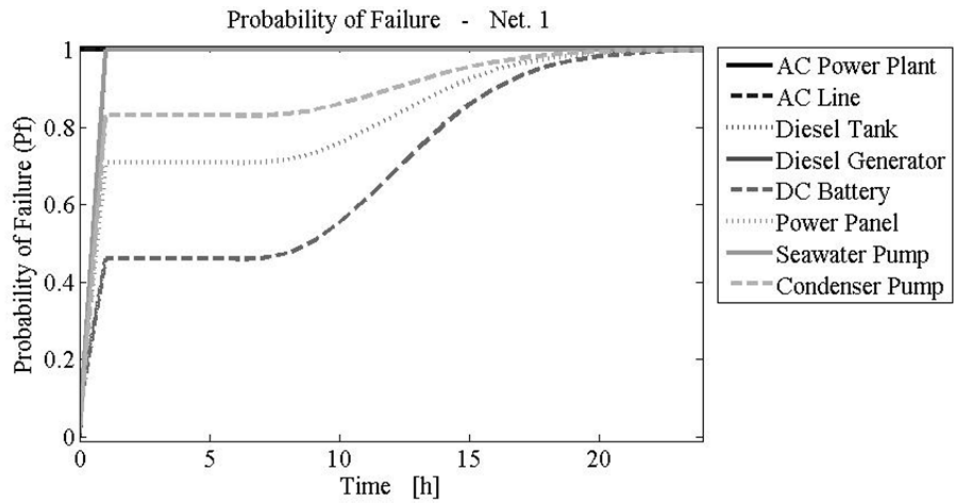

**Figure 8.45** Probability of failure of nodes of the electric network for the simplified model.

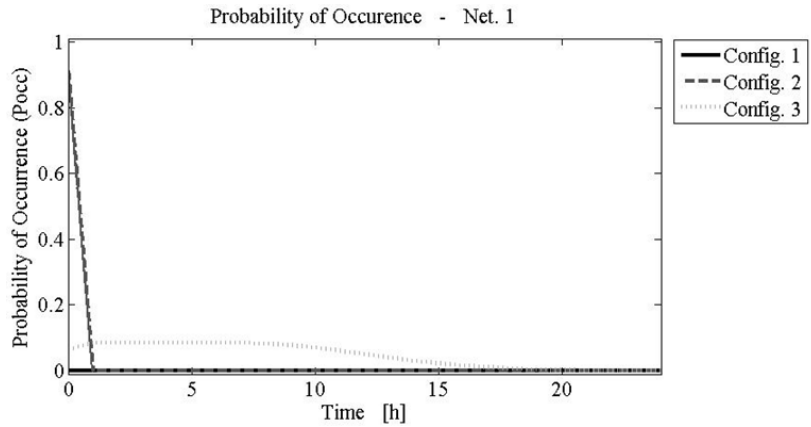

**Figure 8.46** Probability of occurrence of configurations of the electric network for the simplified model.

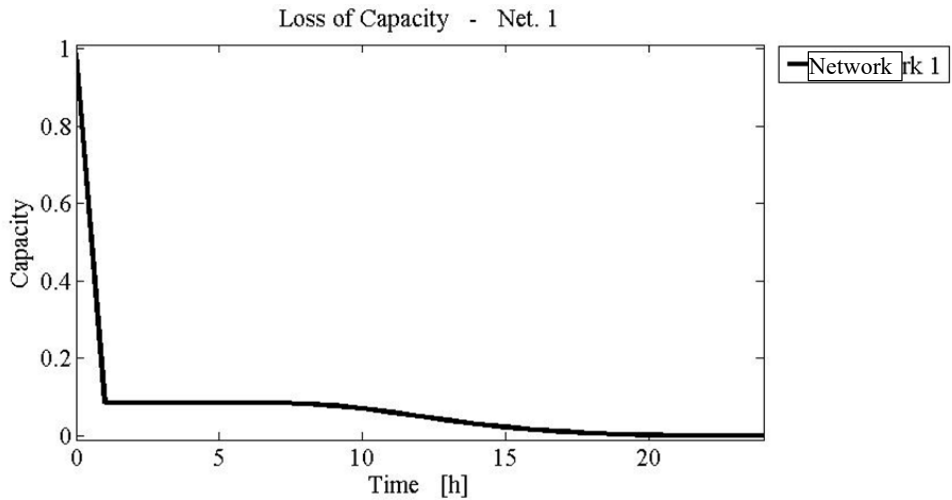

**Figure 8.47** Loss of capacity of the electric network, for the simplified model.

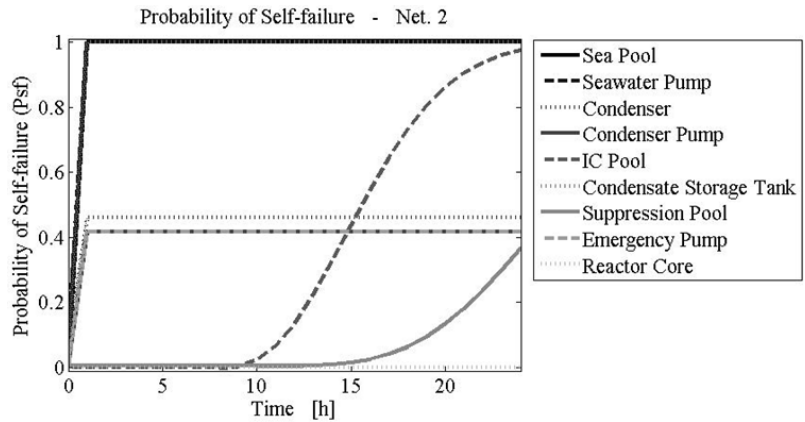

**Figure 8.48** Probability of self-failure of nodes of the water network for the simplified model.

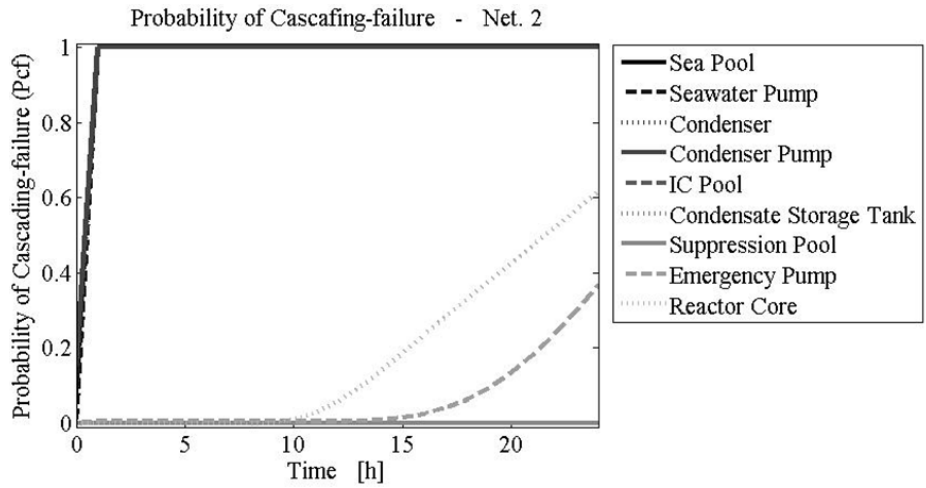

**Figure 8.49** Probability of cascading-failure of nodes of the water network for the simplified model.

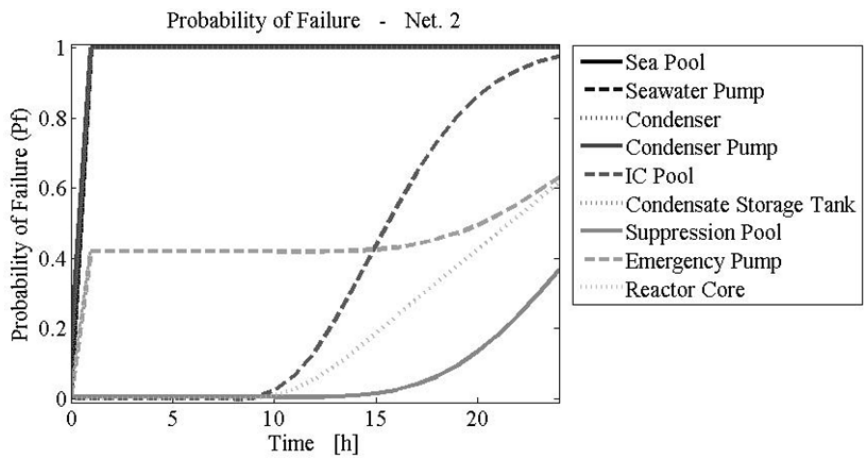

**Figure 8.50** Probability of failure of nodes of the water network for the simplified model.

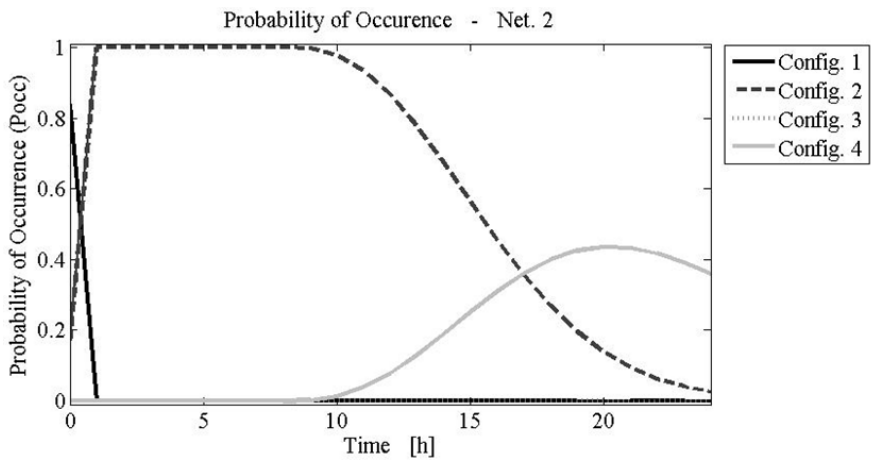

**Figure 8.51** Probability of occurrence of configurations of the water network, for the simplified model.

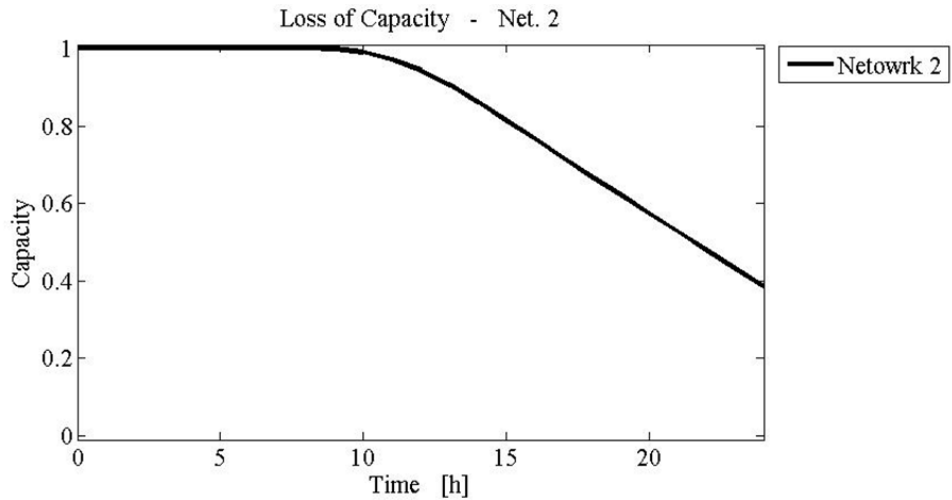

Figure 8.52   Loss of capacity of the water network for the simplified model.

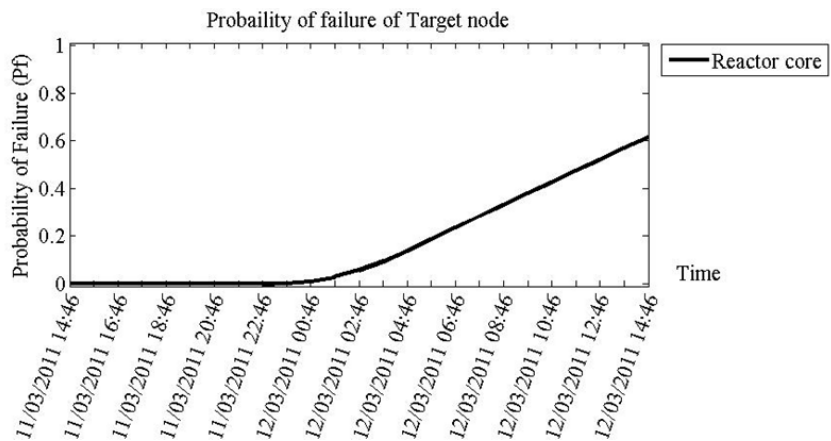

Figure 8.53   Probability of failure of the reactor core cooling for the simplified model.

#### 8.3.3.2  Detailed Regional- to Local-Scale Model

In the detailed model, the earthquake is also responsible for the shutdown of the NPP turbine and of the off-site AC power as AC transmission lines collapse. The loss of off-site AC power propagates the inoperability to the ordinary cooling configuration. Electricity is still provided by diesel generators that feed the RHR system and the control room. Damage caused by the earthquake to emergency cooling systems imply that although all of the first three backup lines have a probability of occurrence $P_{ooc} \neq 0$, the most likely backup to be activated is the RHR.

After the tsunami, diesel tanks, CSTs, and seawater pumps are completely damaged. The access to the NPP is not possible because of the debris deposited by the wave. Rescuers have difficulty in accessing the plant and need time to restore functional access. The water network tries to switch to the IC and HPCI configurations, but to control their valves, DC power is needed. There is a low probability that this is available because the batteries have a relevant probability of failure; therefore, the loss of capacity sharply increases. After three hours, an IC valve is manually opened, and cooling is provided by the IC until its autonomy runs out, resulting in a complete LoC; see Figure 8.64.

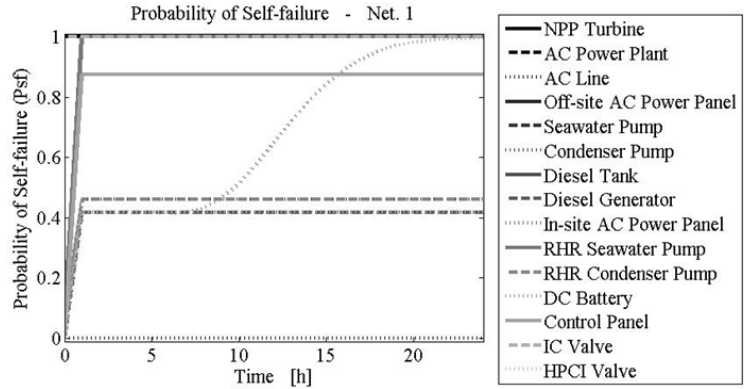

**Figure 8.54** **Probability of self-failure of nodes of the electric network for the detailed model.**

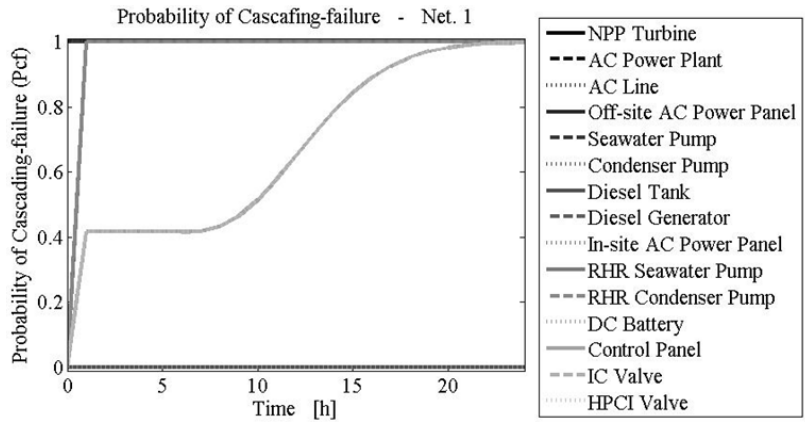

**Figure 8.55** **Probability of cascading failure of nodes of the electric network for the detailed model.**

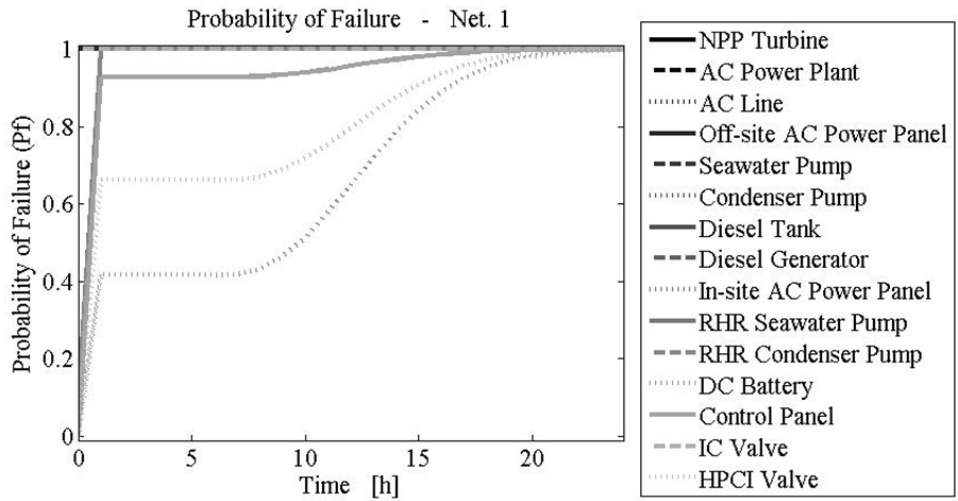

**Figure 8.56** **Probability of failure of nodes of the electric network for the detailed model.**

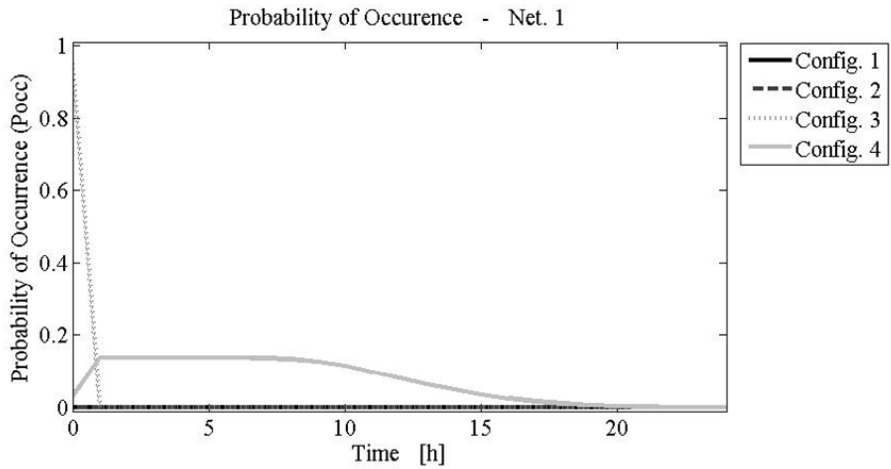

**Figure 8.57** Probability of occurrence of configurations of the electric network for the detailed model.

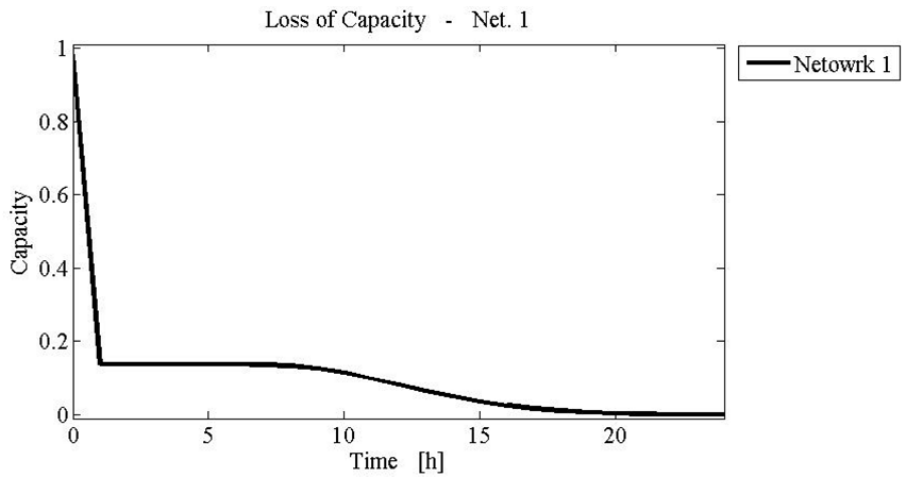

**Figure 8.58** Loss of capacity of the electric network for the detailed model.

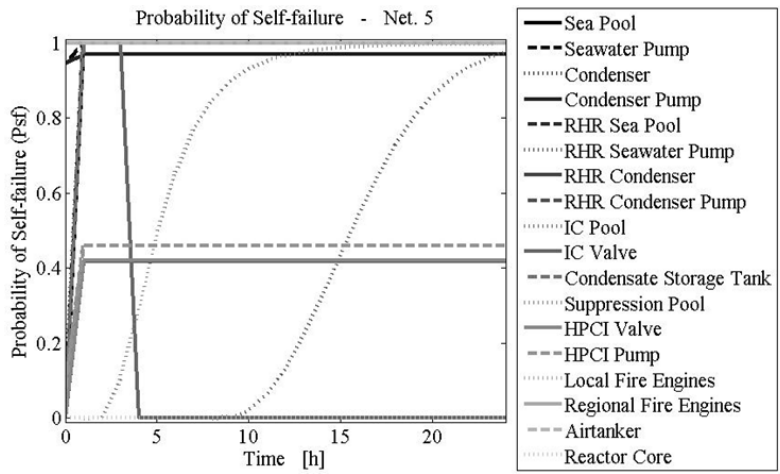

**Figure 8.59** Probability of self-failure of nodes of the water network for the detailed model.

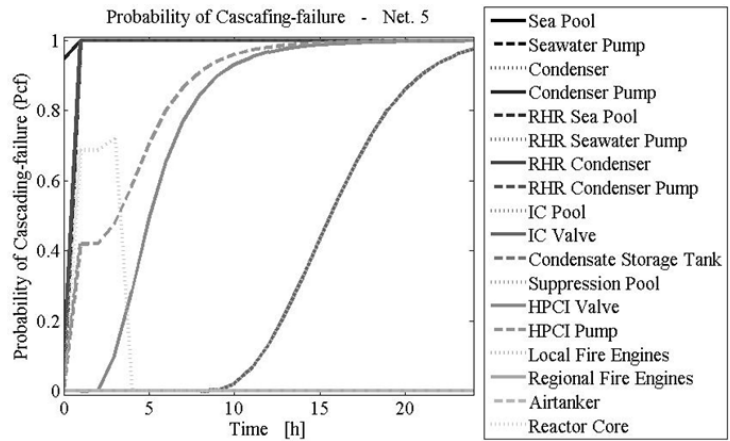

**Figure 8.60** Probability of cascading failure of nodes of the water network for the detailed model

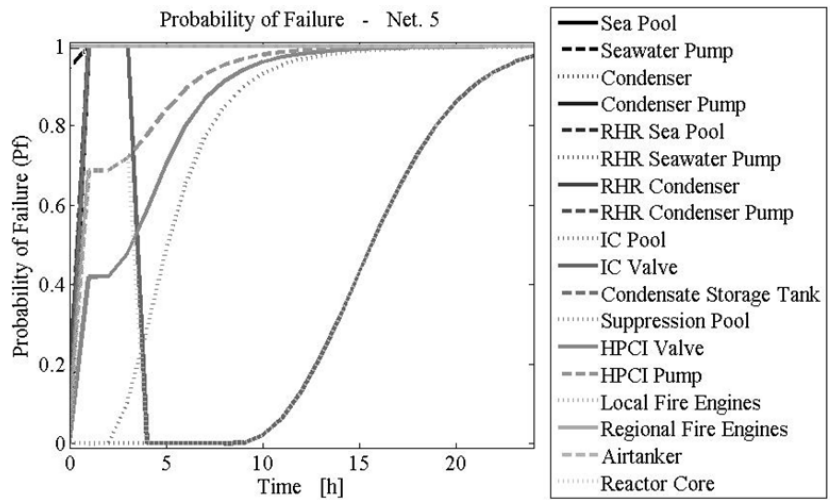

**Figure 8.61** Probability of failure of nodes of the water network for the detailed model.

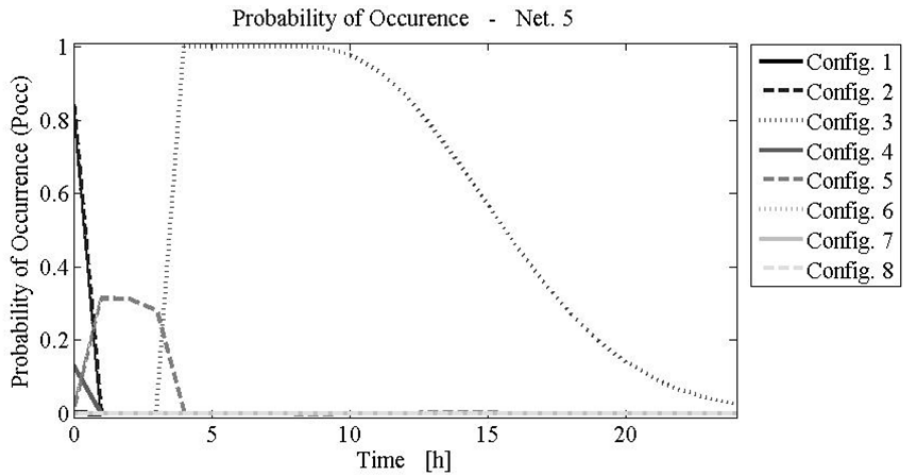

**Figure 8.62** Probability of occurrence of configurations of the water network for the detailed model.

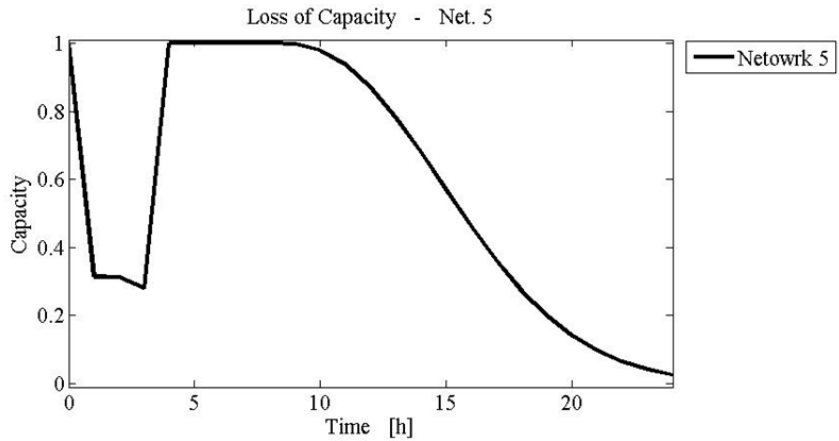

**Figure 8.63** Loss-of-capacity of the water network for the detailed model.

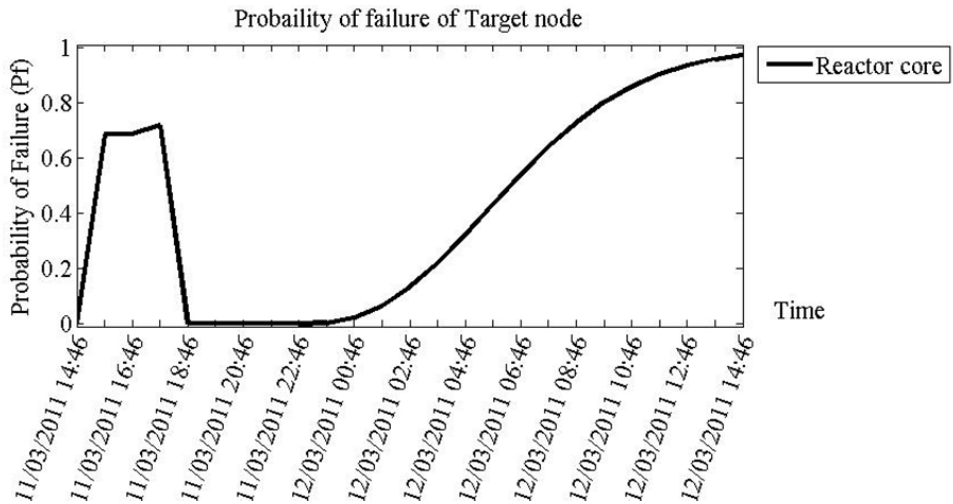

**Figure 8.64** Probability of failure of the reactor core cooling for the detailed model.

### *8.3.3.3 Comparison with the 2011 Fukushima NPP disaster*

Next we ran a sensitivity analysis of results obtained applying the modified IIM to both the simplified- and the regional-scale model. The aim was to evaluate if additional detail would result in a more accurate results than those results obtained with a simple model at a certain time step. Typical in civil engineering design, a preliminary simple analysis is performed, which then is refined based on the results of the simple analysis. Here it is shown that in certain cases the level of detail is critical as it determined the sequence of events. Because the model's topology and parameters were calibrated based on a real case study, it will be shown that in certain cases the level of detail is critical as it determined the sequence of events.

The events that occurred on March 11, 2011, at Unit 1 of the Fukushima Daiichi NNP are represented in the form of event trees. The red line indicates the sequence of events that took place. The probability of occurrence of these configurations is, of course, $P_{ooc} = 1$, because it is the real scenario. For every network, tables compare the $P_{ooc}$ obtained from the simple model and the detailed model with the real case study. The problem is analyzed at three different times; after the earthquake, after the tsunami, and before the arrival of the rescue operation at the NPP.

Figure 8.65 shows the status of the reactor at 14.46. A study of the event trees show how the actual configuration of the power network relied on diesel generators (Figure 8.67), while both the RHR and the IC cooling systems for the water network were active (Figure 8.68). The status of the electric network is well represented by both models (Table 8.2). For the water network, however, the simplified model considers all the power sources feeding the ordinary cooling line and the prediction is incorrect; see Table 8.3.

At 15.35, the second tsunami wave hit the NPP. Both the electric and water network lost their capacity because of damage to their nodes and cascading effects. Although the electric network is well represented by both models, only the detailed model fits results of the water network. The simple model lacks the presence of the IC valve, which in reality was closed and couldn't be remotely activated. Thus, the cooling system of the simple model relies on the IC, while the IC valve in the detailed model was considered closed. Therefore, the probability of a loss of capacity is $LoC = 86,3\%$.

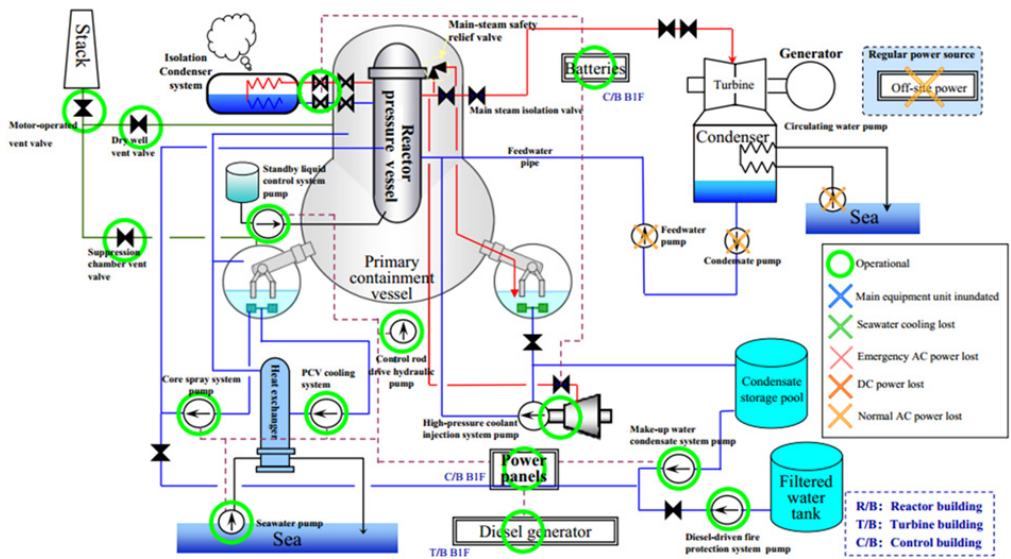

**Figure 8.65** Post-earthquake situation at Unit 1, Fukushima Nuclear Power Plant (source: TEPCO [2012]).

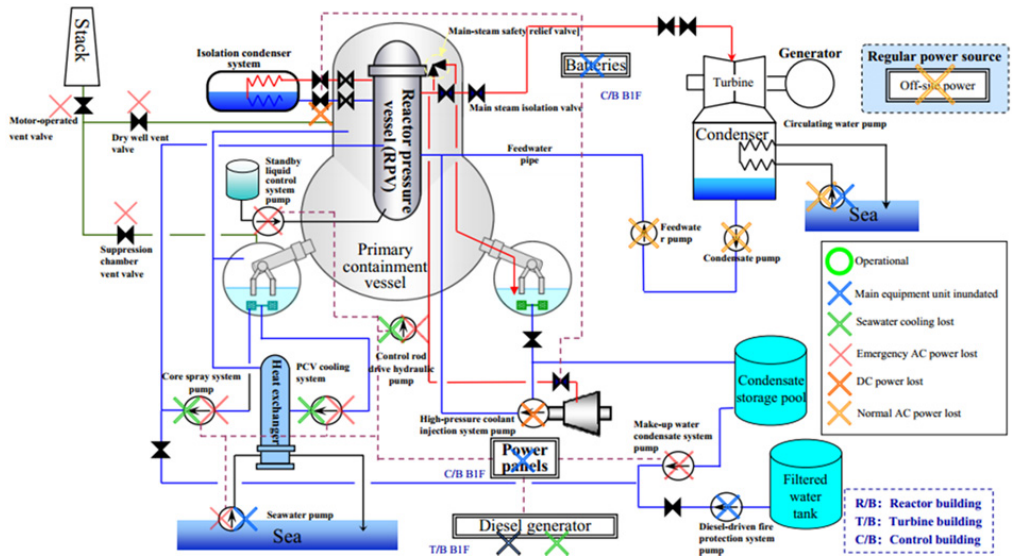

**Figure 8.66** Post-tsunami situation at Unit 1, Fukushima Nuclear Power Plant (source: TEPCO [2012]).

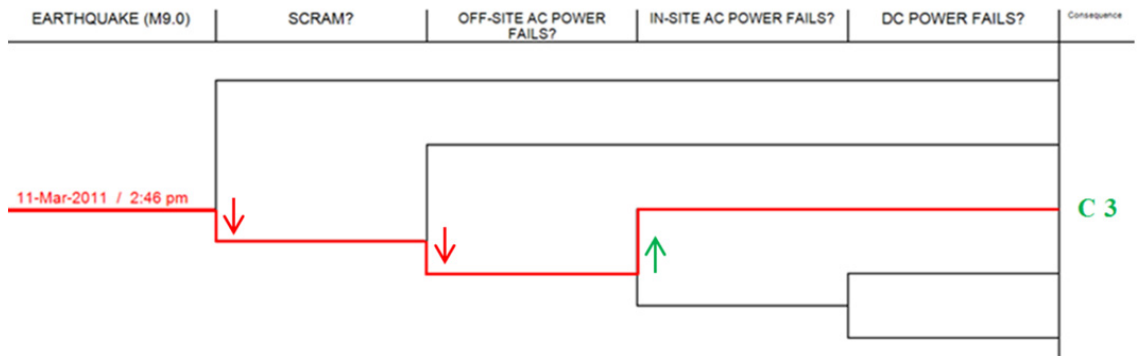

**Figure 8.67** Event tree of the post-earthquake situation for the electric network at Fukushima.

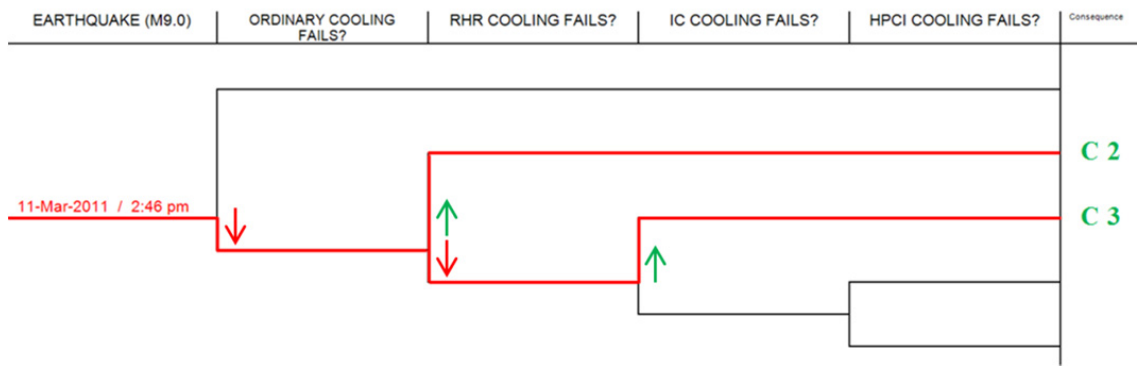

**Figure 8.68** Event tree of the post-earthquake situation for the water network at Fukushima.

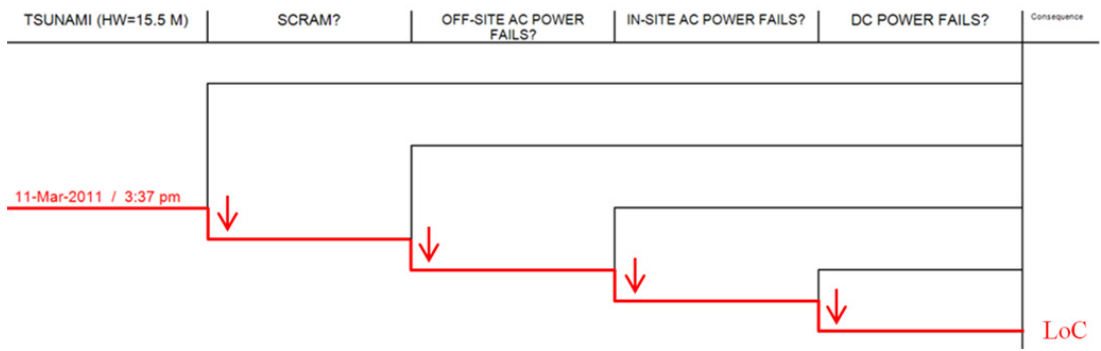

**Figure 8.69** Event tree of the post-tsunami situation for the electric network at Fukushima.

**Table 8.2** Comparison of results for the electric network after the earthquake.

| Configuration | | Probability of occurrence after the earthquake | | |
|---|---|---|---|---|
| | | Fukushima NPP | Simplified model | Detailed model |
| Self-generation | *(-/1) | 0 % | - | 0 % |
| Off-site AC power | * (1/2) | 0 % | 0 % | 0 % |
| In-site AC power | * (2/3) | 100 % | 92,4 % | 95,1 % |
| DC power | * (3/4) | 0 % | 6,5 % | 3,3 % |
| Loss of Capacity | | 0 % | 1,1 % | 1,6 % |

* refers to (simplified model numeration / detailed model numeration)

**Table 8.3** Comparison of results for the water network after the earthquake.

| Configuration | | Probability of occurrence after the earthquake | | |
|---|---|---|---|---|
| | | Fukushima NPP | Configuration | Fukushima NPP |
| Ordinary cooling | *(1/1) | 0 % | **83,9 %** | 0 % |
| RHR cooling | * (-/2) | **100 %** | 16,1 % | **85,6 %** |
| IC cooling | * (2/3) | | | |
| HPCI cooling | * (3+4/4+5) | 0 % | 0 % | 14,1 % |
| Loss of Capacity | | 0 % | 0 % | 0,3 % |

* refers to (simplified model numeration / detailed model numeration)

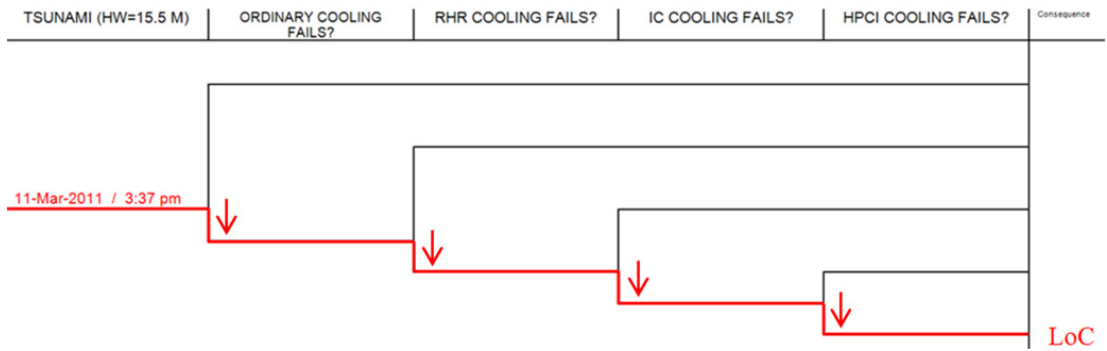

**Figure 8.70** Event tree of the post-tsunami situation for the water network at Fukushima.

**Table 8.4** Comparison of results for the electric network after the tsunami.

| Configuration | | Probability of Occurrence after the ssunami | | |
|---|---|---|---|---|
| | | **Fukushima NPP** | **Simplified model** | **Detailed model** |
| Self-generation | *(-/1) | 0 % | - | 0 % |
| Off-site AC power | * (1/2) | 0 % | 0 % | 0 % |
| In-site AC power | * (2/3) | 0 % | 0 % | 0 % |
| DC power | * (3/4) | 0 % | 8,5 % | 13,7 % |
| Loss of Capacity | | 100 % | 91,5 % | 86,3 % |

* refers to (simplified model numeration / detailed model numeration)

**Table 8.5** Comparison of results for the water network after the tsunami.

| Configuration | | Probability of occurrence after the ssunami | | |
|---|---|---|---|---|
| | | **Fukushima NPP** | **Simplified model** | **Detailed model** |
| Ordinary cooling | *(1/1) | 0 % | 0 % | 0 % |
| RHR cooling | * (-/2) | 0 % | - | 0 % |
| IC cooling | * (2/3) | 0 % | 100 % | 0,1 % |
| HPCI cooling | *(3+4/4+5) | 0 % | 0 % | 31,4 % |
| Loss of Capacity | | 100 % | 0 % | 68,5 % |

* refers to (simplified model numeration / detailed model numeration)

#### 8.3.3.3.1 Before water injection by rescuers

Twenty-four hours after the earthquake, all the systems are down and the core is melting. Temporal effects bring both the simplified and the detail model towards the condition of LoC of both networks. The detailed model is much more precise because it began from a more accurate situation post-tsunami.

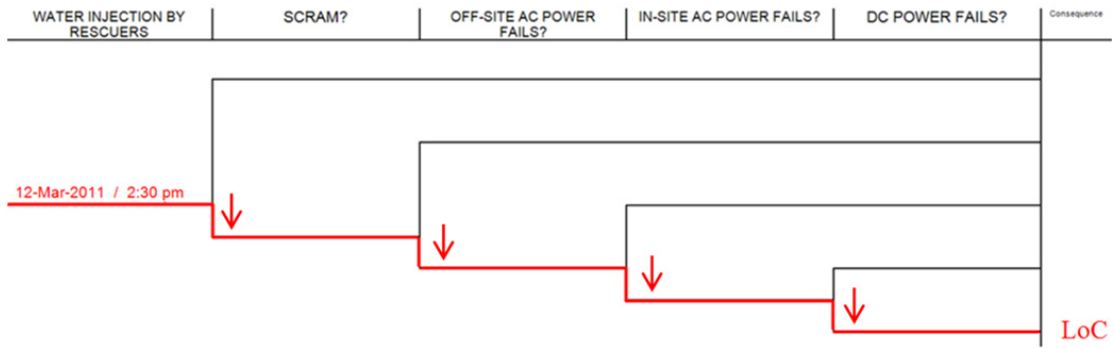

**Figure 8.71** Event tree of the before water injection situation for the electric network at Fukushima.

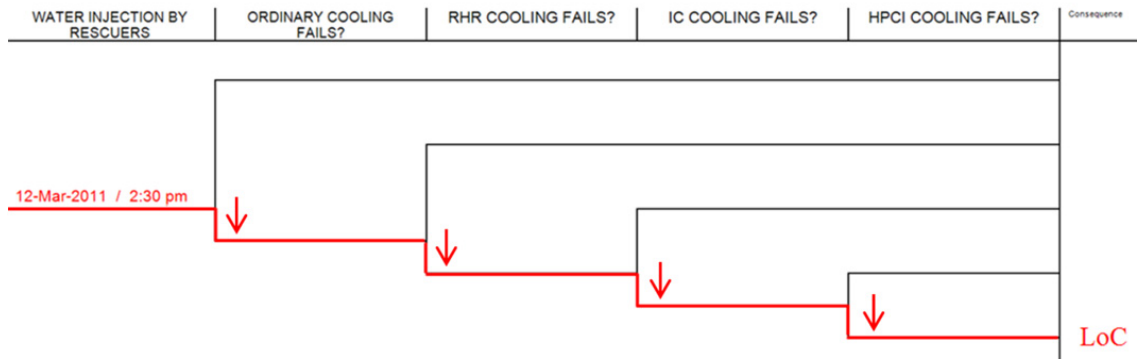

**Figure 8.72** Event tree of the before water injection situation for the water network at Fukushima.

**Table 8.6** Comparison of results for the electric network before water injection.

| Configuration | | Probability of Occurrence before water injection by rescuers | | |
| --- | --- | --- | --- | --- |
| | | Fukushima NPP | Simplified model | Detailed model |
| Self-generation | *(-/1) | 0 % | - | 0 % |
| Off-site AC power | * (1/2) | 0 % | 0 % | 0 % |
| In-site AC power | * (2/3) | 0 % | 0 % | 0 % |
| DC power | * (3/4) | 0 % | 0,4 % | 0,8 % |
| Loss of Capacity | | 100 % | 99,6 % | 99,2 % |

* refers to (simplified model numeration / detailed model numeration)

Table 8.7    Comparison of results for the water network before water injection.

| Configuration | | Probability of occurrence before water injection by rescuer | | |
|---|---|---|---|---|
| | | Fukushima NPP | Simplified model | Detailed model |
| Ordinary cooling | *(1/1) | 0 % | 0 % | 0 % |
| RHR cooling | * (-/2) | 0 % | - | 0 % |
| IC cooling | * (2/3) | 0 % | 4,0 % | 4,1 % |
| HPCI cooling | * (3+4/4+5) | 0 % | 39,2 % | 0 % |
| Loss of Capacity | | 100 % | 56,8 % | 95,9 % |

* refers to (simplified model numeration / detailed model numeration)

### 8.3.3.4  Considerations

The simulations show that although both models effectively reproduced the on-going situation for the electric network, results for the water network were mixed. The only difference between the two models for the electric network were the number of power panels and the split of the supply lines for the target nodes. Target nodes cannot propagate upstream, implying a null effect. In addition to this, the electric network is independent from all of the other networks, so enhancing details of other networks doesn't influence it.

This was not the case for the water network because different sources of power feed different configurations of the water network. The add-on of valves helps to model the inoperability of certain configurations. The simulation with the simplified model proved unreliable. Some aspects of the model, such as the malfunctioning of the HPCI for uncertain reasons, was not well represented, but we should take into consideration the generic data adopted and the difficulty in modeling physical reactions, i.e., the destabilization of reactor bars and pressurization of the core. In conclusion, the detailed system modeled with a greater accuracy the timeline of events that took place at the Fukushima NPP.

## 8.4    REMARKS AND CONCLUSIONS

The operability of lifeline systems post-event is fundamental for the well-being of the community. Disrupting events and major failures of these systems are increasing in number because of the increment of hazards characterizing modern society. To reduce costs related to their failure and recovery, it is important to assess their vulnerability. The first step in addressing this issue to understand that it is not merely a question of analyzing and protecting single networks. Interdependencies and cascading effects play a major role. Interdependencies are not time independent. Various national and federal projects in the U.S. and other parts of the world are addressing this topic to better understand these dynamics and how to mitigate undesirable consequences.

For the analysis of performance, various types of model are commonly employed. Herein, we discussed the strengths and weaknesses of the Input-output Inoperability Method (IIM) developed by Haimes et al. [2005]. After having analyzed its potentialities and limitations, a new methodology was presented. This modified IIM contains three implementations of the traditional IIM. At first, a probabilistic approach is adopted. The modified model incorporates hazard curves, fragility functions, and probabilities of failure, which are quantities much easier to combine and understand than the structure of the original IIM. The second departure was the adoption of a multilayered approach for modeling different interdependent networks. This allows a rapid and intuitive combination of an analysis run on separate networks, and at the same time highlights of the role of interdependencies in the dynamics of failure. The third and major implementation was the adoption of a tensor notation with the aim of taking into account the temporal dimension of the problem. The variable topology and the probability of occurrence of different possible and mutually exclusive configurations enabled the computation of temporal effects, i.e., the autonomy run out of certain nodes, thus providing useful information about the sequence of events disrupting the networks. This modified IIM provides reliability results that are comparable to those obtained performing a probability risk assessment.

The methodology has been applied to a real case study. Lifeline networks serving a nuclear power plant were modeled, from the regional scale to the local scale. To validate the model, the results were compared to the sequence of events that led to the 2011 Fukushima nuclear power plant disaster. A sensitivity analysis of the influence of the level of detail demonstrated that it is important to catch relevant connections among nodes and different networks. In regards to Fukushima, the results obtained using a simplified model were not reliable; however, the detailed model was a reasonably accurate fit with a relatively low approximation. In conclusion, it is possible to affirm that the modified IIM proposed herein provides reliable results in terms of the connectivity of the system. With the addition of tensor notation, it provides the same results as a PRA but begins from topological and fragility data instead of from logical considerations of the system.

# PEER REPORTS

PEER reports are available as a free PDF download from http://peer.berkeley.edu/publications/peer_reports_complete.html. Printed hard copies of PEER reports can be ordered directly from our printer by following the instructions at http://peer.berkeley.edu/publications/peer_reports.html. For other related questions about the PEER Report Series, contact the Pacific Earthquake Engineering Research Center, 325 Davis Hall Mail Code 1792, Berkeley, CA 94720. Tel.: (510) 642-3437; Fax: (510) 642-1655; Email: peer-center@berkeley.edu.

**PEER 2016/08**  *Resilience of Critical Structures, Infrastructure, and Communities*, Gian Paolo Cimellaro, Ali Zamani-Noori, Omar Kamouh, Vesna Terzic, and Stephen A. Mahin, September 2016.

**PEER 2016/07**  *Hybrid Simulation Theory for a Classical Nonlinear Dynamical System*, Paul L. Drazin and Sanjay Govindjee, September 2016.

**PEER 2016/05**  *Ground-Motion Prediction Equations for Arias Intensity Consistent with the NGA-West2 Ground-Motion Models.* Charlotte Abrahamson, Hao-Jun Michael Shi, and Brian Yang. July 2016.

**PEER 2016/04**  *The $M_W$ 6.0 South Napa Earthquake of August 24, 2014: A Wake-Up Call for Renewed Investment in Seismic Resilience Across California.* Prepared for the California Seismic Safety Commission, Laurie A. Johnson and Stephen A. Mahin. May 2016.

**PEER 2016/03**  *Simulation Confidence in Tsunami-Driven Overland Flow.* Patrick Lynett. May 2016.

**PEER 2016/02**  *Semi-Automated Procedure for Windowing time Series and Computing Fourier Amplitude Spectra for the NGA-West2 Database.* Tadahiro Kishida, Olga-Joan Ktenidou, Robert B. Darragh, and Walter J. Silva. May 2016.

**PEER 2016/01**  *A Methodology for the Estimation of Kappa ($\kappa$) from Large Datasets: Example Application to Rock Sites in the NGA-East Database and Implications on Design Motions.* Olga-Joan Ktenidou, Norman A. Abrahamson, Robert B. Darragh, and Walter J. Silva. April 2016.

**PEER 2015/13**  *Self-Centering Precast Concrete Dual-Steel-Shell Columns for Accelerated Bridge Construction: Seismic Performance, Analysis, and Design.* Gabriele Guerrini, José I. Restrepo, Athanassios Vervelidis, and Milena Massari. December 2015.

**PEER 2015/12**  *Shear-Flexure Interaction Modeling for Reinforced Concrete Structural Walls and Columns under Reversed Cyclic Loading.* Kristijan Kolozvari, Kutay Orakcal, and John Wallace. December 2015.

**PEER 2015/11**  *Selection and Scaling of Ground Motions for Nonlinear Response History Analysis of Buildings in Performance-Based Earthquake Engineering.* N. Simon Kwong and Anil K. Chopra. December 2015.

**PEER 2015/10**  *Structural Behavior of Column-Bent Cap Beam-Box Girder Systems in Reinforced Concrete Bridges Subjected to Gravity and Seismic Loads. Part II: Hybrid Simulation and Post-Test Analysis.* Mohamed A. Moustafa and Khalid M. Mosalam. November 2015.

**PEER 2015/09**  *Structural Behavior of Column-Bent Cap Beam-Box Girder Systems in Reinforced Concrete Bridges Subjected to Gravity and Seismic Loads. Part I: Pre-Test Analysis and Quasi-Static Experiments.* Mohamed A. Moustafa and Khalid M. Mosalam. September 2015.

**PEER 2015/08**  *NGA-East: Adjustments to Median Ground-Motion Models for Center and Eastern North America.* August 2015.

**PEER 2015/07**  *NGA-East: Ground-Motion Standard-Deviation Models for Central and Eastern North America.* Linda Al Atik. June 2015.

**PEER 2015/06**  *Adjusting Ground-Motion Intensity Measures to a Reference Site for which $V_{S30}$ = 3000 m/sec.* David M. Boore. May 2015.

**PEER 2015/05**  *Hybrid Simulation of Seismic Isolation Systems Applied to an APR-1400 Nuclear Power Plant.* Andreas H. Schellenberg, Alireza Sarebanha, Matthew J. Schoettler, Gilberto Mosqueda, Gianmario Benzoni, and Stephen A. Mahin. April 2015.

**PEER 2015/04**  *NGA-East: Median Ground-Motion Models for the Central and Eastern North America Region.* April 2015.

**PEER 2015/03**  *Single Series Solution for the Rectangular Fiber-Reinforced Elastomeric Isolator Compression Modulus.* James M. Kelly and Niel C. Van Engelen. March 2015.

**PEER 2015/02**  *A Full-Scale, Single-Column Bridge Bent Tested by Shake-Table Excitation.* Matthew J. Schoettler, José I. Restrepo, Gabriele Guerrini, David E. Duck, and Francesco Carrea. March 2015.

**PEER 2015/01**  *Concrete Column Blind Prediction Contest 2010: Outcomes and Observations.* Vesna Terzic, Matthew J. Schoettler, José I. Restrepo, and Stephen A Mahin. March 2015.

| | |
|---|---|
| **PEER 2012/07** | *Earthquake Engineering for Resilient Communities: 2012 PEER Internship Program Research Report Collection.* Heidi Tremayne (Editor), Stephen A. Mahin (Editor), Collin Anderson, Dustin Cook, Michael Erceg, Carlos Esparza, Jose Jimenez, Dorian Krausz, Andrew Lo, Stephanie Lopez, Nicole McCurdy, Paul Shipman, Alexander Strum, Eduardo Vega. December 2012. |
| **PEER 2012/06** | *Fragilities for Precarious Rocks at Yucca Mountain.* Matthew D. Purvance, Rasool Anooshehpoor, and James N. Brune. December 2012. |
| **PEER 2012/05** | *Development of Simplified Analysis Procedure for Piles in Laterally Spreading Layered Soils.* Christopher R. McGann, Pedro Arduino, and Peter Mackenzie–Helnwein. December 2012. |
| **PEER 2012/04** | *Unbonded Pre-Tensioned Columns for Bridges in Seismic Regions.* Phillip M. Davis, Todd M. Janes, Marc O. Eberhard, and John F. Stanton. December 2012. |
| **PEER 2012/03** | *Experimental and Analytical Studies on Reinforced Concrete Buildings with Seismically Vulnerable Beam-Column Joints.* Sangjoon Park and Khalid M. Mosalam. October 2012. |
| **PEER 2012/02** | *Seismic Performance of Reinforced Concrete Bridges Allowed to Uplift during Multi-Directional Excitation.* Andres Oscar Espinoza and Stephen A. Mahin. July 2012. |
| **PEER 2012/01** | *Spectral Damping Scaling Factors for Shallow Crustal Earthquakes in Active Tectonic Regions.* Sanaz Rezaeian, Yousef Bozorgnia, I. M. Idriss, Kenneth Campbell, Norman Abrahamson, and Walter Silva. July 2012. |
| **PEER 2011/10** | *Earthquake Engineering for Resilient Communities: 2011 PEER Internship Program Research Report Collection.* Heidi Faison and Stephen A. Mahin, Editors. December 2011. |
| **PEER 2011/09** | *Calibration of Semi-Stochastic Procedure for Simulating High-Frequency Ground Motions.* Jonathan P. Stewart, Emel Seyhan, and Robert W. Graves. December 2011. |
| **PEER 2011/08** | *Water Supply in regard to Fire Following Earthquake.* Charles Scawthorn. November 2011. |
| **PEER 2011/07** | *Seismic Risk Management in Urban Areas.* Proceedings of a U.S.-Iran-Turkey Seismic Workshop. September 2011. |
| **PEER 2011/06** | *The Use of Base Isolation Systems to Achieve Complex Seismic Performance Objectives.* Troy A. Morgan and Stephen A. Mahin. July 2011. |
| **PEER 2011/05** | *Case Studies of the Seismic Performance of Tall Buildings Designed by Alternative Means.* Task 12 Report for the Tall Buildings Initiative. Jack Moehle, Yousef Bozorgnia, Nirmal Jayaram, Pierson Jones, Mohsen Rahnama, Nilesh Shome, Zeynep Tuna, John Wallace, Tony Yang, and Farzin Zareian. July 2011. |
| **PEER 2011/04** | *Recommended Design Practice for Pile Foundations in Laterally Spreading Ground.* Scott A. Ashford, Ross W. Boulanger, and Scott J. Brandenberg. June 2011. |
| **PEER 2011/03** | *New Ground Motion Selection Procedures and Selected Motions for the PEER Transportation Research Program.* Jack W. Baker, Ting Lin, Shrey K. Shahi, and Nirmal Jayaram. March 2011. |
| **PEER 2011/02** | *A Bayesian Network Methodology for Infrastructure Seismic Risk Assessment and Decision Support.* Michelle T. Bensi, Armen Der Kiureghian, and Daniel Straub. March 2011. |
| **PEER 2011/01** | *Demand Fragility Surfaces for Bridges in Liquefied and Laterally Spreading Ground.* Scott J. Brandenberg, Jian Zhang, Pirooz Kashighandi, Yili Huo, and Minxing Zhao. March 2011. |
| **PEER 2010/05** | *Guidelines for Performance-Based Seismic Design of Tall Buildings.* Developed by the Tall Buildings Initiative. November 2010. |
| **PEER 2010/04** | *Application Guide for the Design of Flexible and Rigid Bus Connections between Substation Equipment Subjected to Earthquakes.* Jean-Bernard Dastous and Armen Der Kiureghian. September 2010. |
| **PEER 2010/03** | *Shear Wave Velocity as a Statistical Function of Standard Penetration Test Resistance and Vertical Effective Stress at Caltrans Bridge Sites.* Scott J. Brandenberg, Naresh Bellana, and Thomas Shantz. June 2010. |
| **PEER 2010/02** | *Stochastic Modeling and Simulation of Ground Motions for Performance-Based Earthquake Engineering.* Sanaz Rezaeian and Armen Der Kiureghian. June 2010. |
| **PEER 2010/01** | *Structural Response and Cost Characterization of Bridge Construction Using Seismic Performance Enhancement Strategies.* Ady Aviram, Božidar Stojadinović, Gustavo J. Parra-Montesinos, and Kevin R. Mackie. March 2010. |
| **PEER 2009/03** | *The Integration of Experimental and Simulation Data in the Study of Reinforced Concrete Bridge Systems Including Soil-Foundation-Structure Interaction.* Matthew Dryden and Gregory L. Fenves. November 2009. |
| **PEER 2009/02** | *Improving Earthquake Mitigation through Innovations and Applications in Seismic Science, Engineering, Communication, and Response.* Proceedings of a U.S.-Iran Seismic Workshop. October 2009. |

**PEER 2006/11**  *Probabilistic Seismic Demand Analysis Using Advanced Ground Motion Intensity Measures, Attenuation Relationships, and Near-Fault Effects.* Polsak Tothong and C. Allin Cornell. March 2007.

**PEER 2006/10**  *Application of the PEER PBEE Methodology to the I-880 Viaduct.* Sashi Kunnath. February 2007.

**PEER 2006/09**  *Quantifying Economic Losses from Travel Forgone Following a Large Metropolitan Earthquake.* James Moore, Sungbin Cho, Yue Yue Fan, and Stuart Werner. November 2006.

**PEER 2006/08**  *Vector-Valued Ground Motion Intensity Measures for Probabilistic Seismic Demand Analysis.* Jack W. Baker and C. Allin Cornell. October 2006.

**PEER 2006/07**  *Analytical Modeling of Reinforced Concrete Walls for Predicting Flexural and Coupled–Shear-Flexural Responses.* Kutay Orakcal, Leonardo M. Massone, and John W. Wallace. October 2006.

**PEER 2006/06**  *Nonlinear Analysis of a Soil-Drilled Pier System under Static and Dynamic Axial Loading.* Gang Wang and Nicholas Sitar. November 2006.

**PEER 2006/05**  *Advanced Seismic Assessment Guidelines.* Paolo Bazzurro, C. Allin Cornell, Charles Menun, Maziar Motahari, and Nicolas Luco. September 2006.

**PEER 2006/04**  *Probabilistic Seismic Evaluation of Reinforced Concrete Structural Components and Systems.* Tae Hyung Lee and Khalid M. Mosalam. August 2006.

**PEER 2006/03**  *Performance of Lifelines Subjected to Lateral Spreading.* Scott A. Ashford and Teerawut Juirnarongrit. July 2006.

**PEER 2006/02**  *Pacific Earthquake Engineering Research Center Highway Demonstration Project.* Anne Kiremidjian, James Moore, Yue Yue Fan, Nesrin Basoz, Ozgur Yazali, and Meredith Williams. April 2006.

**PEER 2006/01**  *Bracing Berkeley. A Guide to Seismic Safety on the UC Berkeley Campus.* Mary C. Comerio, Stephen Tobriner, and Ariane Fehrenkamp. January 2006.

**PEER 2005/16**  *Seismic Response and Reliability of Electrical Substation Equipment and Systems.* Junho Song, Armen Der Kiureghian, and Jerome L. Sackman. April 2006.

**PEER 2005/15**  *CPT-Based Probabilistic Assessment of Seismic Soil Liquefaction Initiation.* R. E. S. Moss, R. B. Seed, R. E. Kayen, J. P. Stewart, and A. Der Kiureghian. April 2006.

**PEER 2005/14**  *Workshop on Modeling of Nonlinear Cyclic Load-Deformation Behavior of Shallow Foundations.* Bruce L. Kutter, Geoffrey Martin, Tara Hutchinson, Chad Harden, Sivapalan Gajan, and Justin Phalen. March 2006.

**PEER 2005/13**  *Stochastic Characterization and Decision Bases under Time-Dependent Aftershock Risk in Performance-Based Earthquake Engineering.* Gee Liek Yeo and C. Allin Cornell. July 2005.

**PEER 2005/12**  *PEER Testbed Study on a Laboratory Building: Exercising Seismic Performance Assessment.* Mary C. Comerio, Editor. November 2005.

**PEER 2005/11**  *Van Nuys Hotel Building Testbed Report: Exercising Seismic Performance Assessment.* Helmut Krawinkler, Editor. October 2005.

**PEER 2005/10**  *First NEES/E-Defense Workshop on Collapse Simulation of Reinforced Concrete Building Structures.* September 2005.

**PEER 2005/09**  *Test Applications of Advanced Seismic Assessment Guidelines.* Joe Maffei, Karl Telleen, Danya Mohr, William Holmes, and Yuki Nakayama. August 2006.

**PEER 2005/08**  *Damage Accumulation in Lightly Confined Reinforced Concrete Bridge Columns.* R. Tyler Ranf, Jared M. Nelson, Zach Price, Marc O. Eberhard, and John F. Stanton. April 2006.

**PEER 2005/07**  *Experimental and Analytical Studies on the Seismic Response of Freestanding and Anchored Laboratory Equipment.* Dimitrios Konstantinidis and Nicos Makris. January 2005.

**PEER 2005/06**  *Global Collapse of Frame Structures under Seismic Excitations.* Luis F. Ibarra and Helmut Krawinkler. September 2005.

**PEER 2005//05**  *Performance Characterization of Bench- and Shelf-Mounted Equipment.* Samit Ray Chaudhuri and Tara C. Hutchinson. May 2006.

**PEER 2005/04**  *Numerical Modeling of the Nonlinear Cyclic Response of Shallow Foundations.* Chad Harden, Tara Hutchinson, Geoffrey R. Martin, and Bruce L. Kutter. August 2005.

**PEER 2005/03**  *A Taxonomy of Building Components for Performance-Based Earthquake Engineering.* Keith A. Porter. September 2005.

**PEER 2002/04**  Consortium of Organizations for Strong-Motion Observation Systems and the Pacific Earthquake Engineering Research Center Lifelines Program: Invited Workshop on Archiving and Web Dissemination of Geotechnical Data, 4–5 October 2001. September 2002.

**PEER 2002/03**  Investigation of Sensitivity of Building Loss Estimates to Major Uncertain Variables for the Van Nuys Testbed. Keith A. Porter, James L. Beck, and Rustem V. Shaikhutdinov. August 2002.

**PEER 2002/02**  The Third U.S.-Japan Workshop on Performance-Based Earthquake Engineering Methodology for Reinforced Concrete Building Structures. July 2002.

**PEER 2002/01**  Nonstructural Loss Estimation: The UC Berkeley Case Study. Mary C. Comerio and John C. Stallmeyer. December 2001.

**PEER 2001/16**  Statistics of SDF-System Estimate of Roof Displacement for Pushover Analysis of Buildings. Anil K. Chopra, Rakesh K. Goel, and Chatpan Chintanapakdee. December 2001.

**PEER 2001/15**  Damage to Bridges during the 2001 Nisqually Earthquake. R. Tyler Ranf, Marc O. Eberhard, and Michael P. Berry. November 2001.

**PEER 2001/14**  Rocking Response of Equipment Anchored to a Base Foundation. Nicos Makris and Cameron J. Black. September 2001.

**PEER 2001/13**  Modeling Soil Liquefaction Hazards for Performance-Based Earthquake Engineering. Steven L. Kramer and Ahmed-W. Elgamal. February 2001.

**PEER 2001/12**  Development of Geotechnical Capabilities in OpenSees. Boris Jeremić. September 2001.

**PEER 2001/11**  Analytical and Experimental Study of Fiber-Reinforced Elastomeric Isolators. James M. Kelly and Shakhzod M. Takhirov. September 2001.

**PEER 2001/10**  Amplification Factors for Spectral Acceleration in Active Regions. Jonathan P. Stewart, Andrew H. Liu, Yoojoong Choi, and Mehmet B. Baturay. December 2001.

**PEER 2001/09**  Ground Motion Evaluation Procedures for Performance-Based Design. Jonathan P. Stewart, Shyh-Jeng Chiou, Jonathan D. Bray, Robert W. Graves, Paul G. Somerville, and Norman A. Abrahamson. September 2001.

**PEER 2001/08**  Experimental and Computational Evaluation of Reinforced Concrete Bridge Beam-Column Connections for Seismic Performance. Clay J. Naito, Jack P. Moehle, and Khalid M. Mosalam. November 2001.

**PEER 2001/07**  The Rocking Spectrum and the Shortcomings of Design Guidelines. Nicos Makris and Dimitrios Konstantinidis. August 2001.

**PEER 2001/06**  Development of an Electrical Substation Equipment Performance Database for Evaluation of Equipment Fragilities. Thalia Agnanos. April 1999.

**PEER 2001/05**  Stiffness Analysis of Fiber-Reinforced Elastomeric Isolators. Hsiang-Chuan Tsai and James M. Kelly. May 2001.

**PEER 2001/04**  Organizational and Societal Considerations for Performance-Based Earthquake Engineering. Peter J. May. April 2001.

**PEER 2001/03**  A Modal Pushover Analysis Procedure to Estimate Seismic Demands for Buildings: Theory and Preliminary Evaluation. Anil K. Chopra and Rakesh K. Goel. January 2001.

**PEER 2001/02**  Seismic Response Analysis of Highway Overcrossings Including Soil-Structure Interaction. Jian Zhang and Nicos Makris. March 2001.

**PEER 2001/01**  Experimental Study of Large Seismic Steel Beam-to-Column Connections. Egor P. Popov and Shakhzod M. Takhirov. November 2000.

**PEER 2000/10**  The Second U.S.-Japan Workshop on Performance-Based Earthquake Engineering Methodology for Reinforced Concrete Building Structures. March 2000.

**PEER 2000/09**  Structural Engineering Reconnaissance of the August 17, 1999 Earthquake: Kocaeli (Izmit), Turkey. Halil Sezen, Kenneth J. Elwood, Andrew S. Whittaker, Khalid Mosalam, John J. Wallace, and John F. Stanton. December 2000.

**PEER 2000/08**  Behavior of Reinforced Concrete Bridge Columns Having Varying Aspect Ratios and Varying Lengths of Confinement. Anthony J. Calderone, Dawn E. Lehman, and Jack P. Moehle. January 2001.

**PEER 2000/07**  Cover-Plate and Flange-Plate Reinforced Steel Moment-Resisting Connections. Taejin Kim, Andrew S. Whittaker, Amir S. Gilani, Vitelmo V. Bertero, and Shakhzod M. Takhirov. September 2000.

**PEER 2000/06**  Seismic Evaluation and Analysis of 230-kV Disconnect Switches. Amir S. J. Gilani, Andrew S. Whittaker, Gregory L. Fenves, Chun-Hao Chen, Henry Ho, and Eric Fujisaki. July 2000.

**PEER 2000/05**  *Performance-Based Evaluation of Exterior Reinforced Concrete Building Joints for Seismic Excitation*. Chandra Clyde, Chris P. Pantelides, and Lawrence D. Reaveley. July 2000.

**PEER 2000/04**  *An Evaluation of Seismic Energy Demand: An Attenuation Approach*. Chung-Che Chou and Chia-Ming Uang. July 1999.

**PEER 2000/03**  *Framing Earthquake Retrofitting Decisions: The Case of Hillside Homes in Los Angeles.* Detlof von Winterfeldt, Nels Roselund, and Alicia Kitsuse. March 2000.

**PEER 2000/02**  *U.S.-Japan Workshop on the Effects of Near-Field Earthquake Shaking*. Andrew Whittaker, Editor. July 2000.

**PEER 2000/01**  *Further Studies on Seismic Interaction in Interconnected Electrical Substation Equipment.* Armen Der Kiureghian, Kee-Jeung Hong, and Jerome L. Sackman. November 1999.

**PEER 1999/14**  *Seismic Evaluation and Retrofit of 230-kV Porcelain Transformer Bushings*. Amir S. Gilani, Andrew S. Whittaker, Gregory L. Fenves, and Eric Fujisaki. December 1999.

**PEER 1999/13**  *Building Vulnerability Studies: Modeling and Evaluation of Tilt-up and Steel Reinforced Concrete Buildings*. John W. Wallace, Jonathan P. Stewart, and Andrew S. Whittaker, Editors. December 1999.

**PEER 1999/12**  *Rehabilitation of Nonductile RC Frame Building Using Encasement Plates and Energy-Dissipating Devices*. Mehrdad Sasani, Vitelmo V. Bertero, James C. Anderson. December 1999.

**PEER 1999/11**  *Performance Evaluation Database for Concrete Bridge Components and Systems under Simulated Seismic Loads*. Yael D. Hose and Frieder Seible. November 1999.

**PEER 1999/10**  *U.S.-Japan Workshop on Performance-Based Earthquake Engineering Methodology for Reinforced Concrete Building Structures*. December 1999.

**PEER 1999/09**  *Performance Improvement of Long Period Building Structures Subjected to Severe Pulse-Type Ground Motions*. James C. Anderson, Vitelmo V. Bertero, and Raul Bertero. October 1999.

**PEER 1999/08**  *Envelopes for Seismic Response Vectors.* Charles Menun and Armen Der Kiureghian. July 1999.

**PEER 1999/07**  *Documentation of Strengths and Weaknesses of Current Computer Analysis Methods for Seismic Performance of Reinforced Concrete Members*. William F. Cofer. November 1999.

**PEER 1999/06**  *Rocking Response and Overturning of Anchored Equipment under Seismic Excitations*. Nicos Makris and Jian Zhang. November 1999.

**PEER 1999/05**  *Seismic Evaluation of 550 kV Porcelain Transformer Bushings.* Amir S. Gilani, Andrew S. Whittaker, Gregory L. Fenves, and Eric Fujisaki. October 1999.

**PEER 1999/04**  *Adoption and Enforcement of Earthquake Risk-Reduction Measures.* Peter J. May, Raymond J. Burby, T. Jens Feeley, and Robert Wood. August 1999.

**PEER 1999/03**  *Task 3 Characterization of Site Response General Site Categories*. Adrian Rodriguez-Marek, Jonathan D. Bray and Norman Abrahamson. February 1999.

**PEER 1999/02**  *Capacity-Demand-Diagram Methods for Estimating Seismic Deformation of Inelastic Structures: SDF Systems*. Anil K. Chopra and Rakesh Goel. April 1999.

**PEER 1999/01**  *Interaction in Interconnected Electrical Substation Equipment Subjected to Earthquake Ground Motions*. Armen Der Kiureghian, Jerome L. Sackman, and Kee-Jeung Hong. February 1999.

**PEER 1998/08**  *Behavior and Failure Analysis of a Multiple-Frame Highway Bridge in the 1994 Northridge Earthquake*. Gregory L. Fenves and Michael Ellery. December 1998.

**PEER 1998/07**  *Empirical Evaluation of Inertial Soil-Structure Interaction Effects.* Jonathan P. Stewart, Raymond B. Seed, and Gregory L. Fenves. November 1998.

**PEER 1998/06**  *Effect of Damping Mechanisms on the Response of Seismic Isolated Structures*. Nicos Makris and Shih-Po Chang. November 1998.

**PEER 1998/05**  *Rocking Response and Overturning of Equipment under Horizontal Pulse-Type Motions*. Nicos Makris and Yiannis Roussos. October 1998.

**PEER 1998/04**  *Pacific Earthquake Engineering Research Invitational Workshop Proceedings, May 14–15, 1998: Defining the Links between Planning, Policy Analysis, Economics and Earthquake Engineering*. Mary Comerio and Peter Gordon. September 1998.

**PEER 1998/03**  *Repair/Upgrade Procedures for Welded Beam to Column Connections*. James C. Anderson and Xiaojing Duan. May 1998.

# ONLINE PEER REPORTS

The following PEER reports are available by Internet only at http://peer.berkeley.edu/publications/peer_reports_complete.html.

**PEER 2012/103** *Performance-Based Seismic Demand Assessment of Concentrically Braced Steel Frame Buildings.* Chui-Hsin Chen and Stephen A. Mahin. December 2012.

**PEER 2012/102** *Procedure to Restart an Interrupted Hybrid Simulation: Addendum to PEER Report 2010/103.* Vesna Terzic and Bozidar Stojadinovic. October 2012.

**PEER 2012/101** *Mechanics of Fiber Reinforced Bearings.* James M. Kelly and Andrea Calabrese. February 2012.

**PEER 2011/107** *Nonlinear Site Response and Seismic Compression at Vertical Array Strongly Shaken by 2007 Niigata-ken Chuetsu-oki Earthquake.* Eric Yee, Jonathan P. Stewart, and Kohji Tokimatsu. December 2011.

**PEER 2011/106** *Self Compacting Hybrid Fiber Reinforced Concrete Composites for Bridge Columns.* Pardeep Kumar, Gabriel Jen, William Trono, Marios Panagiotou, and Claudia Ostertag. September 2011.

**PEER 2011/105** *Stochastic Dynamic Analysis of Bridges Subjected to Spacially Varying Ground Motions.* Katerina Konakli and Armen Der Kiureghian. August 2011.

**PEER 2011/104** *Design and Instrumentation of the 2010 E-Defense Four-Story Reinforced Concrete and Post-Tensioned Concrete Buildings.* Takuya Nagae, Kenichi Tahara, Taizo Matsumori, Hitoshi Shiohara, Toshimi Kabeyasawa, Susumu Kono, Minehiro Nishiyama (Japanese Research Team) and John Wallace, Wassim Ghannoum, Jack Moehle, Richard Sause, Wesley Keller, Zeynep Tuna (U.S. Research Team). June 2011.

**PEER 2011/103** *In-Situ Monitoring of the Force Output of Fluid Dampers: Experimental Investigation.* Dimitrios Konstantinidis, James M. Kelly, and Nicos Makris. April 2011.

**PEER 2011/102** *Ground-Motion Prediction Equations 1964–2010.* John Douglas. April 2011.

**PEER 2011/101** *Report of the Eighth Planning Meeting of NEES/E-Defense Collaborative Research on Earthquake Engineering.* Convened by the Hyogo Earthquake Engineering Research Center (NIED), NEES Consortium, Inc. February 2011.

**PEER 2010/111** *Modeling and Acceptance Criteria for Seismic Design and Analysis of Tall Buildings.* Task 7 Report for the Tall Buildings Initiative - Published jointly by the Applied Technology Council. October 2010.

**PEER 2010/110** *Seismic Performance Assessment and Probabilistic Repair Cost Analysis of Precast Concrete Cladding Systems for Multistory Buildings.* Jeffrey P. Hunt and Božidar Stojadinovic. November 2010.

**PEER 2010/109** *Report of the Seventh Joint Planning Meeting of NEES/E-Defense Collaboration on Earthquake Engineering.* Held at the E-Defense, Miki, and Shin-Kobe, Japan, September 18–19, 2009. August 2010.

**PEER 2010/108** *Probabilistic Tsunami Hazard in California.* Hong Kie Thio, Paul Somerville, and Jascha Polet, preparers. October 2010.

**PEER 2010/107** *Performance and Reliability of Exposed Column Base Plate Connections for Steel Moment-Resisting Frames.* Ady Aviram, Božidar Stojadinovic, and Armen Der Kiureghian. August 2010.

**PEER 2010/106** *Verification of Probabilistic Seismic Hazard Analysis Computer Programs.* Patricia Thomas, Ivan Wong, and Norman Abrahamson. May 2010.

**PEER 2010/105** *Structural Engineering Reconnaissance of the April 6, 2009, Abruzzo, Italy, Earthquake, and Lessons Learned.* M. Selim Günay and Khalid M. Mosalam. April 2010.

**PEER 2010/104** *Simulating the Inelastic Seismic Behavior of Steel Braced Frames, Including the Effects of Low-Cycle Fatigue.* Yuli Huang and Stephen A. Mahin. April 2010.

**PEER 2010/103** *Post-Earthquake Traffic Capacity of Modern Bridges in California.* Vesna Terzic and Božidar Stojadinović. March 2010.

**PEER 2010/102** *Analysis of Cumulative Absolute Velocity (CAV) and JMA Instrumental Seismic Intensity ($I_{JMA}$) Using the PEER–NGA Strong Motion Database.* Kenneth W. Campbell and Yousef Bozorgnia. February 2010.

**PEER 2010/101** *Rocking Response of Bridges on Shallow Foundations.* Jose A. Ugalde, Bruce L. Kutter, and Boris Jeremic. April 2010.

**PEER 2009/109** *Simulation and Performance-Based Earthquake Engineering Assessment of Self-Centering Post-Tensioned Concrete Bridge Systems.* Won K. Lee and Sarah L. Billington. December 2009.

**PEER 2009/108**  *PEER Lifelines Geotechnical Virtual Data Center.* J. Carl Stepp, Daniel J. Ponti, Loren L. Turner, Jennifer N. Swift, Sean Devlin, Yang Zhu, Jean Benoit, and John Bobbitt. September 2009.

**PEER 2009/107**  *Experimental and Computational Evaluation of Current and Innovative In-Span Hinge Details in Reinforced Concrete Box-Girder Bridges: Part 2: Post-Test Analysis and Design Recommendations.* Matias A. Hube and Khalid M. Mosalam. December 2009.

**PEER 2009/106**  *Shear Strength Models of Exterior Beam-Column Joints without Transverse Reinforcement.* Sangjoon Park and Khalid M. Mosalam. November 2009.

**PEER 2009/105**  *Reduced Uncertainty of Ground Motion Prediction Equations through Bayesian Variance Analysis.* Robb Eric S. Moss. November 2009.

**PEER 2009/104**  *Advanced Implementation of Hybrid Simulation.* Andreas H. Schellenberg, Stephen A. Mahin, Gregory L. Fenves. November 2009.

**PEER 2009/103**  *Performance Evaluation of Innovative Steel Braced Frames.* T. Y. Yang, Jack P. Moehle, and Božidar Stojadinovic. August 2009.

**PEER 2009/102**  *Reinvestigation of Liquefaction and Nonliquefaction Case Histories from the 1976 Tangshan Earthquake.* Robb Eric Moss, Robert E. Kayen, Liyuan Tong, Songyu Liu, Guojun Cai, and Jiaer Wu. August 2009.

**PEER 2009/101**  *Report of the First Joint Planning Meeting for the Second Phase of NEES/E-Defense Collaborative Research on Earthquake Engineering.* Stephen A. Mahin et al. July 2009.

**PEER 2008/104**  *Experimental and Analytical Study of the Seismic Performance of Retaining Structures.* Linda Al Atik and Nicholas Sitar. January 2009.

**PEER 2008/103**  *Experimental and Computational Evaluation of Current and Innovative In-Span Hinge Details in Reinforced Concrete Box-Girder Bridges. Part 1: Experimental Findings and Pre-Test Analysis.* Matias A. Hube and Khalid M. Mosalam. January 2009.

**PEER 2008/102**  *Modeling of Unreinforced Masonry Infill Walls Considering In-Plane and Out-of-Plane Interaction.* Stephen Kadysiewski and Khalid M. Mosalam. January 2009.

**PEER 2008/101**  *Seismic Performance Objectives for Tall Buildings.* William T. Holmes, Charles Kircher, William Petak, and Nabih Youssef. August 2008.

**PEER 2007/101**  *Generalized Hybrid Simulation Framework for Structural Systems Subjected to Seismic Loading.* Tarek Elkhoraibi and Khalid M. Mosalam. July 2007.

**PEER 2007/100**  *Seismic Evaluation of Reinforced Concrete Buildings Including Effects of Masonry Infill Walls.* Alidad Hashemi and Khalid M. Mosalam. July 2007.